%% file: dim_form_fin.tex
\documentclass[a4paper,12pt]{amsart}

\usepackage{amsmath, amsfonts, amsthm, amsfonts,url,eucal,xcolor}

\usepackage{multicol}
\usepackage[T1]{fontenc}
\usepackage[latin1]{inputenc}

\def\]{\textup{\mbox{]\hspace{-.15em}]}}}
\def\[{\textup{\mbox{[\hspace{-.15em}[}}}

\def\got{\mathfrak}

\makeatletter

\@addtoreset{equation}{section}
\makeatother

\newcommand{\Hom}{\mathrm{Hom}}

\newcommand{\pn}{\par \noindent}

\newcommand{\Frob}{\mathrm{Frob}}

\newcommand{\diag}{\mathrm{diag}}
\newcommand{\GL}{\mathrm{GL}}
\newcommand{\SL}{\mathrm{SL}}

\newcommand{\ps}{\par \smallskip}

\newcommand{\N}{\mathbb{N}}
\newcommand{\Z}{\mathbb{Z}}

\newcommand{\Q}{\mathbb{Q}}

\newcommand{\R}{\mathbb{R}}
\newcommand{\C}{\mathbb{C}}

\newcommand{\AAA}{\mathbb{A}}

\newcommand{\isomo}{\overset{\sim}{\rightarrow}}

\newcommand{\PGL}{{\rm PGL}}
\renewcommand{\SL}{{\rm SL}}

\newtheorem{definition}[subsection]{Definition}

\newtheorem{thmconda}[subsection]{Theorem$^\ast$}
\newtheorem{thmcondb}[subsection]{Theorem$^\ast{}^\ast$}
\newtheorem{prop}[subsection]{Proposition}
\newtheorem{propconda}[subsection]{Proposition$^\ast$}
\newtheorem{propcondb}[subsection]{Proposition$^\ast{}^\ast$}
\newtheorem{cor}[subsection]{Corollary}
\newtheorem{corconda}[subsection]{Corollary$^\ast$}
\newtheorem{corcondb}[subsection]{Corollary$^\ast{}^\ast$}
\newtheorem{lemme}[subsection]{Lemma}

\newtheorem{conj}[subsection]{Conjecture}
\newtheorem{prob}[subsection]{Problem}
\newtheorem{rem}[subsection]{Remark}

\newtheorem{propdef}[subsection]{Proposition-Definition}

\newenvironment{pf}
{\medskip\noindent {\it Proof --- \ }}
{\hfill\nobreak $\Box$ \par\bigbreak}

\input xy
\xyoption{all}

\begin{document}

\title{Level one algebraic cusp forms of classical groups of small ranks}
\bigskip
\bigskip
\author{ Ga\"etan Chenevier }
\author{David Renard}

\address{Centre de Math\'ematiques
Laurent Schwartz
\\ Ecole Polytechnique \\ 91128 Palaiseau Cedex \\ FRANCE }

\thanks{Ga\"etan Chenevier is supported by the C.N.R.S. and by the French ANR-10-BLAN 0114 project.}


\begin{abstract} We determine the number of level $1$, polarized,
algebraic regular, cuspidal automorphic representations of ${\rm GL}_n$
over $\Q$ of any given infinitesimal character, for essentially all $n \leq
8$.  For this, we compute the dimensions of spaces of level $1$ automorphic
forms for certain semisimple $\Z$-forms of the compact groups ${\rm SO}_7$,
${\rm SO}_8$, ${\rm SO}_9$ (and ${\rm G}_2$) and determine Arthur's
endoscopic partition of these spaces in all cases.  We also give
applications to the $121$ even lattices of rank $25$ and determinant $2$
found by Borcherds, to level one self-dual automorphic representations of
${\rm GL}_n$ with trivial infinitesimal character, and to vector valued
Siegel modular forms of genus $3$.  A part of our results are conditional to certain expected results in the theory of twisted endoscopy. 
\end{abstract}

\maketitle
\bigskip
\bigskip

\section{Introduction}
\bigskip
\subsection{A counting problem}\label{countingpb} Let $n\geq 1$ be an integer. Consider the cuspidal automorphic representations
$\pi$ of $\GL_n$ over $\Q$ (see \cite[Ch. 3]{GGPS},\cite[\S
4]{boreljacquet},\cite{cogdell}) such that : \begin{itemize} \ps \ps
\item[(a)] (polarization) $\pi^\vee \simeq \pi \otimes |\cdot|^{w}$ for some
$w\in \Z$, \ps
\item[(b)] (conductor 1) $\pi_p$ is unramified for each prime $p$, \ps
\item[(c)] (algebraicity) $\pi_\infty$ is algebraic and regular.\ps
\end{itemize}

Our main aim in this paper is to give for small values of $n$, namely for
$n\leq 8$, the number of such representations as a function of $\pi_\infty$. 
Recall that by the Harish-Chandra isomorphism, the infinitesimal character
of $\pi_\infty$ may be viewed following Langlands as a semisimple conjugacy class in ${\rm
M}_n(\C)$ (see~\S\ref{langparam},\,\S\ref{demialg}). Condition (c) means\footnote{\label{footnotealg}The term {\it algebraic} here is 
in the sense of Borel~\cite[\S 18.2]{borelcorvallis}, and is reminescent to
Weil's notion of Hecke characters of type ${\rm A}_0$ : see~\S\ref{demialg}. Langlands
also uses the term {\it of type Hodge}, e.g.  in~\cite[\S
5]{langlandsedin}. See also~\cite{buzzardgee}, who would employ here the term
{\it ${\rm L}$-algebraic}, 
for a discussion of other notions of algebraicity, as the one used
by Clozel in~\cite{clozel}.} that the eigenvalues of this
conjugacy class are distinct integers. The opposite of these integers will be called the {\it
weights} of $\pi$ and we shall denote them by $k_1 > k_2 >
\cdots > k_n$. When $n \equiv 0 \bmod 4$, we will eventually allow that $k_{n/2}=k_{n/2+1}$
but to simplify we omit this case in the dicussion for the moment. If $\pi$ satisfies (a), the necessarily unique integer $w \in \Z$ such that $\pi^\vee \simeq \pi \otimes |\cdot|^w$ will be called the {\it motivic weight} of $\pi$, and denoted $w(\pi)$.

\begin{prob}\label{mainpb} For any $n\geq 1$, determine the
number ${\rm N}(k_1,k_2,\cdots,k_n)$ of cuspidal automorphic representations 
$\pi$ of $\GL_n$ satisfying {\rm (a)}, {\rm (b)} and {\rm (c)} above, and of weights
$k_1>k_2>\cdots>k_n$. 
\end{prob}

An important finiteness result of Harish-Chandra~(\cite[Thm. 
1.1]{harishchandra}) asserts that this number is indeed finite, even if we
omit assumption (a).  As far as we know, those numbers have been previously
computed only for $n\leq 3$.  For $n=1$, the structure of the id\`eles of
$\Q$ shows that if $\pi$ satisfies (a), (b) and (c) then $w(\pi)=2\, k_1$ is
even and $\pi=|.|^{-k_1}$.  By considering the central character of $\pi$,
this also shows the relation $n \,w(\pi) = 2\sum_{i=1}^n k_i$ for general
$n$.  More interestingly, classical arguments show that ${\rm N}(k-1,0)$
coincides with the dimension of the space of cuspidal modular forms of
weight $k$ for ${\rm SL}_2(\Z)$, whose dimension is well-known (see
e.g.~\cite{serre}) and is about\footnote{We denote by $[x]$ the floor of the real number $x$.} $[k/12]$.  Observe that up to twisting $\pi$ by
$|\cdot|^{k_n}$, there is no loss of generality in assuming that $k_n=0$ in
the above problem. Moreover, condition (a) implies for $i=1,\cdots,n$
the relation $k_i+k_{n+1-i}=w(\pi)$.\ps

\subsection{Motivations} There are several motivations for this problem. A first one is the deep
conjectural relations, due on the one hand to Langlands~\cite{Lgl}, in the
lead of Shimura, Taniyama, and Weil, and on the other hand to Fontaine and
Mazur~\cite{fontainemazur}, that those numbers ${\rm N}(k_1,k_2,\cdots,k_n)$
share with arithmetic geometry and pure motives\footnote{The reader is free
here to choose his favorite definition of a pure motive~\cite{motives}.}
over $\Q$. More precisely, consider the three following type of objects :
\begin{itemize} \ps 

\item[(I)] Pure motives $M$ over $\Q$, of weight $w$ and rank $n$, with coefficients in
$\overline{\Q}$, which are : simple, of conductor $1$, such that $M^\vee
\simeq M(w)$, and whose Hodge numbers satisfy $h^{p,q}(M)=1$ if $(p,q)$ is
of the form $(k_i,w-k_i)$ and $0$ otherwise.\ps

\item[${\rm (II)}_\ell$] Continuous irreducible representations
$\rho : {\rm Gal}(\overline{\Q}/\Q) \rightarrow {\rm GL}_n(\overline{\Q}_\ell)$
which are unramified outside $\ell$, cristalline at $\ell$ with Hodge-Tate numbers $k_1>\cdots > k_n$,  and such that\footnote{Here $\omega_\ell$ denotes the $\ell$-adic cyclotomic character of ${\rm Gal}(\overline{\Q}/\Q)$, and our convention is that its Hodge-Tate number is $-1$. } $\rho^\vee \simeq \rho \otimes
\omega_\ell^w$.\ps

\item[(III)] Cuspidal automorphic representations $\pi$ of
${\rm GL}_n$ over $\Q$ satisfying (a), (b) and (c) above, of weights
$k_1>k_2>\cdots>k_n$, \ps

\end{itemize}

Here $\ell$ is a fixed prime, and $\overline{\Q}$ and $\overline{\Q}_\ell$
are fixed algebraic closures of $\Q$ and $\Q_\ell$. To discuss the
aforementionned conjectures we need to fix a pair
of fields embeddings $\iota_\infty :
\overline{\Q}\rightarrow \C$ and $\iota_\ell: \overline{\Q} \rightarrow
\overline{\Q}_\ell$.  According to Fontaine and Mazur, Grothendieck's
$\ell$-adic \'etale cohomology, viewed with $\overline{\Q}_\ell$
coefficients via $\iota_\ell$, should induce a bijection between isomorphism
classes of motives of type {\rm (I)} and isomorphism classes of Galois
representations of type ${\rm (II)}_\ell$.  Moreover, according to Langlands, the
${\rm L}$-function of the $\ell$-adic realizations of a motive of type {\rm (I)}, which
makes sense via $\iota_\infty$ and $\iota_\ell$, should be the standard
${\rm L}$-function of a unique $\pi$ of type {\rm (III)}, and vice-versa. These conjectural bijections are actually expected
to exist in greater generality (any conductor, any weights, not necessarily
polarized), but we focus on this case as it is the one we really
consider in this paper.  In particular, ${\rm N}(k_1,\cdots,k_n)$ {\it is also
the conjectural number of isomorphism classes of objects of type {\rm (I)}
or ${\rm (II)}_\ell$ for any $\ell$}.  Let us mention that there has been
recently important progresses toward those conjectural bijections.  First of
all, by the works of many authors (including Deligne, Langlands, Kottwitz,
Clozel, Harris, Taylor, Labesse, Shin, Ng\^o and Waldspurger,
see~\cite{grfa},\cite{shin} and \cite{chharris}), if $\pi$ is of type
(III) then there is a unique associated semisimple
representation $\rho_\pi : {\rm Gal}(\overline{\Q}/\Q) \rightarrow
{\rm GL}_n(\overline{\Q}_\ell)$ of type ${\rm (II)}_\ell$ with the same
${\rm L}$-function as $\pi$ (via $\iota_\infty,\iota_\ell$), up to the fact that $\rho_{\pi,\iota}$ is only
known to be irreducible when $n\leq 5$ (see~\cite{calegarigee}). Second,
the advances in modularity results in the lead of Wiles and Taylor,
such as the proof of Serre's conjecture by Khare and Wintenberger (see e.g.
\cite{khare1}), or the recent results~\cite{bggt}, contain striking results toward the
converse statement. \ps

An important source of objects ot type (I) or (II) comes from the cohomology of proper smooth schemes (or
stacks) over $\Z$, about which solving problem~\ref{mainpb} would
thus shed interesting lights. This applies in particular to the moduli
spaces $\overline{\mathcal{M}_{g,n}}$ of stable curves of genus $g$ with $n$-marked points
and to certain spaces attached to the moduli spaces of principally polarized abelian varieties (see e.g.~\cite{BFVdG} and
\cite{faltings}). As an example, the vanishing of some ${\rm
N}(k_1,\cdots,k_n)$ translates to a conjectural non-existence theorem about Galois representations or motives. 
A famous result in this style is the proof by
Abrashkin and Fontaine that there are no abelian scheme over $\Z$ (hence no
projective smooth curve over $\Z$ of nonzero genus), which had been
conjectured by Shaffarevich
(see~\cite{fontaineabelienne},\cite{fontaineprojlisse}).  The corresponding vanishing statement about
cuspidal automorphic forms had been previously checked by Mestre and Serre
(see~\cite{mestre}).  See also Khare's paper~\cite{kh} for other conjectures
in this spirit as well as a discussion about the applications to the
generalized Serre's conjecture. \ps

A second motivation, which is perhaps more exotic, is the well-known
problem of finding an integer $n\geq 1$ such that the cuspidal
cohomology ${\rm H}^\ast_{\rm cusp}({\rm SL}_n(\Z),\Q)$ does not vanish.  It
would be enough to find an integer $n\geq 1$ such that ${\rm
N}(n-1,\cdots,2,1,0) \neq 0$. Results of Mestre~\cite{mestre},
Fermigier~\cite{fermigier} and Miller~\cite{miller} ensure that such an $n$
has to be $\geq 27$ (although those works do not assume the self-duality
condition).  We shall go back to these questions at the end of this
introduction.  \ps

Last but not least, it follows from Arthur's endoscopic classification~\cite{arthur} that the dimensions of
various spaces of modular forms for classical reductive groups over $\Z$
have a "simple" expression in terms of these numbers. Part of this paper is
actually devoted to explain this relation in very precise and concrete terms. This
includes vector valued holomorphic Siegel modular forms for ${\rm Sp}_{2g}(\Z)$ and level
$1$ algebraic automorphic forms for the $\Z$-forms of ${\rm SO}_{p,q}(\R)$ which
are semisimple over $\Z$ (such group schemes exist when $p-q \equiv 0,\pm 1
\bmod 8$). It can be used in both ways : either to deduce the dimensions
of these spaces of modular forms from the knowledge of the integers ${\rm
N}(-)$, or also to compute these last numbers from known dimension formulas. 
We will say much more about this in what follows as this is the main theme of this
paper (see Chapter \ref{sectionarthur}).  \ps

\subsection{The main result} We will now state our main theorem. As many results that we prove in
this paper, it depends on the fabulous work of Arthur in~\cite{arthur}.  As
explained {\it loc.  cit.}, Arthur's results are still conditional to
the stabilization of the twisted trace formula at the moment. All the
results below depending on this assumption will be marked by a simple star
$^\ast$.  We shall also need to use certain results concerning inner forms
of classical groups which have been announced by Arthur (see~\cite[Chap. 
9]{arthur}) but which are not yet available or even precisely stated.  We
have thus formulated the precise general results that we expect in two
explicit Conjectures~\ref{conjectureun} and~\ref{conjecturedeux}.  Those
conjectures include actually a bit more than what has been announced by
Arthur in~\cite{arthur}, namely also the standard expectation that for
Adams-Johnson archimedean Arthur parameters, there is an identification
between Arthur's packets in~\cite{arthur} and the ones of Adams and Johnson
in~\cite{AJ}.  The precise special cases that we need are detailed in
\S\ref{parmultform}.  We state in particular Arthur's multiplicity formula
in a completely explicit way, in a generality that might be useful to
arithmetic geometers.  We will go back to the shape of this formula in~\S\ref{amfintro}. All the results below depending on the
assumptions of~\cite{arthur} as well as on the
assumptions~\ref{conjectureun} and~\ref{conjecturedeux} will be marked by a
double star $^\ast{}^\ast$.  Of course, the tremendous recent progresses in
this area allow some optimism about the future of all these assumptions ! 
Besides Arthur's work, let us mention the following results which play a
crucial role here : the proof by Chaudouard, Laumon, Ng\^o and
Waldspurger of Langlands' fundamental lemma
(\cite{walspond},\cite{ngo},\cite{chaulau}), the works of
Shelstad~\cite{shelstad} and Mezo~\cite{mezopreprint} on endoscopy for real groups, and the recent works
of Labesse and Waldspurger on the twisted trace formula \cite{labwals}.  \ps

\begin{thmcondb}\label{mainthmintro} Assume $n\leq 8$ and $n \neq 7$. There is an
explicit, computable, formula for ${\rm N}(k_1,\cdots,k_n)$.
\end{thmcondb}

Although our formulas are explicit, one cannot write them down here as they
are much too big : see \S\ref{discussionformuleso7} for a discussion of the formula.  Nevertheless, we implemented them on a computer and have a
program which takes $(k_1,k_2,\dots,k_n)$ as input and returns ${\rm
N}(k_1,\dots,k_n)$.  When $k_1-k_n \leq 100$, the computation takes less
than ten minutes on our machine\footnote{Four processors Nortwood Pentium 4, 
2.80 GHz, 5570.56 BogoMIPS.} : see
the website~\cite{homepage} for some data.  We also have some partial
results concerning ${\rm N}(k_1,\dots,k_7)$.  This includes an explicit
upper bounds for these numbers as well as their values modulo $2$, which is
enough to actually determine them in quite a few cases (for instance
whenever $k_1-k_7 \leq 26$).  On the other hand, as we shall see in
Propostion~\ref{corintrosp6} below, these numbers are also closely related
to the dimensions of the spaces of vector valued Siegel modular forms for
${\rm Sp}_6(\Z)$.  In a remarkable recent work, Bergstr\"om, Faber and van
der Geer~\cite{BFVdG} actually found a conjectural explicit formula for
those dimensions.  Their method is completely different from ours : they
study the number of points over finite fields of $\mathcal{M}_{3,n}$ and of 
certain bundles over the moduli space of principally polarized abelian 
varieties of dimension $3$. Fortunately, in the few hundreds of
cases where our work allow to compute this dimension as well, it fits the
results found by the formula of these authors ! Even better, if we assume their formula we obtain in turn a conjectural
explicit formula for ${\rm N}(k_1,\cdots,k_7)$.  \ps
\subsection{Langlands-Sato-Tate groups}\label{lstintro} We not only determine ${\rm N}(k_1,\cdots,k_n)$ for $n\leq 8$ (with the caveat above for $n=7$) but
we give as well the conjectural number of $\pi$ of weights $k_1>\dots>k_n$
having any possible {\it Langlands-Sato-Tate
group}. We refer to the appendix~\ref{parsatotate} for a brief introduction to this
conjectural notion (see also~\cite[Ch. 1, appendix]{serreabelian}). Here are 
certain of its properties. \ps

First, a representation $\pi$ as above 
being given, the Langlands-Sato-Tate group of $\pi$ (or, for short, its {\it Sato-Tate
group}) is a compact Lie group $\mathcal{L}_\pi \subset \GL_n(\C)$, which is
well-defined up to ${\rm GL}_n(\C)$-conjugacy. It is "defined" as the image of the
conjectural Langlands group $\mathcal{L}_\Z$ of $\Z$, that we view as a
topological following Kottwitz, under the hypothetical morphism
$\mathcal{L}_\Z \rightarrow \GL_n(\C)$ attached to $\pi \otimes |\cdot|^{\frac{w(\pi)}{2}}$ (\cite{Lgl},\cite{kottwitz},\cite{arthurconjlan}).
The natural representation of $\mathcal{L}_\pi$ on $\C^n$ is irreducible and self-dual. \ps

The group $\mathcal{L}_\pi$ is equipped with a collection of conjugacy
classes $${\rm Frob}_p \subset \mathcal{L}_\pi$$ which are indexed by the
primes $p$, and such that for each $p$ the $\GL_n(\C)$-conjugacy class of
${\rm Frob}_p$ is the Satake parameter of $\pi_p$ multiplied by the scalar
$p^{\frac{-w(\pi)}{2}}$.  Observe that a necessary condition for this is
that the eigenvalues of the Satake parameter of $\pi_p$ all have absolute
value $p^{\frac{w(\pi)}{2}}$. This is the so-called Ramanujan conjecture
for $\pi$, and it is actually known for each $\pi$ satisfying (a) and (c)
thanks to Deligne's proof of Weil's conjectures and results of
Clozel-Harris-Labesse, Shin and Caraiani (see~\cite{grfa},\cite{shin} and
\cite{caraiani}).  \ps

The volume $1$ Haar measure of $\mathcal{L}_\pi$ induces
a natural measure on its space of conjugacy classes and one of the key
expected properties of $\mathcal{L}_\pi$ is that the ${\rm Frob}_p$ are
equidistributed in this space.  Let us mention that a pleasant consequence
of property (b) of $\pi$ is that $\mathcal{L}_\pi$ is necessarily connected (as ${\rm Spec}(\Z)$ is simply connected !).  A case-by-case argument,
solely based on the fact that $\C^n$ is a self-dual irreducible representation
of the connected compact group $\mathcal{L}_\pi$, shows that the list of all the possible Sato-Tate
groups is rather small when
$n\leq 8$ : see Appendix \ref{parsatotate}.  \ps

\subsection{The symplectic-orthogonal alternative}\label{altsointro}Our main aim now will be to discuss our results in each particular dimension
$n$.  Before doing so it will be convenient to introduce more
notations. We shall make an important use of automorphic representations
$\pi$ of $\GL_n$ satisfying an assumption which is slightly weaker than (c),
that we now have to introduce.  Assume that $\pi$ is a cuspidal automorphic
representation of $\GL_n$ over $\Q$ satisfying property (a) above, so that
the integer $w(\pi)$ still makes sense in particular.  Consider the property
: \begin{itemize} \ps

\item[(c')] The eigenvalues of the infinitesimal character of $\pi_\infty$, viewed as a semisimple conjugacy class in ${\rm M}_n(\C)$, are integers. Moreover, each of these eigenvalues has multiplicity one, except perhaps the eigenvalue $-w(\pi)/2$ which is allowed to have multiplicity $2$ when $n \equiv 0 \bmod 4$. \ps
\end{itemize}

A $\pi$ satisfying (a), (b) and (c') still has weights $k_1 \geq \cdots \geq k_n$ defined as the opposite of the eigenvalues of the infinitesimal character of $\pi_\infty$, counted with multiplicities. When (c) is not satisfied, then $n \equiv 0 \bmod 4$ and $k_1> \cdots > k_{n/2} = k_{n/2+1}  > \cdots > k_n$. It would not be difficult to extend the conjectural picture suggested by the Langlands and Fontaine-Mazur conjectures to those $\pi$'s, but we shall not do so here. Let us simply say that when $\pi$ satisfies (a), (b), (c') but not (c), one still knows how to construct a semisimple continuous Galois representation $\rho_\pi : {\rm Gal}(\overline{\Q}/\Q) \rightarrow \GL_n(\overline{\Q}_\ell)$ unramified outside $\ell$ and with the same ${\rm L}$-function as $\pi$ : see~\cite{goldring}. It is expected but not known that $\rho_\pi$ is cristalline at $\ell$, and that the Ramanujan conjecture holds for $\pi$. However, we know from~\cite{taibi} that $\rho_\pi$ is Hodge-Tate at $\ell$ (with Hodge-Tate numbers the $k_i$)  and that for any complex conjugation $c \in {\rm Gal}(\overline{\Q}/\Q)$, we have ${\rm Trace} \, \rho_\pi(c)=0$.

\ps 

We now consider a quite important property of the $\pi$ satisfying (a), (b) and (c'), namely their {\it orthogonal-symplectic alternative}.

\begin{definition}\label{defos} Let $\pi$ be a cuspidal automorphic
representation of $\GL_n$ over $\Q$ satisfying {\rm (a), (b)} and {\rm (c')} above.  We say that
$\pi$ is {\it symplectic} if $w(\pi)$ is odd, and {\it orthogonal}
otherwise.  \end{definition}

Let $k_1>\cdots>k_n$ denote the weights of $\pi$. The relation $n \,w(\pi)=
2\sum_{i=1}^n k_i$ shows that $\pi$ is necessarily orthogonal if $n$ is odd. 
Definition~\ref{defos} fits with the conjectural picture described above. 
For instance, the main theorem of~\cite{bchsigne} asserts that if $\pi$ is orthogonal (resp. 
symplectic), and if $\rho_\pi$ denotes the Galois representation associated
to $\pi$ and $(\iota_\infty,\iota_\ell)$ as discussed above, then there is a
nondegenerate symmetric (resp.  alternate) ${\rm
Gal}(\overline{\Q}/\Q)$-equivariant pairing $\rho_\pi \otimes \rho_\pi
\rightarrow \omega_\ell^{-w(\pi)}$ (see also~\cite{taibi}
when (c) is not satisfied).  As we shall see in \S\ref{parorthosymp}, the
definition above also fits with the Arthur-Langlands classification of
self-dual cuspidal automorphic representations of $\GL_n$.  This means that
if $\pi$ satisfy (a), (b) and (c'), then $\pi$ is symplectic (resp. 
orthogonal) if and only if the self-dual representation $\pi \otimes
|\cdot|^{w(\pi)/2}$ is so in the sense of Arthur.  \ps

We now come to an important notation that we shall use. Assume that $\pi$
satisfies (a), (b), (c') and is of weights $k_1>\cdots>k_n$.  We will say
that $\pi$ is {\it centered} if $k_n=0$.  Up to twisting $\pi$ if necessary,
we may focus on centered $\pi$'s.  Assume that $\pi$ is centered.  The
symmetry property $k_i+k_{n+1-i}=w(\pi)$ (for $i=1,\dots,n$) shows that
$k_1=w(\pi)$ is at the same time the biggest weight and the motivic weight
of $\pi$.  Set $r=[n/2]$ and introduce the integers $$w_i= 2 k_i - w(\pi)$$
for $i=1,\dots,r$.  Those numbers will be called the {\it Hodge weights} of
$\pi$.  Observe that $w_1=w(\pi)$, $w_i \equiv w(\pi) \bmod 2$ for each $i$, 
and that $w_1>\cdots>w_r \geq 0$.  The $n$ weights $k_i$ of $\pi$ can be recovered from
the $r$ Hodge weights $w_i$ : they are the $2r$ integers $\frac{w_1 \pm w_i}{2}$
when $n=2r$ is even, and the $2r+1$ integers $\frac{w_1 \pm w_i}{2}$ and
$w_1/2$ when $n=2r+1$ is odd.  Observe also that $w_r=0$ if and only if property
(c) is not satisfied.

\begin{definition}\label{defsooetoile} Let $r\geq 1$ be an integer. Let $w_1>\cdots>w_r$ be nonnegative integers which are all congruent mod $2$. \begin{itemize}\ps
\item[(i)] If the $w_i$ are odd, we denote by ${\rm S}(w_1,\cdots,w_r)$ the number of cuspidal automorphic representations $\pi$ of $\GL_{2r}$ satisfying {\rm (a), (b), (c)}, which are symplectic, and with Hodge weights $w_1>\cdots>w_r$. \ps
\item[(ii)] If the $w_i$ are even, we denote by ${\rm O}(w_1,\cdots,w_r)$ {\rm (}resp. ${\rm O}^\ast(w_1,\cdots,w_r)${\rm )} the number of cuspidal automorphic representations $\pi$ of $\GL_{2r}$ (resp. $\GL_{2r+1}$) satisfying {\rm (a), (b), (c')}, which are orthogonal, and with Hodge weights $w_1>\cdots>w_r$. 
\end{itemize}
\end{definition}

It follows from these definitions that if $k_1>\cdots>k_n$ are distinct
integers, if $k_n=0$, and if $w_i=2k_i-k_1$ for $i=1,\dots,[n/2]$, then 
${\rm N}(k_1,\cdots,k_n)$ coincides with : ${\rm S}(w_1,\cdots,w_{n/2})$ if $n$ is even and $k_1$ is
odd, ${\rm O}(w_1,\cdots,w_{n/2})$ if $n$ is even and
$k_1$ is even, ${\rm
O}^\ast(w_1,\cdots,w_{(n-1)/2})$ if $n$ is odd and $k_1$ is even.  \ps

\subsection{Case-by-case description, examples in low motivic weight}\label{casebycaselow}

Let us start with the symplectic cases. As already mentioned, a
standard translation ensures that for each odd integer $w\geq 1$ the number ${\rm
S}(w)$ is the dimension of the space ${\rm S}_{w+1}(\SL_2(\Z))$ of cusp
forms of weight $w+1$ for the full modular group ${\rm SL}_2(\Z)$.  We
therefore have the well-known formula \begin{equation}\label{dimsw} {\rm
S}(w)=\dim {\rm S}_{w+1}(\SL_2(\Z))= [\frac{w+1}{12}]-\delta_{w \equiv 1
\bmod 12}\cdot \delta_{w>1}\end{equation} where $\delta_P$ is $1$ if
property $P$ holds and $0$ otherwise.  The Sato-Tate group of each $\pi$ of
$\GL_2$ satisfying (a), (b) and (c) is necessarily the compact group ${\rm
SU}(2)$.

\ps

The next symplectic case is to give ${\rm S}(w,v)$ for $w>v$ odd positive
integers.  This case, which is no doubt well known to the experts, may be
deduced from Arthur's results~\cite{arthur} and a computation by R. 
Tsushima~\cite{tsushima}.  Let ${\rm S}_{(w,v)}({\rm Sp}_4(\Z))$ be the
space of vector-valued Siegel modular forms of genus $2$ for the coefficient
systems ${\rm Sym}^j \otimes \det^k$ where $j=v-1$ and $k=\frac{w-v}{2}+2$
(we follow the conventions in~\cite[\S 25]{vdg}).  Using the geometry of the
Siegel threefold, Tsushima was able to give an explicit formula for $\dim
{\rm S}_{(w,v)}({\rm Sp}_4(\Z))$ in terms of $(w,v)$. This formula is
already too big to give it here, but see loc.  cit.  Thm.  4. There
is a much simpler closed formula for the Poincar\'e series of the ${\rm
S}(w,1)$ due to Igusa : see~\cite[\S 9]{vdg}. An examination of Arthur's
results~\cite{arthur} for the Chevalley group ${\rm SO}_{3,2}={\rm PGSp}_4$
over $\Z$ shows then that \begin{equation}\label{dimswv}{\rm S}(w,v)=\dim
{\rm S}_{(w,v)}({\rm Sp}_4(\Z)) - \delta_{v=1}\cdot \delta_{w \equiv 1 \bmod
4} \cdot {\rm S}(w).\end{equation} The term which is subtracted is actually
the dimension of the Saito-Kurokawa subspace of ${\rm S}_{(w,v)}({\rm
Sp}_4(\Z))$.  That this is the only term to subtract is explained by
Arthur's multiplicity formula (see~\S\ref{parsp4}).  We refer to
Table~\ref{tableSwv} for the first nonzero values of ${\rm S}(w,v)$.  It
follows for instance that for $w\leq 23$, then ${\rm S}(w,v)=0$ unless
$(w,v)$ is in the following list : $$(19,7), (21,5), (21,9), (21,13),
(23,7), (23,9), (23,13)$$ In all those cases ${\rm S}(w,v)=1$.  Moreover the
first $w$ such that ${\rm S}(w,1)\neq 0$ is $w=37$.  The Sato-Tate group of
each symplectic $\pi$ of $\GL_4$ satisfying (a), (b) and (c) is either the
compact connected Lie group ${\rm Spin}(5) \subset {\rm GL}_4(\C)$ or ${\rm SU}(2) \subset {\rm GL}_4(\C)$ (symmetric cube
of the standard representation).  This latter case
should only occur when $w=3v$ (and ${\rm S}(v)$ times!), so at least for
$w=33$, and thus should not occur in the examples above. This special case
of Langlands functoriality is actually a theorem of Kim and Shahidi~\cite{kimsha}. 
 \ps

Before going further, let us mention that as Bergstr\"om and Faber pointed
out to us, although Tsushima's explicit formula is expected to hold for all
$(w,v) \neq (3,1)$ it is only proved in \cite{tsushima} for $w-v \geq 6$. 
That it holds as well in the remaining cases $w-v=2,4$ (and $v>1$) has actually been
recently proved by Ta\"ibi (at least under assumption $\ast$, see his forthcoming work),
and we shall assume it here to simplify the discussion. This actually would
not matter for the numerical applications that we will discuss in this
introduction. Indeed, the method that we shall describe leads to an
independent upper bound on ${\rm S}(w,v)$ showing that ${\rm S}(w,v)=0$
whenever $w-v \leq 4$ and $w \leq 27$, as predicted by Tsushima's formula\footnote{
When $w-v=2$ (resp. $w-v=4$) this follows from case (vi) of~\S\ref{endocase42}
(resp. case (ii) of~\S\ref{stcasso8}).}. \ps

	Our first serious contribution is the computation of ${\rm S}(w_1,\dots,w_r)$ for $r=3$ and $4$ (and any $w_i$). Our strategy is to compute first the dimension of the spaces of level $1$ automorphic forms for two
certain special orthogonal $\Z$-group  schemes 
${\rm SO}_7$ and ${\rm SO}_9$ which are reductive over $\Z$. These groups are the special orthogonal groups of the root lattice 
${\rm E}_7$ and  ${\rm E}_8 \oplus {\rm A}_1$ respectively. 
They have compact real points ${\rm SO}_n(\R)$ for $n=7$ and $9$ and
both have class number $1$, so that we are reduced to determine the dimension 
of the invariants of their integral points, namely the positive Weyl group of
$E_7$ 
and the Weyl group of $E_8$, in any given finite dimensional irreducible representation 
of the corresponding ${\rm SO}_n(\R)$. We will say more about this computation in~\S\ref{discussionformuleso7} and in Chapter \ref{polinv}. The second important 
step is to rule out all the endoscopic or non-tempered contributions predicted by Arthur's theory for those groups to get the exact values of ${\rm S}(-)$. This is done case-by-case by using the explicit form of Arthur's multiplicity formula that we expect. In the cases of ${\rm SO}_7$ (resp. ${\rm SO}_9$) there are for instance $9$ (resp. $16$) multiplicity formulas to determine. They require in particular the computations of ${\rm S}(-)$, ${\rm O}(-)$ and 
${\rm O}^\ast(-)$ for smaller ranks $n$ :  we refer to Chapters~\ref{chapSO7} and~\ref{chapSO9} for the complete study.\ps

We refer to Tables~\ref{tableSwvu} and~\ref{tableSwvut} for the first non zero values of ${\rm S}(w_1,\cdots,w_r)$ for $r=3,4$, and to the url~\cite{homepage} for much more data.
Here is a small sample of our results.

\begin{corcondb}\label{corintroso7}  \begin{itemize} \item[(i)] ${\rm S}(w_1,w_2,w_3)$ vanishes for $w_1<23$.\ps
\item[(ii)] There are exactly $7$ triples $(w_1,w_2,w_3)$ with $w_1=23$ such that ${\rm S}(w_1,w_2,w_3)$ is nonzero : 
$$(23,13,5), \, \, (23,15,3), \,\, (23,15,7), \, \, (23,17,5),\, \, (23,17,9),\, \, (23,19,3), \,\, (23,19,11),$$
and for all of them ${\rm S}(w_1,w_2,w_3)=1$. \ps
\end{itemize}
\end{corcondb}

\begin{corcondb}\label{corintroso9} \begin{itemize}\item[(i)] ${\rm S}(w_1,w_2,w_3,w_4)$ vanishes for $w_1<25$.\ps
\item[(ii)] There are exactly $33$ triples $(w_2,w_3,w_4)$ such that ${\rm S}(25,w_2,w_3,w_4)\neq 0$
and for all of them ${\rm S}(25,w_2,w_3,w_4)=1$, except ${\rm S}(25,21,15,7)={\rm S}(25,23,11,5)=2$ and ${\rm S}(25,23,15,5)=3$.
\end{itemize}
\end{corcondb}

The conjectural Sato-Tate group of a symplectic $\pi$ of $\GL_6$ (resp.
$\GL_8$) satisfying (a), (b) and (c) is either ${\rm SU}(2)$, ${\rm
SU}(2) \times {\rm SO}(3)$ (resp. ${\rm SU}(2)^3/\{\pm 1\}$), or the 
compact connected, simply connected, Lie group of type ${\rm C}_3$  (resp. ${\rm C}_4$). 
This latter case should occur for each of the $7$ automorphic representations of Cor.~\ref{corintroso7}  (ii) 
(resp. of the $37$ automorphic representations of Cor.~\ref{corintroso9}
(ii)). 
\ps \ps

Let us discuss now the orthogonal case. We start with two general useful
facts.  Although the first one is quite simple, the second one is rather
deep and relies on Arthur's proof that the root number of an orthogonal
$\pi$ is always $1$ (see Proposition~\ref{paramalg} and \S\ref{paragraphepsilon}).

\begin{prop}\label{vanishorodd} If $r$ is odd, then ${\rm O}(w_1,\dots,w_r)=0$ for all $w_1 > w_2 > \cdots > w_r$. 
\end{prop}

\begin{propconda}\label{propepsilonintro} 
If $\frac{1}{2}(\sum_{i=1}^r w_i) \not\equiv [\frac{r+1}{2}] \bmod 2$, then
$${\rm O}(w_1,\dots,w_r)={\rm
O}^\ast(w_1,\dots,w_r)=0.$$ \end{propconda}

Each time we shall write ${\rm O}(w_1,\dots,w_r)$ and ${\rm O}^\ast(w_1,\dots,w_r)$ we shall thus assume from now on that 
$\frac{1}{2}(\sum_i w_i) \equiv [\frac{r+1}{2}] \bmod 2$. Here are those numbers for $r\leq
2$.
\begin{thmcondb} \label{proporthointro}
\begin{itemize}\item[(i)] ${\rm O}^\ast(w)={\rm S}(\frac{w}{2})$, \ps
\item[(ii)] ${\rm O}(w,v)={\rm S}(\frac{w+v}{2})\cdot {\rm S}(\frac{w-v}{2})$ if $v \neq 0$, and ${\rm O}(w,0)= \frac{{\rm S}(w/2)\cdot ({\rm S}(w/2)-1)}{2}$, \ps
\item[(iii)] ${\rm O}^\ast(w,v)={\rm S}(\frac{w+v}{2},\frac{w-v}{2})$.
\end{itemize}
\end{thmcondb}

Part (i) and (ii) actually only rely on assumption $\ast$. These identities correspond to some simple cases of Langlands functoriality 
related to the exceptional isogenies ${\rm SL_2(\C) \rightarrow 
{\rm SO}_3(\C)}$ ({\it symmetric square}), $\SL_2(\C) \times \SL_2(\C) \rightarrow {\rm SO}_4(\C)$ ({\it tensor product}) and ${\rm Sp}_4(\C) 
\rightarrow {\rm SO}_5(\C)$ ({\it reduced exterior square}). As we shall show, they are all consequences of Arthur's work : see Chapter~\ref{chapsmall}. Another tool in our 
proof is a general, elementary, lifting result for isogenies between Chevalley groups over $\Z$. From the point of view of Langlands conjectures, it asserts that the 
Langlands group of $\Z$, a compact connected topological group, is simply connected (see Appendix~\ref{parsatotate}). 

\ps

Our main remaining contribution in the orthogonal case is thus the assertion in Thm.~\ref{mainthmintro}
about ${\rm O}(w_1,w_2,w_3,w_4)$. We argue as before by considering this time the special orthogonal 
group ${\rm SO}_8$ over $\Z$ of the root lattice ${\rm E}_8$. It also has class number $1$ as ${\rm E}_8$ is the unique even unimodular lattice in rank $8$. For some reasons related to twisted endoscopy between ${\rm SO}_8$ and ${\rm Sp}_6$, the precise numbers that we compute are the ${\rm O}(w_1,w_2,w_3,w_4)$ when $w_4\neq0$, as well as the numbers
 $$2 \cdot {\rm O}(w_1,w_2,w_3,0)+  {\rm O}^\ast(w_1,w_2,w_3).$$
When this latter number is $\leq 1$, it  is thus necessarily equal to ${\rm O}^\ast(w_1,w_2,w_3)$, which leads first to the following partial results for the orthogonal representations $\pi$ of in dimension $n=7$. See Table~\ref{tableOwvutbis} and~\cite{homepage} for more results.

\begin{corcondb}\label{corintrosp6}\begin{itemize}\item[(i)] ${\rm O}^\ast(w_1,w_2,w_3)$ vanishes for $w_1<24$. \ps \item[(ii)] There are exactly $8$ triples $(w_1,w_2,w_3)$ with $w_1\leq 26$ such that 
${\rm O}^\ast(w_1,w_2,w_3) \neq 0$, namely 
$$(24,16,8), (26,16,10), (26,20,6), (26,20,10), (26,20,14), $$ $$(26,24,10), (26,24,14), (26,24,18),$$
in which cases ${\rm O}^\ast(w_1,w_2,w_3)=1$. \end{itemize}
\end{corcondb}

Observe that our approach does not allow to tackle this case directly as there is no semisimple $\Z$-group of type ${\rm C}_3$ 
with compact real points (and actually of type ${\rm C}_l$ for any $l\geq 3$). On the other hand, our results allow to compute 
in a number of cases the dimension of the space ${\rm S}_{w_1,w_2,w_3}({\rm
Sp}_6(\Z))$ of vector valued Siegel modular forms of 
whose infinitesimal character, a semisimple element in
$\mathfrak{so}_7(\C)$, has distinct eigenvalues $\pm \frac{w_1}{2},\pm
\frac{w_2}{2},\pm \frac{w_3}{2},0$, where $w_1>w_2>w_3$ are even positive integers.
Indeed, we deduce from Arthur's multiplicity formula that : \ps

\begin{propcondb} \label{formdimgenus3} $\dim {\rm S}_{w_1,w_2,w_3}({\rm Sp}_6(\Z))={\rm O}^\ast(w_1,w_2,w_3)+ {\rm O}(w_1,w_3)\cdot {\rm O}^\ast(w_2)$
$$+\delta_{w_2 \equiv 0 \bmod 4}\cdot ( \delta_{w_2=w_3+2}\cdot {\rm S}(w_2-1)\cdot {\rm O}^\ast(w_1)+ \delta_{w_1=w_2+2}\cdot {\rm S}(w_2+1)\cdot {\rm O}^\ast(w_3)).$$
 \end{propcondb}
\ps

In particular, in turns out that Corollary~\ref{corintrosp6} and
Theorem~\ref{proporthointro} allow to determine the dimension of ${\rm S}_{w_1,w_2,w_3}({\rm Sp}_6(\Z))$ when $w_1\leq 26$
(which makes $140$ cases). We refer to Chapter~\ref{finalapp} for more about this and to the website~\cite{homepage} for some results. We actually explain in this chapter how to compute {\it for any genus $g$} the dimension of the 
space of Siegel cusp forms for ${\rm Sp}_{2g}(\Z)$ of any given regular infinitesimal character in terms of various numbers ${\rm S}(-)$, ${\rm O}(-)$ and ${\rm O}^\ast(-)$. \ps\ps

The problem of the determination of $\dim {\rm S}_{w_1,w_2,w_3}({\rm
Sp}_6(\Z))$ has been solved by Tsuyumine in~\cite{tsuyumine} when
$w_1-w_3=4$ (scalar valued Siegel modular forms of weight
$k=\frac{1}{2}(w_1+2)$).  As already said when we discussed
Theorem~\ref{mainthmintro}, it has also been studied recently in general by
Van der Geer, Bini, Bergstr\"om and Faber, see e.g.~\cite{BFVdG} for the
latest account of their beautiful results. In this last paper, the authors give in particular a (partly conjectural) table for certain values of 
$\dim {\rm S}_{w_1,w_2,w_3}({\rm Sp}_6(\Z))$ : see Table 1 loc. cit. We checked that this table fits our results. 
In turn, their results allow not only to determine conjecturally each ${\rm
O}^\ast(w_1,w_2,w_3)$, but ${\rm O}(w_1,w_2,w_3,0)$ as well by our work. Let us mention that those authors not only compute dimensions but also certain Hecke eigenvalues. 
\ps\ps

The Sato-Tate group of an orthogonal $\pi$ of $\GL_7$ satisfying (a), (b)
and (c) is either  ${\rm SO}(3)$, ${\rm SO}(7)$ or the compact
Lie group of type ${\rm G}_2$.  In order to enumerate the conjectural number of
$\pi$ having this latter group as Sato-Tate group there is a funny game we
can play with the reductive group ${\rm G}_2$ over $\Z$ such that ${\rm
G}_2(\R)$ is compact, namely the automorphism group scheme of the Coxeter
octonions.  We compute in Chapter \ref{chapG2} the dimension of the spaces
of level $1$ automorphic forms for this $\Z$-group ${\rm G}_2$.  Assuming
Langlands and Arthur's conjectures for the embedding of dual groups ${\rm G}_2(\C) \rightarrow {\rm
SO}_7(\C)$ we are able to compute the conjectural
number $${\mathrm G}_2(v,u)$$ of orthogonal $\pi$ of $\GL_7$ with Hodge weights
$v+u>v>u$ and whose Sato-Tate group is either ${\rm SO}(3)$ or ${\rm G}_2$ as an
explicit function of $v$ and $u$.  Of course it is easy to rule out from
this number the contribution of the ${\rm SO}(3)$ case, which occurs 
${\rm S}(u/2)$ times when $v=2u$ (sixth symmetric power of the standard representation). We obtain in particular a conjectural
minoration of ${\rm O}^\ast(v+u,v,u)$ in general which matches beautifully
the results of corollary (ii) above : the first three $\pi$'s should have ${\rm G}_2$ as
Sato-Tate groups, and the five others ${\rm SO}(7)$.  We refer to
Table~\ref{tableStG2} for a sample of results.  This also confirms certain
similar predictions in~\cite{BFVdG}. There is actually another way, still conjectural
but perhaps accessible nowadays, to think about ${\rm G}_2(v,u)$, using twisted endoscopy for a triality automorphism for ${\rm
PGSO}_8$. This concerns triples of weights $(w,v,u)$ with $w=v+u$. With this
theory in mind, it should follow that whenever ${\rm O}(v+u,v,u)=1$ the
unique orthogonal $\pi$ of $\GL_7$ with Hodge weights $v+u>v>u$ should have
${\rm G}_2$ or ${\rm SO}(3)$ as Sato-Tate group. This criterion applies
to the first three $\pi$'s given by the corollary (ii) (the Sato-Tate group
${\rm SO}(3)$ being obviously excluded in these cases) and thus comforts the
previous predictions.  \ps \ps

Modular forms of level one for the Chevalley group  of type ${\rm G}_2$, and
whose archimedean component is a quaternionic discrete series, have been
studied by Gan, Gross and Savin in~\cite{gangrosssavin}.  They define a
notion of Fourier coefficients for those modular forms and give interesting
examples of Eisenstein series and of two exceptional theta series coming
from the modular forms of level $1$ and trivial coefficient of the
anisotropic form of ${\rm F}_4$ over $\Q$.  Table~\ref{tableStG2} suggests that
the first cusp form for this ${\rm G}_2$ whose conjectural transfert to
$\GL_7$ is cuspidal should occur for the weight $k=8$, which is the first
integer $k$ such that ${\rm G}_2(2k,2k-2)\neq 0$.  Modular forms for the
anisotropic $\Q$-form of ${\rm G}_2$ have also been studied by Gross,
Lansky, Pollack and Savin in~\cite{grosssavin}, \cite{grosspollack},
\cite{lanskypollack} and \cite{pollack}, partly in order to find
$\Q$-motives with Galois group of type ${\rm G}_2$, a problem initially
raised by Serre~\cite{serremotives}.  The automorphic forms they consider there are not of level
$1$, but of some prime level $p$ and Steinberg at this prime $p$.  \ps

Let us now give a small sample of results concerning ${\rm O}(w_1,w_2,w_3,w_4)$ for $w_4>0$, see~Table~\ref{tableOwvut} and~\cite{homepage} for more values.
	
\begin{corcondb}\label{corintroso8}  \begin{itemize} \item[(i)] ${\rm O}(w_1,w_2,w_3,w_4)$ vanishes for $w_1<24$. \ps
\item[(ii)] The $(w_1,w_2,w_3,w_4)$ with $0<w_4<w_1\leq 26$ such that ${\rm O}(w_1,w_2,w_3,w_4) \neq 0$ are 
$$ (24,18, 10, 4), (24, 20, 14, 2), (26, 18, 10, 2), (26, 18, 14, 6), (26, 20, 10, 4), (26, 20, 14, 8), $$
$$(26, 22, 10, 6), (26, 22, 14, 2),  (26, 24, 14, 4), (26, 24, 16, 2), (26, 24, 18, 8), (26, 24, 20, 6),$$
and for all of them ${\rm O}(w_1,w_2,w_3,w_4)=1$.\end{itemize}
\end{corcondb}

The Sato-Tate group of an orthogonal $\pi$ of $\GL_8$ satisfying (a), (b)
and (c) can be a priori either ${\rm SU}(2)$, ${\rm SO}(8)$, $({\rm SU}(2)
\times {\rm Spin}(5))/\{\pm 1\}$, ${\rm SU}(3)$ (adjoint representation) or
${\rm Spin}(7)$. The case ${\rm SU}(3)$ should actually never occur
(see Appendix \ref{parsatotate}). Moreover, it is
not difficult to check that the ${\rm Spin}(7)$ case may only occur for Hodge weights $(w_1,w_2,w_3,w_4)$
such that $w_4=|w_1-w_2-w_3|$, in which case it occurs exactly ${\rm
O}^\ast(v_1,v_2,v_3)-{\rm G}_2(v_2,v_3)\cdot \delta_{v_1=v_2+v_3}$ times, where
$(v_1,v_2,v_3)=(w_2+w_3,w_1-w_3,w_1-w_2)$. But for each of the six Hodge weights $(w_1,w_2,w_3,w_4)$ of Corollary~\ref{corintroso8} (ii) such that $w_4=|w_1-w_2-w_3|$, we have $v_1 \neq v_2+v_3$ (i.e. $w_1 \neq w_2+w_3$) and the number ${\rm O}^\ast(v_1,v_2,v_3)+2 \cdot{\rm O}(v_1,v_2,v_3,0)$ that we computed is $1$, so that ${\rm O}^\ast(v_1,v_2,v_3)=1$. It follows that in these six cases the Sato-Tate groups should be ${\rm Spin}(7)$, and thus in the six remaining ones it must be ${\rm SO}(8)$ (it is easy to rule out ${\rm SU}(2)$ and $({\rm SU}(2)
\times {\rm Spin}(5))/\{\pm 1\}$). \ps


\subsection{Generalisations} At the moment, we cannot compute ${\rm
N}(k_1,\cdots,k_n)$ for $n>8$ because we don't know neither the dimensions
of the spaces of vector valued Siegel modular forms for ${\rm Sp}_{2g}(\Z)$
in genus $g \geq 4$, nor the number of level $1$ automorphic representations
$\pi$ of ${\rm SO}_{p,q}$ such that $\pi_\infty$ is a discrete series when
$p+q > 10$.  We actually have in our database the dimensions of the spaces
of level $1$ automorphic forms of the special orthogonal $\Z$-group ${\rm
SO}_{15}$ of the root lattice ${\rm E}_7 \oplus {\rm E}_8$ (note that the
class number is $2$ in this case).  They lead to certain upper bound results
concerning the number of symplectic $\pi$ in dimension $n \leq 14$ that we
won't give here.  However, they contain too many unknowns to give as precise
results as the ones we have described so far for $n \leq 8$, because of the
inductive structure of the dimension formulas.  \ps

\subsection{Methods and proofs} \subsubsection{}\label{discussionformuleso7} We now discuss a bit more the methods and proofs. As already explained, a first important technical ingredient to obtain all the numbers above is to be able to compute, say given a finite subgroup $\Gamma$ of a compact connected Lie group $G$,
and given a finite dimensional irreducible representation $V$ of $G$, the dimension $$\dim V^\Gamma$$
of the subspace of vectors in $V$ which are fixed by $\Gamma$. This general problem is studied in Chapter~\ref{polinv} (which is entirely unconditional). The main result there is an explicit general formula for $\dim V^\Gamma$ as a function of the extremal weight of $V$, which is made explicit in the cases alluded above.
Our approach is to write $$\dim V^\Gamma=\frac{1}{|\Gamma|} \sum_{\gamma \in \Gamma} \chi_V(\gamma)$$ where $\chi_V : G \rightarrow \C$ is the character of $V$. The formula we use for $\chi_V$ is a degeneration of the Weyl character formula
which applies to possibly non regular elements and which was established in~\cite{chcl}. Fix a maximal torus $T$ in $G$, with character group $X$, and a set $\Phi^+ \subset X$ of positive roots for the root system $(G,T)$. Let $V_\lambda$ be the irreducible representation with highest weight $\lambda$. Then $$\dim V_\lambda^\Gamma=\sum_{j \in J} a_j\, e^{\frac{2i\pi}{N} \langle b_j,  w_j(\lambda+\rho)-\rho\rangle}\, P_j(\lambda)$$
where $N$ is the lcm of the orders of the elements of $\Gamma$, $a_j \in \Q(e^{2i\pi/N})$, $b_j$ is a certain cocharacter of $T$,  $w_j$ is a certain element in the Weyl group $W$ of $(G,T)$, 
and $P_j$ is a certain rational polynomial on $X \otimes \Q$
which is a product of at most $|\Phi^+|$ linear forms. For each $\gamma \in T$ let $W_\gamma \subset W$ be the Weyl group of the connected centralizer of $\gamma$ in $G$ with respect to $T$. 
Then $$|J|=\sum_{\gamma} |W/W_\gamma|$$
where $\gamma$ runs over a set of representatives of the $G$-conjugacy classes of elements of $\Gamma$. \ps

In practice, this formula for $\dim V_\lambda^\Gamma$ is quite insane.
Consider for instance $G={\rm SO}_7(\R)$ and $\Gamma={\rm W}^+(E_7)$ the
positive Weyl group of the root system of type $E_7$ : this is the
case we need to compute ${\rm S}(w_1,w_2,w_3)$.  Then $|J|=725$, $N=2520$
and $|\Phi^+|=9$ : it certainly impossible to explicitly write down this
formula in the present paper.  This is however nothing (in this case!) for a
computer and we refer to~\S~\ref{computerprog} for some details about the
computer program we wrote using PARI/GP.  Let us mention that we use in an
important way some tables of Carter~\cite{carter} giving the characteristic
polynomials of all the conjugacy classes of a given Weyl group in its
reflexion representation.  \ps

\subsubsection{}\label{amfintro} The second important ingredient we need is Arthur's multiplicity formula in a various number of cases. Concretely this amounts to determining a quite large collection of signs. This is discussed in details in Chapter~\ref{sectionarthur}, in which we specify Arthur's general results to the case of classical 
semisimple $\Z$-groups $G$. This leads first to a number of interesting properties of the automorphic representations $\pi$ satisfying (a), (b) and (c') of this introduction. Of course, a special attention is given to the groups $G$ with $G(\R)$ compact, hence to the integral theory of quadratic forms. We restrict our study to the representations in the discrete spectrum of $G$ which are unramified at each finite place. At the archimedean place we are led to review some properties of the packets of representations defined by Adams-Johnson in~\cite{AJ}. We explain in particular in the appendix~\ref{appendixadamsjohnson} the parameterization of the elements of these packets by the characters of the dual 
component group in the spirit of Adams paper~\cite{adams} in the discrete series case. For our purposes, we need to apply Arthur's results to a number of classical groups of small rank, namely 
$${\rm SL}_2,\,\,\, {\rm Sp}_4, \,\,\,{\rm Sp}_6, \,\,\,{\rm SO}_{2,2}, \,\,\,{\rm SO}_{3,2}, \,\,\,{\rm SO}_7, \,\,\,{\rm SO}_8\,\,\, {\rm and } \,\,\,{\rm SO}_9.$$
When $G(\R)$ is compact, Arthur's multiplicity formula takes a beautifully simple form, in which the half-sum of the positive roots on the dual side plays an important role. 
Let $G$ be any semisimple $\Z$-group such that $G(\R)$ is compact. We do not assume that $G$ is classical here and state the general conjectural formula. Let $\mathcal{L}_\Z$ denote the Langlands group of $\Z$ and let
$$\psi : \mathcal{L}_\Z \times {\rm SL}_2(\C) \rightarrow \widehat{G}$$ be a global Arthur parameter such that $\psi_\infty$ is an 
Adams-Johnson parameter : see~\cite{arthurunipotent} as well as the appendices~\ref{appendixadamsjohnson} and \ref{parsatotate}. Denote by $\pi_\psi$ the irreducible admissible representation of $G(\AAA)$ which is $G(\widehat{\Z})$-spherical and 
with the Satake parameters and infinitesimal character determined by $\psi$
according to Arthur's recipe\footnote{This means that $\psi$ corresponds to $\pi_\psi$ in the sense of the appendix~\ref{parsatotate}, assumption (L5). This uniquely determines $\pi_\psi$ as $G(\R)$ is compact and connected}. Denote also by $e(\psi)$ the (finite) number of $\widehat{G}$-conjugacy classes of
global Arthur parameters $\psi'$ as above such that
$\pi_{\psi'} \simeq \pi_{\psi}$ (for most $\psi$ we have $e(\psi)=1$). The multiplicity
$m(\pi_\psi)$ of $\pi_\psi$ in ${\rm L}^2(G(\Q)\backslash G(\AAA))$ 
should be given in general by 
\begin{equation}\label{formulamultdefgen} m(\pi_\psi)=\left\{\begin{array}{ll}
e(\psi) & {\rm if}\, \, \,  \rho^\vee_{|{\rm C}_\psi} = \varepsilon_\psi, \\ \\ 0 & {\rm otherwise}.\end{array}\right.\end{equation}
The group ${\rm C}_\psi$ is by definition the centralizer of ${\rm Im}\, \psi$
in $\widehat{G}$, it is a finite group. As explained in the appendix~\ref{appendixadamsjohnson}, it is always an elementary abelian $2$-group. The character $\varepsilon_\psi$ is defined by Arthur in~\cite{arthurunipotent}.  
The character  $\rho^\vee$ is defined as follows. Let
$\varphi_{\psi_\infty} : {\rm W}_\R \longrightarrow \widehat{G}$ the Langlands parameter
associated by Arthur to $\psi_\infty$. First, the centralizer in $\widehat{G}$ of $\varphi_{\psi_\infty}({\rm W}_\C)$ is a maximal torus $\widehat{T}$ of $\widehat{G}$, so that $\varphi_{\psi_\infty}(z)=z^\lambda \overline{z}^{\lambda'}$ for some $\lambda \in \frac{1}{2}{\rm X}_\ast(\widehat{T})$ and all $z \in {\rm W}_\C$,  and $\lambda$ is dominant with respect to a unique Borel subgroup $\widehat{B}$ of $\widehat{G}$ containing $\widehat{T}$. Let $\rho^\vee$ denote the half-sum of the positive 
roots  of $(\widehat{G},\widehat{B},\widehat{T})$. As $G$ is semisimple over $\Z$ and $G(\R)$ is compact, this is actually a character of $\widehat{T}$. By construction, we have  ${\rm C}_\psi \subset \widehat{T}$, and thus formula 
\eqref{formulamultdefgen} makes sense. The second important statement is that any automorphic representation of $G$ which is $G(\widehat{\Z})$-spherical has the form $\pi_\psi$ for
some $\psi$ as above. 
\ps 

\subsection{Application to Borcherds even lattices of rank $25$ and
determinant $2$}\label{applborcherdslat} We end this introduction by discussing two other applications. The first one
is very much in the spirit of the work of the first author and Lannes
\cite{cl}.  It concerns the genus of euclidean lattices $L \subset
\R^{25}$ of covolume $\sqrt{2}$ which are {\it even}, in the sense that $x
\cdot x \in 2\Z$ for each $x \in L$.  A famous computation by Borcherds
in~\cite{borcherdsthese} asserts that there are up to isometry exactly $121$
such lattices.  It follows that there are exactly $121$ level $1$ automorphic
representations of the special orthogonal group ${\rm SO}_{25}$ over $\Z$ of
the root lattice ${\rm E}_8^3 \oplus {\rm A}_1$, for the trivial coefficient. The dual group of ${\rm SO}_{25}$ is ${\rm
Sp}_{24}(\C)$.  \ps

Observe now our tables : we have found exactly $23$ cuspidal automorphic representations 
$\pi$ of $\GL_n$ (for any $n$) satisfying conditions (a), (b) and (c) above,
centered and with motivic weight $\leq 23$, namely : \ps

\begin{itemize}
\item[(a)] the trivial representation of ${\rm GL}_1$, \ps
\item[(b)] $7$ representations of ${\rm GL}_2$, \ps
\item[(c)] $7$ symplectic representations of ${\rm GL}_4$,\ps 
\item[(d)] $7$ symplectic representations of ${\rm GL}_6$,\ps
\item[(e)] The orthogonal representation of ${\rm
GL}_3$ symmetric square of the representation of ${\rm GL}_2$ of
motivic weight $11$ associated to a generator of ${\rm S}_{12}({\rm SL}_2(\Z))$.\ps
\end{itemize}

We refer to \S~\ref{arthurparam} for the notion of global Arthur parameters
for a $\Z$-group such as ${\rm SO}_{25}$. This is the non-conjectural
substitute used by Arthur for the conjugacy classes of morphisms $\psi : \mathcal{L}_\Z \times {\rm SL}_2(\C) \rightarrow
{\rm Sp}_{24}(\C)$ with finite centralizers that we just discussed
in~\S\ref{amfintro} (see Appendix~\ref{parsatotate}). We have now this first crazy coincidence, which is easy to check with a computer. 

\begin{prop}\label{121param} There are exactly $121$ global Arthur parameters for ${\rm SO}_{25}$ which have trivial infinitesimal character 
that one can form using only those $23$ cuspidal automorphic representations. 
\end{prop}

See the table of Appendix \ref{table121} for a list of these parameters,
using notations of~\S\ref{arthurparam}. One
also uses the following notation~: if ${\rm S}(w_1,\dots,w_r)=1$ we denote by
$\Delta_{w_1,\cdots,w_r}$ the twist by $|\cdot|^{w_1/2}$ of the unique
centered $\pi \in \Pi({\rm GL}_{2r})$ satisfying (i) to (iii) and with Hodge
weights $w_1,\cdots,w_r$.  When ${\rm S}(w_1,\dots,w_r)=k>1$ we denote by
$\Delta_{w_1,\dots,w_r}^k$ any of the $k$ representations of ${\rm GL}_{2r}$
with this latter properties.  \ps

The second miracle is that for each of the $121$ parameters $\psi$ that we
found in Proposition~\ref{121param}, the
unique level $1$ automorphic representation $\pi_\psi$ of ${\rm SO}_{25}(\AAA)$ in the
packet $\Pi(\psi)$ (see Def.~\ref{defpipsi}) has indeed a nonzero
multiplicity, that is multiplicity $1$.  In other words, we have the
following theorem.

\begin{thmcondb}\label{thmborcherds} The $121$ level $1$ automorphic
representations of ${\rm SO}_{25}$ with trivial coefficient are the ones
given in Appendix~\ref{table121}.  \end{thmcondb}

The $24$ level $1$ automorphic representations of ${\rm O}_{24}$ with
trivial coefficient ("associated" to the $24$ Niemeier lattices) and the
$32$ level $1$ automorphic representations of ${\rm SO}_{23}$ with trivial
coefficient (associated to the $32$ even lattices of rank $23$ of covolume
$\sqrt{2}$) had been determined in~\cite{cl}.  As in {\it loc.  cit.},
observe that given the shape of Arthur's multiplicity formula, the naive
probability that Theorem~\ref{thmborcherds} be true was close to $0$ (about
$2^{-450}$ here in we take in account the size of ${\rm C}_\psi$ for each
$\psi$), so something quite mysterious seems to occur for these small
dimensions and trivial infinitesimal character.  The miracle in all these
cases is that whenever we can write down some $\psi$, then $\rho^\vee_{|{\rm
C}_\psi}$ is always equal to $\varepsilon_\psi$.  \ps

It is convenient for us to include here the proof of this theorem, although
it uses freely the notations of Chapter~\ref{sectionarthur}. 

\begin{pf} To check that each parameter has multiplicity one, we apply for instance the following
simple claim already observed in~\cite{cl2}. 
Let $\psi=(k,(n_i),(d_i),(\pi_i))$ be a global Arthur parameter for 
${\rm SO}_{8m \pm 1}$ with trivial infinitesimal character. Assume there exists an integer  $1\leq i \leq k$ such that $\pi_i=1$ 
and $\pi_j$ is symplectic if $j \neq i$. Then the unique $\pi \in \Pi(\psi)$ has a nonzero multiplicity 
(hence multiplicity $1$) if and only if for each $j \neq i$ one has 
either $\varepsilon(\pi_j)=1$ or $d_j<d_i$. Indeed, when the infinitesimal
character of $\psi$ is trivial the formula for
$\rho^\vee(s_j)$ given in~\S\ref{mformoddorth} shows that
$\rho^\vee(s_j)=\varepsilon(\pi_j)$ for each $j \neq i$. The claim follows as
$\varepsilon_\psi(s_j)=\varepsilon(\pi_j)^{{\rm Min}(d_i,d_j)}$ by
definition (see~\S\ref{paragraphepsilon}), as $d_i$ is even but $d_j$ is odd.  \ps
 Among the $21$ symplectic $\pi$'s
above of motivic weight $\leq 23$, 
one observes that  exactly $4$ of them have epsilon factor $-1$, namely
$$\Delta_{17},\Delta_{21},\Delta_{23,9}\,\,  {\rm and}\, \, \Delta_{23,13}.$$
A case-by-case check at the list concludes the proof thanks to this claim except for the parameter $${\rm Sym}^2 \Delta_{11} [2] \oplus \Delta_{11}[9]$$
which is the unique parameter which is not of the form above. But it is clear that for such a $\psi$ one has $\varepsilon_\psi=1$ and one observes that $\rho^\vee_{|{\rm C}_\psi}=1$ as well in this case, which concludes the proof 
(see~\S\ref{mformoddorth}).
\end{pf}

\subsection{A level $1$, non-cuspidal, tempered automorphic representation of
$\GL_{28}$ over $\Q$ with weights $0,1,2,\cdots,27$}\label{lfunctionmyrifique}

 None of the $121$ automorphic representations of ${\rm SO}_{25}$ discussed above is tempered. This is clear since none of the $21$ symplectic $\pi$'s above admit the Hodge weight $1$. Two representations in the list are not too far from being tempered however, namely the ones whose Arthur parameters are
$$\Delta_{23, 13, 5}\oplus\Delta_{21,9}\oplus\Delta_{19,7}\oplus\Delta_{17}\oplus\Delta_{15}\oplus\Delta_{11}\oplus[4],$$
$$\Delta_{23, 15, 3}\oplus\Delta_{21,5}\oplus\Delta_{19,7}\oplus\Delta_{17}\oplus\Delta_{11}[3]\oplus[2].$$
\ps
 \noindent It is thus tempting to consider the following problem : for which integers $n$ can we find \begin{itemize}\ps
\item[-] a partition $n=\sum_{i=1}^r n_i$ in integers $n_i \geq 1$, \ps
\item[-] for each $i=1,\cdots,r$, a cuspidal automorphic
representation $\pi_i$ of ${\rm GL}_{n_i}$ satisfying 
assumptions (a), (b) and (c), and of motivic weight $n-1$, \ps
\end{itemize}
such that the parabolically induced representation
$$\pi=\pi_1 \oplus \pi_2 \oplus \cdots \oplus \pi_r$$ of ${\rm
GL}_n(\AAA)$ has the property that the eigenvalues of the infinitesimal character of $\pi_\infty$ are all the integers between $0$ and $n-1$ ? 
\ps

By assumption, the $\pi_i$ are non necessarily centered but share the same motivic weight
$n-1$, so that $\pi$ is esentially tempered. It follows that the ${\rm L}$-function $${\rm L}(\pi,s)=\prod_i {\rm L}(\pi_i,s)$$ of such a
$\pi$ shares much of the analytic properties of the ${\rm L}$-function of
a cuspidal $\pi'$ of ${\rm GL}_n$ satisfying (a), (b) and (c) and with weights $n-1,n-2,\dots,1,0$ : 
they have the same archimedean factors and both satisfy Ramanujan's conjecture. In particular, it
seems that the methods of~\cite{fermigier}, hence his results in~\S 9 {\it loc.  cit.},
apply to these more general ${\rm L}$-functions. They say that such an ${\rm L}$-function
(hence such a $\pi$) does not exist if $1<n<23$, and even if $n=24$ if one assumes the Riemann hypothesis. 
This is fortunately compatible with our previous result !  \ps

On the other hand, our tables allow to show that the above problem has a positive answer for
$n=28$, which leads to a very interesting ${\rm L}$-function in this dimension.  \ps

\begin{thmcondb}\label{thminftriv} There is a non cuspidal automorphic representation of ${\rm GL}_{28}$ over $\Q$ which 
satisfies {\rm (a), (b)} and {\rm (c)}, and whose weights are all the integers between $0$ and $27$, namely the
twist by $|\cdot|^{\frac{1-n}{2}}$ of 
$$\Delta_{27,23,9,1}\oplus \Delta_{25,13,3} \oplus \Delta_{21,5} \oplus \Delta_{19,7} \oplus \Delta_{17} \oplus \Delta_{15} \oplus \Delta_{11}.$$
\end{thmcondb}

It simply follows from the observation that ${\rm S}(27,23,9,1)={\rm
S}(25,13,3)={\rm S}(21,5)={\rm S}(19,7)=1$.  Actually, it is remarkable that
our whole tables only allow to find a single representation with these
properties, and none in rank $1 < n < 28$. It seems quite reasonable to conjecture that this is indeed the only
one in rank $28$ and that there are none in rank $1<n <28$. From the example above,
one easily deduces examples for any even $n\geq 28$. On the other hand,
the first odd $n>1$ for which our tables allow to find a suitable $\pi$ is
$n=31$ (in which case there are several).  \ps \medskip 
\bigskip
\bigskip
{{\sc \bf Acknowledgements}:  We thank Jean Lannes for a number of stimulating discussions about lattices. This paper actually owes a lot to the work~\cite{cl2} in which the authors compute in particular 
the Kneser $p$-neighborhood graph of the Niemeier lattices for each prime $p$. Certain observations here concerning the 
shape of Arthur's multiplicity formula, e.g. formula~\eqref{formulamultdefgen}, were partially realized there, as well as the feeling that it would be possible to obtain Theorem~\ref{thmborcherds}. The first author is also grateful to Colette Moeglin and Olivier Ta\"ibi for useful discussions about Arthur's results. We also thank Xavier Caruso for 
his help in creating tables, as well as Wee Teck Gan, Benedict Gross, Guy Henniart and Jean-Pierre Serre for their remarks. In all our computations, we heavily used the PARI/GP computer system. Last but not least, it is a pleasure to thank J. Arthur for his beautiful conjectures and 
monumental work~\cite{arthur}, without which this paper would not have
existed.}
\newpage 

\tableofcontents

%
%

\newpage 
\section{Polynomial invariants of finite subgroups of compact connected Lie groups}\label{polinv}

\subsection{The setting}\label{settinginv}  Let $G$ be a compact connected Lie group and consider $$\Gamma \subset G$$
a finite subgroup. Let $V$ be a finite dimensional complex continuous representation of $G$. The general problem addressed in this chapter is to compute the dimension
		$$\dim V^\Gamma$$
of the subspace $V^\Gamma=\{ v \in V, \gamma(v)=v \, \, \,\forall \gamma\, \in \Gamma\}$ of $\Gamma$-invariants in $V$. Equivalently,
\begin{equation}\label{scalproddim} \dim V^\Gamma = \frac{1}{|\Gamma|} \sum_{\gamma \in \Gamma} \chi_V(\gamma)\end{equation}
where $\chi_V : G \rightarrow \C$ is the character of $V$. One may of course reduce to the case where $V$ is irreducible and 
we shall most of the time do so. In order to apply
formula~\eqref{scalproddim} it is enough to know : \ps\ps

(a) The value of the character $\chi_V$ on each conjugacy class in $G$, \ps\ps
(b) For each $\gamma \in \Gamma$, a representative of the conjugacy class $c(\gamma)$ of $\gamma$ in $G$. \ps\ps

\noindent Of course, $c(\gamma)$ only depends on the conjugacy class of $\Gamma$, but the induced map $c: {\rm Conj}(\Gamma)  \rightarrow {\rm Conj}(G)$ needs not to be injective in general. Here 
 ${\rm Conj}(H)$ denotes the set of conjugacy classes of the group $H$. \ps
 
 We will be especially interested in cases where $\Gamma \subset G$ are fixed, but with $V$ varying over all the possible irreducible representations of $G$. 
 With this in mind, observe that problem (b) has to be solved once, but problem (a) for infinitely many $V$ whenever $G \neq \{1\}$. \ps

Consider for instance the group $\Gamma \subset {\rm SO}_3(\R)$ of positive isometries of a given regular tetrahedron in the euclidean $\R^3$ with center $0$. Each numbering of the vertices of the tetrahedron defines an isomorphism 
$$\Gamma \simeq \mathfrak{A}_4$$ and we fix one. For each odd integer $n\geq 1$ denote by $V_n$ the $n$-dimensional irreducible representations of ${\rm SO}_3(\R)$. This representation $V_n$ is well-known to be unique up to isomorphism,
and if $g_\theta \in {\rm SO}_3(\R)$ is a non-trivial rotation with angle $\theta$ then $$\chi_{V_n}(g_\theta)=\frac{\sin(n \frac{\theta}{2})}{\sin(\frac{\theta}{2})}.$$ 
The group $\Gamma$ has $4$ conjugacy classes, with representatives $1,(12)(34),(123),(132)$ and respective orders $1,3,4,4$. These representatives act on $\R^3$ as rotations with respective angles $0,\pi,2\pi/3,2\pi/3$. Observe that $(123)$ and $(132)$ are conjugate in ${\rm SO}_3(\R)$ but not in $\Gamma$. Formula~\eqref{scalproddim} thus writes 
$$\dim V_n^\Gamma = \frac{1}{12}( n + 3 \, \frac{\sin(n\pi/2)}{\sin(\pi/2)}+ 8\, \frac{\sin(n\pi/3)}{\sin(\pi/3)})=\left\{\begin{array}{lll}\lceil \frac{n}{12} \rceil & {\rm if}\, \, \, \, n \equiv 1,7,9 \bmod 12, \\ \\ \lfloor \frac{n}{12}\rfloor & {\rm if} \, \, \, \, n \equiv 3,5,11 \bmod 12.\end{array}\right.$$
This formula is quite simple but already possesses some features of the general case.  

\subsection{The degenerate Weyl character formula}\label{degwform} A fundamental ingredient for the above approach is a 
formula for the character $\chi_V(g)$ where $V$ is any irreducible representations of $G$ and $g \in G$ is any element as well. When 
$g$ is either central or regular, such a formula is given by Weyl's dimension formula and Weyl's character formula 
respectively. These formulas have been extended by Kostant to the more general case where the centralizer of $g$ is a Levi subgroup of $G$, 
and by the first author and Clozel in general in~\cite[Prop. 1.9]{chcl}. Let us now recall this last result. \ps

We fix once and for all a maximal torus $T \subset G$ and denote by $$X=X^*(T)={\rm Hom}(T,\mathbb{S}^1)$$ the character group 
of $T$. We denote by $\Phi=\Phi(G,T) \subset X \otimes \R$ the root system of $(G,T)$ and $W=W(G,T)$ its Weyl group. We choose 
$\Phi^+ \subset \Phi$ a system of positive roots, say with base $\Delta$, and we fix as well a $W$-invariant scalar product $\left(,\right)$ on $X \otimes \R$. 
Recall that a dominant weight is an element $\lambda \in X$ such that $\left( \lambda, \alpha \right) \geq 0$ for all $\alpha \in \Delta$. The Cartan-Weyl theory defines a canonical bijection $$\lambda \mapsto V_\lambda$$ between
the dominant weights and the irreducible representations of $G$. The representation $V_\lambda$ is uniquely characterized by the following property. If $V$ is a representation of $G$, denote by $P(V) \subset X$ the subset of $\mu \in X$ appearing in $V_{|T}$. 
If we consider the partial ordering on $X$ defined by 
$\lambda \leq \lambda'$ if and only if $\lambda'-\lambda$ is a finite sum of elements of $\Delta$, then $\lambda$ is the maximal element of $P(V_\lambda)$. One says that $\lambda$ is the highest weight of 
$V_\lambda$.\ps

Let us fix some dominant weight $\lambda \in X$. Recall that the inclusion $T \subset G$ induces a bijection $$W\backslash T \isomo {\rm Conj}(G),$$ 
it is thus enough to determine $\chi_{V_\lambda}(t)$ for any $t \in T$. Fix some $t \in T$ and denote by $$M=C_G(t)^0$$ the neutral component of the centralizer of $t$ in $G$. Of course, $t \in T \subset M$ and $T$ is maximal torus of $M$. 
Set $\Phi_M^+=\Phi(M,T) \cap \Phi^+$ and
consider the set $$W^M=\{ w \in W, w^{-1} \Phi_M^+ \subset \Phi^+\}.$$
Let $\rho$ and $\rho_M \in X \otimes \R$ denote respectively the half-sum of the elements of $\Phi^+$ and of
$\Phi^+_M$. If $w \in W^M$, we set $\lambda_w=w(\lambda+\rho)-\rho_M \in X \otimes
\R$. Observe that $$2 \frac{\left( \alpha, \lambda_w \right)}{\left( \alpha,\alpha \right)} \in \N, \, \, \,\, \, \, \, 
\forall\, \, \alpha \in \Phi_M^+.$$ It follows that $\lambda_w$ is 
a dominant weight for some finite covering of $M$, that we may choose to be the smallest finite covering 
$\widetilde{M} \rightarrow M$ for which $\rho-\rho_M$ becomes a character.
This is possible as $2 \frac{\left( \alpha, \rho-\rho_M \right)}{\left(
\alpha,\alpha \right)} \in \Z,\, \, \, \forall\, \, \alpha \in \Phi_M^+$. 
It follows from the Weyl dimension formula that the dimension of the irreducible representation of $\widetilde{M}$ with highest weight $\lambda_w$ is $P_M(\lambda_w)$ where we set 
$$P_M(v)=\prod_{\alpha \in \Phi^+_M} \frac{\left( \alpha, v + \rho_M\right)}{\left( \alpha,\rho_M\right)} \, \, \, \, \, \, \forall v \in X \otimes \R.$$

We need two last notations before stating the main result. We denote by $\varepsilon : W \rightarrow \{\pm 1\}$ the signature, and for $x \in X$ it will be convenient to write $t^x$ for $x(t)$. It 
is well-known that $w(\mu+\rho)-\rho \in X$ for all $w \in W$ and $\mu \in X$. \ps

\begin{prop}{\it (Degenerate Weyl character formula) } Let $\lambda \in X$ be a dominant weight, $t \in T$ and $M=C_G(t)^0$. Then
 $$\chi_{V_\lambda}(t)=\frac{\sum_{w \in W^M} \varepsilon(w)\cdot  t^{w(\lambda+\rho)-\rho} \cdot P_M(w(\lambda+\rho)-\rho_M)}{\prod_{\alpha \in \Phi^+\backslash \Phi^+_M} (1-t^{-\alpha})}.$$ 

\end{prop}

\begin{pf} This is the last formula in the proof of~\cite[Prop. 1.9]{chcl}. Note that it is unfortunately incorrectly stated in the beginning of that proof that up to replacing $G$ by a finite covering one may assume that $\rho$ and $\rho_M$ are characters. 
It is however not necessary for the proof to make any reduction on the group $G$. Indeed, we rather have to introduce the inverse image $\widetilde{T}$ of $T$ in the covering $\widetilde{M}$ defined above and argue as {\it loc. cit.} but in the Grothendieck group of 
characters of $\widetilde{T}$. The argument given there 
shows that for any element $z \in \widetilde{T}$ whose image in $T$ is $t$, we have
$$\chi_{V_\lambda}(t)=z^{\rho_M-\rho}\frac{\sum_{w \in W^M} \varepsilon(w) z^{\lambda_w} P_M(\lambda_w)}{\prod_{\alpha \in \Phi^+\backslash \Phi_M^+}(1-{t}^{-\alpha})}.$$
We conclude as $\lambda_w+\rho_M-\rho=w(\lambda+\rho)-\rho \in X$, so $z^{\rho_M-\rho} z^{\lambda_w}=t^{w(\lambda+\rho)-\rho}$. 
\end{pf}


\subsection{A computer program}\label{computerprog} We now return to the main problem discussed in~\S\ref{settinginv}. We fix a compact connected Lie group $G$ and a finite subgroup $\Gamma \subset G$. In order to 
enumerate the irreducible representations of $G$ we fix as in the previous paragraph a maximal torus $T \subset G$ and a subset $\Phi^+$ of positive roots for $(G,T)$. For each dominant weight $\lambda$ one thus has a unique irreducible
representation $V_\lambda$ with highest weight $\lambda$, hence a number 
$$\dim(V_\lambda^\Gamma)= \frac{1}{|\Gamma|} \sum_{\gamma \in \Gamma} \chi_{V_\lambda}(t(\gamma)),$$
where for each $\gamma \in \Gamma$ we define $t(\gamma)$  to be any element in $T$ which is conjugate to $\gamma$ in $G$. The last ingredient to be given for the computation is thus a list of these 
elements $t(\gamma) \in T$, which is a slightly more precise form of problem~(b) of~\S\ref{settinginv}. Recall that the elements of $T$ may be described as follows. Denote by $X^\vee={\rm Hom}(\mathbb{S}^1,T)$ the cocharacter group of $T$ and 
$\langle\,,\,\rangle : X \otimes X^\vee \rightarrow \Z$ the canonical perfect pairing. If $\mu \in X^\vee \otimes \C$, denote 
by  $e^{2i\pi \mu}$ the unique element $t \in T$ such that $$\forall \lambda \in X, \, \, \, \lambda(t)=e^{2i\pi \langle \lambda ,\mu\rangle}.$$ 
The map $\mu \mapsto e^{2i\pi\mu}$ defines an isomorphism
$(X^\vee\otimes \C)/X^\vee \isomo T$. \ps\ps

We thus wrote a computer program with the following property. It takes as input : \begin{itemize}\ps\ps

\item[(a)] The based root datum of $(G,T,\Phi^+)$, i.e. the collection $(X,\Phi,\Delta,X^\vee,\Phi^\vee,\langle,\rangle,\iota)$, where $\Phi^\vee \subset X^\vee$ is the set of coroots of $(G,T)$ and $\iota : \Phi \rightarrow \Phi^\vee$ 
is the bijection $\alpha \mapsto \alpha^\vee$. \ps\ps
\item[(b)] A finite set of pairs $(\mu_j,C_j)_{j \in J}$, 
where $\mu_j \in X^\vee \otimes \Q$ and $C_j \in \N$, 
with the property that there exists a partition $\Gamma = \coprod_{j \in J} \Gamma_j$ such that $|\Gamma_j|=C_j$ and each element 
$\gamma \in \Gamma_j$ is conjugate in $G$ to the element $e^{2i\pi \mu_j} \in T$. \ps \ps
\item[(c)] A dominant weight $\lambda \in X$.\ps
\ps\ps
\end{itemize}
It returns $\dim(V_\lambda^\Gamma)=|\Gamma|^{-1}\sum_{j \in J} C_j \, \chi_{V_\lambda}(e^{2i\pi\mu_j})$. \ps\ps

Recall that for $\alpha \in \Phi^+$ and $v \in X \otimes \R$ one has the relation 
$2\frac{\left(v,\alpha\right)}{\left(\alpha,\alpha\right)}=\langle v,\alpha^\vee\rangle$, 
thus (a), (b) and (c) contain indeed everything needed to evaluate the degenerate Weyl character formula. 
Although in theory the Weyl group $W$ of $(G,T)$ may be deduced from (a) we also take it as an input in practice. The program computes 
in particular for each $t_j=e^{2i\pi\mu_j}$ the 
root system of $M_j=C_G(t_j)^0$ and the set $W^{M_i}$. Of course it is often convenient to take $X=X^\vee=\Z^n$ with the canonical pairing. A routine in~PARI/GP may be found at the url~\cite{homepage}. \ps

\subsection{Some numerical applications} \label{appnum}We shall present in this paper four numerical applications of our computations. They concern the respective compact groups 
$$G={\rm SO}_7(\R), \, \,\, {\rm SO}_8(\R), \, \, \, {\rm SO}_9(\R), \, \, \, {\rm and} \,\,\,{\rm G}_2$$ and each time a very specific finite subgroup $\Gamma$. We postpone to~\S~\ref{parinvG2} the discussion of the case ${\rm G}_2$ and concentrate here on the first three cases. The general context is as follows.  \ps

Let $V$ be a finite dimensional vector space over $\R$ and let $R \subset V$ be a reduced root system in the sense of Bourbaki~\cite[Chap. VI \S 1]{bourbaki}. 
Let ${\rm W}(R)$ denote the Weyl group of $R$ and fix a ${\rm W}(R)$-invariant scalar product on $V$, so that $${\rm W}(R) \subset {\rm O}(V).$$
Assume that $R$ is irreducible. Then $V$ is irreducible as a representation of ${\rm W}(R)$ (\cite[Chap. VI \S 2]{bourbaki}). 
Let $\varepsilon : {\rm W}(R) \rightarrow \{\pm 1\}$ the signature of ${\rm W}(R)$, i.e. $\varepsilon(w)=\det(w)$ for each $w \in {\rm W}(R)$, and set 
 	$${\rm W}(R)^+={\rm W}(R) \cap {\rm SO}(V).$$
We are in the general situation of this chapter with $G={\rm SO}(V)$ and $\Gamma={\rm W}(R)^+$. Beware that the root system $\Phi$ of $({\rm SO}(V),T)$ is not the root system $R$ above ! We choose the standard based root datum for $({\rm SO}(V),T)$ as follows. If $l=\lfloor \frac{\dim(V)}{2}\rfloor$ we set
$X=X^\vee=\Z^l$, equipped with the canonical pairing : if $(e_i)$ denotes the canonical basis of $\Z^n$, then $\langle e_i,e_j \rangle = \delta_{i=j}$. There are two cases depending whether $\dim(V)$ is odd or even : \begin{itemize}\ps
\item[(i)] $\dim(V)=2l+1$. Then $\Phi^+=\{e_i, e_i \pm e_j, 1\leq i < j \leq n\}$, $e_i^\vee = 2 e_i$ for all $i$, and $(e_i\pm e_j)^\vee=e_i\pm e_j$ for all $i<j$. \ps
\item[(ii)] $\dim(V)=2l$. Then $\Phi^+=\{e_i \pm e_j, 1\leq i < j \leq n\}$ and $(e_i\pm e_j)^\vee=e_i\pm e_j$ for all $i<j$. \ps
\end{itemize}

The dominant weights are thus the $\lambda=(n_1,\cdots,n_l)=\sum_{i=1}^l n_i e_i \in X$ such that $n_1\geq n_2 \geq \cdots \geq n_l \geq 0$ if $\dim(V)=2l+1$,  and such that 
$n_1\geq n_2 \geq \cdots \geq n_{l-1} \geq |n_l|$ if $\dim(V)=2l$.  \ps

Consider now the input (b) for the program. Recall that at least if $\dim(V)$ is odd, the conjugacy class of any element $g \in {\rm SO}(V)$ is uniquely determined by the characteristic polynomial of $g$ acting on $V$. It turns out that for any 
reduced root system $R$, the characteristic polynomial of each conjugacy class of elements of ${\rm W}(R)$ has been determined by Carter in~\cite{carter}. We make an important use of these results, especially when 
$R$ is of type ${\rm E}_7$ and ${\rm E}_8$ for the applications here, in which case it is given in Tables 10 and 11 {\it loc. cit.} \ps

\subsubsection{Case I : $R$ is of type $E_7$} Then $-1 \in {\rm W}(R)$ and ${\rm W}(R) = {\rm W}(R)^+ \times \{\pm 1\}$, so the conjugacy classes in ${\rm W}(R)^+$ coincide with the conjugacy classes in ${\rm W}(R)$ belonging to ${\rm W}(R)^+$, i.e. with determinant $1$. 
From Table 10 {\it loc.cit.} one sees that  ${\rm W}(R)^+$ has exactly $27$ conjugacy classes $(c_j)$ and for each of them it gives its order $C_j$ and its characteristic polynomial, from which we deduce $\mu_j$ : this is the datum we need for (b). 
For each dominant weight $\lambda=(n_1,n_2,n_3) \in \Z^3$ our computer program then returns  $\dim(V_\lambda)^{{\rm W}(R)^+}$ : see Table~\ref{tableSO7nue} for a sample of results and to the url~\cite{homepage} for much more. \ps

\subsubsection{Case II : $R$ is of type $E_8$} This case presents two little differences compared to the previous one. First the characteristic polynomial of an element $g \in {\rm SO}(V)$ does only determine its ${\rm O}(V)$-conjugacy class as $\dim(V)=8$ is even. It determines its ${\rm SO}(V)$ conjugacy class if and only if $\pm 1$ is an eigenvalue of $g$. Let $C \subset {\rm W}(R)^+$ be a ${\rm W}(R)$-conjugacy class and let $P$ be its characteristic polynomial. If $\pm 1$ is a root of $P$, there is thus a unique conjugacy class in ${\rm SO}(V)$ with this characteristic polynomial. Otherwise, $C$ meets exactly two conjugacy classes in ${\rm SO}(V)$, it follows that $C=C_1 \coprod C_2$ where the $C_i$ are ${\rm W}(R)^+$-conjugacy classes permuted by any element in ${\rm W}(R)\backslash {\rm W}(R)^+$, and in particular $|C_1|=|C_2|$. It follows that the table of Carter gives input (b) as well in this case.  \ps

We refer to Table~\ref{tableSO8nue} for a sample of values of  nonzero $\dim(V_\lambda^{{\rm W}(R)^+})$ for $\lambda=(n_1,n_2,n_3,n_4)$ dominant with\footnote{One easily sees that $\dim V_\lambda^{{\rm W}(R)^+} = \dim V_{\lambda'}^{{\rm W}(R)^+}$ if $\lambda=(n_1,n_2,n_3,n_4)$ and $\lambda'=(n_1,n_2,n_3,-n_4)$. Better, the triality 
$(n_1,n_2,n_3,n_4) \mapsto (\frac{n_1+n_2+n_3+n_4}{2},\frac{n_1+n_2-n_3-n_4}{2},\frac{n_1-n_2+n_3-n_4}{2},\frac{-n_1+n_2+n_3-n_4}{2})$ preserves as well the table. This has a natural explanation when we identify ${\rm W}(R)$ as a certain orthogonal group over $\Z$ as in~\S~\ref{pargpSO8}, see~\cite{grossinv} and the forthoming~\cite{cl2}. } 
$n_4\geq 0$. As $-1 \in {\rm W}(R)^+$, one must have $n_1+n_2+n_3+n_4 \equiv 0 \bmod 2$. \ps

\subsubsection{Case III : the Weyl group of $E_8$ as a subgroup of ${\rm SO}_9(\R)$}\label{defvprime} This case is slightly 
different and we start with some general facts, keeping the setting of the beginning of~\S~\ref{appnum}.
 Consider now the representation of ${\rm W}(R)$ on $V\oplus \R$ defined by $V'=V\oplus \varepsilon$. The map $w \mapsto (w,\varepsilon(w))$ defines an injective group homomorphism 
$${\rm W}(R) \hookrightarrow {\rm SO}(V'),$$
and we are thus again in the general situation of this chapter with this time $G={\rm SO}(V')$ and $\Gamma={\rm W}(R)$. \ps

Consider now the special case of a $R$ of type $E_8$, so that $\dim(V')=9$. Table 11 of Carter gives 
the characteristic polynomials for the action of $V$ of each ${\rm W}(R)$-conjugacy class in ${\rm W}(R)$, from which we immediately deduce the characteristic polynomial for the action of $V'=V \oplus \varepsilon$, hence the associated conjugacy class in ${\rm SO}(V')$ as $\dim(V')=9$ is odd. This is the datum (b) we need for computing $\dim(V_\lambda)^{{\rm W}(R)^+}$ : see Table~\ref{tableSO9nue} for a sample of values.\ps
 
\subsection{Reliability}  Of course, there is some possibility that we have 
made mistakes during the implementation of the program of~\S~\ref{computerprog} or of the characteristic polynomials from Carter's tables. This seems however unlikely due to the very large number of verifications we have made. \ps

The first trivial check is that the sum of characteristic polynomials of all the elements of $\Gamma$ in cases I and II is $$|{\rm W}(R)^+|(X^{\dim(V)}+(-1)^{\dim(V)})$$ as it should be. Indeed, $V$ is an irreducible representation of ${\rm W}(R)$, and even of ${\rm W}(R)^+$ in both cases. \ps

The second check is that our computer program for $\dim(V_\lambda^{\Gamma})$ always returns a positive integer ... and it does in the several hundreds of cases we have tried. As observed in the introduction, a priori each term in the sum of the degenerate Weyl character 
formula is not an integer but an element of the cyclotomic field $\Q(\zeta)$ where $\zeta$ is a $N$-th root of unity ($N=2520$ in both cases, and we indeed computed in this number field with PARI GP). This actually makes a really good check for both the degenerate Weyl character formula and Carter's tables. \ps
We will present two more evidences in the paper. One just below using a specific family of irreducible representations of ${\rm W}(R)^+$ for which one can compute directly the dimension of the ${\rm W}(R)^+$-invariants. The other one will be done much later in Chapters~\ref{chapSO7},~\ref{chapSO9},~\ref{chapSO8}, where we shall check that our computations beautifully confirm the quite intricate Arthur's multiplicity formula in a large number of cases as well.

\subsection{A check : the harmonic polynomial invariants of a Weyl group}\label{reflexionrep} We keep the notations of~\S~\ref{appnum}. 
For each integer $n\geq 0$, let ${\rm Pol}_n(V)$ denote the space homogeneous polynomials on $V$ of degree $n$ and consider 
the two formal power series in $\Z[[t]$ : $$P_R(t)= \sum_{n\geq 0} \dim({\rm Pol}_n(V)^{{\rm W}(R)})\, t^n,$$ $$A_R(t)=\sum_{n\geq 0} \dim (({\rm Pol}_n(V) \otimes \varepsilon)^{{\rm W}(R)})\, t^n.$$
By~\cite[Chap. V \S 6]{bourbaki}, if $l=\dim(V)$ and $m_1,\cdots,m_l$ are the exponents of ${\rm W}(R)$, then 
$$P_R(t)=\prod_{i=1}^l(1-t^{m_i+1})^{-1} \, \, \,\, {\rm and}\, \, \, \, A_{R}(t)=t^{|R|/2} P_R(t).$$

Let $\Delta$ be "the" ${\rm O}(V)$-invariant Laplace operator on $V$. It
induces an ${\rm O}(V)$-equivariant surjective morphism ${\rm Pol}_{n+2}(V)
\rightarrow {\rm Pol}_n(V)$, whose kernel $${\rm H}_n(V) \subset {\rm
Pol}_n(V)$$ is the space of harmonic polynomials of degree $n$ on $V$.  This
is an irreducible representation of ${\rm SO}(V)$ if $\dim(V) \neq 2$,
namely the irreducible representation with highest weight
$ne_1=(n,0,\cdots,0)$ (see e.g.  \cite[\S 5.2.3]{GW}).  One deduces the
following corollary.

\begin{cor} \begin{itemize} \item[(i)]  $\sum_{n\geq 0} \dim( {\rm H}_n(V)^{{\rm W}(R)^+} ) t^n = (1-t^2)(1+t^{|R|/2})P_R(t)$.\ps
\item[(ii)] $\sum_{n\geq 0} \dim( {\rm H}_n(V')^{{\rm W}(R)} ) t^n = (1+t^{1+|R|/2})P_R(t)$. 
\end{itemize}
\end{cor}

\begin{pf} The generating series of $\dim({\rm Pol}_n(V)^{{\rm W}^+(R)})$ is
$P_R(t)+A_R(t)$, thus the first assertion follows from the ${\rm W}(R)$-equivariant
exact sequence $0 \rightarrow {\rm H}_{n+2}(V) \rightarrow {\rm Pol}_{n+2}(V) \rightarrow {\rm Pol}_n(V) \rightarrow
0$. Observe that (i) holds whenever $R$ is irreducible or not. Assertion
(ii) follows then from (i) applied to the root system $R \cup A_1$ in $V'=V \oplus \R$.
\end{pf}

We are not aware of an infinite family $(V_i)$ of irreducible representations of ${\rm SO}(V)$ other than the ${\rm H}_i(V)$ with a simple close formula for $\dim V_i^{{\rm W}(R)^+}$. We end with some examples. Consider for instance the special case where $R$ is of type $E_7$. The exponents of ${\rm W}(R)$ are 
$1$, $5$, $7$, $9$, $11$, $13$, $17$, and $|R|=18 \cdot 7 = 126$. The power series of the corollary (i) thus becomes 
$$ \frac{1+t^{63}}{(1-t^6)(1-t^8)(1-t^{10})(1-t^{12})(1-t^{14})(1-t^{18})}$$
$$=1 + t^6 + t^8 + t^{10} + 2\,t^{12} + 2\,t^{14} + 2\,t^{16} + 4\,t^{18} + 4\,t^{20} + 4\,t^{22} + 7\,t^{24} + 7\,t^{26} + 8\,t^{28}+ o(t^{28})$$
In the case $R$ is of type $E_8$, the exponents of ${\rm W}(R)$ are $1,7,11,13,17,19,23,29$ and $|R|=8\cdot 30= 240$, so the power series of the corollary (i) is
$$\frac{1+t^{120}}{(1-t^8)(1-t^{12})(1-t^{14})(1-t^{18})(1-t^{20})(1-t^{24})(1-t^{30})}$$
$$= 1 + t^8 + t^{12} + t^{14} + t^{16} + t^{18} + 2\, t^{20} + t^{22} + 3 \, t^{24} + 2 \, t^{26} + 3 \, t^{28} + 3 \, t^{30} + 5 \, t^{32} + 3 \, t^{34} + 6 \, t^{36} + o(t^{36})$$
The power series in case (ii) for $R$ of type ${\rm E}_8$ thus starts with 
$$1 + t^2 + t^4 + t^6 + 2\, t^8 + 2 \, t^{10} + 3 \, t^{12} + 4 \, t^{14} + 5 \,t^{16} + 6 \, t^{18} + 8 \, t^{20} + 9 \, t^{22} + 12 \, t^{24} + 14\, t^{26} + 17\,t^{28} + o(t^{28})$$

In the three cases, those numbers turn out to perfectly fit our computations of the previous paragraph with the degenerate Weyl character formula! \ps
\newpage

\section{Automorphic representations of classical groups : review of Arthur's results}\label{sectionarthur}

In this section, we review Arthur's recent results~\cite{arthur} on the endoscopic classification of
discrete automorphic representations of classical groups. Our main aim is to
apply it to the level $1$ automorphic representations of certain very specific
classical groups schemes defined over the ring of integers $\Z$, namely the ones
which are reductive over $\Z$, for which the theory is substantially simpler.

\subsection{Classical semisimple groups over $\Z$}\label{quadform} By a
$\Z$-group we shall mean an affine group scheme over $\Z$ which is of finite
type. Besides the $\Z$-group ${\rm SL}_n$, the symplectic $\Z$-group ${\rm Sp}_{2g}$ and their
respective isogeny classes, we shall mainly focus on a collection of special orthogonal
$\Z$-groups that we shall briefly recall now. We refer to \cite[Ch. IV \& V]{serre}, \cite{grossinv}, \cite[Appendix
C]{conradgp}  and \cite{cl2} for a more complete discussion.
\ps

Let $L$ be a {\it quadratic abelian group} of rank $n$, which means that $L$ 
is a free abelian
group of rank $n$ equipped with a quadratic form, that we will denote by
${\rm q} : L \rightarrow \Z$. We denote by ${\rm O}_L \subset {\rm Aut}_L$ 
the {\it orthogonal group scheme over $\Z$ associated to $L$}. Recall that by definition, if
$A$ is any commutative ring then ${\rm O}_L(A)$ is the subgroup of the
general linear group ${\rm Aut}(L \otimes A)$ consisting of the elements $g$
satisfying ${\rm q}_A \circ g = {\rm q}_A$, where ${\rm q}_A : L \otimes_\Z
A \rightarrow A$ is the extension of scalars of ${\rm q}$. The {\it isometry
group} of $L$ is the group ${\rm O}(L):={\rm O}_L(\Z) \subset {\rm Aut}(L)$. \ps

A quadratic abelian group $L$ has a {\it determinant} $\det(L) \in \Z$, which is by
definition the determinant of the symmetric bilinear form $x \cdot y
= {\rm q}(x+y)-{\rm q}(x)-{\rm q}(y)$ on $L$. We say that $L$ is {\it nondegenerate}
if $\det(L)=\pm 1$ or $\det(L)=\pm 2$. This terminology is
non standard, but will be convenient for us. Note that $x \cdot x
=2 {\rm q}(x) \in 2 \Z$, thus $(x,y) \mapsto x \cdot y$ is alternate on
$L/2L$.  This forces $n$ to be even (resp.  odd) if $\det(L)=\pm 1$
(resp.  $\pm 2$). Define also the {\it signature} of $L$ as the signature $(p,q)$ of ${\rm q}_\R$. \ps

Assume now that $L$ is nondegenerate. If $n$ is even, then ${\rm O}_L$ is
smooth over $\Z$. It has exactly two connected components and we shall
denote by ${\rm SO}_L \subset {\rm O}_L$ the neutral one. This $\Z$-group may be also described as the kernel of the
Dickson-Dieudonn\'e morphism ${\rm O}_L \rightarrow \Z/2\Z$, which refines
the usual homomorphism $\det : {\rm O}_L \rightarrow \mu_2$ defined for any
quadratic abelian group $L$. When $2$ is not a zero divisor in the
commutative ring $A$, it turns out that ${\rm SO}_L(A)=\{g \in {\rm O}_L(A),
\det(g)=1\}$, but this does not hold in general.  If $n$ is
odd, we simply define ${\rm SO}_L \subset {\rm O}_L$ as the
kernel of $\det$, and we have ${\rm O}_L \simeq \mu_2 \times {\rm SO}_L$. 
In all cases,  ${\rm
SO}_L$ is then reductive over $\Z$, and actually semisimple if $n \neq 2$ (see \cite[Appendix
C]{conradgp},\cite[V.23.6]{borellivre}).  We also set ${\rm SO}(L):={\rm
SO}_L(\Z)$.  \ps

Let $L$ be a nondegenerate quadratic abelian group of rank $n$. The following two important properties hold (see~\cite[Ch. V]{serre} when $n$ is even) : \begin{itemize}\ps
\item[(i)]  If $(p,q)$ denotes the signature of $L$, then $p-q \equiv -1,0,1 \bmod 8$.\ps
\item[(ii)] For each prime $\ell$, there is a $\Z_\ell$-basis of $L \otimes \Z_\ell$ in which ${\rm q}_{\Z_\ell}$ has the form 
$$ \left\{ \begin{array}{ll} x_1x_2+x_3x_4+\dots + x_{n-1}x_n & {\rm if}\, \, n \equiv 0 \bmod 2, \\
x_1x_2+x_3x_4+\dots + x_{n-2}x_{n-1} + (-1)^{[n/2]}\frac{1}{2}\det(L) \, \, x_n^2 &  {\rm if}\, \, n \equiv 1 \bmod 2.\end{array}\right .$$
\end{itemize}
In standard terminology, part (ii) implies that the nondegenerate quadratic abelian
groups of given signature and determinant form a single genus. \ps

We now briefly discuss certain aspects of the classification of non
degenerate quadratic abelian groups $L$, starting with the definite case, which
is not only the most important one for us but the most difficult case as
well. A standard reference for this is the book by Conway and
Sloane~\cite{conwaysloane}. Replacing ${\rm q}$ by $-{\rm q}$, there is no loss of generality in restricting to the positive ones (of signature $(n,0)$),
in which case such an $L$ may be viewed as a lattice in the euclidean space
$L \otimes \R$.  Consider thus the standard euclidean space $\R^n$, with 
scalar product $(x_i)\cdot (y_i)= \sum_{i=1}^n x_i y_i$, and denote by
$$\mathcal{L}_n$$ the set of lattices $L \subset \R^n$ such that the map $x
\mapsto \frac{x \cdot x}{2}$ defines a structure of nondegenerate quadratic
abelian group on $L$. It is equivalent to ask that the lattice $L$ is even, 
i.e.  $x \cdot x \in 2\Z$ for each $x \in L$, and that $L$ has covolume $1$
(resp. $\sqrt{2}$) if $n$ is even (resp. odd). The euclidean 
isometry group ${\rm O}(\R^n)$ naturally acts on $\mathcal{L}_n$ and we shall denote by $$\mathcal{X}_n = {\rm O}(\R^n)
\backslash \mathcal{L}_n$$ the quotient set.  The map sending a positive definite quadratic
abelian group $L$ to the isometry class of the euclidean lattice $L$ inside
$L \otimes \R$ defines then a bijection between the set of isomorphism classes of non
degenerate quadratic abelian groups of rank $n$ and $\mathcal{X}_n$.  The set $\mathcal{X}_n$ is a finite by
reduction theory.  Here is what seems to be currently known about
its cardinality $h_n=|\mathcal{X}_n|$, thanks to works of Mordell, Witt,
Kneser, Niemeier, and Borcherds : $$h_1=h_7=h_8=h_9=1, \, \,\, h_{15}=h_{16}=2, \, \, \,
h_{17}=4, \, \, \, h_{23}=32, \, \, \, h_{24}=24, \, \, \, h_{25}=121$$
In all thoses cases explicit representatives of $\mathcal{X}_n$ are known,
and we recall some of them just below~:
see for instance~\cite{conwaysloane} and~\cite{borcherdsthese}.  When $n\geq 31$ then the
Minkowski-Siegel-Smith mass formula shows that $\mathcal{X}_n$ is huge, and
$h_n$ has not been determined in any case. One sometimes need to consider the set $$\widetilde{\mathcal{X}}_n={\rm SO}(\R^n)\backslash \mathcal{L}_n$$ of
direct isometry classes of even lattices $L \in \mathcal{L}_n$. One
has a natural surjective map $\widetilde{\mathcal{X}}_n \rightarrow \mathcal{X}_n$. 
The inverse image of the class of a lattice $L$ has one element if ${\rm O}(L)
\neq {\rm SO}(L)$, and two elements otherwise.  In particular,
$\widetilde{\mathcal{X}}_n \isomo \mathcal{X}_n$ if $n$ is odd.\ps 
Some important even euclidean lattices are related to root systems as
follows.  Let $R \subset \R^n$ be a root system of
rank $n$ such that each $x \in R$ satisfies $x \cdot x =2$. In particular,
the irreducible components of $R$ are of type $A$, $D$ or $E$. The set $R$
generates a lattice of $\R^n$ denoted by ${\rm Q}(R)$ in~\cite[Ch. VI \S 1]{bourbaki}, that we view as a quadratic
abelian group via the quadratic form $x \mapsto \frac{x \cdot
x}{2}$.  It is called the {\it root lattice} associated to $R$.  It contains
exactly the same information as $R$, because of the well known property
$R=\{x \in Q(R), x \cdot x =2\}$.  The Cartan matrix of the root system $R$
is symmetric and is a Gram matrix for the bilinear form of ${\rm Q}(R)$;  its determinant is the index of connexion of $R$. It follows that the root lattices
${\rm A}_1$, ${\rm E}_7$, ${\rm E}_8$ and ${\rm E}_8 \oplus {\rm A}_1$, associated
respectively to root systems of type
$A_1$, $E_7$, $E_8$ and $E_8\coprod A_1$, are nondegenerate of ranks $n=1, 7, 8$ and $9$. Up to
isometry, they are the unique such lattices in these dimensions, and the only ones we
shall really need in this paper. The isometry groups of these lattices, and
more generally of root lattices, are well known.  Indeed, the
isometry group of ${\rm Q}(R)$ is by definition the group denoted ${\rm A}(R)$
in~\cite[Ch. VI \S 1]{bourbaki}. It 
contains the Weyl group ${\rm W}(R)$ as a normal subgroup. Moreover, if $B \subset R$ is
a basis of $R$, and if $\Gamma \subset {\rm A}(R)$ denotes de subgroup
preserving $B$, then $\Gamma$ is isomorphic to the automorphism group of the
Dynkin diagramm of $R$ and the group ${\rm A}(R)$ is a semi-direct product of $\Gamma$ by ${\rm W}(R)$ by~\cite[Ch. VI no 1.5, Prop. 16]{bourbaki}.  It
follows that in the four cases above, we have ${\rm O}({\rm Q}(R))={\rm W}(R)$. \ps

In general dimension $n = 8k + s$ with $s =-1$ (resp. $s=0$, resp. $s=1$). 
We obtain an example of positive definite nondegenerate quadratic abelian
group $\Lambda_n$ by considering the orthogonal direct sum ${\rm E}_7 \oplus {\rm E}_8^{k-1}$ (resp.  ${\rm
E}_8^k$, resp.  ${\rm E}_8^k\oplus {\rm A}_1$). Let us call it the {\it standard} postive definite quadratic abelian group of rank
$n$. We will simply write $${\rm O}_n \, \, \, {\rm  and} \, \, \, {\rm SO}_n$$
for ${\rm O}_{\Lambda_n}$ and ${\rm SO}_{\Lambda_n}$. More generally, if $p \geq q$ are
nonnegative integers, and if $p-q
\equiv -1,0,1 \bmod 8$, the orthogonal direct sum of $q$ hyperbolic planes\footnote{The hyperbolic plane over $\Z$ is the abelian group $\Z^2$ equipped with the quadratic form ${\rm q}(x,y)=xy$.}
over $\Z$ and of $\Lambda_{p-q}$ is a quadratic abelian group of signature $(p,q)$, that we shall call {\it standard} as well for this signature. 
When $q>0$ it turns out to be the only nondegenerate quadratic abelian
group of signature $(p,q)$ up to isometry (see \cite[Ch. V]{serre},\cite{grossinv}), and we shall simply denote by ${\rm
SO}_{p,q}$ its special orthogonal group scheme. When $|p-q|\leq 1$, this is a
Chevalley group. In low dimension, we have
the following exceptional isomorphisms over $\Z$ : $${\rm
SO}_{1,1} \simeq \mathbb{G}_m, \, \, \, \, \, {\rm SO}_{2,1} \simeq {\rm
PGL}_2, \,\,\,\,\,{\rm SO}_{3,2}\simeq {\rm PGSp}_4,$$ as well as a central
isogeny ${\rm SO}_{2,2} \rightarrow {\rm PGL}_2 \times {\rm PGL}_2$.\ps

\begin{rem}\label{remaut}{\rm If $L$ is any of the standard quadratic abelian group defined above, then it always contains elements $\alpha$ such that 
$\alpha \cdot \alpha = 2$. If ${\rm s}_\alpha$ denotes the 
orthogonal symmetry with respect to such an $\alpha$ 
then ${\rm s}_\alpha \in {\rm O}(L) \backslash {\rm SO}(L)$, and the conjugation by ${\rm s}_\alpha$ defines a $\Z$-automorphism of ${\rm SO}_L$.}\end{rem}

\subsection{Discrete automorphic representations}\label{discautrep} Let $G$ by a semisimple $\Z$-group.  Denote by
$\Pi(G)$ the set of isomorphism classes of complex representations $\pi$ of $G(\AAA)$ such that $\pi \simeq \pi_\infty \otimes
\pi_f$, where : \begin{itemize}\ps

\item[(i)] $\pi_f$ is a smooth irreducible complex representation of $G(\AAA_f)$, and $\pi_f$ is {\it unramified}, {\it i.e.} such that $\pi_f^{G(\widehat{\Z})} \neq
0$, \ps
\item[(ii)] $\pi_\infty$ is an irreducible unitary representation of
$G(\R)$. \end{itemize}\ps

Of course $\AAA=\R \times \AAA_f$ denotes the ad\`ele ring of $\Q$
and $\AAA_f=\widehat{\Z} \otimes \Q$ the ring of finite ad\`eles. Denote by $\mathcal{H}(G)$  the complex Hecke-algebra of the pair $(G(\AAA_f),G(\widehat{\Z}))$. 
By well-known results of Satake and Tits (\cite{satake},\cite{titscorvallis}), the ring $\mathcal{H}(G)$ is commutative, so that $\dim \pi_f^{G(\widehat{\Z})}=1$ for each unramified smooth irreducible complex representation $\pi_f$ of $G(\AAA_f)$. \ps

Recall that the homogeneous space $G(\Q) \backslash G(\AAA)$ has a nonzero
$G(\AAA)$-invariant Radon measure (Weil) of finite volume (Borel,
Harish-Chandra, see~\cite[\S 5]{borelfini}).  Consider the Hilbert space $$\mathcal{L}(G)={\rm
L}^2(G(\Q)\backslash G(\AAA)/G(\widehat{\Z}))$$ of square-integrable
functions on $G(\Q)\backslash G(\AAA)$ for this measure which are
$G(\widehat{\Z})$-invariant on the right (\cite[\S 4]{boreljacquet},\cite[Ch. 3]{GGPS}).  This space $\mathcal{L}(G)$ is
equipped with a unitary representation of $G(\R)$ by right translations and
with an action of the Hecke algebra $\mathcal{H}(G)$ commuting with $G(\R)$.  The subspace
$\mathcal{L}_{\rm disc}(G) \subset \mathcal{L}(G)$ is defined as the closure
of the sum of the irreducible closed subspaces for the $G(\R)$-action, it
is stable by $\mathcal{H}(G)$.  A fundamental result of Harish-Chandra~\cite[Ch. 1 Thm. 1]{harishchandra}
asserts that each irreducible representation of $G(\R)$ occurs with finite multiplicity in $\mathcal{L}(G)$.  It follows
that \begin{equation}\label{eqldisc}\mathcal{L}_{\rm disc}(G) =
\underset{\pi \in \Pi(G)}{\overline{\bigoplus}} m(\pi)\, \, \pi_\infty
\otimes \pi_f^{G(\widehat{\Z})},\end{equation} where the integer $m(\pi)\geq
0$ is the multiplicity of $\pi$ as a sub-representation of $\mathcal{L}(G)$.  We
denote by $$\Pi_{\rm disc}(G) \subset \Pi(G)$$ the subset of $\pi$ such that
$m(\pi) \neq 0$ and call them the {\it discrete automorphic representations
of the $\Z$-group $G$}.  A classical result of Gelfand and Piatetski-Shapiro~\cite{GGPS}
asserts that the subspace of cuspforms of $G$, which is stable by $G(\R)$
and $\mathcal{H}(G)$, is included in $\mathcal{L}_{\rm disc}(G)$ and we
denote by $$\Pi_{\rm cusp}(G) \subset \Pi_{\rm disc}(G)$$ the subset of
$\pi$ consisting of cusp forms.  \ps

\subsection{The case of Chevalley and definite semisimple
$\Z$-groups}\label{gpdef} All those automorphic representations have various
models, depending on the specific group $G$ and the kind of $\pi_\infty$ we
are interested in.  We shall content ourselves with the following classical
descriptions.  Consider the class set of $G$ $${\rm Cl}(G)=G(\Q)\backslash
G(\AAA_f)/G(\widehat{\Z}).$$ This is a finite set by~\cite[\S 5]{borelfini} and we set
$h(G)=|{\rm Cl}(G)|$.  A well-known elementary fact (see {\it loc. cit.} \S 2) is that $$h({\rm
SL}_n)=h({\rm PGL}_n)=h({\rm Sp}_n)=h({\rm PGSp}_n)=1,$$ and part of what we said
in~\S~\ref{quadform} amounts to
saying as well that $h({\rm SO}_{p,q})=1$ if $pq \neq 0$.  More generally, the
strong approximation theorem of Kneser ensures that $h(G)=1$ if $G$ is simply
connected and $G(\R)$ has no compact factor (see ~\cite{pr}).  Recall that a {\it
Chevalley group} is a split semisimple $\Z$-group.  We refer to~\cite{sga}
and~\cite{conradgp} for the general theory of Chevalley groups.

\begin{prop}\label{hchev} Let $G$ be a Chevalley group. Then $h(G)=1$ and the inclusion $G(\R) \rightarrow G(\AAA)$ induces a homeomorphism 
$$G(\Z) \backslash G(\R) \rightarrow G(\Q) \backslash G(\AAA)/G(\widehat{\Z}).$$  
Moreover, $G(\R)/G(\Z)$ is connected and ${\rm Z}(G)(\R)={\rm Z}(G)(\Z)$.
\end{prop}

\begin{pf} We refer to~\cite[Exp. XXII~\S 4.2, \S 4.3]{sga} and~\cite[Chap.
6]{conradgp} for central isogenies between semisimple group schemes.  Let $s
: G_{\rm sc} \rightarrow G$ be a central isogeny with $G_{\rm sc}$ simply
connected.  The $\Z$-group $G_{\rm sc}$ is a Chevalley group as well.  Let
$T$ be a maximal $\Z$-split torus in $G$ and let
$T_{\rm sc} \subset G_{\rm sc}$ be the split maximal torus defined as the inverse image of $T$ by $s$. Recall that for any field $k$, we have the following
simple facts from Galois cohomology : \begin{itemize}

\item[(i)] $s(G_{\rm sc}(k))$ is a normal
subgroup of $G(k)$,\ps

\item[(ii)] $T(k)s(G_{\rm sc}(k)) = G(k)$. \ps

\end{itemize}


In particular, $G(\AAA_f)=T(\AAA_f)s(G_{\rm sc}(\AAA_f))G(\widehat{\Z})$.  But
$T(\AAA_f)=T(\Q)T(\widehat{\Z})$ as $T$ is $\Z$-split. It follows that 
\begin{equation}\label{decompggsc}G(\AAA_f)=T(\Q)T(\widehat{\Z})s(G_{\rm sc}(\AAA_f))G(\widehat{\Z})=T(\Q)
s(G_{\rm sc}(\AAA_f))G(\widehat{\Z}),\end{equation}
where the last equality comes from (i) above. But $h(G_{\rm sc})=1$ by the
strong approximation theorem, thus $h(G)=1$ as well by the above identity.
\ps

Observe now that the map of the first statement is trivially injective, and
even surjective as $h(G)=1$.  It is moreover continuous
and open, as $G(\R) \times G(\widehat{\Z})$ is open in $G(\AAA)$, hence a homeomorphism. 
\ps Let us check that $G(\R)/G(\Z)$ is connected. By (ii) again, observe that
$$G(\R)=s(G_{\rm sc}(\R))T(\R).$$ 
But $G_{\rm sc}(\R)$ is connected as
$G_{\rm sc}$ is connected and simply connected, by a classical result of Steinberg.  We conclude as $T(\Z)
\subset G(\Z)$ meets every connected component of $T(\R)$, since $T$ is
$\Z$-split.  The last assertion follows from the following simple fact
applied to $A={\rm Z}(G)$ : if $A$ is a finite multiplicative $\Z$-group scheme,
then the natural map $A(\Z) \rightarrow A(\R)$ is bijective (reduce to the
case $A=\mu_n$).\end{pf}

When $G={\rm PGSp}_{2g}$ or ${\rm Sp}_{2g}$, the cuspidal automorphic representations $\pi$ of $G$ such that $\pi_\infty$ is a holomorphic 
discrete series representation are closely related to vector valued Siegel cuspforms : see e.g.~\cite{asgarischmidt}.\ps
A semisimple $\Z$-group $G$ will be said {\it definite} if $G(\R)$ is compact. This is somewhat the opposite case of Chevalley groups, but a case of great interest in this paper. For instance the semisimple $\Z$-group ${\rm SO}_n$ defined in~\S\ref{quadform} is definite and there is a natural bijection 
\begin{equation}\label{classquad} {\rm Cl}({\rm SO}_n) \isomo \widetilde{X}_n,\end{equation} because the rank $n$ 
definite quadratic forms over $\Z$ form a single genus, as recalled in~\S\ref{quadform} (see~\cite[\S 2]{borelfini}). If $G$ is definite, then 
$$\mathcal{L}(G)=\mathcal{L}_{\rm disc}(G)$$
by the Peter-Weyl theorem, and the discrete automorphic representations of $G$ have very simple models. Automorphic forms for definite semisimple $\Z$-groups are a special case of "algebraic modular forms" in the sense of Gross~\cite{grossalgform}.

\begin{prop}\label{definiteform} Let $G$ be a semisimple definite $\Z$-group and let $(\rho,V)$ be an irreducible continous representation of $G(\R)$. The vector space ${\rm Hom}_{G(\R)}(V,\mathcal{L}(G))$ 
is canonically isomorphic to the space of covariant functions $$M_\rho(G)=\{ f : G(\AAA_f)/G(\widehat{\Z}) \rightarrow V^\ast, \, \, \, f(\gamma g)={}^{t}\!\rho(\gamma)^{-1} f(g) \,\, \, \forall \gamma \in G(\Q), g \in G(\AAA_f)\}.$$
In particular, $\dim(M_\rho(G))=\underset{\pi \in \Pi_{\rm disc}(G), \pi_\infty \simeq
V}{\sum}\, m(\pi)$.
\end{prop}
The canonical bijection of the statement is $\varphi \mapsto (g \mapsto (v \mapsto \varphi(v)(1 \times g)))$, where $\varphi \in  {\rm Hom}_{G(\R)}(V,\mathcal{L}(G))$, $v \in V$ and $g \in G(\AAA_f)$. If $g_1,\cdots,g_{h(G)} \in G(\AAA_f)$ are representatives for the classes in ${\rm Cl}(G)$, the evaluation map $f \mapsto (f(g_i))$ defines thus a bijection $$M_\rho(G) \isomo \prod_{i=1}^{h(G)} (V^\ast)^{\Gamma_i}$$
where $\Gamma_i$ is the finite group $G(\R) \cap g_i^{-1}G(\Q)g_i$. In particular, to compute $M_\rho(G)$ we are reduced to compute invariants of the finite group $\Gamma_i \subset G(\R)$ in the representation $V$, 
what we have already studied in Chapter~\ref{polinv}. Indeed, the compact group $G(\R)$ is always connected by a classical result of Chevalley~\cite[V.24.6 (c) (ii)]{borellivre}. \ps
Of course if $g_i=1$, then $\Gamma_i=G(\Z)$. In the example of the group $G={\rm SO}_n$, if $L_i \in \mathcal{L}_n$ is the lattice corresponding to $g_i$ via the bijection
\eqref{classquad}, then $\Gamma_i={\rm SO}(L_i)$. Later, we will study in details the cases $G={\rm SO}_n$ where $n=7,8$ and $9$, and the definite semisimple $\Z$-group ${\rm G}_2$.

\subsection{Langlands parameterization of $\Pi_{\rm
disc}(G)$}\label{langparam} In this paragraph, we discuss a parameterization 
of the elements of $\Pi(G)$ due to Satake and Harish-Chandra, according to Langlands point of
view~\cite[\S 2]{langlandsyale}. An important role is played by the Langlands dual group
of $G$, for
which a standard reference is Borel's paper~\cite{borelcorvallis}.\ps

Let $G$ be any semisimple $\Z$-group. As observed by Gross in~\cite{grossinv}, the natural action of the absolute Galois group of $\Q$ on the
based root datum of $G_{\overline{\Q}}$ is trivial, as each non trivial number field has a ramified prime. It follows that 
the $\Q$-group $G_\Q$ is an inner form of a split Chevalley group over $\Q$. In particular, the Langlands
dual group of $G$ may simply be defined as a complex semisimple algebraic group $\widehat{G}$
equipped with an isomorphism between the dual based root datum of $\widehat{G}$
and the based root datum of $G_{\overline{\Q}}$ (see \cite{borelcorvallis}). The group $\widehat{G}$
itself is well defined up to inner automorphism. When $G$ is either ${\rm PGL}_n$, 
${\rm Sp}_{2n}$ or of the form ${\rm SO}_L$ for $L$ a nondegenerate quadratic abelian
group of rank $n$, it is well-known that $\widehat{G}$ 
is respectively isomorphic to $${\rm SL}_n(\C),  \, {\rm SO}_{2n+1}(\C), \, {\rm
Sp}_{n-1}(\C)\,\, (n \, \,{\rm odd})\, \,\, \, {\rm or}\, \, \,\,{\rm SO}_n(\C) \,
\, (n\, \, {\rm even}).$$
If $G={\rm SO}_L$, then both $G_\Q$ and the $G(\AAA_f)$-conjugacy class of $G(\widehat{\Z}) \subset G(\AAA_f)$ only depend on 
the signature $L \otimes \R$ by the property of the genus of $L$ discussed in~\S\ref{quadform}. We shall thus loose nothing in assuming once and for all 
that $L$ is the standard quadratic abelian group as defined {\it loc. cit}. \ps

If $H$ is the group of $\C$-points of a complex semisimple algebraic group over $\C$, we shall denote 
by $$\mathcal{X}(H)$$ the set of collections $(c_v)$ indexed by the places $v$ of $\Q$, where each $c_p$ (resp. $c_\infty$) is a semisimple conjugacy class in 
$H$ (resp. ${\rm Lie}_\C(H)$). If $\pi \in \Pi(G)$, then $\pi_f$ is isomorphic to the restricted tensor product over all primes $p$ of irreducible 
smooth representations $\pi_p$ of $G(\Q_p)$ which are well defined up to isomorphism and unramified (i.e. $\pi_p^{G(\Z_p)} \neq 0$). By Langlands' interpretation of the work of Harish-Chandra and Satake, we have a natural parameterization map 
		$$c : \Pi(G) \longrightarrow \mathcal{X}(\widehat{G}), \quad \pi \mapsto (c_p(\pi)),$$
where : \begin{itemize} \item[(i)] for each prime $p$ the semisimple
conjugacy class $c_p(\pi)$ is Satake parameter of $\pi_p$.\ps
\item[(ii)] $c_\infty(\pi)$ is the conjugacy class defined by the infinitesimal character of $\pi_\infty$ 
and the Harish-Chandra isomorphism.  
\end{itemize}

See~\cite[\S 2]{langlandsyale}, \cite{borelcorvallis} and~\cite{grossatake} for a discussion of
Satake's parameterization in those terms. We recall parameterization (ii) for the convenience of the reader (see e.g. Delorme's survey~\cite{delorme} for precise references). Let $\mathfrak{g}$ be the complex Lie algebra of $G(\C)$. If $V$ is a
unitary representation of $G(\R)$, the subspace $V^\infty \subset V$ of
indefinitely differentiable vectors for the action of $G(\R)$ carries an action of the envelopping algebra ${\rm
U}(\mathfrak{g})$ of $\mathfrak{g}$, and is dense in $V$. If $V$ is irreducible, a version of Schur's lemma implies that the center
${\rm Z}({\rm
U}(\mathfrak{g}))$ of ${\rm U}(\mathfrak{g})$ acts by scalars on $V^\infty$, which
thus defines a
$\C$-algebra homomorphism ${\rm Z}({\rm U}(\mathfrak{g})) \rightarrow \C$ called
the {\it infinitesimal character} of $V$. The last point to understand the meaning of (ii) above is
that there is a canonical bijection between ${\rm Hom}_{\C-{\rm alg}}({\rm
Z}({\rm
U}(\mathfrak{g})),\C)$ and the set of semisimple conjugacy classes in
$\widehat{\mathfrak{g}}$, that we now recall. \ps

Let $\mathfrak{t}$ be a Cartan algebra of $\mathfrak{g}$ and let ${\rm W}$ denote the Weyl
group of $(\mathfrak{g},\mathfrak{t})$. The Harish-Chandra
isomorphism is a canonical isomorphism ${\rm Z}({\rm U}(\mathfrak{g})) \isomo 
({\rm Sym}\, \mathfrak{t})^{\rm W}$. It follows that an infinitesimal character may be
viewed as a ${\rm W}$-orbit of elements in the dual vector space ${\rm Hom}(\mathfrak{t},\C)$.
But the ${\rm L}$-group datum defining $\widehat{G}$ naturally identifies ${\rm Hom}(\mathfrak{t},\C)$ with a
Cartan algebra $\widehat{\mathfrak{t}}$ of $\widehat{\mathfrak{g}}$, and $W$ with the Weyl group
of $(\widehat{\mathfrak{g}},\widehat{\mathfrak{t}})$. It follows that the
set of ${\rm
W}$-orbits of
elements in ${\rm Hom}(\mathfrak{t},\C)$ is in canonically bijection with
the set of semisimple conjugacy classes in $\widehat{\mathfrak{g}}$. \ps

A result of Harish-Chandra asserts that up to isomorphism, there are
at most finitely many irreducible unitary representations of $G(\R)$ of any given
infinitesimal character~\cite[Cor. 10.37]{knapplivre}.  When $G(\R)$ is compact, in which case it
is necessarily connected by a theorem of Chevalley already cited in~\S\ref{gpdef}, the situation is
much simpler. Indeed, let $\mathfrak{t}
\subset \mathfrak{b} \subset \mathfrak{g}$ be a Borel subalgebra and let
$\rho \in {\rm Hom}(\mathfrak{t},\C)$ be the half-sum of the positive roots
of $(\mathfrak{g},\mathfrak{t},\mathfrak{b})$. The infinitesimal character of the irreducible
representation $V_\lambda$ of $G(\R)$ of highest weight $\lambda \in {\rm
Hom}(\mathfrak{t},\C)$ relative to $\mathfrak{b}$ is the conjugacy class of the element
$\lambda + \rho$ (\cite[\S 7.4.6]{dixmier}), viewed as an element of $\widehat{\mathfrak{t}}={\rm Hom}(\mathfrak{t},\C)$. In particular, it uniquely determines
$V_\lambda$. \ps

\subsection{Arthur's symplectic-orthogonal alternative} \label{parorthosymp} 
By a {\it classical
semisimple group over $\Z$} we shall mean either ${\rm Sp}_{2g}$ for $g\geq
1$, or ${\rm SO}_L$ for $L$ a standard quadratic form over $\Z$ of rank
$\neq 2$. In particular, ${\rm SO}_L$ is either ${\rm SO}_n$ or ${\rm SO}_{p,q}$ defined in \S\ref{quadform}. The {\it classical Chevalley groups} are the $\Z$-groups ${\rm Sp}_{2g}$, ${\rm SO}_{p,q}$ with $p-q \in \{0,1\}$, and the trivial $\Z$-group ${\rm SO}_1$. The definite classical semisimple groups over $\Z$ are the ${\rm SO}_n$.\ps

Fix $G$ a classical semisimple group over $\Z$. Arthur's classification describes 
$\Pi_{\rm disc}(G)$ in terms of the $\Pi_{\rm cusp}(\PGL_m)$ for various $m$'s, and our aim from now is to 
recall this classification. We shall denote by $${\rm St} :  \widehat{G}
\hookrightarrow \SL_n(\C)$$ the standard representation of its dual group, which defines in particular the integer
$n=n(G)$. For instance $n({\rm Sp}_{2g})=2g+1$ and $n({\rm SO}_m)=2[m/2]$. This group homomorphism defines in particular a natural map
$\mathcal{X}(\widehat{G}) \rightarrow \mathcal{X}(\widehat{\PGL_n})$ that
we shall still denote by ${\rm St}$.

\begin{thmconda}\label{orthosymp} (Arthur) For any $n\geq 1$ and any given self-dual $\pi \in
\Pi_{\rm cusp}(\PGL_n)$ there is a unique classical Chevalley group
$G^\pi$ with $n(G^\pi)=n$ such that there exists 
$\pi' \in \Pi_{\rm disc}(G^\pi)$ satisfying ${\rm St}(c(\pi'))=c(\pi)$. 
\end{thmconda}

 This is~\cite[Thm. 1.4.1]{arthur}. As $n(G^\pi)=n$, the only possibilities for $G^\pi$ are thus $G^\pi={\rm SO}_1$ if $n=1$, $G^\pi={\rm Sp}_{n-1}$ if $n>1$ is odd, $G^\pi={\rm
SO}_{\frac{n}{2},\frac{n}{2}-1}$ or 
${\rm SO}_{\frac{n}{2},\frac{n}{2}}$ if $n$ is even. This last case only exists for $n>2$, which forces $G^\pi={\rm SO}_{2,1}\simeq {\rm
PGL}_2$ if $n=2$. \ps

As self-dual $\pi \in \Pi_{\rm cusp}(\PGL_n)$ will be said {\it orthogonal}
{\rm (resp. {\it symplectic})} if $\widehat{G^\pi}$ is isomorphic to a complex special orthogonal group (resp. symplectic group). For short, we shall define the {\it sign} of $\pi$ $$s(\pi) \in
\{\pm 1\}$$
to be $1$ if $\pi$ is orthogonal, $-1$ otherwise. If $n$ is odd then
$\pi$ is necessarily orthogonal, i.e. $s(\pi)=1$. \ps

\begin{definition}\label{defpbos} Let $n\geq 1$ be an integer. We denote by : \begin{itemize}\ps
\item[-] $\Pi_{\rm cusp}^\bot(\PGL_n) \subset \Pi_{\rm cusp}(\PGL_n)$ the subset of self-dual $\pi$, i.e. such that $\pi^\vee\simeq \pi$,\ps
\item[-] $\Pi^{\rm s}_{\rm cusp}(\PGL_n) \subset \Pi_{\rm cusp}^\bot(\PGL_n)$ the subset of symplectic $\pi$,\ps
\item[-] $\Pi^{\rm o}_{\rm cusp}(\PGL_n) \subset \Pi_{\rm cusp}^\bot(\PGL_n)$ the subset of orthogonal $\pi$.
\end{itemize}
We have $\Pi^\bot_{\rm cusp}(\PGL_n)=\Pi^{\rm s}_{\rm cusp}(\PGL_n) \coprod \Pi_{\rm cusp}^{\rm o}(\PGL_n)$.
\end{definition}

\subsection{The symplectic-orthogonal alternative for polarized algebraic
regular cuspidal automorphic representations of $\GL_n$ over
$\Q$}\label{demialg} Let $\pi$ be a cuspidal automorphic representation of
${\rm GL}_n$ over $\Q$ satisfying (a), (b) and (c') of \S\ref{altsointro}.  The self-dual representation
$\pi':=\pi \otimes |\cdot|^{\frac{w(\pi)}{2}}$ has a trivial central character,
hence defines a self-dual element of $\Pi_{\rm cusp}(\PGL_n)$.  Our aim in
this paragraph is to show that $\pi$ is orthogonal (resp.  symplectic) in
the sense of~\S\ref{altsointro} if and only if $\pi'$ is so in the sense of
Arthur (\S\ref{parorthosymp}).  A key role will be played by condition (c')
on $\pi$.  This forces us to discuss first Langlands parameterization for
$\GL_n(\R)$ in more details than we have done so far.  We refer to
\cite{langlandsreel}, \cite{borelcorvallis}, and especially \cite{knapp},
for more details.  \ps

Recall that the Weil group of $\C$ is the topological group ${\rm W}_\C:=\C^\ast$. The Weil group
of $\R$, denoted ${\rm W}_\R$, is a 
non-split extension of ${\rm Gal}(\C/\R)$ by ${\rm W}_\C$, for the
natural action of ${\rm Gal}(\C/\R)$ on $\C^\times$. The set ${\rm W}_\R \backslash {\rm W}_\C$
contains a unique ${\rm W}_\C$-conjugacy class of elements $j \in {\rm W}_\R \backslash {\rm W}_\C$ such that $j^2=-1$ (as elements of $\C^\ast$), and we fix once and for all such an element. According to Langlands parameterization, if $\pi$ is a cuspidal automorphic representation of ${\rm GL}_n$ over $\Q$ then the unitary representation $\pi_\infty$ of $\GL_n(\R)$ is 
uniquely determined up to isomorphism by its {\it Langlands parameter}. By definition, this is an isomorphism class of 
continuous semisimple representations $${\rm L}(\pi_\infty) : {\rm W}_\R \rightarrow
\GL_n(\C).$$
It refines the infinitesimal character of $\pi_\infty$, viewed as a semisimple conjugacy class in ${\rm M}_n(\C)={\rm Lie}_\C \widehat{\GL_n}$, which may actually be read on the restriction of ${\rm L}(\pi_\infty)$ to ${\rm W}_\C$. Concretely, this restriction is a direct sum of continuous homomorphisms $\chi_i : \C^\ast \rightarrow \C^\ast$ for $i=1,\cdots,n$, which are unique up to reordering. For each $i$, there are unique $\lambda_i,\mu_i \in \C$ such that $\lambda_i-\mu_i \in \Z$,  
satisfying $$\chi_i(z) = z^{\lambda_i} \overline{z}^{\mu_i}$$ for all $z \in \C^\ast$ : the $n$ eigenvalues of the infinitesimal character of $\pi_\infty$ are the elements $\lambda_1,\cdots,\lambda_n \in \C$. We have used Langlands convenient notation : if $\lambda,\mu \in \C$ satisfy $\lambda-\mu \in \Z$, and if $z \in \C^\ast$, then $z^\lambda \overline{z}^\mu$ denotes the element of $\C^\ast$ defined by $(z\overline{z})^{\frac{\lambda+\mu}{2}} (\frac{z}{|z|})^{\lambda-\mu}$.\ps

Langlands parameters are easy to classify, as the irreducible continuous representations of ${\rm W}_\R$ are either one dimensional or induced from a one dimensional representation of ${\rm W}_\C$. The first ones are described thanks to the natural isomorphism ${\rm W}_\R^{\rm ab} \isomo  \R^\ast$ sending $z \in {\rm W}_\C$ to $z \overline{z}$ : they have thus the form $|\cdot|^s$ or $\varepsilon_{\C/\R} |\cdot|^s$, where $s \in \C$ and $\varepsilon_{\C/\R}(x)=x/|x|$ is the sign character. The irreducible continuous $2$-dimensional representations of ${\rm W}_\R$ are the ${\rm I}_w \otimes |\cdot|^s$ for $w>0$ and $s \in \C$, where we have set for any integer $w\geq 0$
$${\rm I}_w={\rm Ind}_{W_\C}^{W_\R} \,z^{-w/2}\overline{z}^{w/2}.$$  \ps

\begin{prop}\label{paramalg} Let $\pi$ be a cuspidal automorphic representation of $\GL_n$ over $\Q$ satisfying {\rm (a)}, {\rm (b)} and {\rm (c')} of the introduction, of Hodge weights $w_i$. Then 
$${\rm L}(\pi_\infty) \otimes |\cdot|^{w(\pi)/2} \simeq \left\{ \begin{array}{ll} \bigoplus_{i=1}^{\frac{n}{2}} {\rm I}_{w_i}  & {\rm if  }\, \, \,  n \equiv 0 \bmod 2, \\ 
\varepsilon_{\C/\R}^{\frac{n-1}{2}} \oplus \bigoplus_{i=1}^{\frac{n-1}{2}} {\rm I}_{w_i} & {\rm if} \, \, \, n \equiv 1 \bmod 2. \end{array} \right.$$
Moreover, if $w$ and $n$ are even then $n \equiv 0 \bmod 4$.
\end{prop}

The main ingredient in the proof of this proposition is the following special case of Clozel's purity lemma~\cite[Lemma 4.9]{clozel}.

\begin{lemme}\label{clopl} (Clozel's purity lemma) Let $\pi$ be a cuspidal automorphic representation of $\GL_n$ over $\Q$. Assume that the eigenvalues of the infinitesimal character of $\pi_\infty$ are in $\frac{1}{2}\Z$. Then there is an element $w \in \Z$ such that ${\rm L}(\pi_\infty) \otimes |\cdot|^{w/2}$ is a direct sum of representations of the form $1$, $\varepsilon_{\C/\R}$, or ${\rm I}_{w'}$ for $w' \in \Z$.
\end{lemme}

\begin{pf}  (of Proposition~\ref{paramalg}) We apply Clozel's purity lemma to $\pi$. Condition (a) on $\pi$ ensures that ${\rm L}(\pi_\infty)^\ast \simeq {\rm L}(\pi_\infty) \otimes |\cdot|^{w(\pi)}$. As both the ${\rm I}_{w'}$, $1$ and $\varepsilon_{\C/\R}$ are self-dual, it follows that the element $w$ given by the purity lemma coincides with $w(\pi)$. Condition (c') on $\pi$, and the relation $-k_i+w(\pi)/2 = -w_i$ for $i=1,\cdots,[n/2]$, concludes the proof when $\frac{w(\pi)}{2}$ is not a weight of $\pi$ (e.g. when $w(\pi)$ is odd). By assumption (c'), if $\frac{w(\pi)}{2}$ is a weight of $\pi$ then it has multiplicity $1$ if $n$ is odd and $2$ if $n \equiv 0 \bmod 4$.  But by condition (b) on $\pi$ and the structure of the id\`eles of $\Q$, the global central character of $\pi$ is $|\cdot|^{-\frac{n \,w(\pi)}{2}}$, so that $\det({\rm L}(\pi_\infty) \otimes |\cdot|^{w(\pi)/2})=1$. Observe that for any $w' \in \Z$ we have
$$\det({\rm I}_{w'}) = \varepsilon_{\C/\R}^{w'+1}.$$ 
The proposition follows when $n$ is odd, as well as when $n \equiv w(\pi) \equiv 0 \bmod 2$ since $${\rm I}_0 \simeq 1 \oplus \varepsilon_{\C/\R}.$$ 
\end{pf}

We are now able to state a strenghtening of Arthur's Theorem~\ref{orthosymp}, which is more precise at the infinite place. Indeed, let $\pi \in \Pi_{\rm cusp}(\PGL_n)$ be self-dual and let
$$\mathrm{L}(\pi_\infty) : {\rm W}_\R \longrightarrow {\rm SL}(n,\C)$$
be the Langlands parameter of $\pi_\infty$. Arthur shows that
$\mathrm{L}(\pi_\infty)$ maybe conjugated into 
${\rm St}(\widehat{G^\pi}) \subset \SL(n,\C)$ (\cite[Thm. 1.4.2]{arthur}). 
Note that the $\widehat{G}^\pi$-conjugacy class of the resulting Langlands parameter 
\begin{equation}\label{arthurliftreel}\widetilde{\mathrm{L}}(\pi_\infty) : {\rm W}_\R \longrightarrow \widehat{G^\pi}\end{equation} is not quite canonical, but so is its ${\rm Out}(\widehat{G^\pi})$-orbit. 

\begin{cor}\label{cororthosymp} Let $\pi$ be a cuspidal automorphic representation of $\GL_n$ over $\Q$ satisfying {\rm (a)}, {\rm (b)} and {\rm (c')} of the introduction. Then $\pi$ is orthogonal (resp. symplectic) in the sense of~\S\ref{altsointro} if and only if $\pi \otimes |\cdot|^{w(\pi)/2}$ is so in the sense of \S\ref{parorthosymp}.
\end{cor}

\begin{pf} Consider the self-dual representation $\pi'=\pi \otimes |\cdot|^{w(\pi)/2}$ in $\Pi_{\rm cusp}(\PGL_n)$. By Proposition~\ref{paramalg}, ${\rm L}(\pi'_\infty) : {\rm W}_\R
\rightarrow {\rm SL}_n(\C)$ is a direct sum of distinct irreducible self-dual
representations of ${\rm W}_\R$. It follows that if ${\rm L}(\pi'_\infty)$ preserves a nondegenerate
pairing on $\C^n$ then each irreducible subspace is nondegenerate as well. Moreover, the
$2$-dimensional representation ${\rm I}_w$ has determinant $\varepsilon_{\C/\R}^{w+1}$
and may be conjugate into ${\rm O}_2(\C)$ if and only if $w$ is even. The
result follows as ${\rm L}(\pi_\infty)$ may be conjugate into ${\rm
St}(\widehat{G^\pi})$ by the aforementioned result of Arthur.
\end{pf}

It will be convenient in the sequel to adopt a slightly different point of
view, although eventually equivalent, on the representations $\pi$ studied
in the introduction.  Consider the cuspidal automorphic representations
$\pi$ of ${\rm GL}_n$ over $\Q$ such that : \begin{itemize} \item[(i)]
(self-dual) $\pi^\vee \simeq \pi$,\ps \item[(ii)] (conductor $1$) $\pi_p$ is
unramified for each prime $p$, \ps \item[(iii)] (regular half-algebraicity)
the representation ${\rm L}(\pi_\infty)$ is multiplicity free and the
eigenvalues of the infinitesimal character of $\pi_\infty$ are in
$\frac{1}{2}\Z$.  \end{itemize}

Such a $\pi$ necessarily has a trivial central character, hence may be viewed as well as an element of $\Pi_{\rm cusp}(\PGL_n)$. 

\begin{prop}\label{classpialg} The map $\pi \mapsto \pi \otimes |\cdot|^{w(\pi)/2}$ defines a bijection between the set of centered cuspidal automorphic representations of $\GL_n$ over $\Q$ satisfying conditions {\rm (a), (b)} and {\rm (c')} (\S\ref{countingpb},\S\ref{altsointro}) and the set of cuspidal automorphic representations of ${\rm GL}_n$ over $\Q$ satisfying {\rm (i), (ii)} and {\rm (iii)} above. \end{prop}

\begin{pf} Proposition~\ref{paramalg} shows that if $\pi$ satisfies {\rm (a), (b)} and {\rm (c')} then $\pi \otimes |\cdot|^{w(\pi)/2}$ satisfies (i), (ii) and {\rm (iii)}. It also shows that the map of the statement is injective. Assume conversely that $\pi$ satisfies  {\rm (i), (ii)} and {\rm (iii)}. Clozel's purity lemma~\ref{clopl} and condition (iii) imply that ${\rm L}(\pi_\infty)$ is a direct sum of non-isomorphic representations of the form $1,\varepsilon_{\C/\R}$ or ${\rm I}_{w'}$ for $w' >0$. As explained in the proof of corollary~\ref{cororthosymp}, it follows from the existence of $\widetilde{{\rm L}}(\pi_\infty)$ that each of these summands has the same symplectic/orthogonal alternative than $\pi$. Recall that ${\rm I}_{w'}$ preserves a nondegenerate symplectic pairing if and only if $w'$ is odd.
\ps

Assume first that $\pi$ is symplectic. Then the representations $1$, $\varepsilon_{\C/\R}$ and ${\rm I}_{w'}$ for $w' \equiv 0 \bmod 2$ do not occur in ${\rm L}(\pi_\infty)$. In other words, $n$ is even and 
\begin{equation}\label{piinfinisymp}{\rm L}(\pi_\infty) \simeq \bigoplus_{i=1}^{\frac{n}{2}} {\rm I}_{w_i}\end{equation}
for some unique odd positive integers $w_1>\cdots > w_{n/2}$. \ps

Assume now that $\pi$ is orthogonal. If $n$ is odd, we have
\begin{equation}\label{piinfiniorthodd}{\rm L}(\pi_\infty) \simeq \chi \oplus \bigoplus_{i=1}^{\frac{n-1}{2}} {\rm I}_{w_i},\end{equation}
where $\chi \in \{1,\varepsilon_{\C/\R}\}$ and for some unique even positive integers $w_1>\cdots > w_{n/2}$. As $\pi$ has trivial central character, we have $\det({\rm L}(\pi_\infty))=1$, thus $\chi=\varepsilon_{\C/\R}^{\frac{n-1}{2}}$ is uniquely determined. If $n$ is even, and if $0$ is not an eigenvalue of the infinitesimal character of $\pi_\infty$, then
\begin{equation}\label{piinfiniortheven}{\rm L}(\pi_\infty) \simeq  \bigoplus_{i=1}^{\frac{n}{2}} {\rm I}_{w_i},\end{equation}
for some unique even positive integers $w_1>\cdots > w_{n/2}$. If $0$ is an eigenvalue of the infinitesimal character of $\pi_\infty$, it has necessarily multiplicity $2$ and the two characters $1$ and $\varepsilon_{\C/\R}$ occur in ${\rm L}(\pi_\infty)$, so that the above isomorphism still holds for some unique even nonnegative integers $w_1>\cdots > w_{n/2}$, and with $w_{n/2}=0$. Note that $n \equiv 0 \bmod 4$ if $n$ is even, as $\det({\rm L}(\pi_\infty))=1$. \ps

In all these cases, we have defined a sequence of nonnegative integers $w_1  > \cdots > w_{[n/2]}$ having the same parity. The cuspidal automorphic representation $\pi \otimes |\cdot|^{-w_1/2}$ of ${\rm GL}_n$ satisfies (a), (b) and (c'), for the motivic weight $w_1$ and the Hodge weights $w_i$ (and is centered).
\end{pf}

\begin{definition}\label{defpialg} We denote by $\Pi_{\rm alg}^\bot(\PGL_n) \subset \Pi_{\rm cusp}^\bot(\PGL_n)$ the subset of $\pi$ satisfying {\rm (i)}, {\rm (ii)} and {\rm (iii)} above. For $\ast= {\rm o}$ or ${\rm s}$ we also set $\Pi_{\rm alg}^\ast(\PGL_n)=\Pi_{\rm cusp}^\ast(\PGL_n) \cap \Pi_{\rm alg}^\bot(\PGL_n)$ (see Definition~\ref{defpbos}).\ps
\end{definition}

\begin{definition}\label{defpihw}
If $\pi \in \Pi_{\rm alg}^\bot(\PGL_n)$, its {\rm Hodge weights} $$w_1>\cdots>w_{[n/2]}$$ are the Hodge weights of the cuspidal automorphic representation $\pi_0$ of $\GL_n$ over $\Q$ such that $\pi \simeq \pi_0\otimes |\cdot|^{w(\pi_0)/2}$ given by Proposition~\ref{classpialg}. \ps
They are odd if $\pi$ is symplectic and even otherwise. They determine $\pi_\infty$ by the formula~\eqref{piinfinisymp} {\rm if $n$ is even (}resp. by the formula \eqref{piinfiniorthodd} if $n$ is odd{\rm )}.
\end{definition}

\subsection{Arthur's classification : global parameters}\label{arthurparam} Let $G$ be a classical semisimple group
over $\Z$ and let $n=n(G)$. Define $s(G) \in \{\pm 1\}$ by $s(G)=1$ if
$\widehat{G}$ is a special orthogonal group, $-1$ otherwise. Denote by $\Psi_{\rm glob}(G)$ the set of quadruples
$$(k,(n_i),(d_i),(\pi_i))$$
where $1\leq k \leq n$ is an integer, where for each $1\leq i \leq k$ then
$n_i \geq 1$ is an integer and $d_i$ is a divisor of $n_i$, and where
$\pi_i \in \Pi_{\rm cusp}^\bot(\PGL_{n_i/d_i})$, such that :
\begin{itemize}
\item[(i)] $\sum_{i=1}^k n_i = n$, \ps
\item[(ii)] for each $i$, $s(\pi_i)(-1)^{d_i+1}=s(G)$,\ps
\item[(iii)] if $i\neq j$ and $(n_i,d_i)=(n_j,d_j)$ then $\pi_i \neq
\pi_j$.\ps
\end{itemize}

The set $\Psi_{\rm glob}(G)$ only depends on $n(G)$ and $s(G)$. Two elements
$(k,(n_i),(d_i),(\pi_i))$ and $(k',(n'_i),(d'_i),(\pi'_i))$ in $\Psi_{\rm glob}(G)$ are
said {\it equivalent} if $k=k'$ and if there exists $\sigma \in \got{S}_k$ such
that $n'_i=n_{\sigma(i)}$, $d'_i=d_{\sigma(i)}$ and $\pi'_i=\pi_{\sigma(i)}$
for each $i$. An element of $\Psi_{\rm glob}(G)$ will be called
a {\it global Arthur parameter} for $G$. The class $\underline{\psi}$ of
$\psi=(k,(n_i),(d_i),(\pi_i))$ will also be denoted symbolically by 
$$\underline{\psi}=\pi_1[d_1]\oplus \pi_2[d_2] \oplus \cdots \oplus
\pi_k[d_k].$$
In the writting above we shall replace the symbol $\pi_i[d_i]$ by $[d_i]$ if
$n_i=d_i$ (as then $\pi_i$ is the trivial representation), and by $\pi_i$ if
$d_i=1$ and $n_i \neq d_i$. \ps

Let $\psi \in \Psi_{\rm glob}(G)$. Recall that for each integer $d\geq 1$, the
$\C$-group ${\rm SL}_2$ has a unique irreducible $\C$-representation
$\nu_d$ of dimension $d$, namely ${\rm Sym}^{d-1}(\C^2)$. Condition (i) on $\psi$
allows to define a morphism 
$$\rho_\psi : \prod_{i=1}^k \SL_{n_i/d_i} \times \SL_2 \longrightarrow
\SL_n$$
(canonical up to conjugation by ${\rm SL}_n(\C)$) obtained as the direct sum of the representations $\C^{n_i/d_i} \otimes
\nu_{d_i}$. One obtains this way a canonical map 
$$ \rho_\psi : \prod_{i=1}^k \mathcal{X}({\rm SL}(n_i/d_i)) \times
\mathcal{X}({\rm SL}_2) \longrightarrow \mathcal{X}(\SL_n).$$
A specific element of $\mathcal{X}(\SL_2)$ plays an important role in Arthur's theory : it
is the element $e=(e_v)$ defined by $$e_p=\diag(p^{1/2},p^{-1/2})$$
(positive square roots) for each prime $p$, and by
$$e_\infty=\diag(1/2,-1/2).$$ As is well known, $e=c(1)$ where $1 \in \Pi_{\rm
disc}(\PGL_2)$ is the trivial representation.  

\begin{thmconda} (Arthur's classification) \label{thmarthurclass}Let $G$ be any classical semisimple
group over $\Z$ and let $\pi \in \Pi_{\rm disc}(G)$. There is a
$\psi(\pi)=(k,(n_i),(d_i),(\pi_i)) \in
\Psi_{\rm glob}(G)$ unique up to equivalence such that $${\rm
St}(c(\pi))=\rho_\psi(\prod_{i=1}^k c(\pi_i) \times e).$$
\end{thmconda}

When $G$ is a Chevalley group this follows from~\cite[Thm. 1.5.2]{arthur}, otherwise it expected to be part of Arthur's treatment of inner forms of quasi-split classical groups over $\Q$ (see the last chapter {\it loc.cit.} ; note that the special case needed here, namely for pure inner forms, is presumably much simpler because none of the difficulties mentionned by Arthur seems to occurs.). The uniqueness of $\psi(\pi)$ up to equivalence is actually due to Jacquet-Shalika \cite{jasha}. The part of the theorem concerning 
the infinitesimal character is a property of Shelstad's transfert : see~\cite[Lemma 15.1]{shelstad},~\cite[Lemma 25]{mezo}.
\ps

\begin{definition} The global Arthur parameter $\psi(\pi)$ will be called the
global Arthur parameter of $\pi$. 
\end{definition}

For instance if $1_G \in \Pi_{\rm disc}(G)$ denotes the trivial
representation of $G$, then it is well-known that the Arthur parameter of
$1_G$ is $[n(G)]$, unless $\widehat{G} \simeq {\rm SO}_{2m}(\C)$ in which case it is
$[1]\oplus [n(G)-1]$.\ps

Let $\psi=(k,(n_i),(d_i),(\pi_i)) \in \Psi_{\rm glob}(G)$. The associated triple $(k,(n_i,d_i))$, taken up to permutations of the $(n_i,d_i)$, will be called the {\it endoscopic type} of $\psi$.
One usually says that $\psi$ is {\it stable} if $k=1$ and {\it endoscopic} otherwise. The generalized Ramanujan conjecture asserts that each $\pi_i$ is tempered.
We shall thus say that $\psi$ is {\it  tempered} if $d_i=1$ for all $i$. If $\psi=\psi(\pi)$, the Ramanujan conjecture asserts then that $\pi$ is tempered if and only if 
$\psi(\pi)$ is. In some important cases, e.g. the special case where $G(\R)$ is compact, this conjecture is actually known in most cases : see Corollary~\ref{corramcomp}.
We will say that $\pi$ is {\it stable}, {\it endoscopic} or {\it formally tempered} if $\psi(\pi)$ is respectively stable, endoscopic or tempered. 
We will also talk about the endoscopic type of a $\pi$ for the endoscopic type of $\psi(\pi)$.

\ps\ps

Our last task is to explain
Arthur's converse to the theorem above, namely to decide whether a given
$\psi \in \Psi_{\rm glob}(G)$ is in the image of the map
$\pi \mapsto \psi(\pi)$. This is the content of the so-called
{\it Arthur's multiplicity formula}. Our aim until the end of this chapter will be to state certain special cases of this formula. \ps

\subsection{The packet $\Pi(\psi)$ of a $\psi \in \Psi_{{\rm glob}}(G)$}\label{defpipsi}
 Fix $\psi=(k,(n_i),(d_i),(\pi_i)) \in \Psi_{\rm glob}(G)$. If $p$ is a prime number, define $$\Pi_p(\psi)$$
as the set of isomorphism classes of $G(\Z_p)$-spherical (i.e. unramified) irreducible smooth representations of $G(\Q_p)$ whose Satake 
parameter $s_p$, a semisimple conjugacy class in $\widehat{G}$, satisfies $${\rm St}(s_p) = \rho_\psi(\prod_{i=1}^k c_p(\pi_i) \times e_p).$$
This relation uniquely determines the ${\rm Out}(\widehat{G})$-orbit of $s_p$. It follows that  $\Pi_p(\psi)$ is a singleton, unless $\widehat{G} \simeq {\rm SO}(2m,\C)$ and ${\rm St}(s_p)$ does not 
possess the eigenvalue $\pm 1$ (which implies that each $n_i$ is even), in which case it has exactly $2$ elements.
\ps
We shall now associate to $\psi$ a ${\rm Out}(\widehat{G})$-orbit of
equivalence classes of {\it archimedean Arthur parameters} for $G_\R$, which will eventually lead in some cases to a definition of a set $\Pi_\infty(\psi)$ of irreducible unitary representations of $G(\R)$. 
Denote by $\Psi(G_\R)$ the set of such parameters, i.e. of 
continuous homomorphisms 
$$\psi_\R: {\rm W}_\R \times {\rm SL}_2(\C) \longrightarrow \widehat{G}$$
which are $\C$-algebraic on the ${\rm SL}_2(\C)$-factor and such that the image of any element of ${\rm W}_\R$ is semisimple. Two such parameters are said {\it equivalent} if they are conjugate under $\widehat{G}$. An important invariant of an equivalence class of parameters $\psi_\R$ is its {\it infinitesimal character} $$z_{\psi_\R}$$
which is a semisimple conjugacy class in $\mathfrak{\widehat{g}}_\C$ given according to a recipe of Arthur : see e.g.~\S\ref{ajparam} for the general definition. It is also the infinitesimal character of the 
Langlands parameter ${\rm W}_\R \rightarrow \widehat{G}$ associated by Arthur to $\psi_\R$. \ps

We now go back to the global Arthur parameter $\psi$. By assumption (ii) on $\psi$, each space $\C^{n_i/d_i} \otimes \nu_{d_i}$ carries a
natural representation of $\widehat{G^{\pi_i}} \times \SL_2(\C)$ which
preserves a nondegenerate bilinear form unique up to scalars, which is
symmetric if $s(\widehat{G})=1$ and antisymmetric otherwise. One thus
obtains a $\C$-morphism 
$$ r_\psi : \prod_{i=1}^k \widehat{G^{\pi_i}} \times \SL_2(\C)
\longrightarrow \widehat{G}.$$
The collection of $\widetilde{{\rm L}}((\pi_i)_\infty) : {\rm W}_\R \rightarrow
\widehat{G^{\pi_i}}$ (see~\eqref{arthurliftreel}) defines by composition with $r_\psi$ a morphism 
$$\psi_\infty : {\rm W}_\R \times {\rm SL}_2(\C) \rightarrow \widehat{G},$$
which is by definition {\it the archimedean Arthur parameter associated to $\psi$}. 
The ${\rm Out}(\widehat{G})$-orbit of the equivalence class of
$\psi_\infty$ only depends on $\underline{\psi}$. In particular, only the ${\rm Out}(\widehat{G})$-orbit of its 
infinitesimal character is well-defined, with this caveat in mind we shall still denote it by $z_{\psi_\infty}$. By definition we have
\begin{equation}\label{formulazpsiglob} {\rm St}(z_{\psi_\infty}) =\rho_\psi(\prod_{i=1}^k c_\infty(\pi_i) \times e_\infty),\end{equation}
which also determines $z_{\psi_\infty}$ uniquely. 
 \ps \ps

Consider the following two properties of an Arthur parameter $\psi_\R \in \Psi(G_\R)$ : \ps
\begin{itemize}
\item[(a)] $z_{\psi_\R}$ is the infinitesimal character of a finite dimensional, irreducible,
$\C$-representation of $G(\C)$, \ps
\item[(b)] ${\rm St} \circ \psi_\R$ is a multiplicity free
representation of $W_\R \times {\rm SL}_2(\C)$.\ps\ps
\end{itemize}

If $\psi$ satisfies (a) (resp (b)) then so does $\tau \circ \psi$ where
$\tau \in {\rm Aut}(\widehat{G})$. In particular, it makes sense to say that
$\psi_\infty$ satisfies (a) if $\psi \in \Psi_{\rm glob}(G)$. \ps

\begin{definition} Let $\Psi_{\rm alg}(G) \subset \Psi_{\rm glob}(G)$
be the subset of $\psi$ such that $\psi_\infty$ satisfies (a).
\end{definition} \ps

The following lemma is mainly a consequence of Clozel's purity lemma~\ref{clopl}.

\begin{lemme}\label{lemmeclass} Let $\psi=(k,(n_i),(d_i),(\pi_i)) \in \Psi_{\rm alg}(G)$.
If $s(G)=1$, assume that $n(G) \not \equiv 2 \bmod 4$. Then : \begin{itemize}
\item[(i)] $\psi_\infty$ satisfies (b). \ps
\item[(ii)] For each $i=1,\dots,k$, we have $\pi_i \in \Pi_{\rm alg}^\bot(\PGL_{n_i/d_i})$, \ps
\item[(iii)] Each $n_i$ is even, except one of them if $n(G)$ is odd, 
and except perhaps exactly two of them if $s(G)=1$ and $n(G) \equiv 0 \bmod
4$. 
Moreover, if $s(G)=1$ and if $n_i$ is even, then $n_i \equiv 0 \bmod 4$. 
\end{itemize}
\end{lemme}

\begin{pf} Let $X$ be the semisimple conjugacy class ${\rm
St}(z_{\psi_\infty}) \subset {\rm M}_n(\C)$.  Of course, we have $X=-X$. 
Property (a) on $\psi_\infty$ is equivalent to the following properties :\ps
\begin{itemize} 
\item[(a1)] the eigenvalues of $X$ are integers if $s(G)=1$, and half odd integers otherwise, \ps 

\item[(a2)] and these eigenvalues are distinct, except if $s(G)=1$, $n(G) \equiv 0 \bmod 2$, and if $0$ is an
eigenvalue of $X$.  In this exceptional case, that we shall call (E), the
eigenvalue $0$ has multiplicity $2$ and the other eigenvalues have
multiplicity $1$.
\ps 
\end{itemize} 

In particular, assertion (i) follows from
(a2) if (E) does not hold. For each $i$, the eigenvalues of $c_\infty(\pi_i) \otimes {\rm
Sym}^{d_i-1}(e_\infty)$ are among those of $X$. Consider the set $I$ of integers $i \in
\{1,\dots,k\}$ such that $c_\infty(\pi_i)$ does not have distinct
eigenvalues.  It follows that $|I|\leq
1$, and if $i \in I$ then $d_i=1$ and $0$ is the only multiple eigenvalue of
${\rm c}_\infty(\pi_i)$ (and has multiplicity $2$). Fix $1\leq i \leq k$. 
The eigenvalues of $c_\infty(\pi_i)$ are in $\frac{1}{2}\Z$ by (a1), thus it
follows from Clozel's purity lemma~\ref{clopl} that ${\rm L}((\pi_i)_\infty)$ 
is a direct sum of representations of ${\rm W}_\R$ of the form ${\rm I}_w$, $1$ or
$\varepsilon_{\C/\R}$. We have $\pi_i \in \Pi_{\rm alg}^\bot$ unless $i
\in I$ and the two characters occuring in ${\rm L}((\pi_i)_\infty)$ are both $1$ or both $\varepsilon_{\C/\R}$. This proves assertion (ii) when (E) does not hold, in which case
$I=\emptyset$. This also shows that if $I=\{i\}$ then $n_i \equiv 0 \bmod 2$. \ps

Observe that the assertion (iii) of the lemma is obvious if $s(G)=-1$. Indeed, for each
$i$ then $s(\pi_i)=(-1)^{d_i}$, so if $d_i$ is odd then $\pi_i$ is
symplectic and thus $n_i/d_i$ is even. It follows that we may assume from
now on that $s(G)=1$. \ps

Let $J \supset I$ be the set of integers $i \in  
\{1,\dots,k\}$ such that $0$ is an eigenvalue of $c_\infty(\pi_i) \otimes {\rm
Sym}^{d_i-1}(e_\infty)$. Then $1\leq |J| \leq 2$, and $|J|=1$ if $n(G)$ is
odd. Let $i \notin J$, we claim that $n_i \equiv 0 \bmod 4$. Indeed, this is
clear if $d_i$ is even as then $\pi_i$ is symplectic. If $d_i$ is odd, then
$\pi_i \in \Pi_{\rm alg}^{\rm o}(\PGL_{n_i/d_i})$ as $i
\notin I$, and $c_\infty(\pi_i)$ does not contain the eigenvalue $0$. It
follows that $n_i/d_i$ is even, in which case $n_i/d_i \equiv 0 \bmod 4$ by Proposition~\ref{vanishorodd}. In particular, we have the congruence
\begin{equation}\label{impcong} \sum_{j \in J} n_j \equiv n(G) \bmod
4.\end{equation}
This proves assertion (iii) of the lemma. \ps

Assume now that $I \neq \emptyset$ and let $i \in I$. Then $J=I=\{i\}$, and also $n(G) \equiv 0 \bmod
4$ by assumption, so $n_i \equiv 0 \bmod 4$. Of course we have $\det({\rm
L}((\pi_i)_\infty))=1$. But ${\rm L}((\pi_i)_\infty)$ is a direct sum of
$n_i/2-1$ non isomorphic representations of the form ${\rm I}_w$ with $w>0$
and $w$ even, and of two characters $\chi_1$ and $\chi_2$ among $1$ and $\varepsilon_{\C/\R}$.
The congruence $n_i \equiv 0 \bmod 4$ implies that
$\chi_1\chi_2=\varepsilon_{\C/\R}$. In other words, $\chi_1 \neq \chi_2$ and thus $\pi_i \in \Pi_{\rm alg}^\bot(\PGL_{n_i/d_i})$. This ends the proof of assertion (ii) of the lemma. \ps

It only remains to prove (i) in case (E). Observe that 
$\psi_\infty$ does not satisfy (b) if and only if the representation 
${\rm St} \circ \psi_\infty$ of ${\rm W}_\R \times {\rm SL}_2(\C)$ contains
either twice the character $1$ or twice the character $\varepsilon_{\C/\R}$
(with trivial action of ${\rm SL}_2(\C)$). This can only happen if either $J=I$ or $|J|=2$, $I=\emptyset$ and $d_i=1$ for all $i \in J$. In the case $I=J$, we
conclude by the previous paragraph. If $|J|=2$, $n_i$ is odd for each
$i \in J$, and the congruence \eqref{impcong} shows that exactly one of the
two $n_i$, $i \in J$, is congruent to $1$ modulo $4$ (resp. to $3$ modulo
$4$). But $\pi_j \in \Pi_{\rm alg}^{\rm o}(\PGL_{n_j})$ for $j \in J$, thus ${\rm L}((\pi_j)_\infty)$ contains $\varepsilon_{\C/\R}$ if $n_j \equiv 1 \bmod 4$ and $1$ otherwise (see~\eqref{piinfiniorthodd}). \end{pf}

This is the first important motivation for the consideration of the properties
(a) and (b). The second is that if $G(\R)$ is compact, and if $\pi \in
\Pi_{\rm disc}(G)$, then $\psi(\pi)_\infty$ obviously satisfies (a),
as well as (b) because $n(G) \equiv -1,0,1 \bmod 8$. 

\begin{corconda}\label{corramcomp} Assume that $n(G) \not \equiv 2  \bmod 4$ if $s(G)=1$. If
$\pi \in \Pi_{\rm disc}(G)$ is such that $\pi_\infty$ has the infinitesimal
character of a finite dimensional irreducible $\C$-representation of $G(\C)$, and if
$\psi(\pi)=(k,(n_i),(d_i),(\pi_i))$, then $\pi_i \in \Pi_{\rm alg}^\bot(\PGL_{n_i/d_i})$ for each $i$.  \ps In particular, each $\pi_i$
satisfies the Ramanujan conjecture, unless perhaps if $s(G)=s(\pi_i)=d_i=1$, $n(G) \equiv
n_i \equiv 0 \bmod 2$, and ${\rm St}(c_\infty(\pi_i))$ contains the eigenvalue $0$. 
\end{corconda}

\begin{pf} It only remains to justify the statement about Ramanujan conjecture, but this follows from Lemma~\ref{lemmeclass} (ii) and the results of Clozel-Harris-Labesse, Shin and Caraiani recalled in~\S\ref{lstintro}.\end{pf}

We shall exclude from now on the particular case $s(G)=1$ and $n(G) \equiv 2
\bmod 4$, i.e.  we assume that $$\widehat{G} \not \simeq {\rm
SO}(4m+2,\C).$$ We already said that $G(\R)$ is an inner form of a split
group.  As $\widehat{G} \simeq {\rm SO}(4m+2,\C)$, it is also an inner form
of a compact group (this is of course obvious if $G(\R)$ is already
compact).  A parameter $\psi_\R \in \Psi(G_\R)$ satisfying conditions (a)
and (b) above is called {\it an Adams-Johnson parameter} for $G_\R$.  The set of
these parameters is denoted by $$\Psi_{\rm AJ}(G_\R) \subset \Psi(G_\R).$$
We refer to the Appendix~\ref{appendixadamsjohnson} for a general discussion about them, and more precisely to Definition~\ref{defappaj} and the discussion that follows.  For $\psi_\R
\in \Psi_{\rm AJ}(G_\R)$, Adams and Johnson have defined in~\cite{AJ} a
finite set $\Pi(\psi_\R)$ of (cohomological) irreducible unitary
representations of $G(\R)$.  In the notations of this appendix, the group
$G(\R)$ is isomorphic to a group of the form $G_t$ for some $t \in
\mathcal{X}_1(T)$.  Recall that up to inner isomorphisms, $G_t$ only depends
on the $W$-orbit of $t{\rm Z}(G)$.  We fix such an isomorphism between $G(\R)$ and
$G_{[t]}$ and set $\Pi(\psi_\R)=\Pi(\psi,G_{[t]})$.  As ${\rm Aut}(G(\R))
\neq {\rm Int}(G(\R))$ in general, this choice of an isomorphism might be
problematic in principle.  However, a simple case-by-case inspection shows
that for any classical semisimple $\Z$-group $G$ the natural map ${\rm
Out}(G) \rightarrow {\rm Out}(G(\R))$ is surjective, so that this choice
virtually plays no role in the following considerations.  We shall say more
about this when we come to the multiplicity formula.  \ps

Let $\psi \in \Psi_{\rm alg}(G)$. If ${\rm Out}(\widehat{G})=1$, or more generally if the ${\rm Out}(\widehat{G})$-orbit of the equivalence class of $\psi_\infty$ has one element, we set 
 $$\Pi_\infty(\psi)=\Pi(\psi_\infty).$$ 
In the remaining case, we define $\Pi_\infty(\psi)$ as the disjoint union of the two sets $\Pi(\psi_\R)$ where $\psi_\R$ is an equivalence class of parameters in the 
${\rm Out}(\widehat{G})$-orbit of $\psi_\infty$. Recall from~\S\ref{parammap} that the isomorphism $G(\R) \rightarrow G_t$ fixed above furnishes a canonical parameterization map 
$$\tau : \Pi_\infty(\psi) \longrightarrow {\rm Hom}({\rm C}_{\rm \psi_\infty},\C^\times).$$
The presence of ${\rm C}_{\psi_\infty}$ in the target, rather than ${\rm S}_{\psi_\infty}$, follows from the fact that $G(\R)$ is a pure inner form of a split group and from Lemma~\ref{factorrho}.
When the ${\rm Out}(\widehat{G})$-orbit of the equivalence class of $\psi_\infty$ has two elements, say $\psi_1,\psi_2$, there is a canonical way of identifying ${\rm C}_{\psi_1}$ and ${\rm C}_{\psi_2}$, thus it 
is harmless to denote them by the same name ${\rm C}_{\psi_\infty}$.
  \ps

\begin{definition} If $\psi \in \Psi_{\rm alg}(G)$ set $\Pi(\psi)=\{ \pi \in \Pi(G), \, \, \pi_v \in \Pi_v(\psi) \, \, \forall
v\}$. \end{definition}

The first conjecture we are in position to formulate is a comparison between the Arthur packet attached to a $\psi_\R \in \Psi_{\rm AJ}(G_\R)$, as defined in his book~\cite[\S 2.2]{arthur} by twisted endoscopy when $G_\R$ is split, and the packet $\Pi(\psi_\R)$ of Adams and Johnson recalled above (in a slightly weak sense in the case $\widehat{G}={\rm SO}_{2r}(\C)$). It seems widely believed that they indeed coincide, although no proof seems to have been given yet. A first consequence would be the following conjecture. Observe that this conjecture is obvious when $G(\R)$ is compact.

\begin{conj}\label{conjectureun} If $\pi \in \Pi_{\rm disc}(G)$ and if $\pi_\infty$ has the infinitesimal character of a finite dimensional irreducible $\C$-representation of $G(\C)$ then $\pi \in \Pi(\psi(\pi))$.
\end{conj}

\ps\ps

So far we have defined for each $\psi \in \Psi_{\rm alg}(G)$ a set $\Pi(\psi)$ as well as a parameterization map $\tau$ of $\Pi_\infty(\psi)$. This set is e.g. a singleton when $G={\rm SO}_n$ with $n$ odd, and it is finite, 
in bijection with $\Pi_\infty(\psi)$,  if ${\rm Out}(\widehat{G})=1$. Arthur's multiplicity formula is a formula for $m(\pi)$ for each $\pi \in \Pi(\psi)$, at least when ${\rm Out}(\widehat{G})=1$. 
This formula contains a last ingredient that we now study.

\subsection{The character $\varepsilon_\psi$ of ${\rm C}_\psi$}
\label{paragraphepsilon} Consider some $\psi=(k,(n_i),(d_i),(\pi_i)) \in
\Psi_{\rm alg}(G)$ and denote by ${\rm C}_{\psi}$ the centraliser of ${\rm
Im}(r_\psi)$ in $\widehat{G}$.  This is an elementary abelian $2$-group that
we may describe as follows. \ps

Observe that ${\rm St} \circ r_\psi$ is a direct
sum of $k$ non-isomorphic irreducible representations of $\prod_{i=1}^k
\widehat{G^{\pi_i}} \times {\rm SL}_2(\C)$, say $\oplus_{i=1}^k V_i$,
where $V_i$ factors through a representation of $\widehat{G^{\pi_i}} \times
{\rm SL}_2(\C)$ whose dimension is $n_i$. Observe that by Lemma~\ref{lemmeclass} (iii) each $n_i$ is even, except perhaps exactly one or two of
them when $\widehat{G}$ is an orthogonal group. If $1\leq i \leq k$ is such
that $n_i$ is even, there is a unique element $$s_i \in \widehat{G}$$
such that ${\rm St}(s_i)$ acts as $-{\rm Id}$ on $V_i$ and as ${\rm Id}$ on
each $V_j$ with $j \neq i$. Of course, we have $s_i^2=1$ and $s_i \in {\rm
C}_\psi$, and the following lemma is clear.

\begin{lemme} ${\rm C}_\psi$ is generated by ${\rm Z}(\widehat{G})$ and by the elements $s_i$, where $i=1,\dots,k$ is such that $n_i$ is even. \end{lemme} \ps

A first important ingredient in Arthur's multiplicity formula is Arthur's character $$\varepsilon_\psi : {\rm C}_\psi \longrightarrow \{\pm 1\}.$$
It has been defined by Arthur in full generality in~\cite{arthurunipotent}. We shall apply formula (1.5.6) of~\cite{arthur}. By definition, $\varepsilon_\psi$ is trivial on ${\rm Z}(\widehat{G}) \subset {\rm C}_\psi$. 
In the special case here, we thus only have to give the $\varepsilon_\psi(s_i)$. As the representation $\nu_a \otimes \nu_b$ of ${\rm SL}_2(\C)$ has exactly ${\rm Min}(a,b)$ irreducible factors, the formula loc. cit. is thus easily seen to be 
\begin{equation}\label{arthurepsilon} \varepsilon_\psi(s_i)=\prod_{j \neq i} \varepsilon(\pi_i \times
\pi_j)^{{\rm
Min}(d_i,d_j) }\end{equation}
where $\varepsilon(\pi_i \times \pi_j)=\pm 1$ is the sign such that  $${\rm L}(1-s,\pi_i \times \pi_j)=\varepsilon(\pi_i \times \pi_j){\rm L}(s,\pi_i \times \pi_j).$$ 
Here ${\rm L}(s,\pi_i \times \pi_j)$ is the completed ${\rm L}$-function of $\pi_i \times \pi_j$, and the functional equation above is due to Jacquet, Shalika and Piatetski-Shapiro : see~\cite[Ch. 9]{cogdell} for a survey. 
An important result of Arthur asserts that $\varepsilon(\pi_i \times \pi_j)=1$ if $s(\pi_i)s(\pi_j)=1$ \cite[Thm. 1.5.3]{arthur}, so that in the product \eqref{arthurepsilon} we may restrict to the $j$ such that $s(\pi_j) \neq s(\pi_i)$. \ps

The cuspidal automorphic representation $\pi_i$ is unramified at each finite place, and also quite specific at the infinite place : it belongs to $\Pi_{\rm alg}^\bot({\rm PGL}_{n_i/d_i})$ by Lemma~\ref{lemmeclass} (ii)). It follows that one has 
an explicit formula for $\varepsilon(\pi_i \times \pi_j)$ in terms of the Hodge weights of $\pi_i$ and $\pi_j$. The precise recipe is as follows. There is a unique collection of complex numbers 
$$\varepsilon(r) \in \{1,i,-1,-i\}$$
defined for all the isomorphism classes of continuous representations 
$r : {\rm W}_\R \rightarrow {\rm GL}_m(\C)$ which are trivial on $\R_{>0} \subset {\rm W}_\C$, such that : \begin{itemize}
\item[(i)] $\varepsilon(r \oplus r')= \varepsilon(r)\varepsilon(r')$ for all $r,r'$,\ps
\item[(ii)] $\varepsilon({\rm I}_w)=i^{w+1}$ for any integer $w\geq 0$,\ps
\item[(iii)] $\varepsilon(1)=1$. \ps
\end{itemize}
As ${\rm I}_0 \simeq 1 \oplus \varepsilon_{\C/\R}$, it follows that $\varepsilon(\varepsilon_{\C/\R})=i$.
For instance, if $w,w'\geq 0$ are integers, then $$\varepsilon({\rm I}_w \otimes {\rm I}_{w'})=(-1)^{1+{\rm Max}(w,w')},$$ as 
${\rm I}_w \otimes {\rm I}_{w'} \simeq {\rm I}_{w+w'} \oplus {\rm I}_{|w-w'|}$. \ps

If $\pi \in \Pi_{\rm alg}^\bot({\rm PGL}_n)$ and $\pi' \in \Pi_{\rm alg}^\bot({\rm PGL}_{n'})$ 
then both ${\rm L}(\pi_\infty)$ and ${\rm L}(\pi'_\infty)$ are trivial on $\R_{>0}$ (see \S\ref{demialg}), and one has 
\begin{equation}\label{formuleeps2}\varepsilon(\pi \times \pi')=\varepsilon({\rm L}(\pi_\infty) \otimes {\rm L}(\pi'_\infty)).\end{equation} 
See~\cite[\S 4]{tate} (the epsilon factor is computed here with respect to $x \mapsto e^{2i\pi x}$),~\cite[\S 1.3]{arthur}, and Cogdell's lectures~\cite[Ch. 9]{cogdell}. This allows to compute the character $\varepsilon_\psi$ in all cases. See~\cite{cl2} for some explicit formulas.  \ps

We are now able to prove Proposition~\ref{propepsilonintro} of the introduction.

\begin{pf} (of Proposition~\ref{propepsilonintro}) Let $\pi \in \Pi_{\rm alg}^{\rm o}(\PGL_n)$ and consider its global epsilon factor $\varepsilon(\pi):=\varepsilon(\pi \times 1)$. Arthur's result \cite[Thm. 1.5.3]{arthur} ensures that $\varepsilon(\pi)=1$ as $\pi$ is orthogonal. On the other hand, if $w_1> \cdots > w_{[n/2]}$ are the Hodge weights of $\pi$ then the formulas~\eqref{piinfiniorthodd} and \eqref{piinfiniortheven} show that
$$\varepsilon(\pi)=\left\{ \begin{array}{cc} (-1)^{\frac{\sum_{j=1}^{[n/2]} (w_j+1)}{2}} & {\rm if }\, \, \, n \not \equiv 3 \bmod 4, \\   
(-1)^{\frac{1+\sum_{j=1}^{[n/2]} (w_j+1)}{2}} & {\rm otherwise.}\end{array}\right.$$ 
\end{pf}

\subsection{Arthur's multiplicity formula}\label{parmultform} Let $G$ be a classical semisimple group over $\Z$ such that $\widehat{G} \neq {\rm SO}(4m+2,\C)$ and let $\psi=(k,(n_i),(d_i),(\pi_i)) \in \Psi_{\rm alg}(G)$. 
Following Arthur, set 
$$m_\psi=\left\{ \begin{array}{lll} 2 & {\rm if}\, \, \, s(\widehat{G})=1 \, \, \, {\rm and}\, \, \, n_i \equiv 0 \bmod 2\, \, \, {\rm for\,\, all}\, \, \,  1 \leq i \leq k, \\ \\ 1 & {\rm otherwise}.\end{array}\right.$$

Consider the following equivalence relation $\sim$ on $\Pi(G)$. 
The relation $\sim$ is trivial (i.e. equality) unless $\widehat{G}$ is an even orthogonal group, in which case one may assume that $G={\rm SO}_L$ is a standard even orthogonal group.  
Consider the outer automorphism $s$ of the $\Z$-group $G$ induced by the conjugation by any $s_\alpha \in {\rm O}(L)$ as in Remark~\ref{remaut}. If $\pi,\pi' \in \Pi(G)$ we define $\pi \sim \pi'$ if $\pi_v \in \{ \pi'_v, \pi'_v \circ s\}$ for each $v$. \ps

For $\pi \in \Pi(G)$, recall that $m(\pi)$ denotes the multiplicity of $\pi$ in ${\rm L}^2_{\rm disc}(G(\Q)\backslash G(\AAA))$. Recall that we have defined a group ${\rm C}_\psi$ in~\S\ref{paragraphepsilon}, as well as a
group ${\rm C}_{\psi_\infty}$ in~\S\ref{defpipsi}. By definition there is a canonical inclusion $${\rm C}_\psi \subset {\rm C}_{\psi_\infty}.$$

\begin{conj} (Arthur's multiplicity formula)\label{conjecturedeux} Let
$\psi=(k,(n_i),(d_i),(\pi_i)) \in \Psi_{\rm alg}(G)$ and let $\pi \in \Pi(\psi)$. Then
$$\sum_{\pi' \in \Pi(\psi), \pi' \sim \pi} m(\pi') = \left\{\begin{array}{lll} 0 & {\rm if}\, \, \tau(\pi_\infty)_{|{\rm C}_\psi} \neq
\varepsilon_\psi, \\ \\
m_\psi & {\rm otherwise.} \end{array}\right.$$ 
\end{conj}

Observe that $\{ \pi' \in \Pi(\psi), \pi \sim \pi'\}$ is the singleton $\{\pi\}$ unless $\widehat{G}$ is an even orthogonal group. \ps

At the moment this multiplicity formula is still conjectural in the form stated here. However, when $G$ is a Chevalley group, it is a Theorem${}^\ast$ by~\cite[Thm. 1.5.2]{arthur} if we replace the parameterized set $(\Pi_\infty(\psi),\tau)$ above by the one abstractly defined by Arthur \cite[Thm. 1.5.1]{arthur}. Actually, an extra subtlety arises in Arthur's work because the archimedean packets he constructs {\it loc. cit.} are a priori multisets rather than sets. The resulting possible extra multiplicities have been neglected here to simplify the exposition, as they are actually expected not to occur according to Arthur (and even more so for the Adams-Johnson packets). Note also that Arthur's formula even holds for all the global parameters $\psi \in \Psi_{\rm glob}(G)$. The case of a general $G$ has also been announced by Arthur : see Chap. 9 loc. cit.  \ps

In this paper, we shall use this conjecture only in the following list of special cases. In each case we will
explicit completely the multiplicity formula in terms of the Hodge weights of the $\pi_i$ appearing in $\psi$. We have already done so for the term $\varepsilon_\psi$ in the previous paragraph. In each case we also discuss the dependence of the multiplicity formula on the choice of the identification of $G(\R)$ that we have fixed in~\S\ref{defpipsi} to define $\tau(\pi_\infty)$. 
\ps\ps

\subsubsection{The definite odd orthogonal group $G={\rm SO}_{2r+1}$} ${}^{}$ \label{mformoddorth}\ps\ps

 In this case $r \equiv  0,3\bmod 4$ and $\widehat{G}={\rm Sp}_{2r}(\C)$. Consider the standard based root datum for $(\widehat{G},\widehat{B},\widehat{T})$ with ${\rm X}^\ast(\widehat{T})=\Z^r$ with canonical basis $(e_i)$ and 
 $$\Phi^+(\widehat{G},\widehat{T})=\{2e_i, 1\leq i \leq r\} \cup \{e_i \pm e_j, 1\leq i < j\leq r\}.$$ We conjugate $r_\psi$ in $\widehat{G}$ so that the centralizer of 
$\varphi_{\psi_\infty}({\rm W}_\C)$ is $\widehat{T}$, and that $\varphi_{\psi_\infty}(z)=z^\lambda \overline{z}^{\lambda'}$ with $\lambda \in \frac{1}{2}{\rm X}_\ast(\widehat{T})$ dominant with respect to $\widehat{B}$. \ps

There is a unique element in $\Pi_\infty(\psi)$, namely the irreducible representation with infinitesimal character 
$z_{\psi_\infty}$. The character $\tau(\pi_\infty)$ is absolutely canonical here as each automorphism of $G(\R)$ is inner and there is a unique choice of strong real form $t$ for $G_t$ (namely $t=1$). By Cor.~\ref{corparamcompact}, this character $\tau(\pi_\infty)$ is $(\rho^\vee)_{|{\rm C}_{\psi_\infty}}$, where $\rho^\vee$ denotes the half-sum of the positive roots of $(\widehat{G},\widehat{B},\widehat{T})$, namely 
$\rho^\vee=r e_1 + (r-1) e_2 + \cdots + e_r$. In particular $\rho^\vee \in {\rm X}^\ast(\widehat{T})$ and it satisfies the congruence $$\rho^\vee \equiv e_r + e_{r-2} + e_{r-4} + \cdots \, \, \, \, \, \bmod 2{\rm X}^\ast(\widehat{T}).$$
Observe that $\rho^\vee(-1)=1$ as $r\equiv 0,3 \bmod 4$, so that $\rho^\vee$ is trivial on ${\rm Z}(\widehat{G})$. \ps

Consider the generators $s_i$ of ${\rm C}_\psi$ introduced in~\S\ref{paragraphepsilon}. We shall now give an explicit formula for the $\rho^\vee(s_i)$. 
Fix some $i \in \{1,\dots,k\}$ and write $n_i=r_id_i$. Assume first that $d_i$
and $r_i$ are even. Then $0$ is not a Hodge weight of $\pi_i$, as otherwise
$z_{\psi_\infty}$ would have twice the eigenvalue $\frac{1}{2}$. The positive
eigenvalues of $z_{\psi_\infty}$ associated to the summand $\pi_i[d_i]$ of
$\psi$ are thus the union of the $d_i$ consecutive
half-integers $$\frac{w_j+d_i-1}{2}, \frac{w_j+d_i-3}{2},
\cdots , \frac{w_j+1-d_i}{2}$$
where $w_j$ runs among the $\frac{r_i}{2}$ Hodge weights of $\pi_i$. It follows that
$$\rho^\vee(s_i)=(-1)^{\frac{d_i}{2}\frac{r_i}{2}}=(-1)^{\frac{n_i}{4}}.$$
If $d_i$ is even and $r_i$ is odd, the positive eigenvalues of 
$z_{\psi_\infty}$ coming from the summand $\pi_i[d_i]$ are of the form above,
plus the $\frac{d_i}{2}$ consecutive half-integers $\frac{d_i-1}{2},\cdots,\frac{3}{2},\frac{1}{2}$. One rather obtains  
$$\rho^\vee(s_i)=\left\{
\begin{array}{ll} -(-1)^{[\frac{r_i}{2}]\frac{d_i}{2}} & {\rm if}\, \,
\frac{d_i}{2} \equiv 1,2 \bmod 4,\\(-1)^{[\frac{r_i}{2}]\frac{d_i}{2}} & {\rm otherwise}.
\end{array}\right.$$
If $d_i$ is odd, in which case $r_i$ is even, the sign $\rho^\vee(s_i)$ depends on the Hodge weights of $\pi_i$. Precisely, denote by  $$w_1 > \cdots > w_r $$ the positive odd integers $w_j$ 
such that the eigenvalues of ${\rm St}(z_{\psi_\infty})$ in ${\rm SL}_{2r}(\C)$ are the $\pm
\frac{w_j}{2}$ (see formula~\eqref{formulazpsiglob}). There is a unique subset $J \subset \{1,\cdots,r\}$ such that the Hodge weights of $\pi_i$ are the $w_j$ for $j \in J$. Denote by $J'$ the subset of $j \in J$ such that $j \equiv r \bmod 2$. It is then clear that 
$$\rho^\vee(s_i)=(-1)^{|J'|}.$$
Although these formulas are explicit, we do not especially recommend to use them in
a given particular case, as usually the determination of ${\rho^\vee}(s_i)$ is pretty immediate by definition from the 
inspection of $\psi$ ! \ps

\subsubsection{The definite even orthogonal group $G={\rm SO}_{2r}$} ${}^{}$ \label{mformevenorth}\ps\ps
 In this case $r \equiv  0\bmod 4$ and $\widehat{G}={\rm SO}_{2r}(\C)$. Consider the standard based root datum for $(\widehat{G},\widehat{B},\widehat{T})$ with ${\rm X}^\ast(\widehat{T})=\Z^r$ with canonical basis $(e_i)$ and 
 $$\Phi^+(\widehat{G},\widehat{T})=\{e_i \pm e_j, 1\leq i < j\leq r\}.$$ We conjugate $r_\psi$ in $\widehat{G}$ as in the odd orthogonal case. \ps 
		 
If the ${\rm Out}(\widehat{G})$-orbit of $\psi_\infty$ consists of only one equivalence class, in which case the ${\rm Out}(\widehat{G})$-orbit of $z_{\psi_\infty}$ is a singleton, then the unique element of $\Pi_\infty(\psi)$
is the representation of $G(\R)$ with infinitesimal character $z_{\psi_\infty}$.  Otherwise, the two elements of $\Pi_\infty(\psi)$, again two finite dimensional irreducible 
representations, have the property that their infinitesimal characters are exchanged by the outer automorphism of $G(\R)$, and both in the ${\rm Out}(\widehat{G})$-orbit of $z_{\psi_\infty}$. Observe that there is still
the possibility that the ${\rm Out}(\widehat{G})$-orbit of $z_{\psi_\infty}$ is a singleton : in this case $\Pi_\infty(\pi)$ consists of two isomorphic representations. However, observe also that by definition all the members of $\Pi(\psi)\subset \Pi(G)$ have the same archimedean component in this case. \ps

Recall we have fixed an isomorphism between $G(\R)$ and $G_t$ for $t=\{\pm 1\} \in {\rm Z}(G)$ as in~\S\ref{strongforms}. Assume first that we actually chosed $t=1$. 
It follows that the one or two elements in $\Pi_\infty(\psi)$ have the same character $\rho^\vee$ by Cor.~\ref{corparamcompact}. Here we have 
$\rho^\vee=(r-1)e_1+(r-2)e_2+\cdots + e_{r-1}$, thus $\rho^\vee \in {\rm X}^\ast(\widehat{T})$ and 
$$\rho^\vee \equiv e_{r-1} + e_{r-3} + e_{r-5} + \cdots  \bmod 2{\rm X}^\ast(\widehat{T}).$$
Observe again that $\rho^\vee(-1)=1$ as $r\equiv 0 \bmod 4$. \ps

Consider the generators $s_i$ of ${\rm C}_\psi$ introduced in~\S\ref{paragraphepsilon}. Fix some $i \in \{1,\dots,k\}$ and write $n_i=r_id_i$. If $d_i$ is even, then $r_i$ is even as well
as $s(\pi_i)=-1$, and we have
			$$\rho^\vee(s_i)=(-1)^{\frac{n_i}{4}}.$$
If $d_i$ is odd, in which case $r_i$ is even as $n_i=d_ir_i$ is even by assumption, the sign $\rho^\vee(s_i)$ depends on the Hodge weights of $\pi_i$. Precisely, denote by  $$w_1 > \cdots > w_r $$ the 
nonnegative even integers $w_j$ such that the eigenvalues of ${\rm St}(z_{\psi_\infty})$ 
in ${\rm SL}_{2r}(\C)$ are the $\pm \frac{w_j}{2}$ (see formula~\eqref{formulazpsiglob}). There is a unique subset $J \subset \{1,\cdots,r\}$ such that the Hodge weights of $\pi_i$ are the $w_j$ for $j \in J$. Denote by $J'$ the subset of $j \in J$ such that $j \equiv {r-1} \bmod 2$. It is then clear that 
$$\rho^\vee(s_i)=(-1)^{|J'|}.$$
\ps
For coherence reasons, we shall check now that the multiplicity formula does not change if we choose to identify $G(\R)$ with $G_{-1}$ or if we modify the fixed isomorphism by the outer automorphism of $G(\R)$. 
This second fact is actually trivial by what we already said, so assume that we identified $G(\R)$ with $G_{-1}$. The effect of this choice is that the one or two elements of $\Pi_\infty(\psi)$ become parameterized by the character $$\rho^\vee+\chi,$$
where $\chi$ is the generator of the group $\mathcal{N}(T)$, by Lemma~\ref{lemmeactz}. As $-1=e^{i\pi\chi}$ we have
$$\chi \equiv \sum_{i=1}^r e_i \bmod 2{\rm X}_\ast(\widehat{T})$$
and we claim that this character is trivial on ${\rm C}_\psi$. Indeed, it follows from Lemma~\ref{lemmeclass} that if $n_i$ is even then $n_i \equiv 0 \bmod 4$, so that $\chi(s_i)=(-1)^{n_i/2}=1$.

\subsubsection{The Chevalley groups ${\rm Sp}_{2g}$, ${\rm SO}_{2,2}$ and ${\rm SO}_{3,2}$} \label{parmultformsp4}${}^{}$ \ps 
	The case of the symplectic groups ${\rm Sp}_{2g}$ will be treated in details in Chapter~\ref{finalapp}, especially in~\S\ref{parholdisc}. We shall only consider there the multiplicity formula for a $\pi$ such that $\pi_\infty$ is a holomorphic discrete series. \ps 
	The cases $G={\rm SO}_{2,2}$ and ${\rm SO}_{3,2}$ will be used in Chapter~\ref{chapsmall}. For ${\rm SO}_{2,2}$ we shall not use that Arthur's packets are the same as the ones of Adams-Johnson. 
For $G={\rm SO}_{3,2}$ we shall need it only in~\S\ref{parsp4}, i.e. to compute ${\rm S}(w,v)$,  for the $\psi \in \Psi_{\rm alg}(G)$  of the form $\pi\oplus [2]$. In this case this is probably not too difficult to check but due to the already substantial length of this paper we decided not to include this twisted character computation here. We hope to do so in the future.

\newpage

\section{Determination of $\Pi_{\rm alg}^\bot(\PGL_n)$ for $n\leq 5$}\label{chapsmall}

In this chapter we justify the formulas for ${\rm S}(w)$ and ${\rm S}(w,v)$ given in the introduction and prove Theorem~\ref{proporthointro} there. We recall that various sets $$\Pi_{\rm alg}^\ast(\PGL_n) \subset \Pi_{\rm cusp}^\ast(\PGL_n) \subset \Pi_{\rm cusp}(\PGL_n)$$ have been introduced in Definitions~\ref{defpbos} and \ref{defpialg}.

\subsection{Determination of $\Pi^\bot_{\rm cusp}(\PGL_2)$} A representation $\pi \in
\Pi_{\rm cusp}(\PGL_2)$ is necessarily self-dual as $g \mapsto {}^tg^{-1}$ is an inner automorphism of $\PGL_2$. It is even symplectic by
Theorem~\ref{orthosymp}, so that $$\Pi_{\rm cusp}(\PGL_2)=\Pi^\bot_{\rm cusp}(\PGL_2)=\Pi_{\rm cusp}^{\rm s}(\PGL_2).$$
If $\pi \in \Pi_{\rm cusp}(\PGL_2)$, the infinitesimal character of $\pi_\infty$ has the form $\diag(\frac{w}{2},-\frac{w}{2}) \in \mathfrak{sl}_2(\C)$ for some integer $w\geq 1$ if and only if $\pi_\infty$ is a discrete series representation, in which case $w$ is odd and determines $\pi_\infty$ (see e.g. \S\ref{demialg} and ~\cite{knapp}, or~\cite[\S 2]{bump}).  \ps

	Let $w\geq 1$ be an odd integer and let $\mathcal{F}_w$ be
the set of $$F=\sum_{m\geq 1}a_m q^m  \in {\rm S}_{w+1}({\rm SL}_2(\Z))$$
which are eigenforms for all the Hecke operators and normalized so that $a_1=1$ : see~\cite{serre}. As is
well-known, and explained by Serre, $\mathcal{F}_w$ is a basis of the complex vector space ${\rm S}_{w+1}({\rm
SL}_2(\Z))$. Moreover, each $F \in \mathcal{F}_w$ generates a $\pi_F \in \Pi_{\rm cusp}(\PGL_2)$, and the map $F \mapsto \pi_F$ is a bijection
between $\mathcal{F}_w$ and the set of $\pi$ in $\Pi_{\rm alg}(\PGL_2)$ such that $\pi_\infty$ has Hodge weight $w$ (see~\cite[\S 3.2]{bump}). In particular $${\rm S}(w)=\dim({\rm S}_{w+1}({\rm SL}_2(\Z)))$$ 
as recalled in the introduction. We shall always identify an $F \in \mathcal{F}_w$ with 
$\pi_F$ in the bijection above, and even write $F \in \Pi_{\rm
alg}(\PGL_2)$. For $w \in \{11, 13, 15, 17, 19, 21\}$ we shall denote by $$\Delta_w \in \Pi_{\rm alg}(\PGL_2)$$
the unique element with Hodge weight $w$, as a reminiscence of the notation $\Delta$ for Jacobi's discriminant function, i.e. $\Delta=\Delta_{11}$.

\subsection{Determination of $\Pi_{\rm alg}^{\rm s}(\PGL_4)$ } ${}^{}$ \ps \label{parsp4}  Fix $w>v$ odd positive integers. Let ${\rm S}_{w,v}({\rm Sp}_4(\Z))$ be the space of Siegel cusp forms of genus $2$ recalled in~\S\ref{casebycaselow} of the introduction. Denote also by $$\Pi_{w,v}({\rm PGSp}_4) \subset \Pi_{\rm cusp}({\rm PGSp}_4)$$
the subset of $\pi \in \Pi_{\rm cusp}({\rm PGSp}_4)$ such that $\pi_\infty$ is the holomorphic discrete series whose infinitesimal character has the 
eigenvalues $\pm \frac{w}{2},\pm \frac{v}{2}$, viewed as a semisimple conjugacy class in ${\mathfrak{sl}}_4(\C)$. It is well-known that to 
each Hecke-eigenform $F$ in  ${\rm S}_{w,v}({\rm Sp}_4(\Z))$ one may associate a unique $\pi_F \in \Pi_{w,v}({\rm PGSp}_4)$, and that the image of the map $F \mapsto \pi_F$ is $\Pi_{w,v}({\rm PGSp}_4)$
 (see e.g.~\cite{asgarischmidt}). \ps

 The semisimple $\Z$-group ${\rm PGSp}_4$ 
is isomorphic to ${\rm SO}(3,2)$ hence we may view it 
as a classical semisimple group over $\Z$. It follows from Arthur's multiplicity formula (\S\ref{parmultform},\S\ref{parmultformsp4}) that the multiplicity of any such $\pi_F$ as above is $1$, so that the Hecke-eigenspace containing a given Hecke-eigenform 
 $F$ is actually not bigger than $\C F$. It follows that if we denote by $\mathcal{F}_{w,v}$ the set of these (one dimensional) Hecke-eigenspaces in ${\rm S}_{w,v}({\rm Sp}_4(\Z))$, then
$$|\mathcal{F}_{w,v}|=\dim {\rm S}_{w,v}({\rm Sp}_4(\Z))=|\Pi_{w,v}({\rm PGSp}_4)|.$$

The following formula was claimed in the introduction.

\begin{propcondb}  For $w>v>0$ odd, ${\rm S}(w,v) = {\rm S}_{w,v}({\rm Sp}_4(\Z)) - \delta_{v=1}\delta_{w \equiv 1 \bmod 4} {\rm S}(w)$. \end{propcondb}


Before starting the proof, recall that if $\varphi$ is a discrete series Langlands parameter for ${\rm PGSp}_4(\R)$, its ${\rm L}$-packet $\Pi(\varphi)$ has two elements 
$\{\pi_{\rm hol},\pi_{\rm gen}\}$ where $\pi_{\rm gen}$ is generic and $\pi_{\rm hol}$ is holomorphic. 
One has moreover
$${\rm C}_{\varphi}={\rm S}_{\varphi} \simeq (\Z/2\Z)^2$$
in the notation of~\S\ref{ajpackets}, and the two Shelstad characters of ${\rm C}_\varphi$ associated to the elements of $\Pi(\varphi)$ are the ones which are trivial 
on the center $Z$ of ${\rm Sp}_4(\C)$. Of course $\tau(\pi_{\rm gen})=1$ 
and so $\tau(\pi_{\rm hol})$ is the unique non-trivial character of ${\rm C}_\varphi$ which is trivial on the center $Z=\{\pm1\}$ of ${\rm Sp}_4(\C)$. \ps

Fix a $\psi \in \Psi_{\rm alg}({\rm PGSp}_4)$ whose infinitesimal character has the eigenvalues 
$\pm \frac{w}{2}, \pm \frac{v}{2}$. One has to determine if $\Pi_\infty(\psi)$ 
contains the holomorphic discrete series and, if it is so, to determine the multiplicity of the unique $\pi \in \Pi(\psi)$ such 
that $\pi_\infty$ is this holomorphic discrete series. Such a $\pi$ is necessarily cuspidal as $\pi_\infty$ is tempered, by a result of Wallach~\cite[Thm. 4.3]{wallach} (as pointed out to us by Wallach, this discrete series case is actually significantly simpler than the general case treated there). We proceed by a case by case argument depending on the global Arthur parameter $\psi$ : \ps\ps

{\bf Case (i)} : (stable tempered case)  $\psi=\pi_1$ where $\pi_1 \in \Pi_{\rm alg}^{\rm s}({\rm PGL}_4)$. In this case $\psi_\infty$ is a discrete series Langlands parameter.
It follows from Arthur's multiplicity formula that $m(\pi)=1$, as ${\rm C}_\psi=Z$. The number of such $\pi$ is the number ${\rm S}(w,v)$ that we want to compute. \ps\ps

{\bf Case (ii)} : $\psi=[4]$. The unique $\pi \in \Pi_{\rm disc}({\rm PGSp}_4)$ with $\psi(\pi)=\psi$ is the trivial representation, for which $\pi_\infty$ is not a discrete series.  \ps \ps

{\bf Case (iii)} : $\psi=\pi_1 \oplus \pi_2$ where $\pi_1,\pi_2 \in \Pi_{\rm alg}(\PGL_2)$ and $\pi_1, \pi_2$ have different Hodge weights. In this case one has 
$${\rm C}_{\psi}={\rm C}_{\psi_{\infty}}=(\Z/2\Z)^2.$$ Moreover, $r_\psi({\rm SL}_2(\C))=1$ so $\varepsilon_\psi$ is trivial and $\psi_\infty$ is a discrete series parameter for ${\rm PGSp}_4(\R)$. If $\pi \in \Pi(\psi)$ is the unique element such that $\pi_\infty$ is holomorphic, Arthur's multiplicity formula thus shows that $m(\pi)=0$ as $\varepsilon_\psi$ is trivial but $\tau(\pi_\infty)$ is not.  \ps\ps

{\bf Case (iv)} :  $\psi=\pi_1 \oplus [2]$ where $\pi_1 \in \Pi_{\rm alg}(\PGL_2)$ with Hodge weight $w\neq 1$ (which is actually automatic as ${\rm S}(1)=0$). Again one has
$${\rm C}_\psi={\rm C}_{\psi_{\infty}}=(\Z/2\Z)^2.$$ 
This time $r_\psi({\rm SL}_2(\C)) \neq 1$, and if $s=s_1$ is the generator of ${\rm C}_\psi/Z$,  then $$\varepsilon_\psi(s)=\varepsilon(\pi_1 \times 1)=\varepsilon(\pi_1)=(-1)^{(w+1)/2}.$$ 
\ps
The Adams-Johnson parameter $\psi_\infty$ has an associated complex Levi subgroup $L$ isomorphic to ${\rm SO}_2(\C) \times {\rm SO}_3(\C)$ (see~\S\ref{ajparam} and \S\ref{ajpackets}). It follows that the set $\Pi_\infty(\psi)$, which has two elements, contains the holomorphic discrete series (associated to the order $2$ element in the center of $L$).  For more details, see Chapter~\ref{finalapp} where the general case ${\rm Sp}_{2g}(\R)$ will be studied. The character of this holomorphic discrete series relative to this $\psi_\infty$ is again the non-trivial character of ${\rm C}_\psi$ trivial on $Z$ by the discrete series case recalled above and Lemma~\ref{reductionsd}. It follows that if $\pi \in \Pi(\psi)$ is the unique element such that $\pi_\infty=\pi_{\rm hol}$, then by Arthur's multiplicity formula we have 
$m(\pi)=0$ if $w \equiv 3 \bmod 4$, and $m(\pi)=1$ if $w \equiv 1 \bmod 4$. \ps

This concludes the proof of the proposition. $\square$

\begin{rem}{\rm By the formula for ${\rm S}(w)$, the first $w$ for which a $\pi$ as in case (iv) exists is for $w=17$, for which $\psi(\pi)=\Delta_{17} \oplus [2]$. The representations $\pi$ occuring in case (iv) have a long history, their existence had been conjectured by Saito and Kurokawa in 1977, and proved independently of this theory by Maass, Andrianov and Zagier. We refer to Arthur's paper~\cite{arthursp4} for a discussion about this (and most of the discussion of this paragraph).}\end{rem} \ps\ps

When ${\rm S}(w,v)=1$ we shall denote by $\Delta_{w,v}$ the unique element of $\pi \in \Pi_{w,v}({\rm PGSp}_4)$ such that $\psi(\pi) \in \Pi_{\rm cusp}(\PGL_4)$. As recalled in~\S\ref{casebycaselow}, an explicit formula for 
$\dim \,{\rm S}_{w,v}({\rm Sp}_4(\Z))$ has been given by T. Tsushima (and by Igusa when $v=1$). See Table~\ref{tableSwv} for a sample of values. For $w<25$, one observes that ${\rm S}(w,v)$ is either $0$ or $1$. For those $w<25$, there are exactly $7$ forms $\Delta_{w,v}$, for the following values $(w,v)$ :
$$(19,7), (21,5), (21,9), (21,13), (23,7), (23,9), (23,13).$$

\medskip
	Contrary to the $\PGL_2$ case where one has simple formulas for the $c_p(\Delta_w)$ thanks to the $q$-expansion of Eisenstein series or the product formula for $\Delta_{11}$, 
much less seems to be known at the moment for the $c_p(\pi)$ where $\pi \in \Pi_{w,v}({\rm PGSp}_4)$, even (say) for $\pi=\Delta_{w,v}$ and $(w,v)$ in the list above. We refer to the recent work~\cite{rrst} 
for a survey on this important problem, as well as some implementation on SAGE. \ps
To cite a few results especially relevant to our purposes here, let us mention first the work of Skoruppa \cite{sko} computing $c_p(\pi)$ for the first $22$ primes $p$ when $\pi$ is any of the $18$ elements in the $\Pi_{w,1}({\rm PGSp}_4)$ for $w \leq 61$. Moreover, works of Faber and Van der Geer (see~\cite[\S 24, \S 25]{vdg}) compute the trace of $c_p(\Delta_{v,w})$ in the standard $4$-dimensional representations 
when $p\leq 37$, and even $c_p(\Delta_{w,v})$ itself when $p\leq 7$, whenever $(w,v)$ is in the list above. In the work~\cite{cl2} of the first author and Lannes, 
the first $4$ of these forms, namely $\Delta_{19,7}$, $\Delta_{21,5}$, $\Delta_{21,9}$ and $\Delta_{21,13}$, appeared in the study of the Kneser $p$-neighbours of the Niemeier lattices. Properties of the Leech lattice also allowed those authors to compute ${\rm Trace}({\rm c}_p(\Delta_{w,v}))$ for those $4$ pairs $(w,v)$ up to $p \leq 79$. \ps

\subsection{An elementary lifting result for isogenies}\label{lzscon}  Consider $\iota : G \rightarrow G'$ a central isogeny between
semisimple Chevalley groups over $\Z$. The morphism $\iota$ is thus a finite flat group scheme homomorphism, 
$Z={\rm Ker} \,\iota \subset {\rm Z}(G)$ is a central multiplicative $\Z$-group scheme and $G'=G/Z$. The following proposition is easy to observe for all the isogenies we shall consider later, but it is perhaps more satisfactory to give a general proof.

\begin{prop}\label{propbijiso} $\iota$ induces a homeomorphism $G(\Q) \backslash G(\AAA)/G(\widehat{\Z}) \isomo G'(\Q) \backslash G'(\AAA)/G'(\widehat{\Z})$. 
\end{prop}

\begin{pf} By Prop.~\ref{hchev}, it is enough to check that the map $$G(\Z)\backslash G(\R) \rightarrow G'(\Z)\backslash G'(\R)$$
induced by $\iota$ is a homeomorphism. As this map is continuous and open it is enough to show it is bijective. As the source and target are connected by Prop.~\ref{hchev}, it is surjective. 
Moreover, it is injective if and only if the inverse image of $G'(\Z)$ in $G(\R)$ coincides with $G(\Z)$, what we check now. The fppf exact sequence defined by $\iota$ leads to the following commutative diagram : 
$$\xymatrix{ 1 \ar@{->}[r]  & {\rm Z}(\R) \ar@{->}[r] &  G(\R) \ar@{->}[r] & G'(\R) \ar@{->}[r]  & H^1(\R,Z) \\
 1 \ar@{->}[r] & {\rm Z}(\Z) \ar@{->}[r]  \ar@{->}[u] &  G(\Z) \ar@{->}[r] \ar@{->}[u] & G'(\Z) \ar@{->}[r]  \ar@{->}[u] & H^1(\Z,Z) \ar@{->}[u] }
$$
The left vertical map is an isomorphism by Prop.~\ref{hchev}. The right vertical one is an isomorphism as well, as so are the natural maps 
$$\Z^\times/(\Z^\times)^n = H^1(\Z,\mu_n) \rightarrow H^1(\R,\mu_n) = \R^\times/(\R^\times)^n$$
for each integer $n\geq 1$. A simple diagram chasing concludes the proof.
\end{pf}

Denote by $\iota^\vee : \widehat{G'} \rightarrow \widehat{G}$ the isogeny dual to $\iota$. We now define a map\footnote{We denote by $\mathcal{P}(X)$ the set of all subsets of $X$.} $$\mathcal{R}_\iota : \Pi(G') \longrightarrow \mathcal{P}(\Pi(G))$$
associated to $\iota$ as follows. If $\pi' =\pi'_\infty \otimes \pi'_f \in \Pi(G')$ we define $\mathcal{R}_\iota(\pi')$ as the set of representations $\pi \in \Pi(G)$ such that : \begin{itemize}
\item[(i)] For each prime $p$ the Satake parameter of $\pi_p$ is $\iota^\vee(c_p(\pi'))$,\ps
\item[(ii)] $\pi_\infty$ is a constituent of the restriction to $G(\R) \rightarrow G'(\R)$ of $\pi'_\infty$.\ps
\end{itemize}
Let $\pi \in \mathcal{R}_\iota(\pi')$. Observe that $\pi_p$ is uniquely determined by (i). Moreover the restriction of $\pi'_\infty$ to $G(\R)$ is a direct sum of finitely many irreducible
representations of same infinitesimal character as $\pi'_\infty$. In particular $\mathcal{R}_\iota(\pi')$ is a finite nonempty set. 
We denote by $[\pi_\infty : \pi'_\infty]$ the multiplicity of $\pi_\infty$ in $(\pi'_\infty)_{|G(\R)}$. If $\pi \in \Pi(H)$ we also write $m_H(\pi)$ for $m(\pi)$ to emphasize the $\Z$-group $H$ (see~\S\ref{discautrep}).

\begin{prop}\label{propisog} If $\pi \in \Pi(G)$ then $$m_G(\pi)=\sum_{\{\pi' \in \Pi(G')\, \, |\, \,  \pi \in \mathcal{R}_\iota(\pi')\}} m_{G'}(\pi') [\pi_\infty,\pi'_\infty].$$ 
In particular, the two following properties hold : \begin{itemize}
\item[(a)] For any $\pi \in \Pi_{\rm disc}(G)$ there exists $\pi' \in \Pi_{\rm disc}(G')$ such that $\pi \in \mathcal{R}_\iota(\pi')$. \ps
\item[(b)] For any $\pi' \in \Pi_{\rm disc}(G')$ then $\mathcal{R}_\iota(\pi') \subset \Pi_{\rm disc}(G)$.
\end{itemize}
\end{prop}

\medskip
Before giving the proof we need to recall certain properties of the Satake isomorphism. Following Satake, consider the $\C$-linear map 
$$\iota^\ast : \mathcal{H}(G) \rightarrow \mathcal{H}(G')$$ 
sending the characteristic function of $G(\widehat{\Z})gG(\widehat{\Z})$ to the one of $G'(\widehat{\Z})\iota(g)G'(\widehat{\Z})$. It follows from~\cite[Prop. 7.1]{satake}, that $\iota^\ast$ is a ring homomorphism.
Indeed, it is enough to check the assumptions there. Let $\iota_p$ be the morphism $G(\Q_p) \rightarrow G'(\Q_p)$ induced by $\iota$. Then $\iota_p(G(\Q_p))$ is a normal open subgroup 
of $G'(\Q_p)$. Moreover $\iota_p^{-1}(G'(\Z_p))=G(\Z_p)$ as this latter group is a maximal compact subgroup of $G(\Q_p)$ by \cite{titscorvallis} and $\iota_p$ is proper. Last but not least, the Cartan decomposition shows that 
$\iota_p$ induces an injection $G(\Z_p)\backslash G(\Q_p)/G(\Z_p) \rightarrow G'(\Z_p) \backslash G'(\Q_p)/G'(\Z_p)$ (see e.g. \cite{grossatake}). \ps

If $V$ is a representation of $G'(\AAA_f)$, it defines by restriction by $\iota$ a representation $V_\iota$ of $G(\AAA_f)$ as well, and $V^{G'(\widehat{\Z})} \subset V_\iota^{G(\widehat{\Z})}$. The following lemma is presumably well-known.

\begin{lemme}\label{lemmasat} Let $V$ be a complex representation of $G'(\AAA_f)$ and let $T \in \mathcal{H}(G)$. The diagram 
$$\xymatrix{ V^{G'(\AAA_f)}  \ar@{->}[d]_{\iota^\ast(T)} \ar@{^{(}->}[r] &  V_\iota^{G(\AAA_f)}  \ar@{->}[d]^{T}  \\
  V^{G'(\AAA_f)}  \ar@{^{(}->}[r] &   V_\iota^{G(\AAA_f)} }$$
is commutative.
\end{lemme}

\begin{pf} We have to show that if $\psi : G'(\AAA_f) \rightarrow \C$ is a locally constant function which is right $G'(\widehat{\Z})$-invariant and with support in  $\iota(G(\AAA_f))G'(\widehat{\Z})$, 
then $$\int_{G'(\AAA_f)} \psi(g) dg = \int_{G(\AAA_f)} \psi(\iota(h))dh.$$ Here the Haar measures $dg$ and $dh$ on $G'(\AAA_f)$ and $G(\AAA_f)$ are normalized so that $G'(\widehat{\Z})$ and $G(\widehat{\Z})$ have respective measure $1$.
But this follows from the already mentionned equality $\iota^{-1}(G'(\widehat{\Z}))=G(\widehat{\Z})$, and form the well-known fact that $\iota(G(\AAA_f))$ is a normal subgroup of $G'(\AAA_f)$.
\end{pf}
 
 A finer property of $\iota^\ast$ is that it commutes with the Satake isomorphism. Recall that if $\mathcal{H}_p(G)$ denotes the Hecke algebra of $(G(\Q_p),G(\Z_p))$, the Satake isomorphism 
is a canonical isomorphism $$S_{G/\Z_p} : \mathcal{H}_p(G) \isomo R(\widehat{G}) $$
where $R(\widehat{G})$ denotes the $\C$-algebra of polynomial class functions on $\widehat{G}$.  Satake shows {\it loc. cit.} that the diagram
$$\xymatrix{ \mathcal{H}_p(G) \ar@{->}^{S_{G/\Z_p}}[r] \ar@{->}[d]_{\iota^\ast} &  R(\widehat{G})  \ar@{->}[d] ^{\iota^\vee} \\
 \mathcal{H}_p(G') \ar@{->}_{S_{G'/\Z_p}}[r]  &  R(\widehat{G'})}
$$
is commutative, where $\iota^\vee : R(\widehat{G}) \rightarrow R(\widehat{G'})$ also denotes the restriction by $\iota^\vee$. \ps

Proposition~\ref{propbijiso} ensures that the map $f(g) \mapsto f(\iota(g))$ defines a $\C$-linear isomorphism 
\begin{equation}\label{eqreslemme} {\rm Res}_\iota : \mathcal{L}(G') \isomo \mathcal{L}(G).\end{equation}
The homomorphism $\iota^\ast$ defines a natural $\mathcal{H}(G)$-module structure on $\mathcal{L}(G')$ 
and Lemma~\ref{lemmasat} ensures that ${\rm Res}_\iota$ is 
$\mathcal{H}(G)$-equivariant for this structure on the left-hand side and the natural structure on the 
right-hand side. The isomorphism ${\rm Res}_\iota$ is obviously $G(\R)$-equivariant as well. As $\iota(G(\R))$ is open of 
finite index in $G'(\R)$ 
we may replace the two $\mathcal{L}$'s in \eqref{eqreslemme} by $\mathcal{L}_{\rm disc}$. We have thus proved the following proposition.

\begin{prop}\label{propisog2} ${\rm Res}_\iota$ induces an isomorphism $\mathcal{L}_{\rm disc}(G') \isomo \mathcal{L}_{\rm disc}(G)$ which commutes with the natural actions of $G(\R)$ and $\mathcal{H}(G)$ on both sides.
\end{prop}

This proposition implies Proposition~\ref{propisog} thanks to formula~\eqref{eqldisc}.

\begin{cor}\label{corisog} Assume that $m_G(\pi)=1$ for each $\pi \in \Pi_{\rm disc}(G)$. Then $m_{G'}(\pi')=1$ for each $\pi' \in \Pi_{\rm disc}(G')$ as well. Moreover, the $\mathcal{R}_\iota(\pi')$ with 
$\pi' \in \Pi_{\rm disc}(G')$ form a partition of $\Pi_{\rm disc}(G)$. 
\end{cor}

This corollary would apply for instance to the isogeny ${\rm Sp}_{2g} \rightarrow {\rm PGSp}_{2g}$ for any $g\geq 1$ by Arthur's multiplicity formula if we knew that the archimedean Arthur packets are sets rather than multisets (see the discussion following Conjecture~\ref{conjecturedeux}). It applies for $g=1$ by the multtiplicity one theorem of Labesse and Langlands~\cite{labesselanglands}.\ps

\begin{cor}\label{m1so} If $G={\rm SO}_{2,2}$ then $m_G(\pi)=1$ for any $\pi \in \Pi_{\rm
disc}(G)$.
\end{cor}

\begin{pf} 
We just recalled that $m_H(\pi)=1$ for any $\pi \in \Pi_{\rm disc}(H)$ when
$H={\rm SL}_2$, hence for $H={\rm SL}_2 \times {\rm SL}_2$ as well. To
conclude we apply Cor.~\ref{corisog} to the central isogeny 
$$({\rm SO}_{2,2})_{\rm sc} \simeq {\rm SL}_2 \times {\rm SL}_2 \rightarrow {\rm
SO}_{2,2}.$$

\end{pf}

\subsection{Symmetric square functoriality and $\Pi^\bot_{\rm cusp}(\PGL_3)$} It follows from Theorem~\ref{orthosymp} that 
$$\Pi_{\rm cusp}^{\rm o}(\PGL_3)=\Pi_{\rm cusp}^\bot(\PGL_3).$$
Recall the $\C$-morphism ${\rm Sym^2} : {\rm SL}_2(\C) \rightarrow {\rm SL}_3(\C)$. 

\begin{propconda}\label{symsq}There is a unique bijection ${\rm Sym}^2 : \Pi_{\rm cusp}(\PGL_2) \rightarrow \Pi_{\rm cusp}^{\bot}(\PGL_3)$ such that for each 
$\pi \in \Pi_{\rm cusp}(\PGL_2)$ we have $c({\rm Sym}^2\pi)={\rm Sym}^2 c(\pi)$. It induces a bijection 
$\Pi_{\rm alg}(\PGL_2) \isomo \Pi_{\rm alg}^{\rm o}(\PGL_3)$.
\end{propconda}

If $\pi \in \Pi_{\rm alg}(\PGL_2)$ has Hodge weight $w$, it follows that ${\rm Sym}^2(\pi)$ has Hodge weight $2w$. The proposition implies thus part (i) of
Thm.~\ref{proporthointro}. Observe in particular that the 
Hodge weight of any $\pi \in \Pi_{\rm alg}^{\rm o}(\PGL_3)$ is $\equiv 2 \bmod 4$, as asserted in general by Prop.~\ref{propepsilonintro}. \ps

\begin{pf} The existence of a unique map ${\rm Sym}^2 : \Pi_{\rm cusp}(\PGL_2) \rightarrow \Pi_{\rm cusp}^\bot(\PGL_3)$  satisfying $c({\rm Sym}^2\pi)={\rm Sym}^2 c(\pi)$ is due to Gelbart and Jacquet~\cite[Thm. 9.3]{GJ} (the assumption in their theorem is satisfied as $\pi$ as conductor $1$). It is however instructive to deduce it as well from Arthur's results, as follows. Consider the isogeny $\iota : {\rm SL}_2 \rightarrow {\rm PGL}_2$. 
Let $\pi \in \Pi_{\rm disc}({\rm PGL}_2)$ and let $\rho \in {\rm Res}_\iota(\pi)$. By Proposition~\ref{propisog} (b), we have $\rho \in \Pi_{\rm disc}({\rm SL}_2)$. By definition, $c(\rho)$ is the image of $c(\pi)$ under the isogeny $\iota^\vee : {\rm SL}_2(\C) \rightarrow {\rm PGL}_2(\C)={\rm SO}_3(\C)$; observe that the composition of $\iota^\vee$ with the standard representation of ${\rm SO}_3(\C)$ is nothing else than the ${\rm Sym}^2$ representation of ${\rm SL}_2(\C)$. In particular, $\psi(\rho)$ does not depend on the choice of $\rho$ in ${\rm Res}_\iota(\pi)$, and it thus makes sense to consider
$$\widetilde{\psi}(\pi)=\psi(\rho) \in \Psi_{\rm glob}({\rm SL}_2).$$ \ps
We have $\widetilde{\psi}(\pi)=[3]$ if and only if $\rho$ is the trivial representation, which can happen only if $\pi$ is trivial as well (see e.g. the decomposition~\eqref{decompggsc}). Otherwise, the only remaining possibilities are that $\pi \in \Pi_{\rm cusp}(\PGL_2)$ and $\widetilde{\psi}(\pi) \in \Pi_{\rm cusp}^{\rm o}(\PGL_3)$. If we set ${\rm Sym}^2 \pi = \widetilde{\psi}(\pi)$, then $c({\rm Sym}^2 \pi)={\rm Sym}^2 c(\pi)$ by construction: this is another definition of the Gelbart-Jacquet map. \ps

The ${\rm Sym}^2$ map is surjective. Indeed, if $\pi' \in \Pi_{\rm cusp}^\bot(\PGL_3)$ there exists $\rho \in \Pi_{\rm disc}({\rm SL}_2)$ such that $\pi'=\psi(\rho)$ by Arthur's Theorem~\ref{orthosymp}. But there exists $\pi \in \Pi_{\rm disc}({\rm PGL}_2)$ such that $\rho \in {\rm Res}_\iota(\pi)$ by Proposition~\ref{propisog} (a). One sees as above that $\pi$ is non-trivial, hence cuspidal (Selberg). It follows that $\pi'={\rm Sym}^2 \pi$. \ps

It only remains to check that ${\rm Sym}^2$ is injective. Let $\pi, \pi' \in \Pi_{\rm cusp}(\PGL_2)$ be such that ${\rm Sym}^2 \pi \simeq {\rm Sym}^2 \pi'$. For each prime $p$, the one or two elements\footnote{Recall that the image of ${\rm SL}_2(\R) \rightarrow {\rm PGL}_2(\R)$ has index $2$.} in ${\rm Res}_\iota(\pi)$, and the one or two elements in ${\rm Res}_\iota(\pi')$, all have the same Satake parameters at $p$. By the multiplicity formula of  Labesse-Langlands~\cite{labesselanglands}, this implies that all these representations are in a same {\it global ${\rm L}$-packet};  this means here that their archimedean components are all conjugate under ${\rm PGL}_2(\R)$. It follows that ${\rm Res}_\iota(\pi)={\rm Res}_\iota(\pi')$. But by Labesse-Langlands~\cite{labesselanglands} again, each element in $\Pi_{\rm disc}(\PGL_2)$ has multiplicity one. It follows that $\pi \simeq \pi'$ by Corollary~\ref{corisog}. 
\end{pf}

\subsection{Tensor product functoriality and $\Pi_{\rm cusp}^{\rm o}(\PGL_4)$} We consider the natural map 
$$\mathcal{X}(\SL_2(\C)) \times \mathcal{X}(\SL_2(\C)) \rightarrow \mathcal{X}(\SL_4(\C))$$
given by the tensor product $(x,y) \mapsto x \otimes y$ of conjugacy classes. If $X$ is a set, we denote by $\Sigma_2 X$ 
the set of all subsets of $X$ with two elements. 

\begin{propconda}\label{tensdeux} There is a unique bijection $\Sigma_2\, \Pi_{\rm cusp}(\PGL_2) \isomo \Pi_{\rm
cusp}^{\rm o}(\PGL_4)$, that we shall denote 
$\{\pi,\pi'\} \mapsto \pi \otimes \pi'$, such that 
for each $\pi \neq \pi' \in \Pi_{\rm cusp}(\PGL_2)$, $$c(\pi \otimes \pi') = c(\pi) \otimes c(\pi').$$  It induces a bijection $\Sigma_2 \,\Pi_{\rm alg}(\PGL_2) \isomo \Pi_{\rm alg}^{\rm o}(\PGL_4)$.
\end{propconda}

Consider the central isogeny  $\iota : {\rm SO}_{2,2} \rightarrow \PGL_2 \times \PGL_2$. 
Let $(\pi,\pi') \in \Pi_{\rm disc}(\PGL_2)^2$ and let $\rho \in {\rm Res}_\iota((\pi,\pi'))$. By Proposition~\ref{propisog} (b), we have $\rho \in \Pi_{\rm disc}({\rm SO}_{2,2})$. 
By definition, $c(\rho)$ is the image of $c(\pi) \times c(\pi')$ under the isogeny $\iota^\vee : {\rm SL}_2(\C)^2 \rightarrow {\rm SO}_4(\C)$. If we compose this latter isogeny with the standard representation of  ${\rm SO}_4(\C)$, we obtain nothing else than the tensor product representation ${\rm SL}_2(\C)^2 \rightarrow {\rm SL}_4(\C)$. In particular, $\psi(\rho)$ does not depend on the choice of $\rho$ in ${\rm Res}_\iota(\pi)$, and it thus makes sense to 
define $$\psi(\pi,\pi')=\psi(\rho) \in \Psi_{\rm glob}({\rm SO}_{2,2}).$$ \ps
It is clear that  $\psi(\pi,\pi')=\psi(\pi',\pi)$. \ps

\begin{propconda}\label{lemmeso4} Let $\pi,\pi' \in \Pi_{\rm disc}(\PGL_2)$. \begin{itemize}
\item[(i)] If $\pi,\pi'$  are both the trivial representation then $\psi(\pi,\pi')=[3]\oplus [1]$, \ps
\item[(ii)] If $\pi'$ is the trivial representation and $\pi$ is cuspidal  then $\psi(\pi,\pi')=\pi[2]$,\ps
\item[(iii)] If $\pi=\pi'$ is cuspidal, then $\psi(\pi,\pi')={\rm Sym}^2 \pi \oplus [1]$,\ps
\item[(iv)] If $\pi,\pi'$ are distinct and cuspidal, then $\psi(\pi,\pi') \in \Pi_{\rm cusp}^{\rm o}(\PGL_4)$. Moreover, $\psi(\pi,\pi')$ determines the pair $\{\pi,\pi'\}$.
\end{itemize}
\end{propconda}

Note that assertion (iii) makes sense by Proposition~\ref{symsq}.

\begin{pf} Assertions (i), (ii) and (iii) follow from an immediate inspection of Satake parameters and from the uniqueness of global Arthur parameters in Theorem~\ref{thmarthurclass}. \ps

Fix distinct $\pi,\pi' \in \Pi_{\rm cusp}(\PGL_2)$. The strong multiplicity one theorem for $\PGL_2$ shows that the global Arthur parameter $\psi(\pi,\pi')$ cannot contain the symbol $[1]$.  Moreover, Jacquet-Shalika's bound shows that $\psi(\pi,\pi')$ cannot have the form $\pi''[2]$ for $\pi'' \in \Pi_{\rm cusp}(\PGL_2)$. The only remaining possibility is that $\psi(\pi,\pi') \in \Pi_{\rm cusp}^{\rm o}(\PGL_4)$. \ps

Fix now $\omega \in \Pi_{\rm cusp}^\bot(\PGL_4)$ of the form $\psi(\pi,\pi')$ for some distinct $\pi,\pi' \in \Pi_{\rm cusp}(\PGL_2)$. We want to show that $\omega$ determines the pair $\{\pi,\pi'\}$. Consider for this the subset $\mathcal{X}({\rm SL}_4(\C))^\bot$ of $\mathcal{X}({\rm SL}_4(\C))$ 
of all the conjugacy classes which are equal to their inverse, and consider the map 
$$t : \mathcal{X}({\rm SL}_4(\C))^\bot \, \, \longrightarrow \, \, \mathcal{X}({\rm SL}_3(\C))^2/\mathfrak{S}_2$$
defined as follows. Start with the standard representation ${\rm O}_4(\C) \rightarrow {\rm SL}_4(\C)$. 
An element $x \in \mathcal{X}({\rm SL}_4(\C))^\bot$ is the image of a unique ${\rm O}_4(\C)$-conjugacy class $y$ in ${\rm SO}_4(\C)$. The image of $y$ via the isogeny ${\rm SO}_4(\C) \rightarrow \PGL_2(\C)^2$ is a well defined element\footnote{If $X$ is a set $X^2/\mathfrak{S}_2$ denotes the quotient of $X^2$ by the equivalence relation $(x,y) \sim (y,x)$. } $z \in \mathcal{X}(\PGL_2(\C))^2/\mathfrak{S}_2$. Set $t(x)={\rm ad}(z)$ where ${\rm ad}: {\rm PGL}_2(\C) \rightarrow {\rm SL}_3(\C)$ is the adjoint representation. Observe that for each prime $p$, we have $$t(c_p(\omega)) \sim ({\rm ad} \circ \mu (c_p(\pi)),{\rm ad} \circ \mu (c_p(\pi')))$$ where $\mu$ is the isogeny ${\rm SL}_2(\C) \rightarrow {\rm PGL}_2(\C)$. But ${\rm ad} \circ \mu (c_p(\pi)) = c_p({\rm Sym}^2 \pi)$, and similarly for $\pi'$. It follows from Jacquet-Shalika's structure theorem for isobaric representation~\cite{jasha} that the pair $\{{\rm Sym}^2\pi, {\rm Sym}^2 \pi'\}$ is uniquely determined by $\omega$. But by Proposition~\ref{symsq} this in turn determines $\{\pi,\pi'\}$. \end{pf}

Let us finally prove the first assertion of Proposition~\ref{tensdeux}. If $\pi,\pi' \in \Pi_{\rm cusp}(\PGL_2)$ are distinct we set
		$$\pi \otimes \pi' = \psi(\pi,\pi'),$$
so that  $c(\pi \otimes \pi') = c(\pi) \otimes c(\pi')$ by definition. It follows from Proposition~\ref{lemmeso4} (iv) that $\{\pi,\pi'\} \mapsto \pi \otimes \pi'$ 
defines an injection $\Sigma_2(\Pi_{\rm cusp}(\PGL_2)) \rightarrow \Pi_{\rm cusp}^{\rm o}(\PGL_4)$. Let us check that it is surjective. If $\omega \in \Pi_{\rm cusp}^{\rm o}(\PGL_4)$, Theorem~\ref{orthosymp} shows the existence of $\rho \in \Pi_{\rm disc}({\rm SO}_{2,2})$ such that $\omega = \psi(\rho)$. Proposition~\ref{propisog} ensures that $\rho$ belongs to ${\rm Res}_\iota((\pi,\pi'))$ for some $(\pi,\pi') \in \Pi_{\rm disc}(\PGL_2)^2$. But then 
$\omega=\psi(\rho)=\psi(\pi,\pi')$, so $\pi,\pi'$ are distinct and cuspidal by Proposition~\ref{lemmeso4}, hence $\omega = \pi \otimes \pi'$. \ps

If $\pi,\pi' \in \Pi_{\rm alg}(\PGL_2)$ have respective Hodge weights $w\geq w'$, the infinitesimal character of 
$\pi \otimes \pi'$ has the eigenvalues $\pm \frac{w+w'}{2},\pm \frac{w-w'}{2}$. This implies that $\pi \otimes \pi'$ is in $\Pi_{\rm alg}^\bot(\PGL_4)$.
Indeed, this is clear if $w \neq w'$. If $w = w'$, Clozel's purity lemma~\ref{clopl} shows that 
${\rm L}((\pi \otimes \pi')_\infty) = {\rm I}_{2w} \oplus \chi_1 \oplus \chi_2$ where $\chi_1,\chi_2 \in \{1,\varepsilon_{\C/\R}\}$. But as $\pi \otimes \pi'$ has a trivial central character, we have $\chi_1\chi_2 = \det {\rm I}_{2w} = \varepsilon_{\C/\R}$, and we are done.\ps
This ends the proof of the proposition, and shows part (ii) of
Thm.~\ref{proporthointro}. $\square$

\subsection{$\Lambda^\ast$ functorality and $\Pi_{\rm cusp}^{\rm o}(\PGL_5)$} If 
$\pi \in \Pi_{\rm cusp}^{\rm s}(\PGL_4)$, there is a unique element $\widetilde{c(\pi)} \in \mathcal{X}({\rm Sp}_4(\C))$ such that ${\rm St}(\widetilde{c(\pi)}) = c(\pi)$ (see~\S~\ref{parorthosymp}). We denote by $\Lambda^\ast$ the irreducible representation ${\rm Sp}_4(\C) \rightarrow {\rm SL}_5(\C)$, so that $\Lambda^2 \C^4 = \Lambda^\ast \oplus 1$. \ps

\begin{propcondb}\label{lambdadeux} There is a unique map $\Pi_{\rm alg}^{\rm s}(\PGL_4) \isomo \Pi_{\rm alg}^{\rm o}(\PGL_5)$, denoted $\pi \mapsto \Lambda^\ast \pi$, such that 
for each $\pi \in \Pi_{\rm alg}^{\rm s}(\PGL_4)$ we have $\Lambda^\ast (\widetilde{c(\pi)})=c(\Lambda^\ast \pi)$.
\end{propcondb}

Note that if $\pi \in \Pi^{\rm s}_{\rm alg}(\PGL_4)$ has Hodge weights $w>v$, then 
$\Lambda^\ast \pi$ has Hodge weights $w+v>w-v$. The proposition implies thus Thm.~\ref{proporthointro} (iii). 
\ps

Consider the central isogeny  $\iota : {\rm Sp}_4 \rightarrow {\rm PGSp}_4={\rm SO}_{3,2}$. 
Let $\pi \in \Pi_{\rm disc}({\rm SO}_{3,2})$ and let $\rho \in {\rm Res}_\iota(\pi)$. By Proposition~\ref{propisog} (b), we have $\rho \in \Pi_{\rm disc}({\rm Sp}_4)$. By definition, $c(\rho)$ is the image of $c(\pi)$ under the isogeny ${\rm Sp}_4(\C) \rightarrow {\rm SO}_5(\C)$. If we compose this latter isogeny with the standard representation of ${\rm SO}_5(\C)$, we obtain nothing else than the $\Lambda^\ast$ representation of ${\rm Sp}_4(\C)$. In particular, $\psi(\rho)$ does not depend on the choice of $\rho$ in ${\rm Res}_\iota(\pi)$, and it thus makes sense to set $$\widetilde{\psi}(\pi)=\psi(\rho) \in \Psi_{\rm glob}({\rm Sp}_4).$$ \ps
Observe that $\widetilde{\psi}(\pi)= \widetilde{\psi}(\pi')$ if $\psi(\pi)=\psi(\pi')$. 
\begin{propcondb}\label{lemmeso5} Let $\pi \in \Pi_{\rm disc}({\rm SO}_{3,2})$. \begin{itemize}
\item[(i)] If $\psi(\pi)=[4]$ then $\widetilde{\psi}(\pi)=[5]$,\ps
\item[(ii)] If $\psi(\pi)=\pi_1\oplus [2]$ with $\pi_1 \in \Pi_{\rm cusp}(\PGL_2)$ then $\widetilde{\psi}(\pi)=\pi_1[2] \oplus [1]$,\ps
\item[(iii)] If $\psi(\pi)=\pi_1\oplus \pi_2$ with distinct $\pi_1,\pi_2 \in \Pi_{\rm cusp}(\PGL_2)$ then $\widetilde{\psi}(\pi)=\pi_1 \otimes \pi_2 \oplus [1]$,\ps
\item[(iv)] If $\psi(\pi) \in \Pi_{\rm alg}^{\rm s}(\PGL_4)$ then $\widetilde{\psi}(\pi) \in \Pi_{\rm alg}^{\rm o}(\PGL_5)$. Moreover, $\widetilde{\psi}(\pi)$ determines $\psi(\pi)$ in this case. \end{itemize}
\end{propcondb}

\begin{pf} Assertions (i), (ii) and (iii) follow from an immediate inspection of Satake
parameters. Assertion (iii) makes sense by Proposition~\ref{tensdeux}. Let us check (iv). Assume that $\omega:=\psi(\pi)$ is cuspidal. Jacquet-Shalika's bound shows that $\widetilde{\psi}(\pi)$ cannot have the form $[5]$ or $\pi'[2] \oplus [1]$ for $\pi' \in \Pi_{\rm cusp}(\PGL_2)$. The only remaining possibility is that $\widetilde{\psi}(\pi)$ is either cuspidal or of the form $\pi_1 \otimes \pi_2 \oplus [1]$ for two distinct $\pi_1,\pi_2 \in \Pi_{\rm cusp}(\PGL_2)$. To rule out this latter case and prove (iv), we shall need to known a certain property of Arthur's archimedean packets that we have not been able to extract from~\cite{arthur}, namely that if $\psi$ an archimedean {\it generic} parameter in his sense, the packet $\widetilde{\Pi}_{\psi}$ defined {\it loc. cit.} contains with multiplicity one each element of the associated Langlands packet having a Whittaker model (see~\cite{arthurunipotent}). This is why we assume from now on that $\psi(\pi) \in \Pi_{\rm alg}^{\rm s}(\PGL_4)$ and rely instead on Conjecture~\ref{conjecturedeux}. \ps

Consider first any $\psi' \in \Psi_{\rm alg}({\rm SO}_{3,2})$ which is either cuspidal or of the form $\pi_1 \oplus \pi_2$ with distinct $\pi_1,\pi_2 \in \Pi_{\rm cusp}(\PGL_2)$. The archimedean Arthur parameter $\psi'_\infty$ is a discrete series Langlands parameter. The associated set of discrete series of ${\rm SO}_{3,2}(\R)$ with infinitesimal character $z_{\psi'_\infty}$ contains a unique element $\pi_{\rm gen}$ having a Wittaker model; its Shelstad character $\tau(\pi_{\rm gen})$ is trivial by definition (\S\ref{parammap},\S\ref{parsp4}). Arthur's multiplicity formula~\ref{conjecturedeux} for ${\rm SO}_{3,2}$ shows thus that the unique element $\rho' \in \Pi_{\rm disc}({\rm SO}_{3,2})$ such that $\psi(\rho')=\psi'$ and $\rho'_\infty \simeq \pi_{\rm gen}$ has multiplicity $1$. This construction applies for instance to $\psi'=\psi(\pi)=\omega$ and gives a representation $\rho'$ that we shall denote by $\varpi$. \ps

Assume now that $\widetilde{\psi}(\pi)$ has the form $\pi_1 \otimes \pi_2 \oplus [1]$. Consider $\psi'=\pi_1\oplus \pi_2 \in \Psi_{\rm glob}({\rm SO}_{3,2})$. It has the same infinitesimal character as $\psi(\pi)$, so that $\psi' \in \Psi_{\rm alg}({\rm SO}_{3,2})$. Let $\rho'$ be the representation associated to $\psi'$ as in the previous paragraph. Then $$\widetilde{\psi}(\varpi)=\widetilde{\psi}(\rho'), \, \, \, \varpi_\infty \simeq \rho'_\infty, \, \, \, {\rm but}\, \, \, \omega=\psi(\varpi) \neq \psi(\rho')=\pi_1 \oplus \pi_2.$$  The first two equalities imply that ${\rm Res}_\iota(\varpi)={\rm Res}_\iota(\rho')$ (these sets actually have two elements because the restriction of $\pi_{\rm gen}$ to ${\rm Sp}_4(\R)$ has two factors). The last one and Proposition~\ref{propisog} imply then that the elements of ${\rm Res}_\iota(\varpi)$ have multiplicity $\geq 2$ in $\mathcal{L}_{\rm disc}({\rm Sp}_4)$. This contradicts Arthur's multiplicity formula~\ref{conjecturedeux} for ${\rm Sp}_4$. \ps

It follows that $\widetilde{\psi}(\pi)$ is cuspidal. By the exact same argument as in the previous paragraph we see that if $\pi' \in \Pi_{\rm disc}({\rm SO}_{3,2})$ is such that $\psi(\pi')$ is cuspidal and satisfy $\widetilde{\psi}(\pi') = \widetilde{\psi}(\pi)$, then $\psi(\pi')= \psi(\pi)$. \end{pf}

Let $\pi \in \Pi_{\rm cusp}^{\rm s}(\PGL_4)$. By Arthur's Theorem~\ref{orthosymp}, we may find a 
$\rho \in \Pi_{\rm disc}({\rm SO}_{3,2})$ such that $\pi=\psi(\rho)$. We set 
$$\Lambda^\ast \pi = \widetilde{\psi}(\rho)$$
It belongs to $\Pi_{\rm cusp}^{\rm o}(\PGL_5)$ by Proposition~\ref{lemmeso5} (iv) and does not depend on the choice of $\rho$ such that $\psi(\rho)=\pi$. The same proposition shows that $\pi \mapsto \Lambda^\ast \pi$ is injective. It only remains to check the surjectivity. 
If $\omega \in \Pi_{\rm cusp}^{\rm o}(\PGL_5)$, Theorem~\ref{orthosymp} shows the existence of $\rho \in \Pi_{\rm disc}({\rm Sp}_4)$ such that $\omega = \psi(\rho)$. Proposition~\ref{propisog} ensures that $\rho$ belongs to ${\rm Res}_\iota(\pi)$ for some $\pi \in \Pi_{\rm disc}({\rm SO}_{3,2})$. But then $\omega=\psi(\rho)=\widetilde{\psi}(\pi)$ so $\psi(\pi)$ is cuspidal by Proposition~\ref{lemmeso5}. This finishes the proof of Proposition~\ref{lambdadeux}. 

\newpage
\section{ $\Pi_{\rm disc}({\rm SO}_7)$ and $\Pi_{\rm alg}^{\rm s}(\PGL_6)$ }\label{chapSO7}

\subsection{The semisimple $\Z$-group ${\rm SO}_7$} Consider the semisimple
classical $\Z$-group $$G={\rm SO}_7={\rm SO}_{{\rm E}_7},$$ i.e.  the
special orthogonal group of the root lattice ${\rm E}_7$ (\S\ref{quadform}).  Let ${\rm W}(E_7)$ denote the Weyl group of the root system of ${\rm E}_7$, let
$\varepsilon : {\rm W}(E_7 )\rightarrow \{\pm 1\}$ be the signature
and ${\rm W}(E_7)^+={\rm Ker} \, \varepsilon$.  As the Dynkin diagram
of $E_7$ has no non-trivial automorphism one has ${\rm O}({\rm
E}_7)={\rm W}(E_7)$ (see~\S\ref{quadform}), thus
	$$G(\Z)={\rm W}(E_7)^+.$$
The group ${\rm W}(E_7)^+$ has order $1451520=7! \cdot 2^5 \cdot 3^2$, 
it is isomorphic via the reduction modulo $2$ to the finite simple group 
$G(\mathbb{F}_2)\simeq {\rm Sp}_6(\mathbb{F}_2)$ 
(\cite[Ch. VI, Ex. 3 \S 4]{bourbaki}). \ps

The class set ${\rm Cl}(G) \simeq {\rm X}_7$ has one element as $X_7=\{{\rm E}_7\}$
(\S~\ref{quadform},\S~\ref{gpdef}). By Arthur's multiplicity formula, each $\pi \in \Pi_{\rm disc}(G)$ has multiplicity $1$. It follows
from Prop.~\ref{definiteform} that 
the number $m(V)$ of $\pi \in \Pi_{\rm disc}(G)$ such that 
$\pi_\infty$ is a given irreducible representation of $G(\R)$ is 
$$m(V)  = \dim V^{{\rm W}({\rm E}_7)^+},$$
which is exactly the number computed in the first chapter~\S~\ref{appnum} Case I. We refer to Table~\ref{tableSO7nue} for a sample of results. \ps

The dual group of ${\rm SO}_7$ is $\widehat{G}={\rm Sp}_6(\C)$.

\subsection{Parameterization by the infinitesimal
character}\label{preliminf} From the point of view of Langlands
parameterization, it is more natural to label the irreducible
representations of $G(\R)$ by their infinitesimal character rather than
their highest weight.  \ps

Let $H$ be a compact connected Lie group, fix $T \subset H$ a  maximal torus
and $\Phi^+ \subset {\rm X}^\ast(T)$ a set of positive roots as
in~\S~\ref{degwform}.  Denote by $\rho \in {\rm X}^\ast(T)[1/2]$ the half
sum of the elements of $\Phi^+$. As recalled in~\S\ref{langparam}, under the Harish-Chandra isomorphism the
infinitesimal character of the irreducible representation $V_\lambda$ of $H$
of highest weight $\lambda$ is the ${\rm W}(H,T)$-orbit of $\lambda+\rho$. 
\ps

For instance if $H={\rm SO}_n(\R)$, and in terms of the standard root data defined in~\S~\ref{appnum}, 
$$ \rho = \left\{ \begin{array}{ll} \frac{2l-1}{2} e_1 + \frac{2l-3}{2} e_2 + \cdots + \frac{1}{2} e_l & \, \, \, {\rm if}\, \, \, n=2l+1, \\ 
 (l-1) e_1 + (l-2) e_2 + \cdots + e_{l-1}  &  \, \, \, {\rm if}\, \, \, n=2l. \end{array} \right.$$
The map $\lambda \mapsto \lambda+\rho = \sum_i \frac{w_i}{2} e_i$ induces thus a bijection between the dominant weights and the collection of $w_1> w_2 > \cdots >  w_l$
where the $w_i$ are odd positive integers when $n=2l+1$, even integers with $w_{l-1}>|w_l|$ when $n=2l$.
\newcommand{\ww}{{\underline{w}}}

\begin{definition}\label{defuw} Let $n\geq 1$ be an integer,  set $l=[n/2]$,
and let $\ww=(w_1,\cdots,w_l)$ where $w_1 > w_2 > \cdots > w_l \geq 0$ are
distinct nonnegative integers all congruent to $n$ modulo $2$.  We denote by
$$U_\ww$$ the finite dimensional irreducible representation $V_\lambda$ of
${\rm SO}_n(\R)$ such that $\lambda+\rho = \sum_i \frac{w_i}{2} e_i$. 
\end{definition}

As an example, observe that if ${\rm H}_m(\R^n)$ is the representation of ${\rm SO}_n(\R)$ defined in \S~\ref{reflexionrep}, then ${\rm H}_m(\R^n) =  U_\ww$ for 
$$w = \left\{ \begin{array}{ll} (2m+n-2,n-4,n-6,\cdots,3,1)  & \, \, \, {\rm if}\, \, \, n \equiv 1 \bmod 2, \\ 
 (2m+n,n-2,n-4,\cdots,2,0)  &  \, \, \, {\rm if}\, \, \, n \equiv 0 \bmod 2. \end{array} \right.$$

The infinitesimal character $\lambda+\rho$ is related to the Langlands
parameterization of $V_\lambda$ as follows.  Assume to simplify that $H$ is
semisimple and that $-1 \in {\rm W}(H,T)$.  This is always the case if
$H=G(\R)$ and $G$ is semisimple over $\Z$, and for $H={\rm SO}_n(\R)$ this
holds if and only if $n \not \equiv 2 \bmod 4$.  Then the Langlands dual
group of $H$ is a connected semisimple complex group $\widehat{H}$.  Recall
that $\widehat{H}$ is equipped with a maximal torus $\widehat{T}$, a set of
positive roots $(\Phi^\vee)^+$ for $(\widehat{H},\widehat{T})$, and an
isomorphism between the dual based root datum of
$(\widehat{H},\widehat{T},(\Phi^\vee)^+)$ and the one of $(H,T,\Phi^+)$.  In
particular, ${\rm X}_\ast(\widehat{T})$ and ${\rm X}^\ast(T)$ are identified
by definition.  The Langlands parameter of $V_\lambda$ is up to
$\widehat{H}$-conjugation the unique continuous homomorphism ${\rm
L}(V_\lambda) : {\rm W}_\R \rightarrow \widehat{H}$ with finite centralizer
and such that in Langlands' notation (see~\S\ref{demialg}) $${\rm
L}(V_\lambda)(z)=(z/\overline{z})^{\lambda+\rho} \in \widehat{T} \, \, \, \,
\, \, \forall z \in \C^\times={\rm W}_\C.$$

When $H={\rm SO}_n(\R)$ and $\ww=(w_1,w_2,\cdots,w_l)$ is as in definition~\ref{defuw}, it follows that in the standard representation ${\rm St} : \widehat{H} \rightarrow {\rm GL}(2l,\C)$ of the classical group $\widehat{H}$, we have
$${\rm St} \circ {\rm L}(U_\ww) \simeq  \bigoplus_{i=1}^l {\rm I}_{w_i}.$$
This is the reason why the normalization above will be convenient. \ps

\begin{definition} Let $G$ be the semisimple classical definite $\Z$-group ${\rm SO}_n$
defined in~\S\ref{quadform}. If $\ww=(w_1,\cdots,w_l)$ is as in Definition~\ref{defuw} we define 
$$\Pi_\ww(G)=\{ \pi \in \Pi_{\rm disc}(G), \pi_\infty \simeq U_\ww\}$$
and set $m(\ww)=|\Pi_\ww(G)|$. \ps
If $\pi \in \Pi_{\rm disc}(G)$, we shall say that $\pi$ has Hodge weights $\ww$ if $\pi \in \Pi_\ww(G)$. 
\end{definition}
\medskip

\subsection{Endoscopic partition of $\Pi_{\rm disc}({\rm SO}_7)$}\label{endopart}

Recall that if $\pi \in \Pi_\ww({\rm SO}_n)$, it has a global Arthur
parameter $$\psi(\pi)=(k,(n_i),(d_i),(\pi_i)) \in \Psi_{\rm glob}({\rm
SO}_n)$$ whose equivalence class in well-defined (\S\ref{arthurparam}).  The associated collection
$(k,(n_i),(d_i))$ will be called the {\it endoscopic type of $\pi$}.  As for
$\psi(\pi)$, the endoscopic type will be called {\it stable} if $k=1$, and
tempered if $d_i=1$ for $i=1,\cdots,k$.  By Lemma~\ref{lemmeclass},
$\psi(\pi)$ is stable and tempered if and only it belongs to $\Pi_{\rm
alg}^\bot(\PGL(2l))$ where $l=[n/2]$.  \ps

So far we have computed $|\Pi_\ww({\rm SO}_7)|$ for any possible Hodge weights $\ww$.  Our next aim will be to compute the number of elements in $\Pi_\ww({\rm SO}_7)$ of each possible endoscopic type. As we shall see, thanks to Arthur's 
multiplicity formula and our previous computation of ${\rm S}(w)$, ${\rm S}(w,v)$ and ${\rm O}^\ast(w)$, we will be able to compute the contribution of each endoscopic type except one, 
namely the stable and tempered type, which is actually ${\rm S}(w_1,w_2,w_3)$. We will in turn obtain this later number from our computation of $|\Pi_w({\rm SO}_7)|$. 
The Corollary~\ref{corintroso7} and Table~\ref{tableSwvu} will follow form these computations.\ps\ps

Fix a triple $\ww=(w_1,w_2,w_3)$. Fix as well once and for all a global Arthur parameter $$\psi=(k,(n_i),(d_i),(\pi_i)) \in \Psi_{\rm glob}({\rm SO}_7)$$ such that the semisimple conjugacy class ${\rm St}(z_{\psi_\infty})$ in $\mathfrak{sl}_6(\C)$  has the 
eigenvalues $$\pm \frac{w_1}{2}, \pm \frac{w_2}{2}, \pm \frac{w_3}{2}.$$ Let us denote by $\pi$ the unique element in $\Pi(\psi)$. We shall make explicit Arthur's multiplicity formula for $m(\pi)$, which is either $0$ or $1$ as $m_\psi=1$, following~\S\ref{mformoddorth}. Recall the important groups
$${\rm C_\psi} \subset {\rm C}_{\psi_\infty}\subset {\rm Sp}_6(\C).$$ \ps
For each $1 \leq i \leq k$ one has a distinguished element $s_i \in {\rm C}_\psi$ (\S~\ref{paragraphepsilon}). Those $k$-elements $s_i$ generate ${\rm C}_\psi\simeq (\Z/2\Z)^k$ and their product
generates the center $Z=\{\pm 1\}$ of ${\rm Sp}_6(\C)$. 

\subsubsection{The stable case}\label{casstableso7} This is the case $k=1$, i.e. ${\rm C}_\psi=Z$, for which the multiplicity formula trivially gives $m(\pi)=1$. Let us describe the different possibilities for $\psi$. One has $\psi(\pi)=\pi_1[d_1]$ with $d_1|6$, $\pi_1 \in \Pi_{\rm alg}^\bot(\PGL(6/d_1))$ and $(-1)^{d_1-1}s(\pi_1)=-1$. \ps \ps

{\bf Case (i)} : $\psi=\pi_1$ where $\pi_1 \in \Pi_{\rm alg}^{\rm
s}(\PGL_6)$, this is the unknown we want to count.  \ps\ps

{\bf Case (ii)} : $\psi=\pi_1[2]$ where $\pi_1 \in \Pi_{\rm alg}^{\rm o}(\PGL_3)$, say of Hodge weight $u>2$ (so $u \equiv 2 \bmod 4$). This occurs if and only if $\ww$ has the form $(u+1,u-1,1)$. Recall that $\pi_1={\rm Sym}^2 \pi'$ for a unique 
$\pi' \in \Pi_{{\rm alg}}(\PGL_2)$ with Hodge weight $u/2$. \ps\ps

{\bf Case (iii)} : $\psi=\pi_1[3]$ where $\pi_1 \in \Pi_{\rm alg}(\PGL_2)$, say of Hodge weight $u>1$ (an odd integer). This occurs if and only if $w$ has the form $(u+2,u,u-2)$.
\ps\ps

{\bf Case (iv)} : $\psi=[6]$. This occurs if and only if $w=(5,3,1)$, and $\pi$ is then the trivial representation of $G$. \ps \ps

\subsubsection{Endoscopic cases of type $(n_1,n_2)=(4,2)$}\label{endocase42}
In this case $k=2$, $$\psi=\pi_1[d_1]\oplus \pi_2[d_2]$$ and ${\rm C}_\psi
\simeq (\Z/2\Z)^2$.  It follows that ${\rm C}_{\psi}$ is generated by $s_1$ and
the center $Z$.  One will have to describe $\rho^\vee(s_1)$ and
$\varepsilon_\psi(s_1)=\varepsilon(\pi_1 \times \pi_2)^{{\rm Min}(d_1,d_2)}$
in each case.  Recall that $\rho^\vee : {\rm C}_{\psi_\infty} \rightarrow \{\pm
1\}$ is the fundamental character defined in~\S~\ref{mformoddorth}.  There
are three cases.\ps\ps

{\bf Case (v)} :  ({\it tempered case}) $d_1=d_2=1$, i.e. $\pi_1 \in \Pi_{\rm alg}^{\rm s}(\PGL_4)$ and $\pi_2 \in \Pi_{\rm alg}(\PGL_2)$. Denote by $a>b$ the Hodge weights of $\pi_1$ and by $c$ the 
Hodge weight of $\pi_2$. One has $\{a,b,c\}=\{w_1,w_2,w_3\}$.  One sees that 
$$\rho^\vee(s_1)=1 \, \, \, {\rm iff}\, \, \, a>c>b$$
But $\varepsilon_\psi(s_1)=1$ as all the $d_i$ are $1$ (tempered case). It follows that 
$$m(\pi)=1 \, \, \, \, \, \Leftrightarrow a>c>b$$
or which is the same, $m(\pi)=1$ if and only if $(w_1,w_2,w_3)=(a,c,b)$. \ps\ps

{\bf Case (vi)} :  $d_1=1$, $d_2=2$, i.e. $\psi=\pi_1\oplus [2]$ where $\pi_1\in \Pi_{\rm alg}^{\rm s}(\PGL_4)$ has Hodge weights $w_1>w_2$ with $w_2>1$. One sees that 
$\rho^\vee(s_1)=-1$. On the other hand
$\varepsilon_\psi(s_1)=\varepsilon(\pi_1)=(-1)^{(w_1+w_2+2)/2}$, it follows that 
$$m(\pi)=1 \, \, \, \, \, \Leftrightarrow \, \, \, \, \, w_1+w_2 \equiv 0 \bmod 4.$$\ps\ps

{\bf Case (vii)} : $d_1=4$, $d_2=1$, i.e. $\psi=[4]\oplus \pi_2$ where $\pi_2 \in \Pi_{\rm alg}(\PGL_2)$ has Hodge weight $w_1$ with $w_1>3$. One sees that 
$\rho^\vee(s_1)=-1$. On the other hand $\varepsilon_\psi(s_1)=\varepsilon(\pi_2)=(-1)^{(w_1+1)/2}$, it follows that 
$$m(\pi)=1 \, \, \, \, \, \Leftrightarrow \, \, \, \, \, w_1\equiv 1  \bmod 4.$$\ps\ps

\subsubsection{Endoscopic cases of type $(n_1,n_2,n_3)=(2,2,2)$} In this case $k=3$, and ${\rm C}_\psi$ is generated by $Z$ and $s_1,s_2$. There are two cases.\ps\ps

{\bf Case (viii)} : ({\it tempered case}) $d_i=1$ for each $i$, i.e. $\psi=\pi_1\oplus \pi_2 \oplus \pi_3$ where each $\pi_i \in \Pi_{\rm alg}(\PGL_2)$ and $\pi_i$ has Hodge weight $w_i$. 
and $\pi_2 \in \Pi_{\rm alg}(\PGL_2)$. Of course $\varepsilon_\psi$ is trivial here, so $m(\pi)=1$ if and only if $\rho^\vee$ is trivial on ${\rm C}_\psi$. But ${\rm C}_\psi={\rm C}_{\psi_\infty}$ and
$\rho^\vee$ is a non-trivial character, so $$m(\pi)=0$$
in all the cases.\ps\ps

{\bf Case (ix)} : $d_1=d_2=1$ and $d_3=2$, i.e. $\psi=\pi_1\oplus \pi_2 \oplus [2]$ where $\pi_1,\pi_2 \in   \Pi_{\rm alg}(\PGL_2)$ have respective Hodge weights $w_1>w_2$, with $w_2>1$. 
One has thus $\rho^\vee(s_1)=-1$ and $\rho^\vee(s_2)=1$. On the other hand for $i=1,2$ one has $\varepsilon_\psi(s_i)=\varepsilon(\pi_i)=(-1)^{(w_i+1)/2}$. It follows that 
$$m(\pi)=1 \, \, \, \, \, \Leftrightarrow \, \, \, \, \, (w_1,w_2)\equiv (1,3)  \bmod 4.$$

\subsection{Conclusions}\label{conclso7}

First, one obtains the value of ${\rm S}(w_1,w_2,w_3)$ as the difference
between $m(w_1,w_2,w_3)$ and the sum of the eight last contributions above. 
For instance, one sees that if $w_1-2>w_2>w_3+2>3$ then $${\rm
S}(w_1,w_2,w_3)=m(w_1,w_2,w_3)-{\rm S}(w_1,w_3)\cdot {\rm S}(w_2).$$ It
turns out that all the formulas for the nine cases considered above
perfectly fit our computations, in the sense that ${\rm S}(w_1,w_2,w_3)$
always returned to us a positive integer.  This is again a substantial
confirmation for both our computer program and for the remarkable precision
of Arthur's results.  This also gives some mysterious significance for the
first non-trivial invariants of the group ${\rm W}(E_7)^+$.  One
deduces in particular Table~\ref{tableSwvu}, and from this table
Corollary~\ref{corintroso7} of the introduction.  \ps

\begin{corcondb}  If $w_1<23$ then ${\rm S}(w_1,w_2,w_3)=0$. There are exactly $7$ triples $(23,w_2,w_3)$ such that ${\rm S}(23,w_2,w_3)\neq 0$, and for each of them ${\rm S}(23,w_2,w_3)=1$.
\end{corcondb}

As far as we know, none of these $7$ automorphic representations (symplectic
of rank $6$) had been discovered before.  As explained in the introduction,
they are related to the $121$ Borcherds even lattices of rank $25$ and
covolume $\sqrt{2}$, in the same way as the $4$ Tsushima's forms
$\Delta_{19,7}, \Delta_{21,5},\Delta_{21,9}$ and $\Delta_{21,13}$ are
related to Niemeier lattices, as discovered in~\cite{cl2}.  It would be
interesting to know more about those forms, e.g.  some of their Satake
parameters. Our tables actually reveals a number of triples $(w_1,w_2,w_3)$ such that
${\rm S}(w_1,w_2,w_3)=1$. \ps

One obtains as well a complete endoscopic description of each $\Pi_\ww({\rm SO}_7)$. For instance Tables~\ref{tableSO7} 
and~\ref{tableSO72} describe entirely the set $\Pi_{w_1,w_2,w_3}({\rm SO}_7)$ for $w_1\leq 25$ whenever it is
non-empty. Recall the following notation already introduced
in~\S\ref{applborcherdslat}: when ${\rm S}(w_1,\cdots,w_r)=1$ we denote by 
$\Delta_{w_1,\cdots,w_r}$ the unique $\pi \in 
\Pi_{\rm alg}^{\rm s}(\PGL_{2r})$ with Hodge weights $w_1>\cdots>w_r$. When
${\rm S}(w_1,\cdots,w_r)=k$ we also denote by $\Delta^k_{w_1,\cdots,w_r}$ any of
the $k$ elements of $\Pi_{\rm alg}^{\rm s}(\PGL_{2r})$ with Hodge weights
$w_1>\cdots>w_r$.\ps

Let us explore some examples. It follows from case (iii) above that the number of $\pi \in \Pi_\ww({\rm SO}_7)$ such that $\psi(\pi)$ has the form $\pi_1[3]$ 
is $\delta_{w_1=w_3+4} \cdot {\rm S}(w_2)$. For instance the first such $\pi$ is $\Delta_{11}[3]$ which thus belongs to $\Pi_{13,11,9}({\rm SO}_7)$. Our computations gives 
$m(13,11,9)=1$ (hence nonzero!) which is not only in accord with Arthur's result but also says that 
$$\Pi_{13,11,9}(G)=\{ \Delta_{11}[3]\}.$$
The triple $\ww=(13,11,9)$ turns out to be the first triple $\neq (5,3,1)$ such that $m(\ww) \neq 0$. Our table even shows that  
$$\forall\,\,  3 \leq u \leq 25, \, \, \, m(u+2,u,u-2)={\rm S}(u),$$
which describes entirely $\Pi_{u+2,u,u-2}({\rm SO}_7)$ for those $u$. One actually has $$m(29,27,25)=4>{\rm S}(27)=2.$$ 
Let us determine $\Pi_{29,27,25}({\rm SO}_7)$. We already found two forms $\Delta_{27}^2[3]$ (there are two elements  in $\Pi_{\rm alg}(\PGL_2)$ of Hodge weight $27$). On the other hand, 
one checks from Tsushima's formula that ${\rm S}(29,25)=1$, so that there is a unique element in $\Delta_{29,25} \in \Pi_{{\rm alg}}^{\rm s}(\PGL_4)$ with Hodge weights $29>25$. The missing two elements are thus the two 
$\Delta_{29,25}\oplus \Delta_{27}^2$. Indeed, we are here in the endoscopic case (v):  $27$ is between $25$ and $29$.\ps

As another example, consider now the $\pi \in \Pi_\ww(G)$ such that $\psi(\pi)$ has the form $\pi_1[2]$ (endoscopic case (ii)). There are exactly $$\delta_{w_3=1} \cdot \delta_{w_1-w_2=2} \cdot \delta_{w_1 \equiv 1 \bmod 4} \cdot {\rm S}(\frac{w_1+1}{2})$$ such $\pi$'s.  The first one  is thus ${\rm Sym}^2\Delta_{11}[2]$ which belongs to $\Pi_{23,21,1}({\rm SO}_7)$. Our computations gives 
$m(23,21,1)=1$, which is not only in accord with Arthur's result but also says that 
$$\Pi_{23,21,1}({\rm SO}_7)=\{ {\rm Sym}^2\Delta_{11}[2]\}.$$

\newpage 
\section{ Description of $\Pi_{\rm disc}({\rm SO}_9)$ and $\Pi_{\rm alg}^{\rm s}(\PGL_8)$}\label{chapSO9}

\subsection{The semisimple $\Z$-group ${\rm SO}_9$}  Consider the semisimple
classical $\Z$-group $$G={\rm SO}_9,$$ i.e.  the special orthogonal group of
the root lattice $L={\rm A}_1 \oplus {\rm E}_8$ (\S\ref{quadform}). 
Let ${\rm W}(E_8)$ denote the Weyl group of the root system of ${\rm
E}_8$ and let $\varepsilon : {\rm W}(E_8) \rightarrow \{\pm 1\}$ be
the signature.  There is a natural homomorphism ${\rm W}(E_8)
\rightarrow G(\Z)$, if we let ${\rm W}(E_8)$ act on $ {\rm A}_1 \oplus
{\rm E}_8$ as $w \mapsto (w,\varepsilon(w))$ (~\S~\ref{appnum} case III). 
One has ${\rm O}(L)=\{\pm 1\} \times {\rm W}(E_8)$ by~\S\ref{quadform}, thus a natural
isomorphism $${\rm W}(E_8) \isomo G(\Z).$$ The group $W$ has order
$$|{\rm W}(E_8)|=8!\cdot 2^5\cdot 3^3 \cdot 5=696729600$$ and the
natural map ${\rm W}(E_8)\rightarrow {\rm SO}_{{\rm
E}_8}(\mathbb{F}_2)$ is surjective with kernel $\{\pm 1\}$
(\cite[Ch. VI, Ex. 1 \S 4]{bourbaki}).  \ps

The class set ${\rm Cl}(G) \simeq {\rm X}_9$ has one element as $X_9=\{{\rm A}_1\oplus{\rm E}_8\}$
(\S~\ref{quadform},\S~\ref{gpdef}). By Arthur's multiplicity formula, each $\pi \in \Pi_{\rm disc}(G)$ has multiplicity $1$. It follows
from Prop.~\ref{definiteform} that 
the number $m(V)$ of $\pi \in \Pi_{\rm disc}(G)$ such that 
$\pi_\infty$ is a given irreducible representation of $G(\R)$ is 
$$m(V)  = \dim V^{{\rm W}(E_8)^+},$$
which is exactly the number computed in the first chapter~\S~\ref{appnum} Case III. We refer to Table~\ref{tableSO9nue} for a sample of results. \ps

The dual group of ${\rm SO}_9$ is $\widehat{G}={\rm Sp}_8(\C)$.

\subsection{Endoscopic partition of $\Pi_w$}

We proceed in a similar way as in~\S~\ref{endopart}. \ps

Fix $\ww=(w_1,w_2,w_3,w_4)$ with $w_1>w_2>w_3>w_4$ odd positive integers. Fix as well once and for all a global Arthur parameter $$\psi=(k,(n_i),(d_i),(\pi_i)) \in \Psi_{\rm glob}(G)$$ such that the semisimple conjugacy class ${\rm St}(z_{\psi_\infty})$ in $\mathfrak{sl}_8(\C)$  has the 
eigenvalues $$\{\pm \frac{w_i}{2}, 1\leq i \leq 4\}.$$ Let us denote by $\pi$ the unique element in $\Pi(\psi)$. We shall make explicit Arthur's multiplicity formula for $m(\pi)$, which is either $0$ or $1$, as in~\S\ref{mformoddorth}. Recall the groups
$${\rm C_\psi} \subset {\rm C}_{\psi_\infty}\subset \widehat{G}={\rm Sp}_8(\C).$$ 

\subsubsection{The stable cases}This is the case $k=1$, i.e. ${\rm C}_\psi=Z$, for which the multiplicity formula trivially gives $m(\pi)=1$. One has $\psi(\pi)=\pi_1[d_1]$ with $d_1|8$, 
$\pi_1 \in \Pi_{\rm alg}^\bot(\PGL(8/d_1))$, and $(-1)^{d_1-1}s(\pi_1)=-1$. \ps \ps

{\bf Case (i)} : ({\it tempered case}) $\psi=\pi_1$ where $\pi_1\in \Pi_{\rm alg}^{\rm s}(\PGL_8)$, this is the unknown we want to
count. \ps\ps

{\bf Case (ii)} : $\psi=\pi_1[2]$ where $\pi_1 \in \Pi_{\rm alg}^{\rm o}(\PGL_4)$, say of Hodge weights $u>v$ (recall $u,v$ even and $u +v \equiv 2 \bmod 4$). This occurs if and only if $w$ has the form $(u+1,u-1,v+1,v-1)$. 
Recall from Proposition~\ref{tensdeux} that $\pi_1=\pi' \otimes \pi''$ for a unique pair 
$\pi',\pi'' \in \Pi_{{\rm alg}}(\PGL_2)$ with respective Hodge weights $(u+v)/2$ and $(u-v)/2$. \ps\ps

{\bf Case (iii)} : $\psi=[8]$. This occurs if and only if $\ww=(7,5,3,1)$, and $\pi$ is then the trivial representation of $G$. \ps \ps

\subsubsection{Endoscopic cases of type $(n_1,n_2)=(6,2)$} In this case
$k=2$, $$\psi=\pi_1[d_1]\oplus \pi_2[d_2]$$ and ${\rm C}_\psi \simeq
(\Z/2\Z)^2$.  It follows that ${\rm C}_{\psi}$ is generated by $s_1$ and the
center $Z$.  One will have to describe $\rho^\vee(s_1)$ and
$\varepsilon_\psi(s_1)=\varepsilon(\pi_1 \times \pi_2)^{{\rm Min}(d_1,d_2)}$
in each case.  Recall that $\rho^\vee : {\rm C}_{\psi_\infty} \rightarrow \{\pm
1\}$ is the fundamental character defined in~\S~\ref{mformoddorth}.  There
are 6 cases.\ps\ps

{\bf Case (iv)} : ({\it tempered case}) $d_1=d_2=1$, i.e. $\pi_1 \in \Pi_{\rm alg}^{\rm s}(\PGL_6)$ and $\pi_2 \in \Pi_{\rm alg}(\PGL_2)$. Denote by $a>b>c$ the Hodge weights of $\pi_1$ and by $d$ the 
Hodge weight of $\pi_2$. One has $\{a,b,c,d\}=\{w_1,w_2,w_3,w_4\}$.  Moreover $\varepsilon_\psi(s_1)=1$ as all the $d_i$ are $1$ (tempered case), so $m(\pi)=1$ if, and only if, $\rho^\vee(s_1)=1$, i.e. if $d \in \{w_1,w_3\}$. In other words,
$$m(\pi)=1 \, \, \, \, \, \Leftrightarrow \,\,\,d > a>b>c \, \,\, {\rm or}\, \,\, a>b>d>c.$$
 \ps\ps

{\bf Case (v)} :  $d_1=1$, $d_2=2$, i.e. $\psi=\pi_1\oplus [2]$ where $\pi_1 \in \Pi_{\rm alg}^{\rm s}(\PGL_6)$ has Hodge weights $w_1>w_2>w_3$, with $w_3>1$. One sees that 
$\rho^\vee(s_1)=-1$. On the other hand $\varepsilon_\psi(s_1)=\varepsilon(\pi_1)=(-1)^{(w_1+w_2+w_3+3)/2}$, it follows that 
$$m(\pi)=1 \, \, \, \, \, \Leftrightarrow \, \, \, \, \, w_1+w_2 + w_3 \equiv 3 \bmod 4.$$\ps\ps

{\bf Case (vi)} : $d_1=2$, $d_2=1$, i.e. $\psi=\pi_1[2]\oplus \pi_2$ where $\pi_1 \in \Pi_{\rm alg}^{\rm o}(\PGL_3)$ and $\pi_2 \in \Pi_{\rm
alg}(\PGL_2)$. Denote by $a$ and $b$ the respective Hodge weights of $\pi_1$ and $\pi_2$, so that $\{w_1,w_2,w_3,w_4\}=\{a+1,a-1,b,1\}$. There are two cases: either $b>a$ or $b<a$. One sees that $\rho^\vee(s_1)=1$ in both cases. On the other hand 
$$\varepsilon_\psi(s_1)=\varepsilon(\pi_1 \times \pi_2)=-(-1)^{\frac{b+1}{2}+{\rm Max}(a,b)}.$$ It follows that 
$$m(\pi)=1 \, \, \, \, \, \Leftrightarrow \, \, \, \, \,\left\{  \begin{array}{ll} b \equiv 3 \bmod 4, & {\rm if }\, \, \, b >a+1, \\ b \equiv 1 \bmod 4, & {\rm if}\, \, \, b<a-1.\end{array}\right.$$\ps\ps

{\bf Case (vii)} : $d_1=3$, $d_2=1$, i.e. $\psi=\pi_1[3]\oplus \pi_2$ where $\pi_1,\pi_2 \in \Pi_{\rm alg}(\PGL_2)$. Denote by $a$ and $b$ the respective Hodge weights of $\pi_1$ and $\pi_2$, so that 
$\{w_1,w_2,w_3,w_4\}=\{a+1,a,a-1,b\}$. One sees that $\rho^\vee(s_1)=1$ if $b>a+1$, $-1$ otherwise. On the other hand 
$\varepsilon_\psi(s_1)=\varepsilon(\pi_1 \times \pi_2)=1$. It follows that 
$$m(\pi)=1 \, \, \, \, \, \Leftrightarrow \, \, \, \, b>a+1.$$\ps\ps

{\bf Case (viii)} : $d_1=3$, $d_2=2$, i.e. $\psi=\pi_1[3]\oplus[2]$ where
$\pi_1 \in \Pi_{\rm alg}(\PGL_2)$.  The Hodge weight of $\pi_1$ is thus
$w_2>3$.  We have $\rho^\vee(s_1)=-1$ and
$\varepsilon_\psi(s_1)=\varepsilon(\pi_1)^2=1$.  It follows that
$$m(\pi)=0$$\ps\ps

{\bf Case (ix)} : $d_1=6$, $d_2=1$, i.e. $\psi=[6]\oplus \pi_2$ where $\pi_2
\in \Pi_{\rm alg}(\PGL_2)$ has Hodge weight $w_1$ with $w_1>5$.  One sees
that $\rho^\vee(s_1)=1$.  On the other hand
$\varepsilon_\psi(s_1)=\varepsilon(\pi_2)=(-1)^{(w_1+1)/2}$, it follows that
$$m(\pi)=1 \, \, \, \, \, \Leftrightarrow \, \, \, \, \, w_1\equiv 3 \bmod
4.$$\ps\ps

\begin{rem}{\rm Observe that the case $d_1=d_2=2$, i.e. $\psi=\pi_1[2] \oplus [2]$ where $\pi_1 \in \Pi_{\rm alg}^{\rm o}(\PGL_3)$, is impossible as it implies $w_2=w_1=1$.}\end{rem}

\subsubsection{Endoscopic cases of type $(n_1,n_2,n_3)=(4,2,2)$} In this case $k=3$, and ${\rm C}_\psi$ is generated by $Z$ (or $s_3$) and $s_1,s_2$. There are three cases.\ps\ps

{\bf Case (x)} :  ({\it tempered case}) $d_i=1$ for each $i$, i.e. $\psi=\pi_1\oplus \pi_2 \oplus \pi_3$ where $\pi_1 \in \Pi_{\rm alg}^{\rm s}(\PGL_4)$ and $\pi_2,\pi_3 \in \Pi_{\rm alg}(\PGL_2)$. Denote by $a>b$ the Hodge weights of $\pi_1$ and by $c$ and $d$ the ones of $\pi_2,\pi_3$, assuming $c>d$.  Of course $\varepsilon_\psi$ is trivial here, so $m(\pi)=1$ if and only if $\rho^\vee$ is trivial on ${\rm C}_\psi={\rm C}_{\psi_\infty}$. One thus obtains
$$m(\pi)=1  \, \, \, \, \, \Leftrightarrow \, \, \, \, \,c>a>d>b.$$
\ps\ps

{\bf Case (xi)} : $d_1=d_2=1$ and $d_3=2$, i.e. $\psi=\pi_1\oplus \pi_2 \oplus [2]$ where $\pi_1\in \Pi_{\rm alg}^{\rm s}(\PGL_4)$, $\pi_2 \in   \Pi_{\rm alg}(\PGL_2)$ have respective Hodge weights $a>b$ and $c$, with $\{w_1,w_2,w_3,w_4\}=\{a,b,c,1\}$. If $a>c>b$ then $\rho^\vee(s_1)=1$ and $\rho^\vee(s_2)=-1$, otherwise $\rho^\vee(s_1)=-1$ and $\rho^\vee(s_2)=1$. On the other hand for $i=1,2$ one has $\varepsilon_\psi(s_i)=\varepsilon(\pi_i)$. It follows that 

$$m(\pi)=1 \, \, \, \, \, \Leftrightarrow \, \, \, \, \,\left\{  \begin{array}{ll} (a+b,c) \equiv (2,1) \bmod 4, & {\rm if }\, \, \, a>c>b \\ (a+b,c) \equiv (0,3) \bmod 4, & {\rm otherwise.} \end{array}\right.$$\ps\ps

{\bf Case (xii)} : $d_1=4$, $d_2=d_3=1$, i.e.  $\psi=[4] \oplus \pi_2 \oplus \pi_3$ where $\pi_2,\pi_3 \in \Pi_{\rm alg}(\PGL_2)$ with respective Hodge weights $w_1$ and $w_2$, with $w_2>3$. One has $\rho^\vee(s_2)=1$, $\rho^\vee(s_3)=-1$, $\varepsilon_\psi(s_i)=\varepsilon(\pi_i)$ for $i=2,3$, thus
$$m(\pi)=1  \, \, \, \, \, \Leftrightarrow \, \, \, \, \,(w_1,w_2) \equiv (3,1) \bmod 4.$$

\subsubsection{Endoscopic cases of type $(n_1,n_2,n_3,n_4)=(2,2,2,2)$}  In this case $k=4$, and ${\rm C}_\psi$ is generated by $Z$ (or $s_4$) and $s_1,s_2,s_3$. There are two cases.\ps\ps

{\bf Case (xiii)} :  ({\it tempered case}) $d_i=1$ for each $i$, i.e. $\psi=\pi_1\oplus \pi_2 \oplus \pi_3\oplus \pi_4$ where $\pi_i \in \Pi_{\rm alg}(\PGL_2)$ has Hodge weight $w_i$. As $\varepsilon_\psi$ is trivial but not $\rho^\vee$ on ${\rm C}_\psi={\rm C}_{\psi_\infty}$ we have in all cases 
$$m(\pi)=0.$$

{\bf Case (xiv)} : $d_4=2$ and $d_1=d_2=d_3=1$, i.e. $\psi= \pi_1 \oplus \pi_2 \oplus \pi_3 \oplus [2]$ where $\pi_i \in \Pi_{\rm alg}(\PGL_2)$ has Hodge weight $w_i$, and $w_3>1$. One has $\rho^\vee(s_1)=\rho^\vee(s_3)=1$ and $\rho^\vee(s_2)=-1$. On the other hand $\varepsilon_\psi(s_i)=\varepsilon(\pi_i)$ for $i=1,2,3$. It follows that 
$$m(\pi)=1  \, \, \, \, \, \Leftrightarrow \, \, \, \, \,(w_1,w_2,w_3) \equiv (3,1,3) \bmod 4.$$

\subsubsection{Endoscopic cases of type $(n_1,n_2)=(4,4)$}  In this case
$k=2$, $$\psi=\pi_1[d_1]\oplus \pi_2[d_2]$$ and ${\rm C}_\psi \simeq
(\Z/2\Z)^2$.  It follows that ${\rm C}_{\psi}$ is generated by $s_1$ and the
center $Z$.  One only has to describe $\rho^\vee(s_1)$ and
$\varepsilon_\psi(s_1)=\varepsilon(\pi_1 \times \pi_2)^{{\rm Min}(d_1,d_2)}$
in each case.\ps\ps

{\bf Case (xv)} : ({\it tempered case}) $d_1=d_2=1$, i.e. $\psi=\pi_1\oplus \pi_2$ with $\pi_1,\pi_2 \in \Pi_{\rm alg}^{\rm s}(\PGL_4)$. Let $a>b$ be the Hodge weight of $\pi_1$ and $c>d$ the ones of $\pi_2$, one may assume that $a>c$, i.e. $a=w_1$. As $\varepsilon_\psi=1$, one sees that 
$$m(\pi)=1  \, \, \, \, \, \Leftrightarrow \, \, \, \, \,a>c>b>d.$$\ps\ps

{\bf Case (xvi)} : $d_1=1$ and $d_2=4$, i.e. $\psi=\pi_1 \oplus [4]$ where $\pi_1 \in \Pi_{\rm alg}^{\rm s}(\PGL_4)$ has Hodge weights $w_1>w_2$ with $w_2>3$. It follows that $\rho^\vee(s_1)=-1$, and as $\varepsilon_\psi(s_1)=\varepsilon(\pi_1)$ one obtains 
$$m(\pi)=1  \, \, \, \, \, \Leftrightarrow \, \, \, \, w_1+w_2 \equiv 0 \bmod 4.$$

\subsection{Conclusions} \label{exso9}The inspection of each case above, and our previous computation of ${\rm S}(w)$, ${\rm S}(w,v)$, ${\rm S}(w,v,u)$, ${\rm O}^\ast(w)$ and ${\rm O}(w,v)$, allow to compute the contribution of each endoscopic type except one, 
namely the stable and tempered type, which is actually ${\rm S}(w_1,w_2,w_3,w_4)$, that we thus deduce from our computation of $m(w_1,w_2,w_3,w_4)$. The Corollary~\ref{corintroso9} and Table~\ref{tableSwvut} follow form these computations.\ps\ps

\begin{corcondb} If $w_1<25$ then ${\rm S}(w_1,w_2,w_3,w_4)=0$. There are $33$ triples $(w_2,w_2,w_4)$ such that ${\rm S}(25,w_2,w_3,w_4)\neq 0$, and in each case ${\rm S}(25,w_2,w_3,w_4)=1$. 
\end{corcondb}

We refer to Table~\ref{tableSO91} for the description of all the nonempty  $\Pi_\ww({\rm SO}_9)$ when $w_1\leq 23$. \ps\ps

For the application to Theorem~\ref{thminftriv},  consider for instance the problem of describing $\Pi_{27,23,9,1}({\rm SO}_9)$. Our program tells us that $$m(27,23,9,1)=5,$$ so that 
$|\Pi_{27,23,9,1}({\rm SO}_9)|=5$. Fix $\pi \in \Pi_{27,23,9,1}({\rm SO}_9)$ and let $\psi(\pi)=(k,(n_i),(d_i),(\pi_i))$. \ps

Assume first that $\psi(\pi)$ is not tempered, i.e. that some $d_i \neq 1$. We may assume that $d_k>1$. One sees that $k>1$, $d_k=2$ and $d_i=1$ for $i<k$. 
As ${\rm S}(9)=0$ we have $k\leq 3$. If $k=2$ then we are in case (v). As $27+23+9 \equiv 3 \bmod 4$ one really has to compute ${\rm S}(27,23,9)$. Our computer program tells us that $m(27,23,9)=4$.  On the other hand one has 
$${\rm S}(27,23,9)=m(27,23,9)-{\rm S}(27,9)\cdot {\rm S}(23)$$
by~\S~\ref{conclso7}. By Tsushima's formula we have ${\rm S}(27,9)=1$. As ${\rm S}(23)=2$ we obtain $${\rm S}(27,23,9)=2.$$ There are thus two representations $\Delta_{27,23,9}^2\oplus[2]$ in $\Pi_{27,23,9,1}({\rm SO}_9)$. \ps
Assume now that $k=3$, so we are in case (xi) and the Hodge weights of $\pi_1$ are $a$ and $9$. As $a+9 \equiv 0 \bmod 4$, the multiplicity formula forces thus $a=23$. Tsushima's formula shows that ${\rm S}(23,9)=1$. As ${\rm S}(27)=2$ there are indeed two parameters $\Delta_{27}^2\oplus \Delta_{23,9} \oplus [2]$ in case (xi), whose associated $\pi$ each have multiplicity $1$ by the multiplicity formula. \ps
Suppose now that $\pi$ is tempered, i.e. $d_i=1$ for all $i$. The multiplicity formula shows that $1$ and $23$ are Hodge weights of a same $\pi_i$, say $\pi_{i_0}$. 
But we already checked that ${\rm S}(23,1)=0$ and ${\rm S}(23,9,1)=0$, and ${\rm S}(9)=0$, it follows that $k=1$, i.e. $\pi$ is stable.

\begin{corcondb}\label{exutilefin} $\Pi_{27,23,9,1}({\rm SO}_9)=\{ \Delta_{27,23,9}^2\oplus [2], \Delta_{27}^2\oplus\Delta_{23,9} \oplus [2], \Delta_{27,23,9,1}\}$. \end{corcondb}
\newpage

\section{Description of $\Pi_{\rm disc}({\rm SO}_8)$ and $\Pi_{\rm alg}^{\rm o}(\PGL_8)$}\label{chapSO8}

\subsection{The semisimple $\Z$-group ${\rm SO}_8$}\label{pargpSO8} Consider the semisimple classical $\Z$-group 
$$G={\rm SO}_8={\rm SO}_{{\rm E}_8},$$ i.e. the special orthogonal group of the root lattice ${\rm
E}_8$. Recall that ${\rm W}(E_8)$ denote the Weyl group of the root system 
of ${\rm E}_8$, that $\varepsilon : {\rm W}(E_8) \rightarrow \{\pm 1\}$ 
is the signature and that ${\rm W}(E_8)^+={\rm Ker} \, \varepsilon$. 
As the Dynkin diagram of ${\rm E}_8$ has no non-trivial automorphism one has ${\rm O}({\rm E}_8)={\rm W}(E_8)$
(\S\ref{quadform}), thus 
	$$G(\Z)={\rm W}(E_8)^+.$$ The class set ${\rm Cl}(G) \simeq
\widetilde{X}_8$ has one element as $X_8=\{{\rm E}_8\}$ and of course ${\rm
O}({\rm E}_8) \neq {\rm SO}({\rm E}_8)$ (\S~\ref{quadform},\S~\ref{gpdef}). 
\ps

We shall consider quadruples $\ww=(w_1,w_2,w_3,w_4)$ where $w_1 > w_2 > w_3
> w_4\geq 0$ are even integers.  It is not necessary to consider the
$(w_1,w_2,w_3,w_4)$ with $w_4<0$ as ${\rm O}({\rm E}_8)={\rm W}(E_8)$ contains root
reflexions. Indeed, fix such a reflexion $s$.  Then $s$ acts by conjugation on
$\mathcal{L}(G)$, hence on $\Pi_{\rm disc}(G)$, with the following property
: if $\pi_\infty$ has the highest weight $(n_1,n_2,n_3,n_4)$, then
$s(\pi)_\infty$ has the highest weight $(n_1,n_2,n_3,-n_4)$.  Moreover
$m(s(\pi))=m(\pi)$.\ps

Consider the number $m'(\ww):==\sum_{\pi \in \Pi_\ww(G)} m(\pi)$. It follows from Prop.~\ref{definiteform} that 
$$m'(\ww)=\dim U_\ww^{{\rm W}(E_8)^+},$$
which is exactly the number computed in the first chapter~\S~\ref{appnum} Case II. We refer to Table~\ref{tableSO8nue} for a sample of results. \ps

By Arthur's multiplicity formula, for each $\pi \in \Pi_\ww(G)$ we have $m(\pi)+m(s(\pi))\leq 2$. In particular, 
if $\ww=(w_1,w_2,w_3,w_4)$ is such that $w_4 \neq 0$, then $m(\pi)=1$. In this case, it follows that 
$$m(\ww)=m'(\ww)=\dim U_\ww^{{\rm W}({\rm E}_8)}.$$
(Recall that $m(\ww)=|\Pi_\ww(G)|$). \ps

The dual group of ${\rm SO}_8$ is $\widehat{G}={\rm SO}_8(\C)$.

\subsection{Endoscopic partition of $\Pi_w$} We proceed again in a similar way as in~\S~\ref{endopart}. \ps

Fix $\ww=(w_1,w_2,w_3,w_4)$ with $w_1>w_2>w_3>w_4\geq 0$ even integers. Fix as well once and for all a global Arthur parameter $$\psi=(k,(n_i),(d_i),(\pi_i)) \in \Psi_{\rm glob}(G)$$ such that the semisimple conjugacy class ${\rm St}(z_{\psi_\infty})$ in $\mathfrak{sl}_8(\C)$  has the 
eigenvalues $$\{\pm \frac{w_i}{2}, 1\leq i \leq 4\}.$$  
 We shall make explicit Arthur's multiplicity formula for the number
$$m'(\psi)=\sum_{\pi \in \Pi(\psi)\cap \Pi_\ww(G)} m(\pi),$$
following~\S\ref{mformevenorth}. It will be convenient to introduce the number 
$$e(\ww)=\left\{\begin{array}{ll} 1 & {\rm if} \, \, \, w_4>0,\\ 2  & {\rm  otherwise}.\end{array}\right.$$
Recall also the groups
$${\rm C_\psi} \subset {\rm C}_{\psi_\infty} \subset \widehat{G}={\rm SO}_8(\C).$$ 
Denote by $J \subset \{1,\cdots,k\}$ the set of integers $j$ such that $n_j \equiv 1 \bmod 2$. It follows from Lemma~\ref{lemmeclass} that: 
\begin{itemize}
\item[(i)] If $j \notin J$ then $n_j \equiv 0 \bmod 4$.\ps
\item[(ii)] $|J|=0$ or $2$, and in this latter case $\sum_{j \in J} n_j \equiv 0 \bmod 4$.
\end{itemize}
\noindent We will say that $\psi$ is {\it even-stable} if $k=1$, and {\it odd-stable} if $k=2$ and $J=\{1,2\}$. 

\subsubsection{The even-stable cases} \label{stcasso8} We have ${\rm C}_\psi=Z$ so the
multiplicity formula trivially gives $m'(\psi)=e(\ww)$.  One has
$\psi(\pi)=\pi_1[d_1]$ with $d_1|8$, $\pi_1 \in \Pi_{\rm
alg}^\bot(\PGL(8/d_1))$, and $(-1)^{d_1-1}s(\pi_1)=-1$.  \ps \ps

{\bf Case (i)} : ({\it tempered case}) $\psi=\pi_1$ where $\pi_1\in \Pi_{\rm
alg}^{\rm o}(\PGL_8)$, this is the first unknown we want to count.  \ps\ps

{\bf Case (ii)} : $\psi=\pi_1[2]$ where $\pi_1 \in \Pi_{\rm alg}^{\rm s}(\PGL_4)$. This occurs if and only if $w_1-w_2=w_3-w_4=2$ and 
$\pi_1$ has Hodge weights $w_1-1,w_3-1$. \ps\ps

{\bf Case (iii)} : $\psi=\pi_1[4]$ where $\pi_1 \in \Pi_{\rm alg}(\PGL_2)$. This occurs if and only if $w_1=w_4+6$ and $\pi_1$ has Hodge weight $w_1-3$. \ps\ps

\subsubsection{The odd-stable cases} We have again ${\rm C}_\psi=Z$ so the multiplicity formula trivially gives $m'(\psi)=1$. One has $\psi(\pi)=\pi_1[d_1]\oplus \pi_2[d_2]$ with $n_1$, $n_2$, $d_1$ and $d_2$ odd. These cases only occur when $w_4=0$. \ps\ps

{\bf Case (iv)} : $d_1=n_2=1$, i.e $\psi=\pi_1\oplus [1]$ where $\pi_1\in
\Pi_{\rm alg}^{\rm o}(\PGL_7)$, which is the second unknown we want to
count.  \ps\ps

{\bf Case (v)} : $d_1=d_2=1$, $n_1=5$, i.e. $\psi=\pi_1\oplus \pi_2$ where $\pi_1 \in \Pi_{\rm alg}^{\rm o}(\PGL_5)$ and $\pi_2 \in \Pi_{\rm alg}^{\rm o}(\PGL_3)$. \ps\ps

{\bf Case (vi)} : $d_1=n_1=5$, $d_2=1$, i.e. $\psi=[5] \oplus \pi_2$ where $\pi_2 \in \Pi_{\rm alg}^{\rm o}(\PGL_3)$. In this case 
$w_2=4$.\ps\ps

{\bf Case (vii)} : $d_1=1$, $n_1=5$, $d_2=3$, i.e. $\psi=\pi_1 \oplus [3]$ where $\pi_1 \in \Pi_{\rm alg}^{\rm o}(\PGL_5)$. In this case 
$w_3=2$. \ps\ps

{\bf Case (viii)} : $d_1=7$, $n_2=1$, i.e. $\psi=[7]\oplus [1]$. This occurs if and only if $\ww=(6,4,2,0)$ and $\pi$ is then the trivial representation of $G$. \ps\ps

\subsubsection{Endoscopic cases of type $(n_1,n_2)=(4,4)$} In this case
$k=2$, $$\psi=\pi_1[d_1]\oplus \pi_2[d_2]$$ and ${\rm C}_\psi \simeq
(\Z/2\Z)^2$.  It follows that ${\rm C}_{\psi}$ is generated by $s_1$ (or $s_2$)
and the center $Z$.  One will have to describe $\rho^\vee(s_1)$ and
$\varepsilon_\psi(s_1)=\varepsilon(\pi_1 \times \pi_2)^{{\rm Min}(d_1,d_2)}$
in each case.  Recall that $\rho^\vee : {\rm C}_{\psi_\infty} \rightarrow
\{\pm 1\}$ is the fundamental character defined in~\S~\ref{mformevenorth}. 
\ps\ps

{\bf Case (ix)} :  ({\it tempered case}) $d_1=d_2=1$, i.e. $\pi_1,\pi_2 \in \Pi_{\rm alg}^{\rm o}(\PGL_4)$. Denote by $a>b$ 
the Hodge weights of $\pi_1$ and by $c>d$ the ones of $\pi_2$. We may assume $a>c$. 
One has $\{a,b,c,d\}=\{w_1,w_2,w_3,w_4\}$.  Moreover $\varepsilon_\psi(s_1)=1$ as all the $d_i$ are $1$ (tempered case), so $m'(\psi)\neq 0$ if, and only if, $\rho^\vee(s_1)=1$, i.e. if $a>c>b>d$. In other words,
$$m'(\psi)=\left\{ \begin{array}{ll} e(\ww) & {\rm if}\, \, a>c>b>d,\\0 & {\rm otherwise}. \end{array}\right.$$
 \ps\ps

{\bf Case (x)} : $d_1=2$, $d_2=1$, i.e. $\psi=\pi_1[2]\oplus \pi_2$ where $\pi_1 \in \Pi_{\rm alg}(\PGL_2)$ and $\pi_2 \in \Pi_{\rm alg}^{\rm o}(\PGL_4)$. If $a$ is the Hodge weight of $\pi_1$ and $b>c$ are the Hodge weights of $\pi_2$ then 
$\{w_1,w_2,w_3,w_4\}=\{a+1,a-1,b,c\}$. One has $$\varepsilon_\psi(s_1)=\varepsilon(\pi_1 \times \pi_2)=(-1)^{{\rm Max}(a,b)+{\rm Max}(a,c)}.$$ On the other hand $\rho^\vee(s_1)=-1$. It follows that 
$$m'(\psi)=\left\{ \begin{array}{ll} e(\ww) & {\rm if}\, \,b>a>c, \\0 & {\rm otherwise}. \end{array}\right.$$
 \ps\ps

{\bf Case (xi)} : $d_1=d_2=2$, i.e. $\psi=\pi_1[2]\oplus \pi_2[2]$ where $\pi_1,\pi_2 \in \Pi_{\rm alg}(\PGL_2)$ have respective Hodge weights $w_1-1$ and $w_3-1$. One has $\varepsilon_\psi(s_1)=\varepsilon(\pi_1 \times \pi_2)=1$ and $\rho^\vee(s_1)=-1$. It follows that 
$$m'(\psi)=0$$
in all cases. \ps\ps

\subsubsection{Endoscopic cases of type $(n_1,n_2,n_3)=(4,3,1)$} In this case $k=3$, $w_4=0$, 
$$\psi=\pi_1[d_1]\oplus \pi_2[d_2]\oplus [1]$$ and ${\rm C}_\psi \simeq (\Z/2\Z)^2$. It follows that
${\rm C}_{\psi}$ is generated by $s_1$ 
and the center $Z$. We have $$\varepsilon_\psi(s_1)=\varepsilon(\pi_1 \times \pi_2)^{{\rm Min}(d_1,d_2)}\varepsilon(\pi_1).$$  
\ps\ps

{\bf Case (xii)} : ({\it tempered case}) $d_1=d_2=1$, i.e. $\pi_1 \in \Pi_{\rm alg}^{\rm o}(\PGL_4)$ and $\pi_2 \in \Pi_{\rm alg}^{\rm o}(\PGL_3)$. 
Denote by $a>b$ the Hodge weights of $\pi_1$ and $c$ the one of $\pi_2$. One has $\{a,b,c\}=\{w_1,w_2,w_3\}$ and $\varepsilon_\psi=1$. The multiplicity is thus nonzero if and only if $\rho^\vee(s_1)=1$, i.e. $a>c>b$: 
$$m'(\psi)=\left\{ \begin{array}{ll} 1 & {\rm if}\, \, a>c>b,\\0 & {\rm otherwise}. \end{array}\right.$$
\ps\ps

{\bf Case (xiii)} : $d_1=2$, $d_2=1$, i.e. $\pi_1 \in \Pi_{\rm alg}(\PGL_2)$ and $\pi_2 \in \Pi_{\rm cusp}^{\rm o}(\PGL_3)$. If $a$ is the Hodge weight of $\pi_1$ and if $b$ is the one of $\pi_2$, then $\{a+1,a-1,b\}=\{w_1,w_2,w_3\}$. One has $\rho^\vee(s_1)=-1$.
On the other hand, $\varepsilon_\psi(s_1)=\varepsilon(\pi_1 \times \pi_2)\varepsilon(\pi_1)=(-1)^{{\rm Max}(a,b)+1}$.  It follows that 
$$m'(\psi)=\left\{ \begin{array}{lll} 1 & {\rm if}\, \,\, b>a, \\
0 & {\rm otherwise}. \end{array}\right.$$
 \ps\ps

{\bf Case (xiv)} : $d_1=1$, $d_2=3$,  i.e. $\psi=\pi_1\oplus [3] \oplus [1]$ with $\pi \in \Pi_{\rm alg}^{\rm o}(\PGL_4)$ of Hodge weights $w_1>w_2$ (here $w_3=2$). We have $\varepsilon_\psi=1$ and $\rho^\vee(s_1)=-1$, so 
$$m'(\psi)=0$$
in all cases.\ps\ps

{\bf Case (xv)} : $d_1=2$, $d_2=3$, i.e. $\psi=\pi_1[2]\oplus [3] \oplus [1]$ with $\pi \in \Pi_{\rm alg}(\PGL_2)$ of Hodge weight $a=w_1-1$.
We have $\varepsilon_\psi(s_1)=\varepsilon(\pi_1)=(-1)^{\frac{a+1}{2}}$ and $\rho^\vee(s_1)=-1$, so 
$$m'(\psi)=\left\{ \begin{array}{ll} 1 & {\rm if}\, \, a \equiv 1 \bmod 4,\\0 & {\rm otherwise}. \end{array}\right.$$

\subsection{Conclusions} \label{conclso8}The inspection of each case above,
and our previous computation of ${\rm S}(w)$, ${\rm S}(w,v)$, ${\rm
S}(w,v,u)$, ${\rm O}^\ast(w)$, ${\rm O}(w,v)$ and ${\rm O}^\ast(w,v)$, allow
to compute the contribution of each endoscopic type except two, namely the
even and odd stable and tempered types.  The contribution of the even-stable
tempered type is exactly $${\rm O}(w_1,w_2,w_3,w_4)$$ when $w_4 \neq 0$, and
$2 \cdot {\rm O}(w_1,w_2,w_3,w_4)$ when $w_4=0$.  The contribution of the
odd-stable tempered type is $${\rm O}^*(w_1,w_2,w_3).$$ This concludes the proof
of Theorem~\ref{mainthmintro}.  The Corollary~\ref{corintroso8} and
Tables~\ref{tableOwvut} and~\ref{tableOwvutbis} follow form these
computations.\ps\ps Let us mention that we also have in our database the
computation of the number of discrete automorphic representations of the
non-connected group ${\rm O}_8$ of any given infinitesimal character.  We
shall not say more about this in this paper however.

\newpage

\section{Description of $\Pi_{\rm disc}({\rm G}_2)$}\label{chapG2}

\subsection{The semisimple definite ${\rm G}_2$ over $\Z$}

Consider the unique semisimple $\Z$-group $G$ of type ${\rm G}_2$ such that
$G(\R)$ is compact, namely the automorphism group scheme over $\Z$ of "the"
ring of Coxeter octonions~(see \cite{coxeter},\cite{blijspringer1},\cite[\S 4]{grossinv}).  We shall simply
write ${\rm G}_2$ for this $\Z$-group $G$.  The reduction map ${\rm G}_2(\Z)
\rightarrow {\rm G}_2(\mathbb{F}_2)$ is an isomorphism and $$|{\rm
G}_2(\Z)|=2^6\cdot 3^3\cdot 7=12096.$$ The $\Z$-group ${\rm G}_2$ admits a
natural homomorphism into the $\Z$-group ${\rm SO}_7$ by its action on the
lattice $L \simeq {\rm E}_7$ of pure Coxeter octonions.  For a well-chosen basis of $L[1/2]$,
it follows from~\cite[\S 4]{cn} that the group ${\rm G}_2(\Z)$ becomes the
subgroup of ${\rm GL}_7(\Z[1/2])$ generated by the two elements
$$\frac{1}{2}\left[\begin{array}{ccccccc} 0& 1& -1& 0& 0 & 1 & -1\\ 0 & -1 &
0 & -1& -1 & 1 & 0 \\ 0 & -1 & 0 & 1 & 1 & 1 & 0 \\ 0 & 1 & 1 & 0 & 0 & 1 &
1 \\
-2 & 0 & 0 & 0 & 0 & 0 & 0 \\ 0 & 0 & 1& 1 & -1 & 0 & -1 \\ 0 & 0 & -1 & 1 & -1 & 0 & 1\end{array}\right], \, \, \, \, \left[\begin{array}{ccccccc} 0& 0& 0& 0& 0& -1& 0\\ 0& 1& 0& 0& 0& 0& 0\\ 0& 0& 1& 0& 0& 0& 0\\ 0& 0& 0& 0& 0& 0& 1\\ 0& 0& 0& 0& 1& 0& 0\\ 1& 0& 0& 0& 0& 0& 0\\ 0& 0& 0& -1& 0& 0& 0\end{array}\right].$$
This allows not only to enumerate (with the computer) all the elements of ${\rm G}_2(\Z) \subset {\rm GL}_7(\Z[1/2])$ but to compute as well their characteristic polynomials. 
The list of the twelve obtained characteristic polynomials, together with the number of elements with that characteristic polynomial, is given in Table~\ref{tableG2charpoly} (we denote by $\Phi_d$ the $d$-th cyclotomic polynomial). One easily checks for instance with this table that 
if $\chi(g,t)$ denotes the characteristic polynomial of $g$ then $$\frac{1}{12096}\sum_{g \in {\rm G}_2(\Z)} \chi(g,t) = t^7 -t^4 +t^3 -1,$$
which is compatible with well-known fact that $\dim(\Lambda^3 L \otimes \C)^{{\rm G}_2(\C)}=1$. 

\begin{table}[htp]

\caption{Characteristic polynomials of the elements of ${\rm G}_2(\Z) \subset {\rm SO}_7(\R)$.}
\begin{tabular}{|c|c|c||c|c|c|}
\hline {\rm Char. Poly.} & $\sharp$ & {\rm Char. Poly.} & $\sharp$ \\
\hline $\Phi_1 \Phi_3^3$ & 56 & $\Phi_1\Phi_3 \Phi_6^2$  & 504 \\
\hline $\Phi_1 \Phi_2^2 \Phi_4^2$ & 378& $\Phi_1^3 \Phi_4^2$ & 378 \\
\hline $\Phi_1^3 \Phi_2^4$ & 315 & $\Phi_1^7$ & 1 \\
\hline $\Phi_1 \Phi_2^2 \Phi_3 \Phi_6$ & 2016 & $\Phi_1 \Phi_3 \Phi_{12}$ & 3024 \\
\hline $\Phi_1 \Phi_2^2 \Phi_8$ & 1512 & $\Phi_1 \Phi_4 \Phi_8$ & 1512 \\
\hline $\Phi_1\Phi_7$ & 1728 & $\Phi_1^3 \Phi_3^2$ & 672 \\
\hline
\end{tabular}
\label{tableG2charpoly}
\end{table}
\ps

\subsection{Polynomial invariants for ${\rm G}_2(\Z) \subset {\rm G}_2(\R)$}\label{parinvG2}
To describe the finite dimensional representations of ${\rm G}_2(\R)$ we fix a maximal torus $T$ and a system of positive roots $\Phi^+$ of $({\rm G}_2(\R),T)$. Let $X={\rm X}^\ast(T)$, $X^\vee={\rm X}_\ast(T)$ and denote by $\langle\cdot,\cdot\rangle$ the canonical perfect pairing between them.\ps 
Let $\alpha,\beta \in X$ the simple roots in $\Phi^+$ where $\alpha$ is short and $\beta$ is long. The positive roots are thus $$\alpha,\beta,\beta+\alpha,\beta+2\alpha,\beta+3\alpha,2\beta+3\alpha,$$ where 
$\alpha,\beta+\alpha$ and $\beta+2\alpha$ are short, and $X=\Z \alpha \oplus \Z \beta$.
The inverse root system is again of type ${\rm G}_2$, with simple positive roots $\alpha^\vee,\beta^\vee \in X^\vee$ with $\alpha^\vee$ long, and where 
$$\langle \alpha, \beta^\vee \rangle = -1\, \,\,{\rm and}\,\,\, \langle \beta, \alpha^\vee \rangle = -3.$$ It follows that the dominant weights are the $a \alpha + b \beta$ where $a,b \in \Z$ satisfy 
$2b \geq a\geq 3b/2$. The fundamental representations with respective fundamental weights $$\omega_1=2\alpha+\beta, \, \, \, \, \,  \omega_2=3\alpha+2\beta$$ will be denoted by $V_7$ and $V_{14}$, because of their respective 
dimension $7$ and $14$. One easily checks that $V_7=L\otimes \C$ and $V_{\rm 14}$ is the adjoint representation. The half-sum of positive roots is $\rho=5 \alpha + 3 \beta=\omega_1+\omega_2$. 

\begin{definition}\label{defww} If $w>v$ are even non-negative integers, we denote by $U_{w,v}$
the irreducible representation of ${\rm G}_2(\R)$ with highest weight $$\frac{w-v-2}{2} \,\omega_1 + \frac{v-2}{2} \, \omega_2.$$ 
We also denote by $\Pi_{w,v}({\rm G}_2)$ the subset of $\pi \in \Pi_{\rm disc}({\rm G}_2)$ such that $\pi_\infty \simeq U_{w,v}$, and set $m(w,v)=\sum_{\pi \in \Pi_{w,v}} m(\pi)$. 
\end{definition}

This curious looking numbering has the following property. If 
$$\varphi : {\rm W}_\R \longrightarrow \widehat{{\rm G}_2}$$ is the Langlands parameter of $U_{w,v}$, and if 
$\rho_7 : \widehat{{\rm G}_2} \rightarrow {\rm SO}_7(\C)$ is the $7$-dimensional irreducible
representation of $\widehat{{\rm G}_2}$, then $\rho_7 \circ \varphi$ is the representation 
$${\rm I}_{w+v} \oplus {\rm I}_w \oplus {\rm I}_v \oplus \varepsilon.$$
Indeed, the weights of $\rho_7$ are $0,\pm \beta^\vee,\pm (\alpha^\vee+\beta^\vee), \pm(\alpha^\vee+2\beta^\vee)$. \ps

Observe that $\rho_7 \circ \varphi$ determines the equivalence class of $\varphi$. This is a special case of the fact that the conjugacy class of any element $g \in {\rm G}_2(\R)$ (resp. of any semisimple element in ${\rm G}_2(\C)$) is 
uniquely determined by its characteristic polynomial in $V_7$. Indeed, this follows from the identity
$$V_{14} \oplus V_7 \simeq \Lambda^2 V_7.$$
This property makes the embedding ${\rm G}_2(\C) \subset {\rm SO}_7(\C)$ quite suitable to study ${\rm G}_2$ and its subgroups. In particular, Table~\ref{tableG2charpoly} leads to a complete determination 
of the semisimple conjugacy classes in ${\rm G}_2(\R)$ of the elements of ${\rm G}_2(\Z)$, which is the ingredient we need to apply the method of \S~\ref{appnum}.\ps

Van der Blij and Springer have shown in~\cite{blijspringer1} that $|{\rm Cl}({\rm G}_2)|=1$
(se also~\cite[\S 5]{grossinv}), it
follows that $$m(w,v)=\dim U_{w,v}^{{\rm G}_2(\Z)}.$$ See
Table~\ref{tableg2nue} for a sample of computations.  As we shall see below,
one should have $m(\pi)=1$ for each $\pi \in \Pi_{\rm disc}(G)$, and we thus
expect that $m(w,v)=|\Pi_{w,v}({\rm G}_2)|$.  \ps

Automorphic forms for the $\Q$-group ${\rm G}_2$ have been previously studied by Gross, Lansky, Pollack and Savin: see~\cite{grosssavin}, \cite{grosspollack}, \cite{lanskypollack} and~\cite{pollack}. Although most of the automorphic forms studied by those authors are Steinberg at one finite place, they may be trivial at the
infinite place. Pollack and Lansky are also able to compute some Hecke eigenvalues in some cases.

\subsection{Endoscopic classification of $\Pi_{\rm disc}({\rm
G}_2)$}\label{endoG2} We recall Arthur's conjectural description of
$\Pi_{\rm disc}({\rm G}_2)$, following his general conjecture
in~\cite{arthurunipotent}.  {\it Most of the results here will thus be
conditional to the existence of the group $\mathcal{L}_\Z$ discussed in
Appendix~\ref{parsatotate} and to these conjectures}, that we will make
explicit.  All the facts stated below about the structure of ${\rm G}_2$ can
be simply checked on its root system.  We refer to~\cite{gangurevich1}
for a complete analysis of Arthur's conjectures for the split groups of type
${\rm G}_2$, in a much greater generality than we actually need here, and
for a survey of the known results.  \ps

A global discrete Arthur parameter for the $\Z$-group ${\rm G}_2$ is a $\widehat{{\rm G}_2}$-conjugacy class 
of morphisms $$\psi : \mathcal{L}_\Z \times {\rm SL}_2(\C) \longrightarrow \widehat{{\rm G}_2}$$
such that: \ps\ps

(a) ${\rm Im}\, \psi$ has a finite centralizer in $\widehat{{\rm G}_2}$,\ps 
(b) $\psi_\infty=\psi_{|{\rm W}_\R \times {\rm SL}_2(\C)}$ is an Adams-Johnson parameter for ${\rm G}_2(\R)$
(see Appendix \ref{appendixadamsjohnson}). \ps\ps
Observe that by property (b) the centralizer ${\rm C}_{\psi_\infty}$ of
${\rm Im}\, \psi_\infty$ in $\widehat{{\rm G}_2}$ is an elementary abelian $2$-group, hence so is the centralizer ${\rm C}_\psi \subset {\rm C}_{\psi_\infty}$ of $\psi$. 
As $\mathcal{L}_\Z$ is connected, observe that the Zariski-closure of ${\rm Im} \psi$ is a connected complex reductive subgroup of $\widehat{{\rm G}_2}$. \ps

This severely limits the possibilities for ${\rm Im}\,\psi$. Up to conjugacy there are exactly $3$ connected complex reductive subgroups 
of $\widehat{{\rm G}_2}$ whose centralizers are elementary abelian $2$-groups: \begin{itemize}

\item[(i)] the group $\widehat{{\rm G}_2}$ itself, with trivial centralizer,\ps
\item[(ii)] a principal $\PGL_2(\C)$ homomorphism, again with trivial centralizer,\ps
\item[(iii)]  the centralizer $H_s \simeq {\rm SO}_4(\C)$ of an element $s$ of order $2$, whose centralizer is the center $\langle s \rangle$ of $H_s$.\ps
\end{itemize}

Recall that up to conjugacy there is a unique element $s$ of order $2$ in $\widehat{{\rm G}_2}$. The isomorphism $H_s \simeq {\rm SO}_4(\C)$ in (iii) is actually canonical up to inner automorphisms as $H_s$ is its own normalizer in $\widehat{{\rm G}_2}$. 
Indeed, one has two distinguished injective homomorphisms ${\rm SL}_2(\C) \rightarrow H_s$, one of which being a short radicial ${\rm SL}_2(\C)$ and the other one being a long radicial ${\rm SL}_2(\C)$ (the long and short roots being orthogonal)\ps

We shall need some facts about the restrictions of $V_7$ and $V_{14}$ to these groups. We denote by $\nu_{\rm long}$ and $\nu_{\rm short}$ the two $2$-dimensional irreducible representations of $H_s$ which are respectively non-trivial on the long and short ${\rm SL}_2(\C)$ inside $H_s$.

\begin{lemme}\label{lemmerepg2} Let $s \in \widehat{{\rm G}_2}$ be an element of order $2$.
 \begin{itemize} \item[(i)] $(V_7)_{|H_s} = \nu_{\rm long} \otimes \nu_{\rm short} \oplus {\rm Sym}^2 \nu_{\rm short}$,\ps 
\item[(ii)] $(V_{14})_{|H_s}={\rm Sym}^2 \nu_{\rm long}\oplus {\rm Sym}^2 \nu_{\rm short} \oplus {\rm Sym}^3 \nu_{\rm short} \otimes \nu_{\rm long}$.\end{itemize}

\noindent Moreover, the restriction of $V_7$ to a principal ${\rm PGL}_2(\C)$ is isomorphic to $\nu_7={\rm Sym}^6(\C^2)$.
\end{lemme}
\ps

If $\psi$ is a global Arthur parameter for ${\rm G}_2$, then $\rho_7 \circ \psi$ actually defines a 
global Arthur parameter
for ${\rm Sp}_6$, that we shall denote $\psi^{\rm SO}$. The previous lemma and discussion show that the equivalent class of $\psi^{\rm SO}$ determines the conjugacy class of
$\psi$. \ps\ps

Fix a global Arthur parameter $\psi$ as above. We denote by $$\pi(\psi) \in \Pi({\rm G}_2)$$
the unique representation $\pi$ such that $c(\pi)$ is associated to $\psi$  
by the standard Arthur recipe. Explicitly, for each prime $p$ we have 
$c_p(\pi)=\psi({\rm Frob}_p \times e_p)$ (see~\S~\ref{arthurparam}), and
 $\pi_\infty$ is the 
unique representation of ${\rm G}_2(\R)$ whose infinitesimal character is the 
one of the Langlands parameter $\varphi_{\psi_\infty}$ (assumption (b) 
on $\psi$). Arthur's conjectures describe $\Pi_{\rm disc}({\rm G}_2)$ as follows. First, any 
$\pi \in \Pi_{\rm disc}({\rm G}_2)$ should be of the form $\pi(\psi)$ for a unique $\psi$ satisfying (a) and (b). 
Second, they describe $m(\pi(\psi))$ for each $\pi$ as follows. \ps\ps

{\bf Case (i)} : ({\it stable tempered cases}) $\psi^{\rm SO} \in \Pi_{\rm
alg}^{\rm o}(\PGL_7)$. This is when ${\rm C}_\psi=1$ and $\psi({\rm
SL}_2(\C))=1$. In this case $$m(\pi)=1.$$
 By
Prop.\ref{propsatotate}, a $\pi \in \Pi_{\rm alg}^{\rm o}(\PGL_7)$ has the form
$\psi^{\rm SO}$ for a stable tempered $\psi$ if and only if 
$c(\pi) \in \rho_7(\mathcal{X}(\widehat{G}_2(\C)))$. It is equivalent to ask that $c(\pi) \times  1$, viewed as an element in $\mathcal{X}({\rm SO}_8(\C))$ is invariant by a triality automorphism. 
Moreover, ${\rm Im}\, \psi$ is either isomorphic to the compact group ${\rm {\rm G}_2}$ or
to ${\rm SO}(3)$. The latter case occurs if and only if $\pi(\psi)_\infty \simeq U_{w,v}$ 
where $v \equiv 2 \bmod 4$ and $w=2v$, in which case it occurs exactly
${\rm S}(v/2)$ times. \ps \ps

{\bf Case (ii)} : ({\it stable non-tempered case}) $\psi^{\rm SO}=[7]$. 
Then $\pi$ is the trivial representation, the unique element in
$\Pi_{4,2}({\rm G}_2)$. \ps\ps

There are three other cases for which ${\rm Im}(\psi)=H_s$. In those cases we have $${\rm
C}_\psi=\langle s \rangle \simeq \Z/2\Z.$$ Arthur's multiplicity formula
requires two ingredients. The first one is the character 
$$\varepsilon_{\psi} : {\rm C}_\psi \rightarrow \C^\times$$
given by Arthur's general recipe~\cite{arthurunipotent}. This character is trivial if  $\psi({\rm SL}_2(\C))=1$. Otherwise there are two distinct cases: \begin{itemize}
\item[(i)] $\nu_{\rm short} \circ  \psi_{|{\rm SL}_2(\C)}=\nu_2$ and $\nu_{\rm long} \circ  \psi_{|\mathcal{L}_\Z} = r(\pi)$ for some $\pi \in \Pi_{\rm alg}(\PGL_2)$ (see Appendix \ref{parsatotate}). 
Then $$\varepsilon_\psi(s)=\varepsilon(\pi)=(-1)^{(w+1)/2}$$ where $w$ is the Hodge weight of $\pi$. \ps
\item[(ii)] $\nu_{\rm long} \circ  \psi_{|{\rm SL}_2(\C)}=\nu_2$ and $\nu_{\rm short} \circ  \psi_{|\mathcal{L}_\Z} = r(\pi)$ for some $\pi \in \Pi_{\rm alg}(\PGL_2)$. Then 
$\varepsilon_\psi(s)=\varepsilon({\rm Sym}^3\pi)$. If $w$ is the Hodge weight of $\pi$, observe that 
$$\varepsilon({\rm Sym}^3\pi )=(-1)^{(w+1)/2+(3w+1)/2}=-1$$
for each $w$, thus $\varepsilon_\psi$ is the non-trivial character in this case.\ps
\end{itemize}
Observe that the a priori remaining case $\psi(\mathcal{L}_\Z)=1$ does not occur as property (b) is not satisfied for such a $\psi$ (the infinitesimal character $z_{\psi_\infty}$ is not regular). \ps\ps

The second ingredient is the restriction to ${\rm C}_\psi$ of the character $\rho^\vee : {\rm C}_{\psi_\infty} \rightarrow \C^\times$. The multiplicity formula will then take the form: 
$m(\pi)=1$ if $\rho^\vee(s)=\varepsilon_\psi(s)$ and $m(\pi)=0$ otherwise. \ps

In order to compute $\rho^\vee(s)$ we fix $\widehat{T}$ a maximal torus in $\widehat{{\rm G}_2}$ such that ${\rm X}^\ast(\widehat{T})=X^\vee$. Observe that the centralizer $T'$ of 
$\rho_7(\widehat{T})$ in ${\rm SO}_7(\C)$ is a maximal torus of the latter group. We consider the standard root system $\Phi'$ for 
$({\rm SO}_7(\C),T')$ recalled in \S~\ref{appnum}, in particular ${\rm X}^\ast(T')=\oplus_{i=1}^3 \Z e_i$. Then $\Phi^{\vee}=\Phi'_{|\widehat{T}}$ is a root system for $(\widehat{T},\widehat{{\rm G}_2})$ with positive roots $(\Phi^\vee)^+=(\Phi^{'})^+_{|\widehat{T}}$: 
up to conjugating $\rho_7$ we may thus assume that $(\Phi^\vee)^+$ is the positive root system of \S~\ref{parinvG2}. 

\begin{lemme} Under the assumptions above, we have $\rho^\vee(s)=e_2(\rho_7(s))$. \end{lemme}

\begin{pf} Under the assumptions above, if $\lambda \in {\rm X}_\ast(\widehat{T})$ is such that $\rho_7(\lambda)$ is $(\Phi')^+$-dominant then $\lambda$ is $(\Phi^\vee)^+$-dominant. 
We have already seen that the respective restriction to $\widehat{T}$ of $e_1,e_2,e_3$ are the elements $2\beta^\vee+\alpha^\vee$, $\beta^\vee+\alpha^\vee$ and $\beta^\vee$.
The lemma follows from the identity $$\rho^\vee = 5 \beta^\vee + 3 \alpha^\vee \equiv (e_2)_{|\widehat{T}} \bmod 2{\rm X}^\ast(\widehat{T}).$$\end{pf}

\noindent We can now make explicit the three remaining multiplicity formulae.\ps\ps

{\bf  Case (iii) :}  ({\it tempered endoscopic case}) $\psi^{\rm SO}=\pi_{\rm long} \otimes \pi_{\rm short} \oplus {\rm Sym}^2 \pi_{\rm short}$ where $\pi_{\rm short},\pi_{\rm long} \in \Pi_{\rm alg}(\PGL_2)$ have respective Hodge weights 
$w_{\rm short},w_{\rm long}$. Of course, ${\rm Sym}^2\pi_{\rm short} \in \Pi_{\rm alg}^{\rm o}(\PGL_3)$ has Hodge weight $2 w_{\rm short}$ and $\pi_{\rm long} \otimes \pi_{\rm short} \in \Pi_{\rm alg}^{\rm o}(\PGL_4)$ has Hodge weights $w_{\rm short}+w_{\rm long}$ and $|w_{\rm short}-w_{\rm long}|$. We also have $\varepsilon_\psi(s)=1$. But by Lemma~\ref{lemmerepg2} (i) we have $e_2(\rho_7(s))=1$ if and only if 
$$w_{\rm long}+w_{\rm short} > 2 w_{\rm short} > w_{\rm long} - w_{\rm short},$$
thus $m(\pi)=1$ in this case and $m(\pi)=0$ otherwise.
\ps\ps

{\bf Case (iv)} : ({\it non-tempered endoscopic case $1$})  $\psi^{\rm SO}=\pi[2] \oplus {\rm Sym^2 \pi}$ where  $\pi \in \Pi_{\rm alg}({\rm PGL}_2)$, say with 
Hodge weight $w$. We have seen that in this case $\varepsilon_\psi(s)=-1$. On the other hand $e_2(\rho_7(s))=-1$ if and only if $w-1<2w$, which is always satisfied as $w>1$. 
Arthur's multiplicity formula tells us that $$m(\pi)=1$$ in all cases. \ps\ps

{\bf Case (v) :} ({\it non-tempered endoscopic case $2$}) $\psi^{\rm SO}=\pi[2] \oplus [3]$ where $\pi \in \Pi_{\rm alg}({\rm PGL}_2)$, say with Hodge weight $w$. 
We have seen that $\varepsilon_\psi(s)=\varepsilon(\pi)=(-1)^{(w+1)/2}$. Observe that $e_2(\rho_7(s))=-1$ as $w-1>3$. 
Arthur's multiplicity formula tells us then that $$m(\pi)=1 \Leftrightarrow w \equiv 1 \bmod 4.$$
\ps
Let us mention that the multiplicity formula for the Arthur's packets appearing in case (v) has been established for the split groups of type ${\rm G}_2$ in~\cite{gangurevich2}.

\subsection{Conclusions}\label{conclg2}The inspection of each case above and the well-known formula for ${\rm S}(w)$ allow to compute the conjectural number 
${\rm G}_2(w,v)$ of $\pi \in \Pi_{\rm alg}^{\rm o}(\PGL_7)$ such that $c(\pi) \in \rho_7(\mathcal{X}(\widehat{{\rm G}_2}))$ and with Hodge weights $w+v>w>v$. Concretely, 
$${\rm G}_2(w,v)=m(w,v)-\delta_{w=4}-{\rm O}^\ast(w)\cdot {\rm O}(w+v,v)$$ $$-\delta_{w-v=2}\cdot {\rm S}(w-1)-\delta_{v=2}\cdot \delta_{w \equiv 0 \bmod 4} \cdot {\rm S}(w+1).$$
See Table~\ref{tableStG2} for 
a sample of results when $w+v\leq 58$. The Sato-Tate group of each of the associated $\pi$ is conjecturally the compact group of type ${\rm G}_2 \subset {\rm SO}_7(\R)$ (rather than ${\rm SO}_3(\R)$ principally embedded in the latter) :  this follows from Prop.~\ref{propsatotate} as the motivic weight of ${\rm Sym}^6 \pi$ for $\pi \in \Pi_{\rm alg}(\PGL_2)$ is at least $66>58$.

\newpage
\section{Application to Siegel modular forms}\label{finalapp}

\subsection{Vector valued Siegel modular forms of level $1$}
We consider in this chapter the classical Chevalley $\Z$-group ${\rm Sp}_{2g}$, whose dual group is ${\rm SO}_{2g+1}(\C)$. Let $$\ww=(w_1,w_2,\cdots,w_g)$$ where
the $w_i$ are even positive integers such that $w_1>w_2>\cdots>w_g$. To such a $\ww$ we may associate a semisimple conjugacy class $z_\ww$ in $\mathfrak{so}_{2g+1}(\C)$,  
namely the class with eigenvalues $\pm \frac{w_i}{2}$ for $i=1,\cdots, g$, and $0$. Recall that for any such $\ww$, 
there is an ${\rm L}$-packet of discrete series with infinitesimal characters $z_\ww$, and that this ${\rm L}$-packet contains two "holomorphic" discrete series which are outer conjugate by ${\rm PGSp}_{2g}(\R)$. 
We make once and for all a choice for the holomorphic ones (hence for the anti-holomorphic as well).\ps\ps

Recall the space
				 $${\rm S}_{\ww}({\rm Sp}_{2g}(\Z))$$
of holomorphic vector valued Siegel modular forms with infinitesimal character $\mathfrak{z}_{\ww}$. If $(\rho,V)$ is the irreducible representation of $\GL_g(\C)$ with standard 
highest weight $m_1 \geq m_2 \geq \cdots\geq  m_g$, and if $m_g>g$, recall that a $(\rho,V)$-valued Siegel modular form has infinitesimal character $z_\ww$ where 
$\ww=(w_i)$ and $w_i=2(m_i-i)$ for each $i=1,\cdots,g$ (see e.g.~\cite[\S 4.5]{asgarischmidt}). Denote also
$$\Pi_{\ww}({\rm Sp}_{2g})$$ the set of $\pi \in \Pi_{\rm disc}({\rm Sp}_{2g})$ such that $\pi_\infty$ is the holomorphic discrete series with infinitesimal character $z_\ww$. Such a $\pi$ is tempered at the infinite place, thus it follows from a result of Wallach~\cite[Thm. 4.3]{wallach} that $\Pi_{\ww}({\rm Sp}_{2g}) \subset \Pi_{\rm cusp}({\rm Sp}_{2g})$. By Arthur's multiplicity formula, $m(\pi)=1$ for each 
$\pi \in \Pi_{\rm disc}({\rm Sp}_{2g})$, it is well-known that this implies 
	$$\dim {\rm S}_{\ww}({\rm Sp}_{2g}(\Z))=|\Pi_{\ww}({\rm Sp}_{2g})|.$$ 
By Lemma~\ref{lemmeclass}, Arthur's multiplicity formula allows to express $|\Pi_{\ww}({\rm Sp}_{2g})|$ in terms of various ${\rm S}(-)$, ${\rm O}(-)$ and ${\rm O}^\ast(-)$. We shall give now the two ingredients needed to make this computation in general and we shall apply them later in the special case $g=3$. \ps

\subsection{Two lemmas on holomorphic discrete series} \label{parholdisc}
Let $$\varphi_\ww : {\rm W}_\R \rightarrow {\rm SO}_{2g+1}(\C)$$ be the discrete series Langlands parameter with infinitesimal character $z_\ww$, and let $$\Pi(\varphi_\ww)$$ be the 
associated ${\rm L}$-packet of discrete series representations of ${\rm Sp}_{2g}(\R)$ with infinitesimal character $z_\ww$. Recall that the centralizer of $\varphi_\ww({\rm W}_\C)$ in ${\rm SO}_{2g+1}(\C)$ is a maximal torus $\widehat{T}$ in ${\rm SO}_{2g+1}(\C)$ and that the centralizer 
${\rm C}_{\varphi_\ww}$ of $\varphi_\ww({\rm W}_\R)$ is the $2$-torsion subgroup of $\widehat{T}$. There is also a unique Borel subgroup $\widehat{B} \supset \widehat{T}$ for which the element $\lambda \in {\rm X}_\ast(\widehat{T})[1/2]$ such that $\varphi_\ww(z)=(z/\overline{z})^\lambda$ for all $z \in {\rm W}_\C$ is dominant with respect to $\widehat{B}$. \ps

We consider the setting and notations of~\S\ref{strongforms} with $G={\rm Sp}_{2g}(\C)$. The strong forms $t \in \mathcal{X}_1(T)$ such that $G_t \simeq {\rm Sp}_{2g}(\R)$ are the ones such that $$t^2=-1,$$ and they form a single $W$-orbit. Let us fix an isomorphism between ${\rm Sp}_{2g}(\R)$ and $G_{[t]}$ for $t$ in this $W$-orbit, which thus identify $\Pi(\varphi_\ww)$ with $\Pi(\varphi_\ww,G_{[t]})$ (\S\ref{ajpackets}). This being done,  Shelstad's parameterization gives a canonical injective map (see~\S\ref{parammap})
$$\tau : \Pi(\varphi_\ww) \rightarrow \Hom({\rm C}_{\varphi_\ww},\C^\times).$$ 
(We may replace ${\rm S}_{\varphi_\ww}$ by ${\rm C}_{\varphi_\ww}$ in the range as ${\rm Sp}_{2g}(\R)$ is split, see Cor.\ref{factorrho}). Our first aim is to determine the image of $\pi_{\rm hol}$ and $\pi_{\rm ahol}$, namely the holomorphic and anti-holomorphic discrete series in $\Pi(\varphi_\ww)$. \ps
 
For a well-chosen $\Z$-basis $(e_i)$ of ${\rm X}^\ast(\widehat{T})$, the positive roots of $({\rm SO}_{2g+1}(\C),\widehat{B},\widehat{T})$ are $\{e_i,i=1,\cdots, g\} \cup \{e_i \pm e_j, 1\leq i< j \leq g\}$ as in~\S~\ref{appnum}.
Let $(e_i^\ast) \in {\rm X}^\ast(T)$ denote the dual basis of $(e_i)$. The set of positive roots of $({\rm Sp}_{2g}(\C),B,T)$ dual to the positive root system above is the set $\{2e_i^\ast, i=1,\cdots,g\} \cup \{e_i^\ast \pm e_j^\ast, 1\leq i<j \leq g\}$. If $t \in T$ we also write $t=(t_i)$ where $t_i=e_i^\ast(t)$ for each $i=1,\cdots,g$.\ps

\begin{lemme}\label{calchol} The Shelstad characters of $\pi_{\rm hol}$ and $\pi_{\rm ahol}$ are the restrictions to ${\rm C}_{\varphi_\ww}$ of the following elements of ${\rm X}^\ast(\widehat{T})$: $$e_1+e_3+e_5+\cdots+e_{2[(g-1)/2]+1} \, \, \, {\rm and} \, \, \, e_2+e_4+e_6+\cdots+e_{2[g/2]}.$$
\end{lemme}

\begin{pf} Let $t \in \mathcal{X}_1(T)$ such that $t^2=-1$. Recall that ${\rm Int}(t)$ is a Cartan involution of $G_t$ and that $K_t$ is the associated maximal compact subgroup of $G_t \simeq {\rm Sp}_{2g}(\R)$. Let $\mathfrak{g}=\mathfrak{k}_t \oplus \mathfrak{p}$ the Cartan decomposition relative to ${\rm Int}(t)$. We have 
$\mathfrak{p}=\mathfrak{p}_+\oplus \mathfrak{p}_-$ where $\mathfrak{p}_{\pm} \subset \mathfrak{p}$ 
are two distinct irreducible $K_t$-submodules for the adjoint action. As a general fact, the representation $\pi_t(\lambda)$ is a holomorphic or anti-holomorphic  discrete series of $G_t$ if and only if $\mathfrak{b}$ is included in either $\mathfrak{k}_t\oplus \mathfrak{p}_{+}$ or in $\mathfrak{k}_t\oplus \mathfrak{p}_{-}$. In those cases, $\mathfrak{k}_t$ is thus a standard Levi subalgebra of $(\mathfrak{g},\mathfrak{b},\mathfrak{t})$ isomorphic to $\mathfrak{gl}_g$. There is a unique such Lie algebra, namely the one with positive roots the $e_i^\ast-e_j^\ast$ for $i<j$.  It follows that $\pi_t(\lambda)$ is a holomorphic discrete series if and only if the positive roots of $T$ in $\mathfrak{k}_t$ are the $e_i^\ast-e_j^\ast$ for each $1\leq i<j \leq g$, i.e. if $t_i=t_j$ for $j\neq i$. As $t^2=-1$, the two possibilities are thus the elements 
$$t_+=(i,i,\dots,i) \, \, \, \, {\rm and} \, \, \, \, t_-=(-i,-i,\dots,-i).$$
We have $t_b=e^{i\pi \rho^\vee}=(i^{2g-1},\dots,i,-i,i)$ (see~\S\ref{parammap}), so $t_\pm t_b^{-1}=\pm(\dots, -1,1,-1,1)$. Let 
$\mu=e_1+e_3+\dots$ and $\mu'=e_2+e_4+\dots$ be the two elements of ${\rm X}_\ast(T)={\rm X}^\ast(\widehat{T})$ given in the statement. One concludes as
$$e^{i\pi\mu}=(-1,1,-1,\dots) \,\,\, \, {\rm and}\, \, \, \, e^{i\pi \mu'}=(1,-1,1,\dots).$$
\end{pf}

The second ingredient is to determine which Adams-Johnson packets $\Pi(\psi)$ of ${\rm Sp}_{2g}(\R)$ contains a holomorphic discrete series.

\begin{lemme}\label{crithol} Let $\psi : {\rm W}_\R \times {\rm SL}_2(\C) \rightarrow {\rm SO}_{2g+1}(\C)$ be an Adams-Johnson parameter for ${\rm Sp}_{2g}(\C)$. Then $\widetilde{\Pi}(\psi)$ contains discrete series of ${\rm Sp}_{2g}(\R)$ if and only if the underlying representation of ${\rm W}_\R \times {\rm SL}_2(\C)$ on $\C^{2g+1}$ does not contain any $1 \otimes \nu_q$ or $\varepsilon \otimes \nu_q$ where $q>1$. \ps
Furthermore, if $\psi$ has this property then 
the holomorphic and anti-holomorphic discrete series of ${\rm Sp}_{2g}(\R)$ belong to $\Pi(\psi)$. 
\end{lemme}

\begin{pf} Let $T,B,L$ and $\lambda$ be attached to $\psi$ as in \S\ref{ajparam} and~\S\ref{ajpackets}, recall that $L \subset {\rm Sp}_{2g}(\C)$ is a Levi factor of a parabolic subgroup. From the last example of~\S\ref{ajparam}, from which we take the notations, we have 
$$L \simeq {\rm Sp}_{d-1}(\C) \times \prod_{i \neq 0} {\rm GL}_{d_i}(\C).$$
Moreover, the underlying representation of ${\rm W}_\R \times {\rm SL}_2(\C)$ on $\C^{2g+1}$ does not contain any $1 \otimes \nu_q$ or $\varepsilon \otimes \nu_q$ where $q>1$ if and only if $d=1$. On the other hand,  
$\Pi(\psi)$ contains a discrete series $\pi_{t}(\lambda)$ of ${\rm Sp}_{2g}(\R)$ if and only if there is a $t \in \mathcal{X}_1(T)\cap {\rm Z}(L)$ such that $t^2=-1$, by Lemma~\ref{caraccompact}. As ${\rm Sp}(d-1,\C)$ does not contain
any element of square $-1$ in its center for $d>1$, the first assertion follows. \ps
Assume now that $L$ is a product of general complex linear groups. It is equivalent to ask that the positive roots of $L$ with respect to 
$(B,T)$ are among the $e_i^\ast-e_j^\ast$ for $i<j$. In particular, the element $t_0=\pm (i,i,\cdots,i) \in \mathcal{X}_1(T)$ is in the center of $L$, thus 
$\pi_{t_0,B}(\lambda) \in \widetilde{\Pi}(\psi)$ by Lemma~\ref{caraccompact} and Lemma~\ref{reductionsd}. But we have seen in the proof of Lemma~\ref{calchol} that this is a holomorphic/anti-holomorphic discrete series.
\end{pf}

The difference between $\pi_{\rm hol}$ and $\pi_{\rm ahol}$ is not really meaningful four our purposes, and we 
will not need to say exactly which of the two characters in Lemma~\ref{calchol} corresponds to e.g. $\pi_{\rm hol}$ (of course this would be possible if we had defined $\pi_{\rm hol}$ more carefully). More importantly, 
let $$\chi=\sum_{i=1}^g e_i$$ be the sum of the two elements of the statement of Lemma~\ref{calchol}. Fix $\psi=(k,(n_i),(d_i),(\pi_i)) \in \Psi_{\rm glob}({\rm Sp}_{2g})$ with infinitesimal character $z_\ww$, one has  canonical embeddings $${\rm C}_\psi \subset {\rm C}_{\psi_\infty} \subset  {\rm C}_{\varphi_\ww}$$ 
by~\S\ref{defpipsi} and~\S\ref{parammap}, as $(\psi_{\infty})_{\rm disc}=\varphi_\ww$. Recall that ${\rm C}_\psi$ is generated by elements $s_i$ as in~\S\ref{paragraphepsilon}. 

\begin{lemme}\label{verifsympl} For each $i=1,\cdots,k$ such that $n_i$ is even we have $\chi(s_i)=1$. \end{lemme}

\begin{pf} Indeed, it follows from Lemma~\ref{lemmeclass} that if $n_i$ is even then $n_i \equiv 0 \bmod 4$. 
\end{pf}

Assume now that $\Pi_\infty(\psi)$ contains holomorphic discrete series, i.e. that $n_i \neq d_i$ for each $i$ such that $n_i>1$ by Lemma~\ref{crithol}. In this case it follows from Lemma~\ref{reductionsd} that the characters of 
$\pi_{\rm hol}$ and $\pi_{\rm ahol}$ viewed as elements of $\Pi_\infty(\psi)$ are again the two characters of Lemma~\ref{calchol}. But it follows then from the Lemma~\ref{verifsympl} above
that the multiplicity formula is the same for the two $\pi \in \Pi(\psi)$ such that $\pi_\infty$ is either holomorphic or anti-holomorphic.  \ps

To conclude this paragraph let us say a word about the choice of the isomorphism that we fixed between ${\rm Sp}_{2g}(\R)$ and $G_{[t]}$ (for $t^2=-1$), which allowed to fix the parameterization $\tau$.
Consider for this the order $2$ outer automorphism of ${\rm Sp}_{2g}(\R)$ obtained as the conjugation by any element of  ${\rm GSp}_{2g}(\R)$ with similitude factor $-1$. It defines in particular element defines an involution 
of $\Pi(\varphi_\ww,G_{[t]})$ and we want to check the effect of this involution on Shelstad's parameterization. The next lemma shows that it is quite benign.

\begin{lemme}\label{chgtiso} If $\pi \in \Pi(\varphi_\ww)$, then $\tau(\pi \circ \theta)=\tau(\pi)+\chi$.
\end{lemme}

\begin{pf} Fix some $t$ such that $t^2=-1$ and view $\theta$ as an outer automorphism of 
$G_{t}$. A suitable representant of $\theta$ in ${\rm Aut}(G_t)$ preserves $(K_t,T_c)$, and the automorphism of $T_c$ obtained this way is well-defined up to ${\rm W}(K_t,T_c)$. It is a simple exercise
to check that it coincides here with the class of the inversion $t \mapsto t^{-1}$ of $T_c$. As $-1 \in {\rm W}(G,T)$, it follows that $\pi_t(\lambda) \circ \theta = \pi_{t^{-1}}(\lambda)$. In other words, 
$\tau_0(\pi \circ \theta)=-\tau_0(\pi)$. As $\tau(\pi)=\tau_0(\pi)-\rho^\vee$, it follows that $$\tau(\pi \circ \theta)=-\tau(\pi)-2\rho^\vee.$$
But observe that $2\rho^\vee = \chi \bmod 2{\rm X}^\ast(\widehat{T})$. As ${\rm C}_{\varphi_\ww}$ is an elementary abelian $2$-group, the lemma follows. \end{pf}

It follows then from Lemma~\ref{verifsympl} that the choice of our isomorphism has no effect on the multiplicity formula for the $\pi \in \Pi(\psi)$ such that $\pi_\infty$ is either holomorphic or anti-holomorphic.

\subsection{An example: the case of genus $3$} 

We shall now describe the endoscopic classification of $\Pi_{\ww}({\rm Sp}_6)$ for any $\ww=(w_1,w_2,w_3)$. As an application, we will deduce in particular the following proposition stated in the introduction.

\begin{propcondb} $\dim {\rm S}_{w_1,w_2,w_3}({\rm Sp}_6(\Z))={\rm O}^\ast(w_1,w_2,w_3)+ {\rm O}(w_1,w_3)\cdot {\rm O}^\ast(w_2)$
$$+\delta_{w_2 \equiv 0 \bmod 4}\cdot ( \delta_{w_2=w_3+2}\cdot {\rm S}(w_2-1)\cdot {\rm O}^\ast(w_1)+ \delta_{w_1=w_2+2}\cdot {\rm S}(w_2+1)\cdot {\rm O}^\ast(w_3)).$$
\end{propcondb}
\ps\ps

Let us fix a $$\psi=(k,(n_i),(d_i),(\pi_i)) \in \Psi_{\rm alg}({\rm Sp}_6)$$ with infinitesimal character $\mathfrak{z}_\ww$. We have to determine first whether or not $\Pi(\psi_\infty)$ contains a holomorphic discrete series. Lemma~\ref{crithol} ensures that it is the case if and only if 
for each $i$ such that $\pi_i= 1$ then $d_i=1$. We thus assume that this property is satisfied and we denote by $\pi$ the unique element in $\Pi(\psi)$ such that $\pi_\infty \simeq \pi_{\rm hol}$. 
We want then to determine $m(\pi)$. By Lemma~\ref{calchol} and the remark that follows, we have 
$$\tau(\pi)_{|C_{\rm \psi}}={e_2}_{|{\rm C}_\psi}.$$
By Lemma~\ref{lemmeclass} (iii), if some $n_i$ is even then $n_i \equiv 0 \bmod 4$, and there is exactly one integer $i$ such that $n_i$ is odd. It follows that either $k=1$ (stable case) or $k=2$ and (up to equivalence) $(n_1,n_2)=(4,3)$. In this latter case 
${\rm C}_\psi=\langle s_1 \rangle \simeq \Z/2\Z$ and $\pi_2 \neq 1$.
 \ps

{\bf Case (i)} : ({\it stable tempered case}) $\psi=\pi_1 \in \Pi_{\rm alg}^{\rm o}(\PGL_7)$. Then $\psi_\infty$ is a discrete series Langlands parameter (hence indeed $\pi_{\rm hol} \in \Pi(\psi_\infty)$) and $m(\pi)=1$ by the multiplicity formula.
The number of such $\pi$ is thus ${\rm O}^\ast(w_1,w_2,w_3)$. \ps\ps

{\bf Case (ii)} : ({\it endoscopic tempered case}) $k=2$, $d_1=d_2=1$, $\psi=\pi_1\oplus \pi_2$ where $\pi_1 \in \Pi_{\rm alg}^{\rm o}({\rm PGL}_4)$ and $\pi_2 \in \Pi_{\rm alg}^{\rm o}(\PGL_3)$. Say $\pi_1$ has Hodge weights $a>b$ and $\pi_2$ has Hodge weight 
$c$. Then again $\psi_\infty$ is a discrete Langlands parameters (hence contains $\pi_{\rm hol}$). In particular $\varepsilon_\psi(s_1)=1$. But $e_2(s_1)=1$ if and only if $a>c>b$, thus 
$m(\pi)=1$ if and only if $a>c>b$. The number of such $\pi$ is thus ${\rm O}^*(w_2)\cdot {\rm O}(w_1,w_3)$.
 \ps\ps

{\bf Case (iii)} :  ({\it endoscopic non-tempered case}) $k=2$, $d_1=2$ and $d_2=1$, i.e. $\psi=\pi_1[2]\oplus \pi_2$ where $\pi_1 \in \Pi_{\rm alg}(\PGL_2)$ and $\pi_2 \in \Pi_{\rm alg}^{\rm o}(\PGL_3)$. Say $\pi_1$ has Hodge weight $a$ and $\pi_2$ has Hodge weight $b$. This time $\psi_\infty$ is not tempered and $$\varepsilon_\psi(s_1)=\varepsilon(\pi_1 \times \pi_2)=(-1)^{1+{\rm Max}(a,b)+\frac{a+1}{2}}.$$
But $e_2(s_1)=-1$, so

$$m(\pi)=\left\{\begin{array}{lll} 1 & {\rm if}\, \, \, \, b>a \, \, \, \, {\rm and}\, \, \, \, a \equiv 3 \bmod 4,\\ 1 & {\rm if}\, \, \, \, a>b \, \, \, \, {\rm and}\, \, \, \, a \equiv 1 \bmod 4, \\ 0 & {\rm otherwise}.\end{array}\right.$$
This concludes the proof of the proposition above. \ps\ps

We remark that we excluded three kinds of parameters thanks to Lemma~\ref{crithol}, namely $[7]$, $\pi_1 \oplus [3]$ and $\pi_1[2]\oplus [3]$. 
Alternatively, we could also have argued directly using only Lemma~\ref{calchol}. Indeed, in those three cases we obviously have $\varepsilon_\psi=1$, and we see also that $e_2(s_1)=-1$. \ps

\newpage 

\appendix

\section{Adams-Johnson packets}
\label{appendixadamsjohnson}

\subsection{Strong inner forms of compact connected real Lie groups}\label{strongforms} Let $K$ be a compact connected semisimple 
Lie group and let $G$ be its complexification. It is a complex semisimple algebraic group equipped with an 
anti-holomorphic group involution $\sigma : g \mapsto \overline{g}$ such that $K=\{g \in G, \overline{g}=g\}$. As is well-known, $K$ is a maximal compact subgroup of $G$.\ps

Let $T_c$ be a maximal torus of $K$ and denote by $T \subset G$ the unique maximal torus of $G$ with maximal compact subgroup $T_c$. Following J. Adams in~\cite{adams}, 
consider the group $$\mathcal{X}_1(T)=\{t \in T, t^2 \in {\rm Z}(G)\}.$$
An element of $\mathcal{X}_1(T)$ will be called a {\it strong inner form} of $K$ (relative to $(G,T)$). As  $K$ is a maximal compact subgroup of $G$, we have ${\rm Z}(G)={\rm Z}(K) \subset T_c$ and thus $\mathcal{X}_1(T) \subset T_c$.
A strong inner form $t \in \mathcal{X}_1(T)$ of $K$ is said {\it pure} if $t^2=1$.\ps

If $t \in \mathcal{X}_1(T)$ we denote by $\sigma_t$ the group automorphism ${\rm Int}(t) \circ \sigma$ of $G$. We have 
$\sigma_t^2={\rm Id}$. It follows that the real linear algebraic Lie group
	$$G_t=\{g \in G, \sigma_t(g)=g\}$$ is an inner form of $G_1=K$ in the usual sense. Observe that $T_c \subset G_t$ and that $G_t$ is stable by $\sigma$. The polar decomposition 
of $G$ relative to $K$ shows then that the group $$K_t = K \cap G_t,$$ which is also the centralizer of $t$ in $K$, is a maximal compact subgroup of $G_t$. The torus $T_c$ is thus a common maximal torus of all the $G_t$.  
Any involution of $G$ of the form ${\rm Int}(g) \circ \sigma$ with $g \in G$ is actually of the form ${\rm Int}(h) \circ \sigma_t \circ {\rm Int}(h)^{-1}$ for some $t \in \mathcal{X}_1(T)$ and some $h \in G$ by~\cite[\S 4.5]{serregc}. In particular, every inner form of $K$ inside $G$ is $G$-conjugate to some $G_t$.\ps

Consider the Weyl group $$W={\rm W}(G,T)={\rm W}(K,T_c).$$ 
It obviously acts on the group $\mathcal{X}_1(T)$, and two strong real forms $t,t' \in \mathcal{X}_1(T)$ are said {\it equivalent} if they are in a same $W$-orbit. If $w \in W$, observe that ${\rm Int}(w)$ defines an isomorphism $G_t \rightarrow G_{w(t)}$ which is well-defined up to inner isomorphisms by $T_c$, so that the group $G_t$ is canonically defined up to inner isomorphisms by the equivalence class of $t$. This is however not the unique kind of redundancy among the groups $G_t$ in general, as for instance $G_t=G_{tz}$ whenever $z \in {\rm Z}(G)$. We shall denote by $[t] \in W\backslash \mathcal{X}_1(T)$ the equivalence class of $t \in \mathcal{X}_1(T)$ and 
by $G_{[t]}$ the group $G_t$ "up to inner automorphisms". It makes sense in particular to talk about representations of $G_{[t]}$.  
\ps

As a classical example, consider the case of the even special orthogonal group $$G={\rm SO}_{2r}(\C)=\{g \in {\rm SL}_{2r}(\C), {}^t g g ={\rm Id}\}$$ with the coordinate-wise complex conjugation $\sigma$, i.e. $K={\rm SO}_{2r}(\R)=G \cap {\rm SL}_{2r}(\R)$. Consider the maximal torus $$T={\rm SO}_2(\C)^r \subset G$$ preserving each plane $P_i=\C e_{2i-1} \oplus \C e_{2i}$ for $i=1,\cdots, r$. Here $(e_i)$ is the canonical basis of $\C^{2r}$. Any 
$t \in \mathcal{X}_1(T)$ is ${\rm W}(G,T)$-equivalent to either a unique element $t_j$, $0\leq j \leq r$, where $t_j$ acts by $-1$ on $P_i$ if $i\leq j$, and by $+1$ otherwise, or to exactly one of the two element $t^\ast_\pm \in T_c$ sending each $e_{2i}$ 
on  $-e_{2i-1}$ for $i<r$ and $e_{2r}$ on $\pm e_{2r-1}$. We have $t_j^2=1$ (pure inner forms) and $(t^\ast_\pm)^2=-1$. We see that $$K_{t_j}=S({\rm O}(2j) \times {\rm O}(2r-2j))$$
and $G_{t_j} \simeq {\rm SO}(2j,2r-2j)$. In particular, $G_{t_j} \simeq G_{t_{j'}}$ if and only if $j=j'$ or $j+j'=r$. Moreover, $K_{t^\ast_\pm}$ is isomorphic to the unitary group in $r$ variables and $G_{t^\ast_\pm}$ is the real Lie group sometimes denoted by ${\rm SO}_{2r}^\ast$. Observe that the only quasi-split group among the $G_{t_j}$ and $G_{t^\ast_\pm}$ is ${\rm SO}(r+1,r-1)$ if $r$ is odd, ${\rm SO}(r,r)$ if $r$ is even. In particular, the split group 
${\rm SO}(r,r)$ is a pure inner form of $K$ if and only if $r$ is even.\ps

We leave as an exercise to the reader to treat the similar cases $G={\rm Sp}_{2g}(\C)$ and $G={\rm SO}_{2r+1}(\C)$ which are only easier. When $G={\rm Sp}_{2g}(\C)$, each twisted form of $K$ is actually 
inner as ${\rm Out}(G)=1$. In this case the (inner) split form ${\rm Sp}_{2g}(\R)$ is not a pure inner form of $K$, it corresponds to the single equivalence class of $t$ such that $t^2=-1$. 
When $G={\rm SO}_{2r+1}(\C)$, then ${\rm Z}(G)={\rm Out}(G)=1$, and the equivalence classes of strong inner forms of $K$ 
are in bijection with the isomorphism classes of inner forms of $K$, namely the real special orthogonal groups ${\rm SO}(2j,2r+1-2j)$ of signature $(2j,2r+1-2j)$ for $j=0,\cdots,r$.

\subsection{Adams-Johnson parameters} \label{ajparam} We refer to Kottwitz' exposition in \cite[p. 195]{kottwitz}
and to Adams paper~\cite{adams}, from which the presentation below is very much inspired.\ps

We keep the assumptions of~\S\ref{strongforms} and we assume from now on that the set of strong real forms of $K$ contains a split real group. It is equivalent to ask that 
the center of the simply connected covering $G_{\rm sc}$ of $G$ is an elementary abelian $2$-group, i.e. $G$ has no 
factor of type $E_6$, or type $A_n$ or $D_{2n-1}$ for $n>1$. We may view the Langlands dual group of $G$ as a complex connected semisimple algebraic group $\widehat{G}$, omitting the trivial Galois action. We fix once and for all a Borel subgroup $B$ containing $T$,  which gives rise to a dual Borel pair $(\widehat{T},\widehat{B})$ in $\widehat{G}$.  \ps\ps

Denote by $\Psi(G)$ the set of Arthur parameters of the inner forms of $K$. This is the set of continuous homomorphisms 
$${\rm W}_\R \times {\rm SL}_2(\C) \longrightarrow \widehat{G}$$
which are $\C$-algebraic on the ${\rm SL}_2(\C)$-factor and such that the image of any element of ${\rm W}_\R$ is semisimple. 
Two such parameters are said equivalent if they are conjugate under $\widehat{G}$. Fix $\psi \in \Psi(G)$.
Let $\widehat{L}$ be the centralizer in $\widehat{G}$ of $\psi(W_\C)$,
which is a Levi subgroup of some parabolic subgroup of $\widehat{G}$ ; as ${\rm W}_\C$ is commutative
$$\psi(W_\C \times {\rm SL}_2(\C)) \subset \widehat{L}.$$  Let us denote by 
${\rm C}_\psi$ the centralizer of ${\rm Im}(\psi)$ in $\widehat{G}$ and consider the 
following two properties of a $\psi \in \Psi(G)$. \ps\ps

\pn (a)  $\psi({\rm SL}_2(\C))$ contains a regular unipotent element of
$\widehat{L}$.\ps

\pn (b) ${\rm C}_\psi$ is finite.\ps\ps

Property (a)  forces in particular the centralizer of $\psi( {\rm
SL}_2(\C))$ in $\widehat{L}$ to be ${\rm Z}(\widehat{L})$, thus under (a) we have  
$${\rm C}_\psi={\rm Z}(\widehat{L})^{\theta}$$ where $\theta={\rm Int}(\psi(j))$ (recall that $j \in {\rm W}_\R \backslash {\rm W}_\C$ satisfies $j^2=-1$, see~\S\ref{demialg}).  Moreover, if one assumes (a) then property (b) is equivalent to the assertion that the
involution $\theta$ acts as the inversion on ${\rm Z}(\widehat{L})^0$. If $A$ is
an abelian group, we denote by $A[2]$ the subgroup of elements $a \in A$
such that $a^2=1$.

\begin{lemme}\label{lemmacent} If $\psi \in \Psi(G)$ satisfies (a) and (b) then ${\rm C}_\psi={\rm Z}(\widehat{L})[2]$. \end{lemme}

\noindent \begin{pf} Indeed, as a general fact one has ${\rm Z}(\widehat{L})={\rm Z}(\widehat{G}){\rm Z}(\widehat{L})^0$, because the character group of the diagonalisable group ${\rm Z}(\widehat{L})/{\rm Z}(\widehat{G})$ is free, being the quotient of the root lattice of $\widehat{G}$ by the root lattice of $\widehat{L}$. We also obviously have ${\rm Z}(\widehat{G}) \subset
{\rm C}_\psi$ ($\theta$ acts trivially on ${\rm Z}(\widehat{G})$), thus 
${\rm C}_\psi = {\rm Z}(\widehat{G}) ({\rm Z}(\widehat{L})^0)^\theta$. By assumption on $G$ one has ${\rm Z}(\widehat{G})={\rm Z}(\widehat{G})[2]$. As
$\theta$ acts as the inversion on ${\rm Z}(\widehat{L})^0$ one obtains
${\rm Z}(\widehat{L})[2]={\rm Z}(\widehat{G})({\rm Z}(\widehat{L})^0[2])={\rm C}_\psi$. 
\end{pf}
\ps\ps

To any $\psi \in \Psi(G)$ one may attach following Arthur a Langlands parameter $$\varphi_\psi
: {\rm W}_\R \rightarrow \widehat{G}$$ defined by restricting $\psi$ along the homomorphism 
$${\rm W}_\R \rightarrow {\rm W}_\R \times {\rm SL}_2(\C)$$ which is the identity on the first 
factor and the representation 
$|\cdot|^{1/2} \oplus |\cdot|^{-1/2}$ on the second factor. Here, $|\cdot| : {\rm W}_\R \rightarrow \R_{>0}$
is the norm homomorphism, sending $j$ to $1$ and $z \in W_\C$ to $z\overline{z}$. Up to $\widehat{G}$-conjugation of $\psi$, 
we may assume that $$\varphi_\psi({\rm W}_\C) \subset \widehat{T}.$$ We
follow Langlands notation\footnote{Recall that if $z \in \C^\times$, and if $a,b \in \C$
satisfy $a-b \in \Z$, we set $z^a\overline{z}^b=e^{ax+b\overline{x}}$ where $x
\in \C$ is any element such that $z=e^x$. The element $z^\lambda\overline{z}^\mu \in \widehat{T}$ is uniquely defined by the formula $\eta(z^\lambda\overline{z}^\mu)=z^{\langle \eta,\lambda\rangle}\overline{z}^{\langle \eta,\mu\rangle}$ for all $\eta \in {\rm X}^\ast(\widehat{T})$. } and write 
\begin{equation}\label{deflambda}\varphi_\psi(z)=z^\lambda \overline{z}^\mu \end{equation} for any $z \in {\rm W}_\C$, 
where $\lambda,\mu \in {\rm X}_\ast(\widehat{T}) \otimes \C={\rm Lie}_\C(\widehat{T})$ and 
$\lambda-\mu \in {\rm X}_\ast(\widehat{T})$. The $\widehat{G}$-conjugacy class of $\lambda$ in $\mathfrak{\widehat{g}}$ 
is called the infinitesimal character of $\psi$ (and $\varphi_\psi$) and will be denoted by $z_\psi$. The last 
condition we shall consider is : \ps\ps

\pn (c) $z_\psi$ is the infinitesimal character of a finite dimensional
$\C$-algebraic representation of $G$. \ps\ps
Under assumption (c), it
follows that $\widehat{T}$ is the centralizer of $\varphi_\psi({\rm W}_\C)$ in $\widehat{G}$, 
and up to conjugating $\psi$ by the Weyl group of $(\widehat{G},\widehat{T})$ one may assume that $\lambda$ is dominant 
with respect to $\widehat{B}$. 

\begin{definition}\label{defappaj} The subset of $\psi \in \Psi(G)$ satisfying (a), (b) and (c) will be denoted by $\Psi_{\rm AJ}(G)$. \end{definition} \ps

When $G$ is classical group, that is either ${\rm SO}_r(\C)$ or ${\rm Sp}_{2g}(\C)$, then so is $\widehat{G}$. A parameter $\psi \in \Psi(G)$ is an Adams-Johnson parameter if and only if it satisfies (c) and ${\rm St} \circ \psi$ is a multiplicity free representation of ${\rm W}_\R \times {\rm SL}_2(\C)$, where ${\rm St}$ denotes the standard representation of $\widehat{G}$.  \ps

As an example, consider the group $G={\rm Sp}_{2g}(\C)$, so that $\widehat{G}={\rm SO}_{2g+1}(\C)$. Let ${\rm St} : \widehat{G} \rightarrow {\rm GL}_{2r+1}(\C)$ be the standard 
representation of $\widehat{G}$. Let $\psi \in \Psi(G)$. Then $\psi \in \Psi_{\rm AJ}(G)$ if and only if 
$${\rm St} \circ \psi \simeq \varepsilon^s \otimes \nu_{d_0} \oplus \bigoplus_{i \neq 0} {\rm I}_{w_i} \otimes \nu_{d_i}$$
for some positive integers $w_i$ and $d_i$ with $(-1)^{w_i+d_i-1}=1$ for each $i$, with the convention $w_0=0$ and where $w_i>0$ if $i\neq 0$, such that the $2g+1$ even integers 
$$\pm w_i+d_i-1,\pm w_i+d_i-3,\cdots,\pm w_i-d_i+1$$
are distinct. The integer $s$ is congruent mod $2$ to the number of $i\neq 0$ such that $w_i$ is even. 
Moreover, the equivalence class of $\psi$ is uniquely determined by the isomorphism class of ${\rm St}\circ \psi$. If $\psi$ is as above, then 
$\widehat{L} \simeq {\rm SO}_{d_0}(\C) \times \prod_{i\neq 0} \GL_{d_i}(\C)$ and ${\rm C}_\psi \simeq \prod_{i \neq 0} \{\pm 1\}$. \ps

The case $G={\rm SO}_{2r+1}(\C)$ is quite similar, one simply has to replace the condition $(-1)^{w_i+d_i-1}=1$ by $(-1)^{w_i+d_i-1}=-1$, and there is no more restriction on $s \bmod 2$.
The case $G=\widehat{G}={\rm SO}_{2r}(\C)$ is slightly different but left as an exercise to the reader. 

\subsection{Adams-Johnson packets}\label{ajpackets}

In the paper~\cite{AJ}, J. Adams and J. Johnson associate to any $\psi \in \Psi_{\rm AJ}(G)$, and to any equivalence class of strong inner forms 
$t \in \mathcal{X}_1(T)$, a finite set of (usually non-tempered) irreducible unitary representations of $G_{[t]}$ satisfying certain predictions of Arthur. 
More precisely to any $t \in \mathcal{X}_1(T)$ they 
associate an isomorphism class $$\pi_t(\lambda)$$ 
of unitary irreducible representations of $G_t$, where $\lambda$ is as in~\eqref{deflambda}. Let us recall briefly their definition.\ps

Let $\widehat{P}$ be the parabolic subgroup of $\widehat{G}$ containing $\widehat{B}$ and 
with Levi subgroup $\widehat{L}$, let $P \subset G$ be parabolic subgroup dual to $\widehat{P}$, 
and let $L \subset P$ be the Levi subgroup dual to $\widehat{L}$. Concretely, as $$\psi({\rm W}_\C) \subset {\rm Z}(\widehat{L})^0 \subset \widehat{T}$$ there is a unique 
$\lambda_0 \in {\rm X}_\ast(\widehat{T})[1/2]$ such that $\psi(z)=(z/\overline{z})^{\lambda_0}$ for all $z \in {\rm W}_\C$, and $L$ is the stabilizer in $G$ of $\lambda_0 \in \mathfrak{g}^\ast$. 
Of course we have $T \subset L$. As $\sigma_t(\lambda_0)=-\lambda_0$, it follows that $\sigma_t$ preserves $L$, and thus the real Lie group $$L_t=G_t \cap L$$ is a real form of $L$ containing $T_c$. The real group $L_t$ is even an inner 
form of $K \cap L$ (a maximal compact subgroup of $L$). Moreover, the Cartan involution ${\rm Int}(t)$ of $G_t$ preserves $P$ as $T \subset P$, and defines a Cartan involution of $L_t$ as well. Assume that $L_t$ is connected to simplify (see loc. cit. for the general case). 
There is a unique one-dimensional unitary character $\chi_\lambda$ of $L_t$ whose restriction to $T_c$ is $\lambda-\rho$ where $\rho$ denotes the half-sum of the positive roots of $(G,T)$ with respect to $B$.
Adams and Johnson define $\pi_t(\lambda)$ as the cohomological induction 
relative to $P$ from the $(\mathfrak{l},K_t \cap L_t)$-module $\chi_\lambda$ to $(\mathfrak{g},K_t)$.  To emphasize the dependence on $P$ in this construction, we shall sometimes write $$\pi_{t,P}(\lambda)$$ rather 
than $\pi_t(\lambda)$.
\ps
The isomorphisms ${\rm Int}(w) : G_t \rightarrow G_{w(t)}$, for $w \in W$, allow to consider the collection of representations $\pi_{w(t)}(\lambda)$ as representations of $G_{[t]}$. The set of such 
representations is the Adams-Johnson packet of $G_{[t]}$ attached to $\psi$, and we shall denote it by $$\Pi(\psi,G_{[t]}).$$
It turns out that for $t,t' \in \mathcal{X}_1(T)$ in a same $W$-orbit, then $\pi_t(\lambda) \simeq \pi_{t'}(\lambda)$ if and only if $t$ and $t'$ are in a same ${\rm W}(L,T)$-orbit.
Observe also that for $t\in\mathcal{X}_1(T)$ we have $$\{w \in W, w(t)=t\}={\rm W}(K_t,T_c).$$ It follows that $\Pi(\psi,G_{[t]}))$ is in natural bijection with ${\rm W}(L,T)\backslash W/{\rm W}(K_t,T_c)$ and in particular that 
$|\Pi(\psi,G_{[t]})|$ is the number of such double cosets.\ps

\begin{lemme}\label{caraccompact} The representation $\pi_{t,P}(\lambda)$ is a discrete series representation if and only if $t \in {\rm Z}(L)$. \end{lemme}

\begin{pf} Indeed, as recalled loc. cit., $\pi_{t,P}(\lambda)$ is a discrete series representation if and only if $L_t$ is compact. The result follows as ${\rm Int}(t)$ is a Cartan involution of $L_t$. 
Note that for such a $t$ the group $L_t$ is of course always connected as so is $L$. 
\end{pf}

In the special case $t \in {\rm Z}(G)$, i.e. $G_t$ is compact, it follows that $\Pi(\psi,G_{[t]})$ is the singleton made of the unique irreducible representations of highest weight $\lambda-\rho$ relative to $B$. 
A more important special case is the one with $\psi({\rm SL}_2(\C))=\{1\}$. In this case $\psi$ is 
nothing more than a discrete series parameter in the sense of Langlands. Here (a) is automatic, (b) implies (c), $\varphi_\psi=\psi$ and $\widehat{L}=\widehat{T}$. Then $\pi_\lambda(t)$ is the 
discrete series representation with Harish-Chandra parameter $\lambda$, and $\Pi(\psi,G_{[t]})$ is simply the set of isomorphism classes of discrete series representations of $G_{[t]}$ with infinitesimal character $z_\psi$. \ps


\subsection{Shelstad's parameterization map}\label{parammap}What follows is again much inspired from \cite[p. 195]{kottwitz}
and~\cite{adams}. We fix a $\psi \in \Psi_{\rm AJ}(G)$ and keep the assumptions and notations of the previous paragraphs. We denote by $$\widetilde{\Pi}(\psi)$$
the disjoint union of the sets $\Pi(\psi,G_{[t]})$ where $[t]$ runs over the equivalence classes of strong real forms of $K$. As already explained in the previous paragraph, the map 
$\mathcal{X}_1(T) \rightarrow \widetilde{\Pi}(\psi), t \mapsto \pi_t(\lambda),$ induces a bijection 
\begin{equation}\label{basicparam} {\rm W}(L,T) \backslash \mathcal{X}_1(T) \isomo \widetilde{\Pi}(\psi). \end{equation} \ps

Define ${\rm S}_\psi$ as the inverse image of ${\rm C}_\psi$ under the simply
connected covering $$p : \widehat{G}_{\rm sc} \rightarrow \widehat{G}.$$
Following~\cite{shelun},~\cite{shelstad}, Langlands, Arthur, ~\cite{AJ},~\cite{abv} and~\cite{kottwitz}, the
set $\widetilde{\Pi}(\psi)$ is equipped with a natural map $$\tau_0 : \widetilde{\Pi}(\psi)
\rightarrow {\rm Hom}({\rm S}_\psi,\C^\times)$$
that we shall now describe in the style of Adams in~\cite{adams}. Observe first that ${\rm S}_\psi$
is the inverse image of ${\rm C}_\psi$ in $\widehat{T}_{\rm
sc}=p^{-1}(\widehat{T})$, hence it is an abelian group. \ps\ps

\begin{lemme}\label{lemmecentdeux}  ${\rm S}_\psi \subset (p^{-1}(\widehat{T}[2]))^{{\rm W}(L,T)}$.\end{lemme}

\begin{pf} By Lemma~\ref{lemmacent} and the inclusion ${\rm C}_\psi \subset {\rm C}_{\varphi_\psi}=\widehat{T}$ described in~\S\ref{ajparam}, one obtains a canonical 
inclusion ${\rm C}_\psi \subset \widehat{T}[2]$. Moreover, $\widehat{L}_{\rm sc}:=p^{-1}(\widehat{L})$ is a Levi subgroup of 
$\widehat{G}$ containing $\widehat{T}_{\rm sc}$ and thus $p^{-1}({\rm Z}(\widehat{L}))={\rm Z}(\widehat{L}_{\rm sc})$ and ${\rm W}(L,T)={\rm W}(\widehat{L},\widehat{T})={\rm W}(\widehat{L}_{\rm sc},\widehat{T}_{\rm sc})$. In particular, ${\rm W}(L,T)$ acts trivially on $p^{-1}({\rm Z}(\widehat{L}))$,
hence trivially on ${\rm S}_\psi$. 
\end{pf}\ps\ps

On the other hand, there is a natural perfect $W$-equivariant pairing $$\mathcal{X}_1(T) \times p^{-1}(\widehat{T}[2]) \rightarrow \C^\times.$$
Indeed, if $P^\vee(T)$ denotes the co-weight lattice of $T$ we have natural identifications 
 $$\mathcal{X}_1(T)=\frac{1}{2}P^\vee(T)/{\rm X}_\ast(T) \, \, \, {\rm
and}\, \, \, \, p^{-1}(\widehat{T}[2])= \frac{1}{2}{\rm
X}_\ast(\widehat{T})/{\rm X}_\ast(\widehat{T}_{\rm sc})$$ via $\mu \mapsto e^{2i\pi \mu}$. 
The pairing alluded above is then $(\mu,\mu') \mapsto e^{i\pi \langle \mu,\mu'\rangle}$, where $\langle,\rangle$ is the canonical perfect pairing 
${\rm X}_\ast(T)\otimes \Q \times {\rm X}_\ast(\widehat{T})\otimes \Q \rightarrow \Q$.  The resulting pairing is perfect as well as ${\rm X}_\ast(\widehat{T}_{\rm sc})$ is canonically identified by $\langle,\rangle$ with the root lattice of $T$. \ps

One then defines
$\tau_0$ as follows. Fix $\pi \in \widetilde{\Pi}(\psi)$. By the bijection~\eqref{basicparam}, there is an element $t \in \mathcal{X}_1(T)$, whose 
$W(L,T)$-orbit is canonically defined, such that $\pi \simeq \pi_t(\lambda)$. The perfect pairing above
associates to $t$ a unique character $p^{-1}(\widehat{T}[2]) \rightarrow
\C^\ast$, whose restriction to ${\rm S}_\psi$ only depends on the
${\rm W}(L,T)={\rm W}(\widehat{L},\widehat{T})$-orbit of $t$ by Lemma~\ref{lemmecentdeux}~: define $\tau_0(\pi)$ as this
character of ${\rm S}_\psi$. \ps

This parameterization is discussed in details in~\cite{adams} in the discrete series case, i.e. when $\widehat{L}=\widehat{T}$. 
It follows from the previous discussion that $\tau_0$ is a bijection in this
case, as ${\rm S}_\psi=p^{-1}(\widehat{T}[2])$. The following simple lemma shows that the determination of the parameterization of discrete series in $\widetilde{\Pi}(\psi)$ 
for general $\psi$ reduces to this latter case.\ps

Observe following~\cite{kottwitz} that for any $\psi \in \Psi_{\rm AJ}(G)$ there is a unique discrete series parameter $\psi_{\rm disc} \in \Psi_{{\rm AJ}}(G)$ such that the centralizers of $\varphi_\psi({\rm W}_\C)$ and $\psi_{\rm disc}({\rm W}_\C)$ coincide, and such that the parameters $\lambda$ for $\psi$ and $\psi_{\rm disc}$ defined by~\eqref{deflambda} coincide as well. In particular, 
$\psi$ and $\psi_{\rm disc}$ have the same infinitesimal character. Of course, $\psi_{\rm disc} \neq \varphi_\psi$ if $\psi \neq \psi_{\rm disc}$. If $\psi$ is normalized as before, we have canonical inclusions
$${\rm C}_\psi \subset {\rm C}_{\psi_{\rm disc}}=\widehat{T}[2] \, \, \, \, {\rm and}\, \, \, \, {\rm S}_\psi \subset {\rm S}_{\psi_{\rm disc}}=p^{-1}(\widehat{T}[2]).$$
The discrete series representations belonging to $\widetilde{\Pi}(\psi)$ are exactly the elements of $\widetilde{\Pi}(\psi)\cap \widetilde{\Pi}(\psi_{\rm disc})$. It will be important to distinguish in the next lemma the 
parameterization maps $\tau_0$ of $\widetilde{\Pi}(\psi)$ and $\widetilde{\Pi}(\psi_{\rm disc})$, so we shall denote them respectively by $\tau_{0,\psi}$ and $\tau_{0,\psi_{\rm disc}}$. Recall form Lemma~\ref{caraccompact}
that $\pi_t(\lambda)$ is a discrete series representation if and only if $t \in {\rm Z}(L)$. The following lemma is a variant of an observation by Kottwitz in~\cite{kottwitz}.

\begin{lemme}\label{reductionsd} Let $\psi \in \Psi_{\rm AJ}(G)$ and let $\pi \in \widetilde{\Pi}(\psi)\cap \widetilde{\Pi}(\psi_{\rm disc})$. Then 
$$\tau_{0,\psi}(\pi)={\tau_{0,\psi_{\rm disc}}}(\pi)_{|{\rm S}_{\psi}}.$$
\end{lemme}

\begin{pf} We have  $\pi \simeq \pi_{t,B}(\lambda) \in \widetilde{\Pi}(\psi_{\rm disc})$ for a unique $t \in \mathcal{X}_1(T)$ and we also have $\pi \simeq \pi_{t',P}(\lambda) \in \widetilde{\Pi}(\psi)$ for a unique element $t' \in \mathcal{X}_1(T)\cap {\rm Z}(L)$ by Lemma~\ref{caraccompact} (note that $t'$ is fixed by $W(L,T)$). 
Applying the "transitivity of cohomological induction"  via the compact connected group $L_t$ (use e.g.~\cite[Cor. 11.86 (b)]{knappvogan}, here $q_0=0$), we have 
$\pi_{t,B}(\lambda) \simeq \pi_{t,P}(\lambda)$. It follows that $t=t'$, which concludes the proof. 
\end{pf}

Observe that for $t \in \mathcal{X}_1(T)$, $G_t$ is compact if and only if $t \in {\rm Z}(G)$, in which case it coincides with its equivalence class (and it is fixed by ${\rm W}(L,T)$). 
The associated representation $\pi_{t,P}(\lambda)$ is the unique finite dimensional representation of $G_t$ with infinitesimal character $z_\psi$. It occurs $|{\rm Z}(G)|$ times in $\widetilde{\Pi}(\psi)$, once for each $t \in {\rm Z}(G)$,  and these representations are perhaps the most obvious elements in $\widetilde{\Pi}(\psi)\cap \widetilde{\Pi}(\psi_{\rm disc})$. To understand their characters we have to describe the image 
 $$\mathcal{N}(T)$$ of ${\rm Z}(G)$ under the homomorphism $\mathcal{X}_1(T) \rightarrow {\rm Hom}(p^{-1}(\widehat{T}[2]),\C^\ast)$ induced by the canonical pairing. 
 Observe that ${\rm Z}(G)=\{t^2, t \in \mathcal{X}_1(T)\}$. The following lemma follows.

\begin{lemme}\label{lemmeactz} The subgroup $\mathcal{N}(T) \subset {\rm Hom}(p^{-1}(\widehat{T}[2]),\C^\ast)$ is the subgroup of squares, or equivalently of characters which are trivial on $\widehat{T}_{\rm sc}[2]$.
\end{lemme} 

The parameterization $\tau_0$ of $\widetilde{\Pi}(\psi)$ introduced so far is the one we shall need up to a translation by a certain character $b_\psi$ of ${\rm S}_\psi$ (or "base point of $\psi$"). Write again temporarily $\tau_{0,\psi}$ for $\tau_0$ in order to emphasize its dependence on $\psi$ and we write character groups additively. The map $$\tau_{\psi}  = \tau_{0,\psi} - b_\psi$$
has to satisfy the following two conditions :\ps\ps

(i) Lemma~\ref{reductionsd} holds with $\tau_{0,\psi}$ and $\tau_{0,\psi_{\rm disc}}$ replaced respectively by $\tau_{\psi}$ and $\tau_{\psi_{\rm disc}}$. \ps
(ii) If $\psi$ is a discrete series parameter, and if $\pi=\pi_t(\lambda) \in \widetilde{\Pi}(\psi)$ satisfies $\tau_\psi(\pi)=1$, i.e. $\tau_{0,\psi}(\pi)=b_\psi$, then $G_t$ is 
a split real group and $\pi$ is generic with respect to some Whittaker functional. \ps

Normalize $\psi$ as in~\S\ref{ajparam}. Following~\cite{adams}, consider the element $$ t_b = e^{i\pi \rho^\vee} \in \mathcal{X}(T)$$ where
$\rho^\vee\in {\rm X}_\ast(T)$ is the half-sum of the positive coroots with respect
to $(G,B,T)$.  Under the identification $\mathcal{X}_1(T)=\frac{1}{2}P^\vee(T)/{\rm
X}_\ast(T)$, $t_b$ is the class of $\frac{1}{2}\rho^\vee$. In particular, under the
canonical pairing between ${\rm S}_{\psi_{\rm disc}}$ and $\mathcal{X}_1(T)$ the element $t_b$ corresponds
to the restriction to ${\rm S}_{\psi_{\rm disc}}$ of the character $\rho^\vee \in {\rm
X}^\ast(\widehat{T}_{\rm sc})$. The characteristic property of $t_b$ is that for any $t$ in the coset 
${\rm Z}(G) t_b \subset \mathcal{X}_1(T)$, then $\pi_{t}(\lambda)$ is
a generic (or "large" in the sense of Vogan) discrete series of the split group $G_t$. To fulfill the conditions (i) and (ii) one simply set 
$b_\psi = \rho^\vee$. 

\begin{definition}\label{defparamcan} If $\psi \in \Psi_{\rm AJ}(G)$, the canonical parameterization $$\tau : \widetilde{\Pi}(\psi) \rightarrow {\rm Hom}({\rm S}_\psi,\C^\ast)$$ is defined by $\tau=\tau_0-{\rho^\vee}_{|{\rm S}_\psi}$
where $\widehat{T}$ is the centralizer of $\varphi_\psi({\rm W}_\C)$, $\widehat{B}$ is the unique Borel subgroup of $\widehat{G}$ containing $\widehat{T}$ with respect to which
 the element $\lambda$ defined by~\eqref{deflambda} is dominant, and $\rho^\vee$ is the half-sum of the positive roots of $(\widehat{G},\widehat{B},\widehat{T})$. 
 \end{definition}

\begin{cor}\label{corparamcompact} If $\pi \in \widetilde{\Pi}(\psi)$ is a finite dimensional representation, then $\tau(\pi) \in \mathcal{N}(T)-\rho^\vee$. \end{cor}

We end this paragraph by collecting a couple of well-known and simple facts we used in the paper.  For $t,t' \in \mathcal{X}_1(T)$, $G_t$ and $G_{t'}$ are pure inner forms if and only if $t^2=(t')^2$. 

\begin{cor} $K$ is a pure inner form of a split group if and only if $\rho^\vee \in {\rm X}^\ast(\widehat{T})$. \end{cor}

Indeed, $G_t$ is a pure inner form of a split group if and only if $t^2=t_b^2=(-1)^{2\rho^\vee}$. \ps

\begin{cor}\label{factorrho} Let $t \in \mathcal{X}_1(T)$. Then $G_t$ is a pure inner form of a split group if and only if the character of $p^{-1}(\widehat{T}[2])$ associated to $tt_b^{-1}$ under the canonical pairing 
factors through $\widehat{T}[2]$.
\end{cor}


\newpage

\section{The Langlands group of $\Z$ and Sato-Tate groups}\label{parsatotate}

In this brief appendix, we discuss a conjectural topological group that
might be called the {\it Langlands group of $\Z$}, and that we shall denote $\mathcal{L}_\Z$. We will define $\mathcal{L}_\Z$ as a suitable quotient of the conjectural Langlands group $\mathcal{L}_\Q$ of the field of rational numbers $\Q$, originally introduced by Langlands
in \cite{Lgl}. The group $\mathcal{L}_\Z$ is especially relevant to understand the level $1$ automorphic representations of reductive groups over $\Z$. As we shall explain, and following~\cite[Ch. II \S 3.6]{chhdr}, it also offers a plausible point of view on the Sato-Tate groups of automorphic representations and motives. Let us stress once and for all that {\it most of this appendix is purely hypothetical}. Nevertheless, we hope it might be a useful and rather precise guide to the understand the philosophy, due to Langlands and Arthur, behind the results of this paper.  \ps

We shall view $\mathcal{L}_\Q$ as a (Hausdorff) locally compact topological group following Kottwitz' point of view in \cite[\S 12]{kottwitzctt}. We refer to Arthur's paper~\cite{arthurconjlan} for a thorough discussion of the expected properties of $\mathcal{L}_\Q$ and for a description of a candidate for this group as well. 

\subsection{The locally compact group $\mathcal{L}_\Z$}

If $p$ is a prime, recall that the group $\mathcal{L}_\Q$ is equipped with a conjugacy class of continuous homomorphisms $\eta_p : I_p \times {\rm SU}(2) \rightarrow \mathcal{L}_\Q$, where $I_p$ is the inertia group of the absolute Galois group of $\Q_p$. {\it We define $\mathcal{L}_\Z$ as the quotient of $\mathcal{L}_\Q$ by the closed normal subgroup generated by the union, over all primes $p$, of ${\rm Im}\, \, \eta_p$}. It is naturally equipped with : \begin{itemize}\ps
\item[-] (Frobenius elements) a conjugacy class ${\rm Frob_p} \subset \mathcal{L}_\Z$ for each prime $p$, \ps
\item[-] (Hodge morphism) a conjugacy class of continuous group homomorphisms $$h : {\rm W}_\R \rightarrow \mathcal{L}_\Z,$$\ps
\end{itemize}
which inherit from $\mathcal{L}_\Q$ a collection of axioms that we partly describe below. \ps

As $\widehat{\GL_n}=\GL_n(\C)$ we have a parameterization map $c : \Pi(\GL_n) \rightarrow \mathcal{X}(\GL_n(\C))$ as in~\S\ref{langparam}. Denote by ${\rm Irr}_n(\Z)$ the set of isomorphism classes of irreducible continuous representations $\mathcal{L}_\Z \rightarrow {\rm GL}_n(\C)$.  \ps
\begin{itemize}\ps
\item[(L1)] ({\it Langlands conjecture}) For any $n\geq 1$ and any $\pi \in \Pi_{\rm cusp}(\GL_n)$, there exists $r_\pi \in {\rm Irr}_n(\Z)$ such that $c_p(\pi)$ is conjugate to $r_\pi(\Frob_p)$ for each prime $p$, and ${\rm L}(\pi_\infty) \simeq r_\pi \circ \mu$ (see~\S\ref{demialg}). Moreover, $\pi \mapsto r_\pi$ defines a bijection $\Pi_{\rm cusp}(\GL_n) \isomo {\rm Irr}_n(\Z)$. \ps
\end{itemize}

\noindent Let $\R_{>0} \subset \R^\times$ be the multiplicative subgroup of positive numbers. The adelic norm $$|\cdot |:  \Q^\times \backslash \AAA^\times / {\widehat{\Z}}^\times \rightarrow \R_{>0}$$ is an isomorphism, thus $\Pi_{\rm cusp}(\GL_1)=\{ |\cdot|^s, s \in \C\}$. Set $|\cdot|_\Z=r_{|\cdot|} \in {\rm Irr}_1(\Z)$. We have $|{\rm Frob}_p|_\Z=p^{-1}$ for any prime $p$ and $|\cdot|_\Z \circ h : {\rm W}_\R \rightarrow \R_{>0}$ coincides with the homomorphism recalled in~\S\ref{demialg}. If $\mathcal{D} H \subset H$ denotes the closed subgroup generated by the commutators of the topological group $H$, and $H^{\rm ab}=H/\mathcal{D}H$, it is natural to ask that : \ps
\begin{itemize}
\item[(L2)] ({\it Class field theory}) $|\cdot|_\Z$ induces a topological isomorphism $\mathcal{L}_\Z^{\rm ab} \isomo \R_{>0}$.\ps
\end{itemize}

\noindent Let $\mathcal{L}_\Z^1=\mathcal{D}\mathcal{L}_\Z$ be the kernel of $|\cdot|_\Z$. \ps\ps
 
\begin{itemize} 
\item[(L3)] ({\it Ramanujan conjecture}) $\mathcal{L}^1_\Z$ is compact. \ps
\end{itemize}\ps\ps

\noindent Properties (L2) and (L3) have the following consequence on the structure of $\mathcal{L}_\Z$. \ps\medskip
{\bf Fact 1 : } ({\it Polar decomposition}) If  $\mathcal{C} \subset \mathcal{L}_\Z$ denotes the neutral component of the center of $\mathcal{L}_\Z$, then 
$\mathcal{L}_\Z = \mathcal{C} \times \mathcal{L}^1_\Z$
and $|\cdot|_\Z$ induces an isomorphism $\mathcal{C} \isomo \R_{>0}$. In particular,  $\mathcal{L}^1_\Z = \mathcal{D} \mathcal{L}^1_\Z$.\ps

{\footnotesize \begin{pf} Let $\mathcal{Z} \subset \mathcal{L}_\Z$ be the centralizer of the compact normal subgroup $\mathcal{L}^1_\Z$. As $\mathcal{L}^1_\Z$ is a compact normal subgroup of $\mathcal{L}_\Z$, and as $\mathcal{L}_\Z/\mathcal{L}^1_\Z$ is connected, a classical result of Iwasawa \cite[\S 1]{iwasawa} ensures that $\mathcal{L}_\Z = \mathcal{Z} \, \mathcal{L}^1_\Z$. The subgroup $\mathcal{Z}^1=\mathcal{Z} \cap \mathcal{L}^1_\Z$ is central in $\mathcal{Z}$, and $|\cdot|_\Z$ induces an isomorphism $\mathcal{Z}/\mathcal{Z}^1 \isomo \R_{>0}$ by (L2) and the open mapping theorem. As any central extension of $\Q$ is abelian, the Hausdorff topological group $\mathcal{Z}$ is abelian : it thus coincides with the center of $\mathcal{L}_\Z$. In particular, $\mathcal{L}^1_\Z=\mathcal{D} \mathcal{L}^1_\Z$. The center of a compact Lie group $H$ such that $\mathcal{D} H=H$ is finite, so the center $\mathcal{Z}^1$ of $\mathcal{L}^1_\Z$ is profinite. The structure theorem of locally compact abelian groups concludes $\mathcal{Z}=\mathcal{C} \times \mathcal{Z}^1$. 
\end{pf}}
It will be convenient to identify once and for all $\mathcal{C}$ and $\R_{>0}$ via $|\cdot|_\Z$. In other words, we view $|\cdot|_\Z$ has a homomorphism $\mathcal{L}_\Z \rightarrow \mathcal{C}$, and write $g = |g| \cdot (g/|g|)$ for the polar decomposition of an element $g \in \mathcal{L}_\Z$. Observe that for any continuous representation $r : \mathcal{L}_\Z \rightarrow \GL_n(\C)$, the elements of the compact group $r(\mathcal{L}^1_\Z)$ are semisimple and all their eigenvalues have norm $1$. Moreover, if $r$ is irreducible then $r(\mathcal{C})$ acts by scalars by Schur's lemma, so there exists $s \in \C$ such that $r(x)=|x|_\Z^s$ for all $x \in \mathcal{C}$. If $r=r_\pi$ for $\pi \in \Pi_{\rm cusp}(\GL_n)$, observe that $\det(r_\pi)=|.|_\Z^{ns}=r_{\omega_\pi}$, where $\omega_\pi \in \Pi_{\rm cusp}(\GL_1)$ denotes the central character of $\pi$. If we consider the image of an element in the conjugacy class ${\rm Frob}_p$ and property (L1) we recover the classical Ramanujan conjecture on the $c_p(\pi)$ for $\pi \in \Pi_{\rm cusp}(\GL_n)$. This "explains" as well Clozel's purity lemma~\ref{clopl}. \ps\medskip

{\bf Fact 2} : ({\it Generalized Minkowski theorem}) $\mathcal{L}^1_\Z$ is connected. The conjugacy class of $h(U(1))$, where $U(1)$ is the maximal compact subgroup of ${\rm W}_\C$, generates a dense subgroup of $\mathcal{L}^1_\Z$. \ps \ps

{\footnotesize \begin{pf} One of the axioms on $\mathcal{L}_\Q$ is that its group of connected components is naturally isomorphic to ${\rm Gal}(\overline{\Q}/\Q)$. The first part of Fact 2 follows then from Minkowski's theorem asserting that any non-trivial number field admits at least a ramified prime. Here is another proof. Assume that $\mathcal{L}_\Z^1$ admits a non-trivial finite quotient $\Gamma$ and choose a non trivial irreducible representation of $\Gamma$, say of dimension $n\geq 1$, that we view as an element $r \in {\rm Irr}_n(\Z)$ trivial on $\mathcal{C}$, of finite image.  Let $\pi \in \Pi_{\rm cusp}(\GL_n)$ be such that $r=r_\pi$. As ${\rm Im}\, r$ is finite, ${\rm L}(\pi_\infty) = r \circ h$ is trivial on the connected subgroup ${\rm W}_\C$. To conclude the proof (of the second statement as well) it is thus enough to show that such a $\pi$ is necessarily the trivial representation of $\GL_1$. But it follows indeed from Weil's explicit formulas that the ${\rm L}$-function of a non-trivial such $\pi$, which is entire and of conductor $1$, does not exist : see~\cite[\S 3]{mestre}.\end{pf}}

We end this paragraph by a definition of the {\it motivic Langlands group of $\Z$}, that we shall define as a certain quotient $\mathcal{L}_\Z^{\rm mot}$ of $\mathcal{L}_\Z$.  \ps

Recall that a $\pi \in \Pi_{\rm cusp}(\GL_n)$ is said {\it algebraic} if the restriction of ${\rm L}(\pi_\infty)$ to ${\rm W}_\C$ is a direct sum of characters of the form $z \mapsto z^{a_i}\overline{z}^{b_i}$ where $a_i,b_i \in \Z$, for $i=1,\dots,n$ : see the footnote~\ref{footnotealg} of the introduction for references about this notion. Clozel's purity lemma ensures that $a_i+b_i$ is independent of $i$, or which is the same, that ${\rm Z}({\rm W}_\R)$ acts as scalars in ${\rm L}(\pi_\infty)$.  Here ${\rm Z}({\rm W}_\R)$ denotes the center of ${\rm W}_\R$, namely the subgroup $\R^\times \subset {\rm W}_\C$. For any $n\geq 1$
define $\Pi_{\rm mot}(\GL_n)$ as the subset of $\pi \in \Pi_{\rm
cusp}(\GL_n)$ such that ${\rm Z}({\rm W}_\R)$ acts as scalars in ${\rm L}(\pi_\infty)$.  Let $\pi \in \Pi_{\rm cusp}(\GL_n)$.  It is now a simple exercise to check that $\pi \in \Pi_{\rm mot}(\GL_n)$ if and only if there exists $s \in \C$ such that $\pi \otimes |\cdot|^s$ is algebraic. For instance, $$\Pi_{\rm alg}^\bot(\PGL_n) \subset \Pi_{\rm mot}(\GL_n)$$ (see Definition~\ref{defpialg}).  \ps

{\it We define $\mathcal{L}_\Z^{\rm mot}$ as the quotient of $\mathcal{L}_\Z$ by the closed normal subgroup generated by the $xyx^{-1}y^{-1}$ where $x \in h({\rm Z}({\rm W}_\R))$ and $y \in \mathcal{L}_\Z$}. By definition, if $\pi \in \Pi_{\rm cusp}(\GL_n)$ then $\pi \in \Pi_{\rm mot}(\GL_n)$ if and only if $r_\pi$ factors through $\mathcal{L}_\Z^{\rm mot}$. The locally compact group $\mathcal{L}_\Z^{\rm mot}$ inherits from $\mathcal{L}_\Z$ all the properties considered so far. Better, the subgroup $h({\rm Z}({\rm W}_\R))$ is a central subgroup of $\mathcal{L}_\Z^{\rm mot}$, so that the polar decomposition is even simpler to understand for $\mathcal{L}_\Z^{\rm mot}$ as $h_{|\R_{>0}}$ defines a central section of $|\cdot|_\Z : \mathcal{L}_\Z^{\rm mot} \rightarrow \R_{>0}$. \ps

\subsection{Sato-Tate groups}

Serre's point of view in~\cite[Ch. 1, appendix]{serreabelian} and \cite[\S 12]{serremotives} suggests the following universal form of the Sato-Tate conjecture (here, in the level $1$ case).

\begin{itemize}\ps
\item[(L4)] ({\it General Sato-Tate conjecture}) The conjugacy classes $\frac{{\rm Frob}_p}{|{\rm Frob}_p|} \subset \mathcal{L}^1_\Z$ are
equidistributed in the compact group $\mathcal{L}^1_\Z$ equipped with its Haar measure of mass $1$. \ps
\end{itemize}
\ps\ps

Note in particular that the union of the conjugacy classes $\frac{{\rm Frob}_p}{|{\rm Frob}_p|}$ is dense in $\mathcal{L}^1_\Z$ ({\it Cebotarev property}), which "explains" the strong multiplicity one theorem for $\GL_n$ by (L1). \ps

\begin{propdef} If $\pi \in \Pi_{\rm cusp}(\GL_n)$, define its Sato-Tate group as $$\mathcal{L}_\pi:=r_\pi(\mathcal{L}^1_\Z).$$ \ps
It is a compact connected subgroup of ${\rm SL}_n(\C)$ well-defined up to ${\rm SL}_n(\C)$-conjugacy, which acts irreducibly on $\C^n$, and such that $\mathcal{L}_\pi^{\rm ab}=1$. The Satake parameters of the $\pi_p$ have well-defined representatives in $\mathcal{L}_\pi$, namely the conjugacy classes $r_\pi({\rm Frob}_p/|{\rm Frob}_p|) \subset \mathcal{L}_\pi$, which are equidistributed for a Haar measure of $\mathcal{L}_\pi$.
\end{propdef}\ps

Remark that if $\pi \in \Pi_{\rm cusp}(\PGL_n)$ then $r(\mathcal{C})=1$ and so $\mathcal{L}_\pi=r(\mathcal{L}_\Z)$. The last property of $\mathcal{L}_\Z$ we would like to discuss is the general Arthur-Langlands conjecture. This is first especially helpful in order to understand the results of Arthur recalled in~\S\ref{sectionarthur} (see also~\cite{arthurunipotent} and the introduction of~\cite{arthur}). This will also give another way to think about $\mathcal{L}_\pi$ when $\pi \in \Pi_{\rm cusp}(\GL_n)$.\ps

Fix $G$ a semisimple group scheme over $\Z$. Following Arthur, define a {\it global Arthur parameter for $G$} as a $\widehat{G}$-conjugacy class of continuous group homomorphisms $$\psi : \mathcal{L}_\Z \times {\rm SL}_2(\C) \rightarrow \widehat{G}$$ such that $\psi_{|{\rm SL}_2(\C)}$ is algebraic, and such that the centralizer ${\rm C}_\psi$ of ${\rm Im}\, \, \psi$ in $\widehat{G}$ is finite. As an example, suppose that $G$ is a classical group and let ${\rm St} : \widehat{G} \rightarrow {\rm GL}_n(\C)$ denote the standard representation. If $\psi$ is a global Arthur parameter for $G$, the finiteness of ${\rm C}_\psi$ ensures that the representation ${\rm St} \circ \psi$ is a direct sum of pairwise non-isomorphic irreducible representations of ${\mathcal{L}}_\Z \times {\rm SL}_2(\C)$, say of dimension $n_i$, hence of the form $r_i \otimes {\rm Sym}^{d_i-1} \,\C^2$ where $d_i | n_i$ and $r_i \simeq r_i^\ast  \in {\rm Irr}_{n_i/d_i}(\Z)$ : via (L1) this "explains" the definition of a global Arthur parameter in~\S\ref{arthurparam}, except property (ii) {\it loc. cit.}  at the moment. (Observe that for any $\pi \in \Pi_{\rm cusp}(\PGL_n)$, we have $r_\pi^\ast = r_{\pi^\vee}$ by (L1)). \ps
\ps

Recall Arthur's morphism $a : \mathcal{L}_\Z \rightarrow \mathcal{L}_\Z \times {\rm SL}_2(\C), \, \, \, g \mapsto (\, g\, ,\, \diag(\, |g|_\Z^{1/2}\,,\,|g|_\Z^{-1/2}\,))$. If $\psi$ is a global Arthur parameter for $G$, then $\varphi_\psi := \psi \circ a$ is a well-defined conjugacy class of continuous homomorphisms ${\mathcal L}_\Z \rightarrow \widehat{G}$. Moreover $\psi_\infty:=\psi \circ h$ is also an Arthur parameter in the sense of~\S\ref{ajparam}. In particular, it possesses an infinitesimal character $z_{\psi_\infty} \subset \widehat{\mathfrak{g}}$ as defined {\it loc. cit.}\ps\ps

\ps \begin{itemize}
\item[(L5)] ({\it Arthur-Langlands conjecture}) For any $\pi \in \Pi_{\rm disc}(G)$, there is a global Arthur parameter $\psi$ for $G$ {\it associated to $\pi$} in the following sense : $\varphi_\psi({\rm Frob}_p)$ is conjugate to $c_p(\pi)$ for each prime $p$ and $z_{\psi_\infty}$ is conjugate to $c_\infty(\pi)$. Conversely, if $\psi$ is a global Arthur parameter for $G$, and if $\Pi(\psi)$ is the finite set of $\pi \in \Pi(G)$ associated to $\psi$, then there is a formula for $\sum_{\pi \in \Pi(\psi)} m(\pi)$. 
\end{itemize}
\ps\ps

Recall that $m(\pi)$ denotes the multiplicity of $\pi$ in $\mathcal{L}_{\rm disc}(G)$. Let us warn that 
there may be in general several $\psi$ associated to a given $\pi$, because there are examples of continuous morphisms  $H_1 \rightarrow H_2$, say between two compact connected Lie groups, which are point-wise conjugate but non conjugate (try for instance $H_1={\rm SU}(3)$ and $H_2={\rm SO}(8)$) : see~\cite{arthurunipotent} and \cite{arthurconjlan} for more about this problem. Let us also mention that Langlands originally considered only the tempered $\pi \in \Pi_{\rm cusp}(G)$ and conjectured the existence of a $\psi$ as in (L5) but trivial on the ${\rm SL}_2(\C)$ factor. \ps

If we consider again the example of classical groups, the first part of (L5) "explains" Arthur's Theorem~\ref{thmarthurclass}. Of course, its second part is too vague as stated here : see~\cite{arthurconjlan} and~\cite{arthurunipotent} for more informations about this quite delicate point  called {\it Arthur's multiplicity formula}. See also~\S\ref{amfintro} for an explicit formula when $G(\R)$ is compact, and to \S\ref{parmultform} for certain explicit special cases for classical groups. Observe also that in this latter case, the group ${\rm C}_\psi$ defined there following Arthur fortunately co\"incides with the group ${\rm C}_\psi$ defined here. \ps

An important situation where we can say more about the second part of (L5) is when ${\rm C}_\psi$ coincides with the center of $\widehat{G}$. In this case, and if $\pi \in \Pi(\psi)$, we have $m(\pi) \neq 0$ if $\pi_\infty$ belongs to Arthur's conjectural set (or "packet") of unitary representations of $G(\R)$ associated to $\psi_\infty$. When $\psi({\rm SL}_2(\C))=1$, this packet is the set of unitary representations of $G(\R)$ associated by Langlands to $\varphi_\psi \circ h$ in \cite{langlandsreel}; it is never empty if $G$ is a Chevalley group. \ps

This last paragraph "explains" for instance Theorem~\ref{orthosymp} (as well as the results in~\S\ref{demialg}). Indeed, a finite-dimensional selfdual irreducible representation of any group preserves a {\it unique} non-degenerate pairing up to scalars, either symmetric or anti-symmetric. This explains as well condition (ii) in the definition of a global Arthur parameter in~\S\ref{arthurparam}. \ps

This also leads to another way of thinking about $\mathcal{L}_\pi$ when $\pi \in \Pi_{\rm cusp}(\GL_n)$, which involves all the semisimple groups over $\Z$. Indeed, fix $\pi \in \Pi_{\rm cusp}(\PGL_n)$ and let $G_\pi$ be the Chevalley group such that $\widehat{G}_\pi$ is a complexification of the compact connected semisimple Lie group $\mathcal{L}_\pi$. By definition, we may factor $r_\pi$ through a homomorphism $$\widetilde{r_\pi} : \mathcal{L}_\Z \rightarrow \widehat{G_\pi}$$ such that ${\rm C}_{\widetilde{r_\pi}}$ is the center of $\widehat{G_\pi}$. We thus obtain {\it \`a la Langlands} a non-empty finite set of representations $\pi' \in \Pi_{\rm disc}(G_\pi)$ associated to $\widetilde{r_\pi}$. This explains for instance the discussion that we had about the group ${\rm G}_2$ in the introduction. From this point of view, the results that we proved in~\S\ref{lzscon} imply the following: \ps \medskip

{\bf Fact 3 :} $\mathcal{L}^1_\Z$ is simply connected. \ps\medskip

Arthur has a similar prediction for $\mathcal{L}_\Q$ in~\cite{arthurconjlan}, although this does not seem to directly imply that $\mathcal{L}_\Z$ should be simply connected as well. As a consequence of Fact 3, it follows that $\mathcal{L}_\Z$ is a direct product of $\R_{>0}$ and of countably many semisimple, connected, simply connected, compact Lie groups. The same property holds for $\mathcal{L}_\Z^{\rm mot}$ by construction. It is a natural question to ask which semisimple, connected, simply connected compact Lie group appear as a direct factor of $\mathcal{L}_\Z^{\rm mot}$ or of $\mathcal{L}_\Z$. The results of this paper show that this is indeed the case (for $\mathcal{L}_\Z^{\rm mot}$) for each such group whose simple factors are of type $A_1$, $B_2$, $G_2$, $B_3$, $C_3$, $C_4$ or $D_4$.  Let us mention that in their work~\cite{ashpollack}, Ash and Pollack did seach factors of $\mathcal{L}_\Z^{\rm mot}$ of type $A_2$ by computing cuspidal cohomology of ${\rm SL}_3(\Z)$ for a quite large number of coefficients : they did not find any.
 

\subsection{A list in rank $n\leq 8$}
Our goal in this last paragraph is to determine the possible Sato-Tate
groups of a $\pi \in \Pi_{\rm alg}^\bot(\PGL_n)$ when $n\leq 8$. For such a $\pi$, define $\mathcal{A}_\pi$ as the compact symplectic group of rank
$n/2$ if $s(\pi)=-1$, the compact special orthogonal group ${\rm SO}(n)$
otherwise.  By Arthur's theorem~\ref{orthosymp}, $\mathcal{L}_\pi$ is isomorphic to a subgroup
of $\mathcal{A}_\pi$.  \ps 

\begin{prop}\label{propsatotate} Assume the existence of $\mathcal{L}_\Z$ satisfying the axioms (L1)--(L5).  Let $\pi \in \Pi_{\rm alg}^\bot({\rm
PGL}_n)$ and assume $n\leq 8$.  Then $\mathcal{L}_\pi \simeq
\mathcal{A}_\pi$ unless : \begin{itemize}

\item[(i)]  $s(\pi)=(-1)^{n+1}$ and there exists a $\pi' \in \Pi_{\rm
alg}({\rm PGL}_2)$ such that $r_\pi \simeq {\rm Sym}^{n-1} r_{\pi'}$.  In this
case $\mathcal{L}_\pi\simeq {\rm SU}(2)$ if $n$ is even, ${\rm SO}(3)$ if
$n$ is odd.\ps

\item[(ii)]  $n=6$, $s(\pi)=-1$, and there exists two distinct $\pi',\pi''
\in \Pi_{\rm alg}({\rm PGL}_2)$ such that $r_{\pi}\simeq r_{\pi'} \otimes {\rm
Sym}^2 r_{\pi''}$.  In this case $\mathcal{L}_\pi \simeq {\rm SU}(2) \times
{\rm SO}(3)$.\ps

\item[(iii)] $n=7$ and $\mathcal{L}_\pi$ is the compact simple group of type ${\rm G}_2$.\ps

\item[(iv)] $n=8$, $s(\pi)=1$ and there exists $\pi' \in \Pi_{\rm alg}({\rm
PGL}_2)$, $\pi'' \in \Pi^{\rm s}_{\rm alg}({\rm PGL}_4)$, such that $r_{\pi}
\simeq r_{\pi'} \otimes r_{\pi''}$.  In this case $\mathcal{L}_\pi$ is the
quotient of ${\rm SU}(2) \times {\rm Spin}(5)$ by the diagonal central $\{\pm
1\}$.  \ps

\item[(v)] $n=8$, $s(\pi)=1$ and there exists two distinct $\pi',\pi'' \in
\Pi_{\rm alg}({\rm PGL}_2)$ such that $r_\pi \simeq r_{\pi'} \otimes {\rm Sym}^3
r_{\pi''}$.  In this case $\mathcal{L}_\pi$ is the quotient of ${\rm SU}(2)
\times {\rm SU}(2)$ by the diagonal central $\{\pm 1\}$.\ps

\item[(vi)] $n=8$, $s(\pi)=1$ and $\mathcal{L}_\pi \simeq {\rm Spin}(7)$. This occurs if and only if 
there exists $\pi' \in \Pi_{\rm alg}^{\rm o}({\rm PGL}_7)$ such that 
$$\mathcal{L}_{\pi'} \simeq {\rm SO}(7) \, \, \,{\rm and}\, \, \,   \rho \circ \xi  \simeq r_\pi,$$
where $\xi : \mathcal{L}_\Z \rightarrow {\rm Spin}(7)$ denotes the unique lift of $\widetilde{r_{\pi'}} :  \mathcal{L}_\Z \rightarrow {\rm SO}(7)$, and where $\rho$ denotes the Spin representation of ${\rm Spin}(7)$. \ps

\item[(vii)] $n=8$, $s(\pi)=-1$ and there exists distinct $\pi',\pi'',\pi'''
\in \Pi_{\rm alg}({\rm PGL}_2)$ such that $r_{\pi} \simeq r_{\pi'}\otimes r_{\pi''}
\otimes r_{\pi'''}$.  In this case $\mathcal{L}_\pi$ is the quotient of ${\rm
SU}(2)^3$ by the central subgroup $\{ (\epsilon_i)\in \{\pm 1\}^3,
\epsilon_1\epsilon_2\epsilon_3=1\}$.

\end{itemize}
\end{prop}

\begin{pf}We first observe that the only simply connected quasi-simple compact Lie groups having a
faithful self-dual finite dimensional irreducible representation of
dimension $\leq 8$ are in types : $A_1$ in each dimension, $B_2=C_2$ in dimensions $4$ and $5$, $C_3$ in dimension $6$,
$G_2$ and $B_3$ in dimension $7$, $A_2$ and $B_3$ in
dimension $8$, and $C_4$ in dimension $8$ (three representations
permuted by triality). \ps

The case $\mathcal{L}_\pi \simeq {\rm SU}(3)$ (type $A_2$),
equipped with its $8$-dimensional adjoint representation, does not occur. Indeed, if $r: {\rm W}_\R \rightarrow {\rm SU}(3)$ is a continuous
$3$-dimensional representation trivial on $\R_{>0} \subset {\rm W}_\C$, then the adjoint representation of $r$ on
${\rm Lie}({\rm SU}(3))$ is never multiplicity free, which contradicts $\pi \in \Pi_{\rm alg}^\bot(\PGL_8)$.  \ps
We conclude the proof by a case-by-case inspection.
\end{pf}

\newpage 
\section{Tables}\label{tables}

\begin{table}[htp]
\begin{center} $G={\rm SO}_7(\R)$, $\Gamma={\rm W}^+({\rm E}_7)$. \end{center}
\include{C_dim_SO7}
\caption{{\small The nonzero $d(\lambda)=\dim V_\lambda^\Gamma$ for $\lambda=(n_1,n_2,n_3)$ with $n_1 \leq 11$.}}
\label{tableSO7nue}
\end{table}
\clearpage
\newpage

\begin{table}[htp]
\begin{center} $G={\rm SO}_8(\R)$, $\Gamma={\rm W}^+({\rm E}_8)$. \end{center}
\include{C_dim_SO8}
\caption{{\small The nonzero $d(\lambda)=\dim V_\lambda^\Gamma$ for $\lambda=(n_1,n_2,n_3,n_4)$ with $n_1 \leq 11$.}}
\label{tableSO8nue}
\end{table}

\newpage \clearpage

\begin{table}[htp]
\begin{center} $G={\rm SO}_9(\R)$, $\Gamma={\rm W}({\rm E}_8)$. \end{center}

\include{C_dim_SO9}
\caption{{\small The nonzero $d(\lambda)=\dim V_\lambda^\Gamma$ for $\lambda=(n_1,n_2,n_3,n_4)$ with $n_1 \leq 9$.}}
\label{tableSO9nue}
\end{table}

\clearpage
\newpage

\begin{table}[htp]
\begin{center} $G={\rm G}_2(\R)$, $\Gamma={\rm G}_2(\Z)$. \end{center}
\include{C_dimG2}
\caption{{\small The nonzero $m(w,v)=\dim U_{w,v}^\Gamma$ for $v+w\leq 56$.}}
\label{tableg2nue}
\end{table}

\newpage
\clearpage

\begin{table}[htp]
\caption{{\small  The nonzero ${\rm S}(\underline{w})$ for $\underline{w}=(w_1,w_2)$ and $w_1\ \leq 43$,
using Tsushima's formula~\cite{tsushima}.}}
\include{C_dim_StSp4}
\label{tableSwv}
\end{table}

\newpage
\clearpage

\begin{table}[htp]
\caption{{\small The nonzero ${\rm S}(\underline{w})$ for $\underline{w}=(w_1,w_2,w_3)$ and $w_1\ \leq 29$.}}
\include{C_dim_StSO7}
\label{tableSwvu}
\end{table}

\newpage
\clearpage

\begin{table}[htp]
\caption{{\small The nonzero ${\rm S}(\underline{w})$ for $\underline{w}=(w_1,w_2,w_3,w_4)$ and $w_1\ \leq 27$.}}
\include{C_dim_StSO9}
\label{tableSwvut}
\end{table}

\newpage
\clearpage

\begin{table}[htp]
\caption{{\small The nonzero ${\rm O}(\underline{w})$ for $\underline{w}=(w_1,w_2,w_3,w_4)$ and $0<w_4<w_1\ \leq 30$.}}
\include{C_dim_StSO8_nonnul}
\label{tableOwvut}
\end{table}

\newpage
\clearpage

\begin{table}[htp]
\caption{{\small The nonzero ${\rm O}(\ww)=2\cdot {\rm O}(w_1,w_2,w_3,0)+{\rm O}^\ast(w_1,w_2,w_3)$ for  $w_1 \leq 34$.}}
\include{C_dim_StSO8_nul}
\label{tableOwvutbis}
\end{table}

\begin{table}[htp]
\caption{{\small The nonzero ${\rm G}_2(\ww)$ for $\ww=(w,v)$ and $w+v \leq 58$.}}
\include{C_dimStG2}
\label{tableStG2}
\end{table}

\newpage\clearpage
\begin{table}[htp]
\ps\ps
\renewcommand{\arraystretch}{1.5}
\caption{{\small The nonempty $\Pi_{w_1,w_2,w_3}({\rm SO}_7)$ for $w_1\leq 23$}}
\begin{tabular}{|c|c||c|c|}
\hline  $(w_1,w_2,w_3)$ & $\Pi_{w_1,w_2,w_3}({\rm SO}_7)$ & $(w_1,w_2,w_3)$ & $\Pi_{w_1,w_2,w_3}({\rm SO}_7)$ \\
\hline (5,3,1) & $[6]$ & (21,19,17) & $\Delta_{19}[3]$\\
\hline (13,11,9) & $\Delta_{11}[3]$ &  (23,9,1) & $\Delta_{23,9} \oplus [2]$\\

\hline (17,3,1) & $\Delta_{17}\oplus [4]$ &(23,11,7) & $\Delta_{23,7} \oplus \Delta_{11}$\\

\hline (17,11,1) & $\Delta_{17}\oplus \Delta_{11} \oplus [2]$ & (23,11,9) & $\Delta_{23,9} \oplus \Delta_{11}$\\

\hline (17,15,1) & $\Delta_{17}\oplus \Delta_{15} \oplus [2]$ & (23,13,1) & $\Delta_{23,13} \oplus [2]$\\

\hline (17,15,13) & $\Delta_{15}[3]$& (23,13,5) & $\Delta_{23,13,5}$\\

\hline (19,11,7) & $\Delta_{19,7}\oplus \Delta_{11}$& (23,15,3) & $\Delta_{23,15,3}$\\
\hline (19,15,7) & $\Delta_{19,7} \oplus \Delta_{15}$&  (23,15,7) & $\Delta_{23,7} \oplus \Delta_{15}$, $\Delta_{23,15,7}$\\
\hline (19,17,7) & $\Delta_{19,7} \oplus \Delta_{17}$ &(23,15,9) & $\Delta_{23,9} \oplus \Delta_{15}$\\

\hline (19,17,15) & $\Delta_{17}[3]$& (23,15,13) & $\Delta_{23,13} \oplus \Delta_{15}$\\
\hline (21,3,1) & $\Delta_{21}\oplus [4]$&  (23,17,5) & $\Delta_{23,17,5} $\\
\hline (21,11,1) & $\Delta_{21}\oplus \Delta_{11}\oplus [2]$&(23,17,7) & $\Delta_{23,7} \oplus \Delta_{17}$\\

\hline (21,11,5) & $\Delta_{21,5}\oplus \Delta_{11}$& (23,17,9) & $\Delta_{23,9} \oplus \Delta_{17}$, $\Delta_{23,17,9}$\\

\hline (21,11,9) & $\Delta_{21,9}\oplus \Delta_{11}$&  (23,17,13) & $\Delta_{23,13} \oplus \Delta_{17}$\\
\hline (21,15,1) & $\Delta_{21} \oplus \Delta_{15} \oplus [2]$&(23,19,3) & $\Delta_{23,19,3}$\\
\hline (21,15,5) & $\Delta_{21,5} \oplus \Delta_{15}$ & (23,19,7) & $\Delta_{23,7} \oplus \Delta_{19}$\\
\hline (21,15,9) & $\Delta_{21,9} \oplus \Delta_{15}$ & (23,19,9) & $\Delta_{23,9} \oplus \Delta_{19}$\\
\hline (21,15,13) & $\Delta_{21,13} \oplus \Delta_{15}$&  (23,19,11) & $\Delta_{23,19,11}$\\

\hline (21,17,5)&  $\Delta_{21,5} \oplus \Delta_{17}$ &(23,19,13) & $\Delta_{23,13} \oplus \Delta_{19}$\\

\hline (21,17,9) & $\Delta_{21,9} \oplus \Delta_{17}$&  (23,21,1) & ${\rm Sym}^2 \Delta_{11}[2]$\\
\hline (21,17,13) & $\Delta_{21,13}\oplus \Delta_{17}$&  (23,21,7) & $\Delta_{23,7} \oplus \Delta_{21}$\\

\hline (21,19,1) & $\Delta_{21}\oplus \Delta_{19} \oplus [2]$&(23,21,9) & $\Delta_{23,9} \oplus \Delta_{21}$\\
\hline (21,19,5) & $\Delta_{21,5} \oplus \Delta_{19}$& (23,21,13) & $\Delta_{23,13}\oplus \Delta_{21}$\\

\hline (21,19,9) & $\Delta_{21,9} \oplus \Delta_{19}$& (23,21,19) & $\Delta_{21}[3]$\\

\hline  (21,19,13) & $\Delta_{21,13} \oplus \Delta_{19}$ & &\\ \hline

\end{tabular}
\label{tableSO7}
\end{table}
\newpage 

\begin{table}[htp]
\ps\ps
\renewcommand{\arraystretch}{1.5}
\caption{{\small The nonempty $\Pi_{25,w_2,w_3}({\rm SO}_7)$ }}
\ps\ps
\begin{tabular}{|c|c||c|c|}
\hline  $(w_1,w_2,w_3)$ & $\Pi_{w_1,w_2,w_3}({\rm SO}_7)$ & $(w_1,w_2,w_3)$ & $\Pi_{w_1,w_2,w_3}({\rm SO}_7)$ \\
\hline (25,3,1) & $\Delta_{25}\oplus [4]$ & (25,19,9) & $\Delta_{25,9}^2 \oplus \Delta_{19}$, $\Delta_{25,19,9}^2$ \\
\hline (25,7,1) & $\Delta_{25,7} \oplus [2]$ & (25,19,11) & $\Delta_{25,11} \oplus \Delta_{19}$ \\
\hline (25,11,1) & $\Delta_{25,11}\oplus [2]$, $\Delta_{25}\oplus\Delta_{11}\oplus [2]$ &  (25,19,13) & $\Delta_{25,13}^2 \oplus \Delta_{19}$, $\Delta_{25,19,13}$ \\
\hline (25,11,5) & $\Delta_{25,5} \oplus \Delta_{11}$ & (25,19,15) & $\Delta_{25,15}\oplus \Delta_{19}$ \\
\hline (25,11,7) & $\Delta_{25,7} \oplus \Delta_{11}$ & (25,19,17) & $\Delta_{25,17} \oplus \Delta_{19}$ \\
\hline (25,11,9) & $\Delta_{25,9}^2\oplus \Delta_{11}$ & (25,21,3) & $\Delta_{25,21,3}^2$ \\
\hline (25,13,3) & $\Delta_{25,13,3}$ & (25,21,5) & $\Delta_{25,5} \oplus \Delta_{21}$ \\
\hline (25,13,7) & $\Delta_{25,13,7}$ & (25,21,7) & $\Delta_{25,7} \oplus \Delta_{21}$, $\Delta_{25,21,7}^2$ \\
\hline (25,15,1) & $\Delta_{25,15}\oplus[2]$, $\Delta_{25}\oplus \Delta_{15} \oplus [2]$ & (25,21,9) & $\Delta_{25,9}^2 \oplus \Delta_{21}$ \\
\hline (25,15,5) & $\Delta_{25,5} \oplus \Delta_{15}$, $\Delta_{25,15,5}$ & (25,21,11) & $\Delta_{25,11} \oplus \Delta_{21}$, $\Delta_{25,21,11}^2$ \\
\hline (25,15,7) & $\Delta_{25,7} \oplus \Delta_{15}$ & (25,21,13) & $\Delta_{25,13}^2 \oplus \Delta_{21}$ \\
\hline (25,15,9) & $\Delta_{25,9}^2\oplus \Delta_{15}$, $\Delta_{25,15,9}$ & (25,21,15) & $\Delta_{25,15} \oplus \Delta_{21}$, $\Delta_{25,21,15}$ \\
\hline (25,15,11) & $\Delta_{25,11} \oplus \Delta_{15}$ & (25,21,17) & $\Delta_{25,17} \oplus \Delta_{21}$ \\
\hline (25,15,13) & $\Delta_{25,13}^2 \oplus \Delta_{15}$ & (25,21,19) & $\Delta_{25,19} \oplus \Delta_{21}$ \\
\hline (25,17,3) & $\Delta_{25,17,3}^2$ & (25,23,1) & $\Delta_{25}\oplus \Delta_{23}^2 \oplus [2]$ \\
\hline (25,17,5) & $\Delta_{25,5} \oplus \Delta_{17}$ & (25,23,5) & $\Delta_{25,5} \oplus \Delta_{23}^2$ \\
\hline (25,17,7) & $\Delta_{25,7}\oplus \Delta_{17}$, $\Delta_{25,17,7}^2$ & (25,23,7) & $\Delta_{25,7} \oplus \Delta_{23}^2$\\
\hline (25,17,9) & $\Delta_{25,9}^2 \oplus \Delta_{17}$ & (25,23,9) & $\Delta_{25,9}^2\oplus \Delta_{23}^2$ \\
\hline (25,17,11) & $\Delta_{25,11} \oplus \Delta_{17}$, $\Delta_{25,17,11}$ & (25,23,11) & $\Delta_{25,11}\oplus \Delta_{23}^2$ \\
\hline (25,17,13) & $\Delta_{25,13}^2 \oplus \Delta_{17}$ & (25,23,13) & $\Delta_{25,13}^2\oplus \Delta_{23}^2$ \\
\hline (25,17,15) & $\Delta_{25,15} \oplus \Delta_{17}$ & (25,23,15) & $\Delta_{25,15}\oplus \Delta_{23}^2$ \\
\hline (25,19,1) & $\Delta_{25,19} \oplus [2]$, $\Delta_{25} \oplus \Delta_{19} \oplus [2]$, $\Delta_{25,19,1}$ & (25,23,17) & $\Delta_{25,17}\oplus \Delta_{23}^2$ \\
\hline (25,19,5) & $\Delta_{25,5} \oplus \Delta_{19}$, $\Delta_{25,19,5}^2$ & (25,23,19) & $\Delta_{25,19}\oplus \Delta_{23}^2$ \\
\hline (25,19,7) & $\Delta_{25,7}\oplus \Delta_{19}$ & (25,23,21) & $\Delta_{23}^2[3]$ \\
\hline
\end{tabular}
\label{tableSO72}
\end{table}

\newpage
\clearpage

\begin{table}[htp]\ps\ps
\renewcommand{\arraystretch}{1.2}
\caption{{\small The nonempty $\Pi_{w_1,w_2,w_3,w_4}({\rm SO}_9)$ for $w_1 \leq 23$ }}
\ps\ps
\begin{tabular}{|c|c||c|c|}
\hline  $(w_1,w_2,w_3,w_4)$ & $\Pi_{w_1,w_2,w_3,w_4}({\rm SO}_9)$ & $(w_1,w_2,w_3,w_4)$ & $\Pi_{w_1,w_2,w_3,w_4}({\rm SO}_9)$ \\
\hline (7, 5, 3, 1) & $[8] $ & (23, 17, 15, 5) & $\Delta_{23,17,5} \oplus \Delta_{15}$ \\

\hline (11, 5, 3, 1) & $\Delta_{11} \oplus [6]$& (23, 17, 15, 9) & $\Delta_{23,17,9} \oplus \Delta_{15}$ \\
\hline (15, 5, 3, 1) & $\Delta_{15} \oplus [6]$& (23, 17, 15, 13) & $\Delta_{23}^2 \oplus \Delta_{15}[3]$ \\
\hline (15, 13, 11, 9) & $\Delta_{15} \oplus \Delta_{11}[3]$& (23, 19, 9, 7) & $\Delta _{23,9} \oplus \Delta_{19,7}$\\

\hline (17, 13, 11, 9) & $\Delta_{17} \oplus \Delta_{11}[3]$& (23, 19, 11, 3) & $\Delta_{23,19,3} \oplus \Delta_{11}$\\

\hline (19, 5, 3, 1) & $\Delta_{19} \oplus [6]$& (23, 19, 11, 7) & $\Delta_{23}^2 \oplus \Delta_{19,7} \oplus \Delta_{11}$\\
\hline (19, 13, 11, 9) & $\Delta_{19} \oplus \Delta_{11}[3]$& (23, 19, 13, 7) & $\Delta_{23,13} \oplus \Delta_{19,7}$ \\

\hline (19, 17, 3, 1) & $\Delta_{19}\oplus \Delta_{17} \oplus [4]$& (23, 19, 15, 3) & $\Delta_{23,19,3} \oplus
\Delta_{15}$ \\
\hline (19, 17, 7, 1) & $\Delta_{19,7}\oplus \Delta_{17} \oplus [2]$& (23, 19, 15, 7) & $\Delta_{23}^2\oplus \Delta_{19,7} \oplus \Delta_{15}$\\

\hline (19, 17, 11, 1) & $\Delta_{19}\oplus \Delta_{17}\oplus \Delta_{11} \oplus [2]$& (23, 19, 15, 11) & $\Delta_{23,19,11} \oplus \Delta_{15}$\\
\hline (19, 17, 15, 1) & $\Delta_{19}\oplus \Delta_{17}\oplus \Delta_{15} \oplus [2]$& (23, 19, 17, 3) & $\Delta_{23,19,3} \oplus \Delta_{17}$ \\
\hline (19, 17, 15, 13) & $\Delta_{19} \oplus \Delta_{15}[3]$& (23, 19, 17, 7) & $\Delta_{23}^2\oplus \Delta_{19,7} \oplus \Delta_{17}$\\
\hline (21, 13, 11, 9) & $\Delta_{21} \oplus \Delta_{11}[3]$& (23, 19, 17, 11) & $\Delta_{23,19,11} \oplus \Delta_{17}$\\

\hline (21, 17, 5, 1) & $\Delta_{21,5} \oplus \Delta_{17} \oplus [2]$& (23, 19, 17, 15) & $\Delta_{23}^2\oplus \Delta_{17}[3]$\\

\hline (21, 17, 9, 1) & $\Delta_{21,9} \oplus \Delta_{17} \oplus [2]$& (23, 21, 3, 1) & $\Delta_{23}^2\oplus \Delta_{21} \oplus [4]$\\

\hline (21, 17, 13, 1) & $\Delta_{21,13} \oplus \Delta_{17} \oplus [2]$& (23, 21, 7, 1) & $\Delta_{23,7}\oplus \Delta_{21} \oplus [2]$\\
\hline (21, 17, 15, 13) & $\Delta_{21} \oplus \Delta_{15}[3]$& (23, 21, 7, 5) & $\Delta_{23,7}\oplus \Delta_{21,5}$\\
\hline (21, 19, 9, 7) & $\Delta_{21,9} \oplus \Delta_{19,7}$& (23, 21, 9, 5) & $\Delta_{23,9}\oplus \Delta_{21,5}$\\

\hline (21, 19, 11, 7) & $\Delta_{21} \oplus \Delta_{19,7} \oplus \Delta_{11}$& (23, 21, 11, 1) & $\Delta_{23}^2 \oplus \Delta_{21} \oplus \Delta_{11} \oplus [2]$\\
\hline  (21, 19, 13, 7) & $\Delta_{21,13} \oplus \Delta_{19,7}$& (23, 21, 11, 5) & $\Delta_{23}^2\oplus \Delta_{21,5}\oplus \Delta_{11}$\\
\hline (21, 19, 15, 7) & $\Delta_{21} \oplus \Delta_{19,7} \oplus \Delta_{15}$& (23, 21, 11, 9) & $\Delta_{23}^2\oplus \Delta_{21,9}\oplus \Delta_{11}$\\

\hline (21, 19, 17, 7) & $\Delta_{21} \oplus \Delta_{19,7} \oplus \Delta_{17}$& (23, 21, 13, 5) & $\Delta_{23,13}\oplus \Delta_{21,5}$\\

\hline (21, 19, 17, 15) & $\Delta_{21} \oplus \Delta_{17}[3]$& (23, 21, 13, 9) & $\Delta_{23,13}\oplus \Delta_{21,9}$\\

\hline (23, 5, 3, 1) & $\Delta_{23}^2 \oplus [6]$& (23, 21, 15, 1) & $\Delta_{23}^2\oplus \Delta_{21} \oplus \Delta_{15} \oplus [2]$\\

\hline (23, 9, 3, 1) & $\Delta_{23,9} \oplus [4]$& (23, 21, 15, 5) & $\Delta_{23}^2\oplus\Delta_{21,5}\oplus \Delta_{15}$\\

\hline (23, 13, 3, 1) & $\Delta_{23,13}\oplus [4]$& (23, 21, 15, 9) & $\Delta_{23}^2\oplus \Delta_{21,9}\oplus \Delta_{15}$\\

\hline (23, 13, 11, 1) & $\Delta_{23,13} \oplus \Delta_{11} \oplus [2]$& (23, 21, 15, 13) & $\Delta_{23}^2 \oplus \Delta_{21,13} \oplus \Delta_{15}$\\
\hline (23, 13, 11, 5) & $\Delta_{23,13,5} \oplus \Delta_{11}$& (23, 21, 17, 1) & ${\rm Sym}^2 \Delta_{11}[2] \oplus \Delta_{17}$\\
\hline (23, 13, 11, 9) & $\Delta_{23}^2\oplus \Delta_{11}[3]$& (23, 21, 17, 5) & $\Delta_{23}^2 \oplus \Delta_{21,5} \oplus \Delta_{17}$\\
\hline (23, 15, 11, 3) & $\Delta_{23,15,3} \oplus \Delta_{11}$& (23, 21, 17, 9) & $\Delta_{23}^2 \oplus \Delta_{21,9} \oplus \Delta_{17}$\\

\hline (23, 15, 11, 7) & $\Delta_{23,15,7} \oplus \Delta_{11}$& (23, 21, 17, 13) & $\Delta_{23}^2 \oplus \Delta_{21,13} \oplus \Delta_{17}$\\
\hline (23, 17, 3, 1) & $\Delta_{23}^2\oplus \Delta_{17}\oplus  [4]$& (23, 21, 19, 1) & $\Delta_{23}^2\oplus \Delta_{21} \oplus \Delta_{19}\oplus [2]$\\
\hline (23, 17, 7, 1) & $\Delta_{23,7} \oplus \Delta_{17} \oplus [2]$& (23, 21, 19, 5) & $\Delta_{23}^2\oplus \Delta_{21,5}\oplus \Delta_{19}$\\
\hline (23, 17, 11, 1) & $\Delta_{23}^2\oplus \Delta_{17}\oplus \Delta_{11}\oplus [2]$& (23, 21, 19, 9) & $\Delta_{23}^2\oplus \Delta_{21,9}\oplus \Delta_{19}$\\

\hline (23, 17, 11, 5) & $\Delta_{23,17,5}\oplus \Delta_{11}$& (23, 21, 19, 13) & $\Delta_{23}^2 \oplus \Delta_{21,13}\oplus \Delta_{19}$\\
\hline (23, 17, 11, 9) & $\Delta_{23,17,9} \oplus \Delta_{11}$& (23, 21, 19, 17) & $\Delta_{23}^2 \oplus \Delta_{19}[3]$\\

\hline (23, 17, 15, 1) & $\Delta_{23}^2\oplus \Delta_{17}\oplus \Delta_{15}\oplus [2]$& (25, 7, 3, 1)& $\Delta_{25,7} \oplus [4]$ \\
\hline
\end{tabular}
\label{tableSO91}
\end{table}

\newpage
\clearpage

\bigskip
\bigskip

\section{The $121$ level $1$ automorphic representations of ${\rm SO}_{25}$ with trivial coefficients}
\label{table121}
\bigskip
\bigskip
{\setlength{\baselineskip}{1.5\baselineskip}
\begin{multicols}{2}
\include{C_listeformesautomorphes_25}

\end{multicols}
\par}

\end{document}

%% file: C_dim_SO7.tex
\begin{tabular}{|c|c||c|c||c|c||c|c||c|c|}
\hline  $\lambda$ & $d(\lambda)$  & $\lambda$ & $d(\lambda)$  & $\lambda$ & $d(\lambda)$  & $\lambda$ & $d(\lambda)$  & $\lambda$ & $d(\lambda)$ \\ 
\hline
(0, 0, 0) & 1& (9, 6, 3) & 2& (10, 7, 2) & 1& (10, 10, 10) & 2& (11, 9, 0) & 2 \\
\hline
(4, 4, 4) & 1& (9, 6, 4) & 1& (10, 7, 3) & 3& (11, 3, 0) & 1& (11, 9, 1) & 1 \\
\hline
(6, 0, 0) & 1& (9, 6, 6) & 1& (10, 7, 4) & 2& (11, 3, 2) & 1& (11, 9, 2) & 4 \\
\hline
(6, 4, 0) & 1& (9, 7, 2) & 1& (10, 7, 5) & 2& (11, 4, 1) & 1& (11, 9, 3) & 4 \\
\hline
(6, 6, 0) & 1& (9, 7, 3) & 1& (10, 7, 6) & 2& (11, 4, 3) & 2& (11, 9, 4) & 5 \\
\hline
(6, 6, 6) & 1& (9, 7, 4) & 2& (10, 7, 7) & 1& (11, 4, 4) & 1& (11, 9, 5) & 4 \\
\hline
(7, 4, 3) & 1& (9, 7, 6) & 1& (10, 8, 0) & 3& (11, 5, 0) & 2& (11, 9, 6) & 5 \\
\hline
(7, 6, 3) & 1& (9, 8, 1) & 1& (10, 8, 2) & 3& (11, 5, 2) & 2& (11, 9, 7) & 3 \\
\hline
(7, 7, 3) & 1& (9, 8, 3) & 1& (10, 8, 3) & 1& (11, 5, 3) & 1& (11, 9, 8) & 2 \\
\hline
(7, 7, 7) & 1& (9, 8, 4) & 1& (10, 8, 4) & 4& (11, 5, 4) & 1& (11, 9, 9) & 1 \\
\hline
(8, 0, 0) & 1& (9, 8, 5) & 1& (10, 8, 5) & 1& (11, 6, 1) & 2& (11, 10, 1) & 3 \\
\hline
(8, 4, 0) & 1& (9, 8, 6) & 1& (10, 8, 6) & 3& (11, 6, 2) & 1& (11, 10, 2) & 3 \\
\hline
(8, 4, 2) & 1& (9, 9, 0) & 1& (10, 8, 7) & 1& (11, 6, 3) & 4& (11, 10, 3) & 5 \\
\hline
(8, 4, 4) & 1& (9, 9, 3) & 1& (10, 8, 8) & 1& (11, 6, 4) & 2& (11, 10, 4) & 4 \\
\hline
(8, 6, 0) & 1& (9, 9, 4) & 1& (10, 9, 1) & 2& (11, 6, 5) & 2& (11, 10, 5) & 6 \\
\hline
(8, 6, 2) & 1& (9, 9, 6) & 1& (10, 9, 2) & 1& (11, 6, 6) & 2& (11, 10, 6) & 5 \\
\hline
(8, 6, 4) & 1& (9, 9, 9) & 1& (10, 9, 3) & 3& (11, 7, 0) & 1& (11, 10, 7) & 5 \\
\hline
(8, 6, 6) & 1& (10, 0, 0) & 1& (10, 9, 4) & 2& (11, 7, 1) & 1& (11, 10, 8) & 3 \\
\hline
(8, 7, 2) & 1& (10, 2, 0) & 1& (10, 9, 5) & 3& (11, 7, 2) & 4& (11, 10, 9) & 2 \\
\hline
(8, 7, 4) & 1& (10, 4, 0) & 2& (10, 9, 6) & 2& (11, 7, 3) & 3& (11, 10, 10) & 2 \\
\hline
(8, 7, 6) & 1& (10, 4, 2) & 1& (10, 9, 7) & 2& (11, 7, 4) & 4& (11, 11, 1) & 1 \\
\hline
(8, 8, 0) & 1& (10, 4, 3) & 1& (10, 9, 8) & 1& (11, 7, 5) & 3& (11, 11, 2) & 2 \\
\hline
(8, 8, 2) & 1& (10, 4, 4) & 2& (10, 9, 9) & 1& (11, 7, 6) & 3& (11, 11, 3) & 3 \\
\hline
(8, 8, 4) & 1& (10, 5, 1) & 1& (10, 10, 0) & 2& (11, 7, 7) & 2& (11, 11, 4) & 2 \\
\hline
(8, 8, 6) & 1& (10, 5, 3) & 1& (10, 10, 2) & 2& (11, 8, 1) & 3& (11, 11, 5) & 3 \\
\hline
(8, 8, 8) & 1& (10, 6, 0) & 2& (10, 10, 3) & 2& (11, 8, 2) & 2& (11, 11, 6) & 3 \\
\hline
(9, 3, 0) & 1& (10, 6, 2) & 2& (10, 10, 4) & 4& (11, 8, 3) & 5& (11, 11, 7) & 3 \\
\hline
(9, 4, 3) & 1& (10, 6, 3) & 1& (10, 10, 5) & 2& (11, 8, 4) & 4& (11, 11, 8) & 2 \\
\hline
(9, 4, 4) & 1& (10, 6, 4) & 3& (10, 10, 6) & 4& (11, 8, 5) & 5& (11, 11, 9) & 1 \\
\hline
(9, 5, 0) & 1& (10, 6, 5) & 1& (10, 10, 7) & 2& (11, 8, 6) & 4& (11, 11, 10) & 1 \\
\hline
(9, 5, 2) & 1& (10, 6, 6) & 2& (10, 10, 8) & 2& (11, 8, 7) & 3& (11, 11, 11) & 1 \\
\hline
(9, 6, 1) & 1& (10, 7, 1) & 2& (10, 10, 9) & 2& (11, 8, 8) & 1& (12, 0, 0) & 2 \\
\hline
\end{tabular}

%% file: C_dim_SO8.tex
\begin{tabular}{|c|c||c|c||c|c||c|c|}
\hline  $\lambda$ & $d(\lambda)$  & $\lambda$ & $d(\lambda)$  & $\lambda$ & $d(\lambda)$  & $\lambda$ & $d(\lambda)$ \\ 
\hline
(0, 0, 0, 0) & 1& (10, 9, 1, 0) & 1& (11, 8, 5, 2) & 1& (11, 11, 3, 3) & 2 \\
\hline
(4, 4, 4, 4) & 1& (10, 9, 4, 3) & 1& (11, 8, 6, 1) & 1& (11, 11, 4, 4) & 1 \\
\hline
(6, 6, 0, 0) & 1& (10, 9, 5, 0) & 1& (11, 8, 6, 3) & 1& (11, 11, 5, 1) & 1 \\
\hline
(6, 6, 6, 6) & 1& (10, 9, 6, 1) & 1& (11, 8, 7, 2) & 1& (11, 11, 5, 5) & 2 \\
\hline
(7, 7, 3, 3) & 1& (10, 9, 7, 0) & 1& (11, 8, 7, 4) & 1& (11, 11, 6, 2) & 1 \\
\hline
(7, 7, 7, 7) & 1& (10, 9, 9, 2) & 1& (11, 8, 8, 3) & 1& (11, 11, 6, 6) & 2 \\
\hline
(8, 0, 0, 0) & 1& (10, 10, 0, 0) & 1& (11, 9, 2, 0) & 1& (11, 11, 7, 1) & 2 \\
\hline
(8, 4, 4, 0) & 1& (10, 10, 2, 2) & 1& (11, 9, 3, 1) & 1& (11, 11, 7, 7) & 2 \\
\hline
(8, 6, 6, 0) & 1& (10, 10, 3, 3) & 1& (11, 9, 4, 2) & 1& (11, 11, 8, 0) & 2 \\
\hline
(8, 7, 7, 0) & 1& (10, 10, 4, 0) & 1& (11, 9, 5, 1) & 1& (11, 11, 8, 4) & 1 \\
\hline
(8, 8, 0, 0) & 1& (10, 10, 4, 4) & 2& (11, 9, 5, 3) & 1& (11, 11, 8, 8) & 1 \\
\hline
(8, 8, 2, 2) & 1& (10, 10, 5, 5) & 1& (11, 9, 6, 0) & 2& (11, 11, 9, 3) & 1 \\
\hline
(8, 8, 4, 4) & 1& (10, 10, 6, 0) & 1& (11, 9, 6, 4) & 1& (11, 11, 9, 9) & 1 \\
\hline
(8, 8, 6, 6) & 1& (10, 10, 6, 2) & 1& (11, 9, 7, 1) & 1& (11, 11, 10, 2) & 1 \\
\hline
(8, 8, 8, 0) & 1& (10, 10, 6, 6) & 2& (11, 9, 7, 3) & 1& (11, 11, 10, 10) & 1 \\
\hline
(8, 8, 8, 8) & 1& (10, 10, 7, 1) & 1& (11, 9, 7, 5) & 1& (11, 11, 11, 3) & 1 \\
\hline
(9, 6, 3, 0) & 1& (10, 10, 7, 7) & 1& (11, 9, 8, 2) & 1& (11, 11, 11, 11) & 1 \\
\hline
(9, 7, 4, 2) & 1& (10, 10, 8, 0) & 1& (11, 9, 9, 3) & 1& (12, 0, 0, 0) & 1 \\
\hline
(9, 8, 6, 1) & 1& (10, 10, 8, 4) & 1& (11, 10, 1, 0) & 1& (12, 4, 0, 0) & 1 \\
\hline
(9, 9, 3, 3) & 1& (10, 10, 8, 8) & 1& (11, 10, 3, 2) & 1& (12, 4, 4, 0) & 1 \\
\hline
(9, 9, 4, 4) & 1& (10, 10, 9, 3) & 1& (11, 10, 4, 1) & 1& (12, 4, 4, 4) & 1 \\
\hline
(9, 9, 6, 6) & 1& (10, 10, 9, 9) & 1& (11, 10, 4, 3) & 1& (12, 5, 3, 2) & 1 \\
\hline
(9, 9, 9, 9) & 1& (10, 10, 10, 2) & 2& (11, 10, 5, 0) & 1& (12, 6, 0, 0) & 1 \\
\hline
(10, 4, 0, 0) & 1& (10, 10, 10, 10) & 2& (11, 10, 5, 2) & 1& (12, 6, 2, 0) & 1 \\
\hline
(10, 4, 4, 2) & 1& (11, 4, 4, 3) & 1& (11, 10, 5, 4) & 1& (12, 6, 4, 0) & 1 \\
\hline
(10, 6, 0, 0) & 1& (11, 5, 2, 0) & 1& (11, 10, 6, 1) & 2& (12, 6, 4, 2) & 1 \\
\hline
(10, 6, 4, 0) & 1& (11, 6, 3, 0) & 1& (11, 10, 6, 5) & 1& (12, 6, 6, 0) & 2 \\
\hline
(10, 6, 6, 2) & 1& (11, 6, 4, 3) & 1& (11, 10, 7, 0) & 3& (12, 6, 6, 4) & 1 \\
\hline
(10, 7, 4, 1) & 1& (11, 6, 6, 3) & 1& (11, 10, 7, 2) & 1& (12, 7, 3, 0) & 1 \\
\hline
(10, 7, 6, 3) & 1& (11, 7, 3, 1) & 1& (11, 10, 7, 6) & 1& (12, 7, 3, 2) & 1 \\
\hline
(10, 7, 7, 2) & 1& (11, 7, 4, 0) & 1& (11, 10, 8, 1) & 1& (12, 7, 4, 1) & 1 \\
\hline
(10, 8, 2, 0) & 1& (11, 7, 5, 1) & 1& (11, 10, 8, 3) & 1& (12, 7, 4, 3) & 1 \\
\hline
(10, 8, 4, 0) & 1& (11, 7, 6, 2) & 1& (11, 10, 9, 2) & 1& (12, 7, 5, 2) & 1 \\
\hline
(10, 8, 4, 2) & 1& (11, 7, 7, 3) & 1& (11, 10, 9, 4) & 1& (12, 7, 6, 1) & 1 \\
\hline
(10, 8, 6, 0) & 1& (11, 8, 3, 0) & 1& (11, 10, 10, 3) & 2& (12, 7, 6, 3) & 1 \\
\hline
(10, 8, 6, 4) & 1& (11, 8, 4, 1) & 1& (11, 11, 1, 1) & 1& (12, 7, 7, 0) & 1 \\
\hline
(10, 8, 8, 2) & 1& (11, 8, 5, 0) & 1& (11, 11, 2, 2) & 1& (12, 7, 7, 4) & 1 \\
\hline
\end{tabular}

%% file: C_dim_SO9.tex
\begin{tabular}{|c|c||c|c||c|c||c|c||c|c|}
\hline  $\lambda$ & $d(\lambda)$  & $\lambda$ & $d(\lambda)$  & $\lambda$ & $d(\lambda)$  & $\lambda$ & $d(\lambda)$  & $\lambda$ & $d(\lambda)$ \\ 
\hline
(0, 0, 0, 0) & 1& (8, 6, 6, 2) & 1& (9, 3, 0, 0) & 1& (9, 7, 7, 3) & 2& (9, 8, 8, 8) & 1 \\
\hline
(2, 0, 0, 0) & 1& (8, 6, 6, 4) & 1& (9, 4, 4, 1) & 1& (9, 7, 7, 4) & 2& (9, 9, 3, 1) & 1 \\
\hline
(4, 0, 0, 0) & 1& (8, 6, 6, 6) & 2& (9, 4, 4, 3) & 1& (9, 7, 7, 6) & 1& (9, 9, 3, 3) & 2 \\
\hline
(4, 4, 4, 4) & 1& (8, 7, 3, 3) & 1& (9, 4, 4, 4) & 1& (9, 7, 7, 7) & 1& (9, 9, 4, 0) & 1 \\
\hline
(5, 4, 4, 4) & 1& (8, 7, 4, 1) & 1& (9, 5, 0, 0) & 1& (9, 8, 1, 0) & 1& (9, 9, 4, 2) & 2 \\
\hline
(6, 0, 0, 0) & 1& (8, 7, 4, 3) & 2& (9, 5, 4, 0) & 1& (9, 8, 2, 2) & 1& (9, 9, 4, 3) & 2 \\
\hline
(6, 4, 4, 4) & 1& (8, 7, 5, 3) & 1& (9, 5, 4, 2) & 1& (9, 8, 3, 0) & 2& (9, 9, 4, 4) & 2 \\
\hline
(6, 6, 0, 0) & 1& (8, 7, 6, 1) & 1& (9, 5, 4, 4) & 1& (9, 8, 3, 2) & 2& (9, 9, 5, 1) & 1 \\
\hline
(6, 6, 2, 0) & 1& (8, 7, 6, 3) & 2& (9, 6, 1, 0) & 1& (9, 8, 4, 1) & 2& (9, 9, 5, 2) & 1 \\
\hline
(6, 6, 4, 0) & 1& (8, 7, 6, 5) & 1& (9, 6, 3, 0) & 2& (9, 8, 4, 2) & 2& (9, 9, 5, 3) & 3 \\
\hline
(6, 6, 6, 0) & 1& (8, 7, 7, 1) & 1& (9, 6, 3, 2) & 1& (9, 8, 4, 3) & 3& (9, 9, 5, 4) & 2 \\
\hline
(6, 6, 6, 6) & 1& (8, 7, 7, 3) & 2& (9, 6, 4, 1) & 2& (9, 8, 4, 4) & 2& (9, 9, 6, 0) & 1 \\
\hline
(7, 4, 4, 4) & 1& (8, 7, 7, 5) & 1& (9, 6, 4, 3) & 2& (9, 8, 5, 0) & 2& (9, 9, 6, 1) & 1 \\
\hline
(7, 6, 1, 0) & 1& (8, 7, 7, 7) & 2& (9, 6, 5, 0) & 2& (9, 8, 5, 2) & 3& (9, 9, 6, 2) & 3 \\
\hline
(7, 6, 3, 0) & 1& (8, 8, 0, 0) & 2& (9, 6, 5, 2) & 1& (9, 8, 5, 3) & 1& (9, 9, 6, 3) & 3 \\
\hline
(7, 6, 5, 0) & 1& (8, 8, 2, 0) & 1& (9, 6, 6, 1) & 2& (9, 8, 5, 4) & 2& (9, 9, 6, 4) & 3 \\
\hline
(7, 6, 6, 6) & 1& (8, 8, 2, 2) & 1& (9, 6, 6, 3) & 2& (9, 8, 6, 1) & 3& (9, 9, 6, 5) & 1 \\
\hline
(7, 7, 3, 3) & 1& (8, 8, 3, 2) & 1& (9, 6, 6, 5) & 1& (9, 8, 6, 2) & 3& (9, 9, 6, 6) & 2 \\
\hline
(7, 7, 4, 3) & 1& (8, 8, 4, 0) & 2& (9, 6, 6, 6) & 1& (9, 8, 6, 3) & 4& (9, 9, 7, 1) & 1 \\
\hline
(7, 7, 5, 3) & 1& (8, 8, 4, 2) & 2& (9, 7, 0, 0) & 1& (9, 8, 6, 4) & 3& (9, 9, 7, 2) & 2 \\
\hline
(7, 7, 6, 3) & 1& (8, 8, 4, 4) & 2& (9, 7, 3, 1) & 1& (9, 8, 6, 5) & 2& (9, 9, 7, 3) & 3 \\
\hline
(7, 7, 7, 3) & 1& (8, 8, 5, 2) & 1& (9, 7, 3, 3) & 2& (9, 8, 6, 6) & 2& (9, 9, 7, 4) & 3 \\
\hline
(7, 7, 7, 7) & 1& (8, 8, 5, 4) & 1& (9, 7, 4, 0) & 2& (9, 8, 7, 0) & 1& (9, 9, 7, 5) & 1 \\
\hline
(8, 0, 0, 0) & 2& (8, 8, 6, 0) & 2& (9, 7, 4, 2) & 3& (9, 8, 7, 1) & 2& (9, 9, 7, 6) & 2 \\
\hline
(8, 2, 0, 0) & 1& (8, 8, 6, 2) & 2& (9, 7, 4, 3) & 2& (9, 8, 7, 2) & 3& (9, 9, 8, 1) & 1 \\
\hline
(8, 4, 0, 0) & 1& (8, 8, 6, 4) & 2& (9, 7, 4, 4) & 2& (9, 8, 7, 3) & 3& (9, 9, 8, 2) & 1 \\
\hline
(8, 4, 4, 0) & 1& (8, 8, 6, 6) & 2& (9, 7, 5, 1) & 1& (9, 8, 7, 4) & 3& (9, 9, 8, 3) & 2 \\
\hline
(8, 4, 4, 2) & 1& (8, 8, 7, 0) & 1& (9, 7, 5, 2) & 1& (9, 8, 7, 5) & 2& (9, 9, 8, 4) & 2 \\
\hline
(8, 4, 4, 4) & 2& (8, 8, 7, 2) & 2& (9, 7, 5, 3) & 3& (9, 8, 7, 6) & 2& (9, 9, 8, 5) & 1 \\
\hline
(8, 5, 4, 1) & 1& (8, 8, 7, 4) & 2& (9, 7, 5, 4) & 1& (9, 8, 7, 7) & 1& (9, 9, 8, 6) & 2 \\
\hline
(8, 5, 4, 3) & 1& (8, 8, 7, 6) & 2& (9, 7, 6, 0) & 2& (9, 8, 8, 1) & 2& (9, 9, 9, 3) & 1 \\
\hline
(8, 6, 0, 0) & 2& (8, 8, 8, 0) & 2& (9, 7, 6, 2) & 3& (9, 8, 8, 2) & 2& (9, 9, 9, 4) & 1 \\
\hline
(8, 6, 2, 0) & 1& (8, 8, 8, 2) & 2& (9, 7, 6, 3) & 2& (9, 8, 8, 3) & 2& (9, 9, 9, 6) & 1 \\
\hline
(8, 6, 4, 0) & 2& (8, 8, 8, 4) & 2& (9, 7, 6, 4) & 2& (9, 8, 8, 4) & 2& (9, 9, 9, 9) & 1 \\
\hline
(8, 6, 4, 2) & 1& (8, 8, 8, 6) & 2& (9, 7, 6, 6) & 1& (9, 8, 8, 5) & 2& (10, 0, 0, 0) & 2 \\
\hline
(8, 6, 4, 4) & 1& (8, 8, 8, 8) & 2& (9, 7, 7, 0) & 1& (9, 8, 8, 6) & 2& (10, 2, 0, 0) & 1 \\
\hline
(8, 6, 6, 0) & 2& (9, 1, 0, 0) & 1& (9, 7, 7, 2) & 2& (9, 8, 8, 7) & 1& (10, 4, 0, 0) & 2 \\
\hline
\end{tabular}

%% file: C_dimG2.tex
\begin{tabular}{|c|c||c|c||c|c||c|c||c|c|}
\hline  $ (w,v) $ & $ m(w,v)$  & $ (w,v) $ & $ m(w,v) $  & $ (w,v) $ & $ m(w,v) $  & $ (w,v) $ & $ m(w,v) $  & $ (w,v) $ & $ m(w,v) $ \\ 
\hline
(4, 2) & 1& (34, 2) & 2& (26, 16) & 4& (30, 18) & 17& (42, 12) & 54 \\
\hline
(16, 2) & 1& (32, 4) & 3& (24, 18) & 6& (28, 20) & 15& (40, 14) & 60 \\
\hline
(20, 2) & 1& (30, 6) & 4& (22, 20) & 3& (26, 22) & 3& (38, 16) & 45 \\
\hline
(12, 10) & 1& (28, 8) & 6& (42, 2) & 5& (48, 2) & 14& (36, 18) & 47 \\
\hline
(16, 8) & 1& (26, 10) & 3& (40, 4) & 11& (46, 4) & 22& (34, 20) & 38 \\
\hline
(24, 2) & 1& (24, 12) & 4& (38, 6) & 13& (44, 6) & 31& (32, 22) & 24 \\
\hline
(22, 4) & 1& (22, 14) & 3& (36, 8) & 15& (42, 8) & 31& (30, 24) & 15 \\
\hline
(20, 6) & 1& (36, 2) & 4& (34, 10) & 16& (40, 10) & 37& (28, 26) & 13 \\
\hline
(16, 10) & 1& (34, 4) & 4& (32, 12) & 17& (38, 12) & 37& (54, 2) & 20 \\
\hline
(24, 4) & 1& (32, 6) & 8& (30, 14) & 12& (36, 14) & 32& (52, 4) & 39 \\
\hline
(22, 6) & 1& (30, 8) & 6& (28, 16) & 11& (34, 16) & 28& (50, 6) & 51 \\
\hline
(20, 8) & 1& (28, 10) & 8& (26, 18) & 9& (32, 18) & 29& (48, 8) & 60 \\
\hline
(18, 10) & 1& (26, 12) & 6& (24, 20) & 2& (30, 20) & 15& (46, 10) & 66 \\
\hline
(28, 2) & 3& (24, 14) & 4& (44, 2) & 10& (28, 22) & 12& (44, 12) & 72 \\
\hline
(24, 6) & 2& (22, 16) & 2& (42, 4) & 14& (26, 24) & 5& (42, 14) & 64 \\
\hline
(22, 8) & 2& (20, 18) & 4& (40, 6) & 18& (50, 2) & 13& (40, 16) & 64 \\
\hline
(20, 10) & 1& (38, 2) & 3& (38, 8) & 20& (48, 4) & 27& (38, 18) & 60 \\
\hline
(16, 14) & 2& (36, 4) & 7& (36, 10) & 25& (46, 6) & 33& (36, 20) & 45 \\
\hline
(30, 2) & 1& (34, 6) & 7& (34, 12) & 17& (44, 8) & 41& (34, 22) & 37 \\
\hline
(28, 4) & 2& (32, 8) & 9& (32, 14) & 20& (42, 10) & 44& (32, 24) & 30 \\
\hline
(26, 6) & 2& (30, 10) & 9& (30, 16) & 17& (40, 12) & 42& (30, 26) & 10 \\
\hline
(24, 8) & 2& (28, 12) & 7& (28, 18) & 11& (38, 14) & 41& (56, 2) & 29 \\
\hline
(22, 10) & 2& (26, 14) & 6& (26, 20) & 6& (36, 16) & 41& (54, 4) & 48 \\
\hline
(20, 12) & 2& (24, 16) & 6& (24, 22) & 6& (34, 18) & 30& (52, 6) & 63 \\
\hline
(32, 2) & 3& (22, 18) & 2& (46, 2) & 9& (32, 20) & 26& (50, 8) & 74 \\
\hline
(30, 4) & 3& (40, 2) & 8& (44, 4) & 16& (30, 22) & 20& (48, 10) & 88 \\
\hline
(28, 6) & 3& (38, 4) & 8& (42, 6) & 21& (28, 24) & 6& (46, 12) & 82 \\
\hline
(26, 8) & 3& (36, 6) & 12& (40, 8) & 28& (52, 2) & 23& (44, 14) & 87 \\
\hline
(24, 10) & 5& (34, 8) & 13& (38, 10) & 25& (50, 4) & 29& (42, 16) & 83 \\
\hline
(22, 12) & 2& (32, 10) & 12& (36, 12) & 27& (48, 6) & 45& (40, 18) & 72 \\
\hline
(20, 14) & 2& (30, 12) & 11& (34, 14) & 26& (46, 8) & 52& (38, 20) & 63 \\
\hline
(18, 16) & 1& (28, 14) & 13& (32, 16) & 19& (44, 10) & 54& (36, 22) & 58 \\
\hline
\end{tabular}

%% file: C_dim_StSp4.tex
\begin{tabular}{|c|c||c|c||c|c||c|c||c|c|}
\hline  $\underline{w}$ & $S(\underline{w})$  & $\underline{w}$ & $S(\underline{w})$  & $\underline{w}$ & $S(\underline{w})$  & $\underline{w}$ & $S(\underline{w})$  & $\underline{w}$ & $S(\underline{w})$ \\ 
\hline
(19, 7) & 1& (29, 21) & 2& (35, 13) & 5& (39, 15) & 10& (43, 5) & 3 \\
\hline
(21, 5) & 1& (29, 25) & 1& (35, 15) & 6& (39, 17) & 8& (43, 7) & 9 \\
\hline
(21, 9) & 1& (31, 3) & 2& (35, 17) & 5& (39, 19) & 11& (43, 9) & 7 \\
\hline
(21, 13) & 1& (31, 5) & 1& (35, 19) & 7& (39, 21) & 10& (43, 11) & 11 \\
\hline
(23, 7) & 1& (31, 7) & 3& (35, 21) & 6& (39, 23) & 10& (43, 13) & 11 \\
\hline
(23, 9) & 1& (31, 9) & 2& (35, 23) & 5& (39, 25) & 10& (43, 15) & 15 \\
\hline
(23, 13) & 1& (31, 11) & 3& (35, 25) & 5& (39, 27) & 9& (43, 17) & 13 \\
\hline
(25, 5) & 1& (31, 13) & 4& (35, 27) & 3& (39, 29) & 7& (43, 19) & 17 \\
\hline
(25, 7) & 1& (31, 15) & 4& (35, 29) & 2& (39, 31) & 6& (43, 21) & 14 \\
\hline
(25, 9) & 2& (31, 17) & 3& (35, 31) & 1& (39, 33) & 4& (43, 23) & 16 \\
\hline
(25, 11) & 1& (31, 19) & 4& (37, 1) & 1& (39, 35) & 1& (43, 25) & 16 \\
\hline
(25, 13) & 2& (31, 21) & 3& (37, 5) & 4& (39, 37) & 1& (43, 27) & 16 \\
\hline
(25, 15) & 1& (31, 23) & 2& (37, 7) & 3& (41, 1) & 1& (43, 29) & 14 \\
\hline
(25, 17) & 1& (31, 25) & 2& (37, 9) & 7& (41, 3) & 1& (43, 31) & 14 \\
\hline
(25, 19) & 1& (33, 5) & 3& (37, 11) & 5& (41, 5) & 6& (43, 33) & 11 \\
\hline
(27, 3) & 1& (33, 7) & 2& (37, 13) & 9& (41, 7) & 4& (43, 35) & 8 \\
\hline
(27, 7) & 2& (33, 9) & 5& (37, 15) & 6& (41, 9) & 9& (43, 37) & 7 \\
\hline
(27, 9) & 1& (33, 11) & 2& (37, 17) & 9& (41, 11) & 6& (43, 39) & 3 \\
\hline
(27, 11) & 2& (33, 13) & 6& (37, 19) & 8& (41, 13) & 13& (45, 1) & 2 \\
\hline
(27, 13) & 2& (33, 15) & 4& (37, 21) & 10& (41, 15) & 10& (45, 3) & 1 \\
\hline
(27, 15) & 2& (33, 17) & 6& (37, 23) & 7& (41, 17) & 13& (45, 5) & 8 \\
\hline
(27, 17) & 1& (33, 19) & 5& (37, 25) & 9& (41, 19) & 11& (45, 7) & 6 \\
\hline
(27, 19) & 1& (33, 21) & 5& (37, 27) & 6& (41, 21) & 14& (45, 9) & 13 \\
\hline
(27, 21) & 1& (33, 23) & 3& (37, 29) & 5& (41, 23) & 11& (45, 11) & 9 \\
\hline
(29, 5) & 2& (33, 25) & 4& (37, 31) & 4& (41, 25) & 15& (45, 13) & 17 \\
\hline
(29, 7) & 1& (33, 27) & 2& (37, 33) & 2& (41, 27) & 11& (45, 15) & 13 \\
\hline
(29, 9) & 3& (33, 29) & 1& (39, 3) & 3& (41, 29) & 11& (45, 17) & 19 \\
\hline
(29, 11) & 1& (35, 3) & 2& (39, 5) & 2& (41, 31) & 9& (45, 19) & 17 \\
\hline
(29, 13) & 4& (35, 5) & 1& (39, 7) & 7& (41, 33) & 8& (45, 21) & 21 \\
\hline
(29, 15) & 2& (35, 7) & 5& (39, 9) & 5& (41, 35) & 4& (45, 23) & 16 \\
\hline
(29, 17) & 3& (35, 9) & 4& (39, 11) & 8& (41, 37) & 3& (45, 25) & 22 \\
\hline
(29, 19) & 2& (35, 11) & 5& (39, 13) & 8& (43, 3) & 5& (45, 27) & 18 \\
\hline
\end{tabular}

%% file: C_dim_StSO7.tex
\begin{tabular}{|c|c||c|c||c|c||c|c|}
\hline  $\underline{w}$ & $S(\underline{w})$  & $\underline{w}$ & $S(\underline{w})$  & $\underline{w}$ & $S(\underline{w})$  & $\underline{w}$ & $S(\underline{w})$ \\ 
\hline
(23, 13, 5) & 1& (27, 17, 11) & 1& (29, 11, 5) & 1& (29, 21, 19) & 1 \\
\hline
(23, 15, 3) & 1& (27, 17, 13) & 1& (29, 13, 3) & 1& (29, 23, 1) & 1 \\
\hline
(23, 15, 7) & 1& (27, 19, 3) & 2& (29, 13, 5) & 1& (29, 23, 3) & 2 \\
\hline
(23, 17, 5) & 1& (27, 19, 5) & 2& (29, 13, 7) & 3& (29, 23, 5) & 5 \\
\hline
(23, 17, 9) & 1& (27, 19, 7) & 3& (29, 13, 9) & 1& (29, 23, 7) & 5 \\
\hline
(23, 19, 3) & 1& (27, 19, 9) & 3& (29, 15, 1) & 1& (29, 23, 9) & 6 \\
\hline
(23, 19, 11) & 1& (27, 19, 11) & 3& (29, 15, 5) & 3& (29, 23, 11) & 7 \\
\hline
(25, 13, 3) & 1& (27, 19, 13) & 2& (29, 15, 7) & 2& (29, 23, 13) & 5 \\
\hline
(25, 13, 7) & 1& (27, 19, 15) & 1& (29, 15, 9) & 3& (29, 23, 15) & 5 \\
\hline
(25, 15, 5) & 1& (27, 21, 1) & 1& (29, 15, 13) & 1& (29, 23, 17) & 3 \\
\hline
(25, 15, 9) & 1& (27, 21, 5) & 4& (29, 17, 3) & 3& (29, 23, 19) & 1 \\
\hline
(25, 17, 3) & 2& (27, 21, 7) & 2& (29, 17, 5) & 1& (29, 25, 3) & 3 \\
\hline
(25, 17, 7) & 2& (27, 21, 9) & 4& (29, 17, 7) & 6& (29, 25, 5) & 3 \\
\hline
(25, 17, 11) & 1& (27, 21, 11) & 2& (29, 17, 9) & 3& (29, 25, 7) & 7 \\
\hline
(25, 19, 1) & 1& (27, 21, 13) & 3& (29, 17, 11) & 3& (29, 25, 9) & 4 \\
\hline
(25, 19, 5) & 2& (27, 21, 15) & 1& (29, 17, 13) & 1& (29, 25, 11) & 7 \\
\hline
(25, 19, 9) & 2& (27, 21, 17) & 1& (29, 19, 1) & 1& (29, 25, 13) & 4 \\
\hline
(25, 19, 13) & 1& (27, 23, 3) & 1& (29, 19, 3) & 1& (29, 25, 15) & 5 \\
\hline
(25, 21, 3) & 2& (27, 23, 5) & 3& (29, 19, 5) & 6& (29, 25, 17) & 3 \\
\hline
(25, 21, 7) & 2& (27, 23, 7) & 1& (29, 19, 7) & 3& (29, 25, 19) & 2 \\
\hline
(25, 21, 11) & 2& (27, 23, 9) & 2& (29, 19, 9) & 7& (29, 25, 21) & 1 \\
\hline
(25, 21, 15) & 1& (27, 23, 11) & 2& (29, 19, 11) & 4& (29, 27, 1) & 1 \\
\hline
(27, 9, 5) & 1& (27, 23, 13) & 1& (29, 19, 13) & 5& (29, 27, 5) & 1 \\
\hline
(27, 13, 5) & 2& (27, 23, 15) & 1& (29, 19, 15) & 1& (29, 27, 7) & 2 \\
\hline
(27, 13, 7) & 1& (27, 23, 17) & 1& (29, 19, 17) & 1& (29, 27, 9) & 3 \\
\hline
(27, 13, 9) & 1& (27, 25, 5) & 2& (29, 21, 3) & 5& (29, 27, 11) & 1 \\
\hline
(27, 15, 3) & 1& (27, 25, 7) & 1& (29, 21, 5) & 1& (29, 27, 13) & 2 \\
\hline
(27, 15, 5) & 1& (27, 25, 9) & 1& (29, 21, 7) & 10& (29, 27, 15) & 1 \\
\hline
(27, 15, 7) & 2& (27, 25, 11) & 1& (29, 21, 9) & 4& (29, 27, 17) & 1 \\
\hline
(27, 15, 9) & 1& (27, 25, 13) & 1& (29, 21, 11) & 8& (29, 27, 19) & 1 \\
\hline
(27, 17, 5) & 4& (27, 25, 15) & 1& (29, 21, 13) & 4& (31, 9, 5) & 1 \\
\hline
(27, 17, 7) & 1& (27, 25, 17) & 1& (29, 21, 15) & 5& (31, 11, 3) & 1 \\
\hline
(27, 17, 9) & 3& (29, 9, 7) & 1& (29, 21, 17) & 1& (31, 11, 7) & 1 \\
\hline
\end{tabular}

%% file: C_dim_StSO9.tex
\begin{tabular}{|c|c||c|c||c|c||c|c|}
\hline  $\underline{w}$ & $S(\underline{w})$  & $\underline{w}$ & $S(\underline{w})$  & $\underline{w}$ & $S(\underline{w})$  & $\underline{w}$ & $S(\underline{w})$ \\ 
\hline
(25, 17, 9, 5) & 1& (27, 17, 13, 7) & 2& (27, 21, 19, 7) & 1& (27, 23, 21, 9) & 1 \\
\hline
(25, 17, 13, 5) & 1& (27, 19, 9, 5) & 1& (27, 21, 19, 9) & 1& (27, 25, 9, 3) & 2 \\
\hline
(25, 19, 9, 3) & 1& (27, 19, 11, 3) & 2& (27, 21, 19, 11) & 1& (27, 25, 11, 1) & 1 \\
\hline
(25, 19, 11, 5) & 1& (27, 19, 11, 5) & 1& (27, 23, 7, 3) & 2& (27, 25, 11, 3) & 1 \\
\hline
(25, 19, 13, 3) & 1& (27, 19, 13, 1) & 1& (27, 23, 9, 1) & 1& (27, 25, 11, 5) & 2 \\
\hline
(25, 19, 13, 5) & 1& (27, 19, 13, 3) & 1& (27, 23, 9, 5) & 2& (27, 25, 13, 3) & 5 \\
\hline
(25, 19, 13, 7) & 1& (27, 19, 13, 5) & 4& (27, 23, 11, 3) & 5& (27, 25, 13, 5) & 1 \\
\hline
(25, 19, 13, 9) & 1& (27, 19, 13, 7) & 1& (27, 23, 11, 5) & 1& (27, 25, 13, 7) & 4 \\
\hline
(25, 19, 15, 5) & 1& (27, 19, 13, 9) & 3& (27, 23, 11, 7) & 4& (27, 25, 13, 9) & 1 \\
\hline
(25, 21, 11, 7) & 1& (27, 19, 15, 3) & 2& (27, 23, 13, 1) & 4& (27, 25, 15, 1) & 3 \\
\hline
(25, 21, 13, 5) & 1& (27, 19, 15, 5) & 1& (27, 23, 13, 3) & 1& (27, 25, 15, 3) & 2 \\
\hline
(25, 21, 13, 7) & 1& (27, 19, 15, 7) & 1& (27, 23, 13, 5) & 6& (27, 25, 15, 5) & 5 \\
\hline
(25, 21, 15, 3) & 1& (27, 19, 15, 9) & 1& (27, 23, 13, 7) & 3& (27, 25, 15, 7) & 3 \\
\hline
(25, 21, 15, 5) & 1& (27, 19, 17, 5) & 1& (27, 23, 13, 9) & 6& (27, 25, 15, 9) & 5 \\
\hline
(25, 21, 15, 7) & 2& (27, 19, 17, 9) & 1& (27, 23, 15, 3) & 7& (27, 25, 15, 11) & 1 \\
\hline
(25, 21, 15, 9) & 1& (27, 21, 9, 3) & 2& (27, 23, 15, 5) & 3& (27, 25, 17, 3) & 7 \\
\hline
(25, 21, 17, 5) & 1& (27, 21, 9, 7) & 1& (27, 23, 15, 7) & 7& (27, 25, 17, 5) & 2 \\
\hline
(25, 21, 17, 7) & 1& (27, 21, 11, 3) & 1& (27, 23, 15, 9) & 4& (27, 25, 17, 7) & 7 \\
\hline
(25, 21, 17, 9) & 1& (27, 21, 11, 5) & 2& (27, 23, 15, 11) & 5& (27, 25, 17, 9) & 4 \\
\hline
(25, 23, 9, 3) & 1& (27, 21, 11, 7) & 2& (27, 23, 15, 13) & 1& (27, 25, 17, 11) & 5 \\
\hline
(25, 23, 11, 1) & 1& (27, 21, 13, 3) & 5& (27, 23, 17, 1) & 5& (27, 25, 17, 13) & 1 \\
\hline
(25, 23, 11, 5) & 2& (27, 21, 13, 5) & 2& (27, 23, 17, 3) & 2& (27, 25, 19, 1) & 3 \\
\hline
(25, 23, 13, 3) & 1& (27, 21, 13, 7) & 6& (27, 23, 17, 5) & 6& (27, 25, 19, 3) & 2 \\
\hline
(25, 23, 13, 7) & 1& (27, 21, 13, 9) & 2& (27, 23, 17, 7) & 5& (27, 25, 19, 5) & 5 \\
\hline
(25, 23, 15, 1) & 1& (27, 21, 15, 1) & 1& (27, 23, 17, 9) & 7& (27, 25, 19, 7) & 3 \\
\hline
(25, 23, 15, 5) & 3& (27, 21, 15, 3) & 2& (27, 23, 17, 11) & 3& (27, 25, 19, 9) & 6 \\
\hline
(25, 23, 15, 9) & 1& (27, 21, 15, 5) & 4& (27, 23, 17, 13) & 4& (27, 25, 19, 11) & 3 \\
\hline
(25, 23, 15, 11) & 1& (27, 21, 15, 7) & 4& (27, 23, 19, 3) & 5& (27, 25, 19, 13) & 3 \\
\hline
(25, 23, 17, 3) & 1& (27, 21, 15, 9) & 4& (27, 23, 19, 5) & 1& (27, 25, 21, 3) & 4 \\
\hline
(25, 23, 17, 5) & 1& (27, 21, 15, 11) & 2& (27, 23, 19, 7) & 6& (27, 25, 21, 7) & 4 \\
\hline
(25, 23, 17, 7) & 1& (27, 21, 17, 3) & 5& (27, 23, 19, 9) & 2& (27, 25, 21, 9) & 2 \\
\hline
(25, 23, 17, 11) & 1& (27, 21, 17, 7) & 6& (27, 23, 19, 11) & 3& (27, 25, 21, 11) & 3 \\
\hline
(25, 23, 19, 5) & 1& (27, 21, 17, 9) & 2& (27, 23, 19, 13) & 1& (27, 25, 21, 13) & 1 \\
\hline
(27, 17, 9, 3) & 1& (27, 21, 17, 11) & 3& (27, 23, 19, 15) & 1& (27, 25, 21, 15) & 1 \\
\hline
(27, 17, 9, 7) & 1& (27, 21, 19, 3) & 1& (27, 23, 21, 1) & 1& (27, 25, 23, 3) & 1 \\
\hline
(27, 17, 13, 3) & 2& (27, 21, 19, 5) & 1& (27, 23, 21, 5) & 1& (27, 25, 23, 9) & 1 \\
\hline
\end{tabular}

%% file: C_dim_StSO8_nonnul.tex
\begin{tabular}{|c|c||c|c||c|c||c|c|}
\hline  $\underline{w}$ & $O(\underline{w})$  & $\underline{w}$ & $O(\underline{w})$  & $\underline{w}$ & $O(\underline{w})$  & $\underline{w}$ & $O(\underline{w})$ \\ 
\hline
(24, 18, 10, 4) & 1& (28, 24, 14, 2) & 2& (30, 22, 14, 2) & 2& (30, 26, 16, 12) & 1 \\
\hline
(24, 20, 14, 2) & 1& (28, 24, 14, 10) & 1& (30, 22, 14, 6) & 3& (30, 26, 18, 2) & 3 \\
\hline
(26, 18, 10, 2) & 1& (28, 24, 16, 4) & 1& (30, 22, 16, 4) & 2& (30, 26, 18, 6) & 2 \\
\hline
(26, 18, 14, 6) & 1& (28, 24, 16, 12) & 1& (30, 22, 16, 8) & 1& (30, 26, 18, 10) & 1 \\
\hline
(26, 20, 10, 4) & 1& (28, 24, 18, 2) & 1& (30, 22, 18, 2) & 1& (30, 26, 18, 14) & 1 \\
\hline
(26, 20, 14, 8) & 1& (28, 24, 18, 6) & 1& (30, 22, 18, 6) & 1& (30, 26, 20, 4) & 3 \\
\hline
(26, 22, 10, 6) & 1& (28, 24, 20, 4) & 1& (30, 22, 18, 10) & 1& (30, 26, 20, 8) & 1 \\
\hline
(26, 22, 14, 2) & 1& (28, 24, 20, 8) & 1& (30, 24, 8, 2) & 1& (30, 26, 22, 2) & 1 \\
\hline
(26, 24, 14, 4) & 1& (28, 26, 12, 2) & 1& (30, 24, 10, 4) & 3& (30, 26, 22, 6) & 2 \\
\hline
(26, 24, 16, 2) & 1& (28, 26, 14, 4) & 1& (30, 24, 12, 2) & 2& (30, 26, 22, 10) & 1 \\
\hline
(26, 24, 18, 8) & 1& (28, 26, 16, 2) & 2& (30, 24, 12, 6) & 2& (30, 28, 10, 4) & 1 \\
\hline
(26, 24, 20, 6) & 1& (28, 26, 18, 8) & 1& (30, 24, 14, 4) & 2& (30, 28, 10, 8) & 1 \\
\hline
(28, 16, 10, 6) & 1& (28, 26, 20, 6) & 1& (30, 24, 14, 8) & 3& (30, 28, 12, 2) & 1 \\
\hline
(28, 18, 8, 2) & 1& (28, 26, 22, 4) & 1& (30, 24, 16, 2) & 3& (30, 28, 14, 4) & 3 \\
\hline
(28, 18, 12, 2) & 1& (30, 14, 8, 4) & 1& (30, 24, 16, 6) & 2& (30, 28, 14, 12) & 1 \\
\hline
(28, 18, 14, 4) & 1& (30, 16, 10, 4) & 1& (30, 24, 16, 10) & 2& (30, 28, 16, 2) & 2 \\
\hline
(28, 20, 10, 2) & 1& (30, 18, 8, 4) & 1& (30, 24, 18, 4) & 4& (30, 28, 16, 6) & 1 \\
\hline
(28, 20, 12, 4) & 1& (30, 18, 10, 2) & 1& (30, 24, 18, 8) & 1& (30, 28, 18, 4) & 1 \\
\hline
(28, 20, 14, 2) & 1& (30, 18, 10, 6) & 1& (30, 24, 18, 12) & 2& (30, 28, 18, 8) & 2 \\
\hline
(28, 20, 14, 6) & 1& (30, 18, 12, 4) & 1& (30, 24, 20, 2) & 2& (30, 28, 20, 2) & 1 \\
\hline
(28, 20, 16, 4) & 1& (30, 18, 14, 2) & 1& (30, 24, 20, 6) & 2& (30, 28, 20, 6) & 3 \\
\hline
(28, 20, 16, 8) & 1& (30, 18, 14, 6) & 1& (30, 24, 20, 10) & 1& (30, 28, 20, 10) & 1 \\
\hline
(28, 22, 8, 2) & 1& (30, 20, 6, 4) & 1& (30, 24, 20, 14) & 1& (30, 28, 22, 4) & 4 \\
\hline
(28, 22, 10, 4) & 1& (30, 20, 10, 4) & 1& (30, 24, 22, 8) & 1& (30, 28, 22, 8) & 1 \\
\hline
(28, 22, 12, 2) & 1& (30, 20, 10, 8) & 1& (30, 24, 22, 16) & 1& (30, 28, 22, 12) & 1 \\
\hline
(28, 22, 12, 6) & 1& (30, 20, 12, 2) & 1& (30, 26, 6, 2) & 1& (32, 16, 8, 4) & 1 \\
\hline
(28, 22, 14, 8) & 1& (30, 20, 14, 4) & 3& (30, 26, 8, 4) & 1& (32, 16, 10, 2) & 1 \\
\hline
(28, 22, 16, 2) & 1& (30, 20, 14, 8) & 1& (30, 26, 10, 2) & 1& (32, 16, 10, 6) & 1 \\
\hline
(28, 22, 16, 6) & 1& (30, 20, 14, 12) & 1& (30, 26, 10, 6) & 2& (32, 18, 8, 2) & 1 \\
\hline
(28, 22, 16, 10) & 1& (30, 20, 16, 2) & 1& (30, 26, 12, 4) & 2& (32, 18, 8, 6) & 1 \\
\hline
(28, 22, 18, 4) & 1& (30, 20, 16, 6) & 1& (30, 26, 12, 8) & 1& (32, 18, 10, 4) & 2 \\
\hline
(28, 24, 8, 4) & 1& (30, 20, 18, 8) & 1& (30, 26, 14, 2) & 3& (32, 18, 10, 8) & 1 \\
\hline
(28, 24, 10, 2) & 1& (30, 22, 8, 4) & 1& (30, 26, 14, 6) & 1& (32, 18, 12, 2) & 1 \\
\hline
(28, 24, 10, 6) & 1& (30, 22, 10, 2) & 2& (30, 26, 14, 10) & 2& (32, 18, 12, 6) & 1 \\
\hline
(28, 24, 12, 4) & 1& (30, 22, 10, 6) & 1& (30, 26, 16, 4) & 2& (32, 18, 14, 4) & 1 \\
\hline
(28, 24, 12, 8) & 1& (30, 22, 12, 4) & 1& (30, 26, 16, 8) & 1& (32, 18, 14, 8) & 1 \\
\hline
\end{tabular}

%% file: C_dim_StSO8_nul.tex
\begin{tabular}{|c|c||c|c||c|c||c|c|}
\hline  $\underline{w}$ & $O(\underline{w})$  & $\underline{w}$ & $O(\underline{w})$  & $\underline{w}$ & $O(\underline{w})$  & $\underline{w}$ & $O(\underline{w})$ \\ 
\hline
(24, 16, 8, 0) & 1& (30, 28, 10, 0) & 2& (32, 30, 10, 0) & 2& (34, 28, 26, 0) & 2 \\
\hline
(26, 16, 10, 0) & 1& (30, 28, 14, 0) & 3& (32, 30, 14, 0) & 4& (34, 30, 4, 0) & 2 \\
\hline
(26, 20, 6, 0) & 1& (30, 28, 18, 0) & 5& (32, 30, 18, 0) & 6& (34, 30, 8, 0) & 2 \\
\hline
(26, 20, 10, 0) & 1& (30, 28, 26, 0) & 1& (32, 30, 26, 0) & 6& (34, 30, 12, 0) & 7 \\
\hline
(26, 20, 14, 0) & 1& (32, 12, 8, 0) & 1& (34, 12, 6, 0) & 1& (34, 30, 16, 0) & 14 \\
\hline
(26, 24, 10, 0) & 1& (32, 14, 10, 0) & 1& (34, 14, 8, 0) & 1& (34, 30, 20, 0) & 6 \\
\hline
(26, 24, 14, 0) & 1& (32, 16, 4, 0) & 1& (34, 16, 6, 0) & 1& (34, 30, 24, 0) & 7 \\
\hline
(26, 24, 18, 0) & 1& (32, 16, 8, 0) & 1& (34, 16, 10, 0) & 3& (34, 32, 2, 0) & 1 \\
\hline
(28, 14, 6, 0) & 1& (32, 16, 12, 0) & 1& (34, 16, 14, 0) & 1& (34, 32, 6, 0) & 2 \\
\hline
(28, 16, 8, 0) & 1& (32, 18, 6, 0) & 1& (34, 18, 4, 0) & 1& (34, 32, 10, 0) & 6 \\
\hline
(28, 18, 10, 0) & 1& (32, 18, 10, 0) & 1& (34, 18, 8, 0) & 1& (34, 32, 14, 0) & 8 \\
\hline
(28, 20, 8, 0) & 1& (32, 18, 14, 0) & 3& (34, 18, 12, 0) & 3& (34, 32, 18, 0) & 13 \\
\hline
(28, 20, 12, 0) & 1& (32, 20, 4, 0) & 1& (34, 20, 6, 0) & 3& (34, 32, 22, 0) & 3 \\
\hline
(28, 22, 14, 0) & 2& (32, 20, 8, 0) & 2& (34, 20, 10, 0) & 3& (34, 32, 26, 0) & 14 \\
\hline
(28, 24, 4, 0) & 1& (32, 20, 12, 0) & 2& (34, 20, 14, 0) & 8& (36, 12, 8, 0) & 1 \\
\hline
(28, 24, 12, 0) & 1& (32, 20, 16, 0) & 3& (34, 20, 18, 0) & 2& (36, 14, 6, 0) & 1 \\
\hline
(28, 24, 16, 0) & 3& (32, 22, 6, 0) & 1& (34, 22, 4, 0) & 1& (36, 14, 10, 0) & 1 \\
\hline
(28, 26, 18, 0) & 2& (32, 22, 10, 0) & 4& (34, 22, 8, 0) & 3& (36, 16, 4, 0) & 1 \\
\hline
(30, 16, 6, 0) & 1& (32, 22, 14, 0) & 1& (34, 22, 12, 0) & 3& (36, 16, 8, 0) & 3 \\
\hline
(30, 16, 10, 0) & 1& (32, 22, 18, 0) & 3& (34, 22, 16, 0) & 5& (36, 16, 12, 0) & 2 \\
\hline
(30, 16, 14, 0) & 1& (32, 24, 4, 0) & 1& (34, 24, 6, 0) & 3& (36, 18, 6, 0) & 2 \\
\hline
(30, 18, 8, 0) & 1& (32, 24, 8, 0) & 5& (34, 24, 10, 0) & 11& (36, 18, 10, 0) & 5 \\
\hline
(30, 20, 6, 0) & 1& (32, 24, 12, 0) & 5& (34, 24, 14, 0) & 7& (36, 18, 14, 0) & 4 \\
\hline
(30, 20, 10, 0) & 4& (32, 24, 16, 0) & 4& (34, 24, 18, 0) & 12& (36, 20, 4, 0) & 2 \\
\hline
(30, 20, 14, 0) & 1& (32, 24, 20, 0) & 6& (34, 24, 22, 0) & 2& (36, 20, 8, 0) & 4 \\
\hline
(30, 20, 18, 0) & 1& (32, 26, 6, 0) & 2& (34, 26, 4, 0) & 1& (36, 20, 12, 0) & 6 \\
\hline
(30, 22, 8, 0) & 1& (32, 26, 10, 0) & 4& (34, 26, 8, 0) & 6& (36, 20, 16, 0) & 6 \\
\hline
(30, 22, 12, 0) & 1& (32, 26, 14, 0) & 8& (34, 26, 12, 0) & 9& (36, 22, 6, 0) & 3 \\
\hline
(30, 24, 6, 0) & 2& (32, 26, 18, 0) & 3& (34, 26, 16, 0) & 7& (36, 22, 10, 0) & 6 \\
\hline
(30, 24, 10, 0) & 2& (32, 26, 22, 0) & 5& (34, 26, 20, 0) & 12& (36, 22, 14, 0) & 10 \\
\hline
(30, 24, 14, 0) & 5& (32, 28, 4, 0) & 2& (34, 26, 24, 0) & 2& (36, 22, 18, 0) & 6 \\
\hline
(30, 24, 18, 0) & 2& (32, 28, 8, 0) & 2& (34, 28, 6, 0) & 6& (36, 24, 4, 0) & 3 \\
\hline
(30, 26, 8, 0) & 1& (32, 28, 12, 0) & 5& (34, 28, 10, 0) & 8& (36, 24, 8, 0) & 8 \\
\hline
(30, 26, 12, 0) & 1& (32, 28, 16, 0) & 9& (34, 28, 14, 0) & 16& (36, 24, 12, 0) & 13 \\
\hline
(30, 26, 16, 0) & 3& (32, 28, 20, 0) & 4& (34, 28, 18, 0) & 11& (36, 24, 16, 0) & 16 \\
\hline
(30, 28, 2, 0) & 1& (32, 28, 24, 0) & 3& (34, 28, 22, 0) & 14& (36, 24, 20, 0) & 12 \\
\hline
\end{tabular}

%% file: C_dimStG2.tex
\begin{tabular}{|c|c||c|c||c|c||c|c||c|c|}
\hline  $ (w,v) $ & $ G_2(\ww)$  & $ (w,v) $ & $ G_2(\ww) $  & $ (w,v) $ & $ G_2(\ww) $  & $ (w,v) $ & $ G_2(\ww) $  & $ (w,v) $ & $ G_2(\ww) $ \\ 
\hline
(16, 8) & 1& (30, 8) & 4& (44, 2) & 7& (28, 22) & 12& (44, 12) & 72 \\
\hline
(20, 6) & 1& (28, 10) & 8& (42, 4) & 13& (26, 24) & 4& (42, 14) & 61 \\
\hline
(16, 10) & 1& (26, 12) & 6& (40, 6) & 18& (50, 2) & 11& (40, 16) & 64 \\
\hline
(24, 4) & 1& (24, 14) & 4& (38, 8) & 18& (48, 4) & 27& (38, 18) & 58 \\
\hline
(20, 8) & 1& (20, 18) & 3& (36, 10) & 25& (46, 6) & 29& (36, 20) & 45 \\
\hline
(18, 10) & 1& (38, 2) & 2& (34, 12) & 15& (44, 8) & 41& (34, 22) & 34 \\
\hline
(28, 2) & 1& (36, 4) & 7& (32, 14) & 20& (42, 10) & 42& (32, 24) & 30 \\
\hline
(24, 6) & 2& (34, 6) & 5& (30, 16) & 15& (40, 12) & 42& (30, 26) & 7 \\
\hline
(22, 8) & 1& (32, 8) & 9& (28, 18) & 11& (38, 14) & 39& (56, 2) & 25 \\
\hline
(20, 10) & 1& (30, 10) & 8& (26, 20) & 6& (36, 16) & 41& (54, 4) & 44 \\
\hline
(16, 14) & 1& (28, 12) & 7& (24, 22) & 4& (34, 18) & 27& (52, 6) & 63 \\
\hline
(28, 4) & 2& (26, 14) & 6& (46, 2) & 7& (32, 20) & 26& (50, 8) & 72 \\
\hline
(26, 6) & 2& (24, 16) & 6& (44, 4) & 16& (30, 22) & 18& (48, 10) & 88 \\
\hline
(24, 8) & 2& (40, 2) & 5& (42, 6) & 19& (28, 24) & 6& (46, 12) & 76 \\
\hline
(22, 10) & 1& (38, 4) & 6& (40, 8) & 28& (52, 2) & 19& (44, 14) & 87 \\
\hline
(20, 12) & 2& (36, 6) & 12& (38, 10) & 23& (50, 4) & 27& (42, 16) & 81 \\
\hline
(32, 2) & 1& (34, 8) & 12& (36, 12) & 27& (48, 6) & 45& (40, 18) & 72 \\
\hline
(30, 4) & 2& (32, 10) & 12& (34, 14) & 24& (46, 8) & 48& (38, 20) & 60 \\
\hline
(28, 6) & 3& (30, 12) & 9& (32, 16) & 19& (44, 10) & 54& (36, 22) & 58 \\
\hline
(26, 8) & 3& (28, 14) & 13& (30, 18) & 15& (42, 12) & 52& (34, 24) & 29 \\
\hline
(24, 10) & 5& (26, 16) & 4& (28, 20) & 15& (40, 14) & 60& (32, 26) & 26 \\
\hline
(20, 14) & 2& (24, 18) & 6& (26, 22) & 3& (38, 16) & 42& (30, 28) & 6 \\
\hline
(34, 2) & 1& (42, 2) & 3& (48, 2) & 11& (36, 18) & 47& (58, 2) & 25 \\
\hline
(32, 4) & 3& (40, 4) & 11& (46, 4) & 18& (34, 20) & 36& (56, 4) & 54 \\
\hline
(30, 6) & 3& (38, 6) & 12& (44, 6) & 31& (32, 22) & 24& (54, 6) & 69 \\
\hline
(28, 8) & 6& (36, 8) & 15& (42, 8) & 29& (30, 24) & 12& (52, 8) & 93 \\
\hline
(26, 10) & 3& (34, 10) & 14& (40, 10) & 37& (28, 26) & 11& (50, 10) & 92 \\
\hline
(24, 12) & 4& (32, 12) & 17& (38, 12) & 35& (54, 2) & 16& (48, 12) & 104 \\
\hline
(22, 14) & 2& (30, 14) & 10& (36, 14) & 32& (52, 4) & 39& (46, 14) & 102 \\
\hline
(36, 2) & 2& (28, 16) & 11& (34, 16) & 26& (50, 6) & 49& (44, 16) & 96 \\
\hline
(34, 4) & 3& (26, 18) & 9& (32, 18) & 29& (48, 8) & 60& (42, 18) & 89 \\
\hline
(32, 6) & 8& (24, 20) & 2& (30, 20) & 12& (46, 10) & 62& (40, 20) & 88 \\
\hline
\end{tabular}

%% file: C_listeformesautomorphes_25.tex
$$[24]$$
$$\Delta_{15}[9]\oplus[6]$$
$$\Delta_{17}[7]\oplus[10]$$
$$\Delta_{19}[5]\oplus[14]$$
$$\Delta_{21}[3]\oplus[18]$$
$$\Delta_{23}^2\oplus[22]$$
$$\Delta_{23}^2\oplus\Delta_{11}[11]$$
$${\rm Sym}^2\Delta_{11}[2]\oplus\Delta_{11}[9]$$
$$\Delta_{19}[5]\oplus\Delta_{11}[3]\oplus[8]$$
$$\Delta_{21}[3]\oplus\Delta_{11}[7]\oplus[4]$$
$$\Delta_{21}[3]\oplus\Delta_{15}[3]\oplus[12]$$
$$\Delta_{21}[3]\oplus\Delta_{17}\oplus[16]$$
$$\Delta_{23}^2\oplus\Delta_{15}[7]\oplus[8]$$
$$\Delta_{23}^2\oplus\Delta_{17}[5]\oplus[12]$$
$$\Delta_{23}^2\oplus\Delta_{19}[3]\oplus[16]$$
$$\Delta_{23}^2\oplus\Delta_{21}\oplus[20]$$
$$\Delta_{21,9}[3]\oplus\Delta_{15}[3]\oplus[6]$$
$$\Delta_{21,13}[3]\oplus\Delta_{17}\oplus[10]$$
$$\Delta_{23,7}\oplus\Delta_{15}[7]\oplus[6]$$
$$\Delta_{21}[3]\oplus\Delta_{15}[3]\oplus\Delta_{11}\oplus[10]$$
$$\Delta_{21}[3]\oplus\Delta_{17}\oplus\Delta_{11}[5]\oplus[6]$$
$$\Delta_{21}[3]\oplus\Delta_{17}\oplus\Delta_{15}\oplus[14]$$
$$\Delta_{23}^2\oplus\Delta_{17}[5]\oplus\Delta_{11}\oplus[10]$$
$$\Delta_{23}^2\oplus\Delta_{19}[3]\oplus\Delta_{11}[5]\oplus[6]$$
$$\Delta_{23}^2\oplus\Delta_{19}[3]\oplus\Delta_{15}\oplus[14]$$
$$\Delta_{23}^2\oplus\Delta_{21}\oplus\Delta_{11}[9]\oplus[2]$$
$$\Delta_{23}^2\oplus\Delta_{21}\oplus\Delta_{15}[5]\oplus[10]$$
$$\Delta_{23}^2\oplus\Delta_{21}\oplus\Delta_{17}[3]\oplus[14]$$
$$\Delta_{23}^2\oplus\Delta_{21}\oplus\Delta_{19}\oplus[18]$$
$$\Delta_{23}^2\oplus\Delta_{21,9}\oplus\Delta_{15}[5]\oplus[8]$$
$$\Delta_{23}^2\oplus\Delta_{21,13}\oplus\Delta_{17}[3]\oplus[12]$$
$$\Delta_{23,7}\oplus\Delta_{21,9}\oplus\Delta_{15}[5]\oplus[6]$$
$$\Delta_{23,9}\oplus\Delta_{17}[5]\oplus\Delta_{11}\oplus[8]$$
$$\Delta_{23,9}\oplus\Delta_{21}\oplus\Delta_{15}[5]\oplus[8]$$
$$\Delta_{23,13}\oplus\Delta_{19}[3]\oplus\Delta_{15}\oplus[12]$$
$$\Delta_{23,13}\oplus\Delta_{21}\oplus\Delta_{17}[3]\oplus[12]$$
$$\Delta_{23, 19, 3}\oplus\Delta_{21}\oplus\Delta_{11}[7]\oplus[2]$$
$$\Delta_{23, 19, 11}\oplus\Delta_{21}\oplus\Delta_{15}[3]\oplus[10]$$
$$\Delta_{23, 19, 11}\oplus\Delta_{21,9}\oplus\Delta_{15}[3]\oplus[8]$$
$$\Delta_{23, 15, 7}\oplus\Delta_{19}[3]\oplus\Delta_{11}[3]\oplus[6]$$
$$\Delta_{21}[3]\oplus\Delta_{17}\oplus\Delta_{15}\oplus\Delta_{11}[3]\oplus[8]$$
$$\Delta_{23}^2\oplus\Delta_{19}[3]\oplus\Delta_{15}\oplus\Delta_{11}[3]\oplus[8]$$
$$\Delta_{23}^2\oplus\Delta_{21}\oplus\Delta_{17}[3]\oplus\Delta_{11}[3]\oplus[8]$$
$$\Delta_{23}^2\oplus\Delta_{21}\oplus\Delta_{19}\oplus\Delta_{11}[7]\oplus[4]$$
$$\Delta_{23}^2\oplus\Delta_{21}\oplus\Delta_{19}\oplus\Delta_{15}[3]\oplus[12]$$
$$\Delta_{23}^2\oplus\Delta_{21}\oplus\Delta_{19}\oplus\Delta_{17}\oplus[16]$$
$$\Delta_{23}^2\oplus\Delta_{21,13}\oplus\Delta_{17}[3]\oplus\Delta_{11}\oplus[10]$$
$$\Delta_{21,5}[3]\oplus\Delta_{17}\oplus\Delta_{15}\oplus\Delta_{11}[3]\oplus[2]$$
$$\Delta_{23,7}\oplus\Delta_{19}[3]\oplus\Delta_{15}\oplus\Delta_{11}[3]\oplus[6]$$
$$\Delta_{23,7}\oplus\Delta_{21}\oplus\Delta_{17}[3]\oplus\Delta_{11}[3]\oplus[6]$$
$$\Delta_{23,7}\oplus\Delta_{21,5}\oplus\Delta_{17}[3]\oplus\Delta_{11}[3]\oplus[4]$$
$$\Delta_{23,9}\oplus\Delta_{21,13}\oplus\Delta_{17}[3]\oplus\Delta_{11}\oplus[8]$$
$$\Delta_{23,13}\oplus\Delta_{19}[3]\oplus\Delta_{15}\oplus\Delta_{11}\oplus[10]$$
$$\Delta_{23,13}\oplus\Delta_{21}\oplus\Delta_{17}[3]\oplus\Delta_{11}\oplus[10]$$
$$\Delta_{23,13}\oplus\Delta_{19,7}[3]\oplus\Delta_{15}\oplus\Delta_{11}\oplus[4]$$
$$\Delta_{23,13}\oplus\Delta_{21,9}\oplus\Delta_{17}[3]\oplus\Delta_{11}\oplus[8]$$
$$\Delta_{23, 17, 5} \oplus\Delta_{21}\oplus\Delta_{19}\oplus\Delta_{11}[5]\oplus[4]$$
$$\Delta_{23, 19, 3}\oplus\Delta_{21,5}\oplus\Delta_{17}\oplus\Delta_{11}[5]\oplus[2]$$
$$\Delta_{23, 19, 11}\oplus\Delta_{21,13}\oplus\Delta_{17}\oplus\Delta_{15}\oplus[10]$$
$$\Delta_{23}^2\oplus\Delta_{21}\oplus\Delta_{19}\oplus\Delta_{15}[3]\oplus\Delta_{11}\oplus[10]$$
$$\Delta_{23}^2\oplus\Delta_{21}\oplus\Delta_{19}\oplus\Delta_{17}\oplus\Delta_{11}[5]\oplus[6]$$
$$\Delta_{23}^2\oplus\Delta_{21}\oplus\Delta_{19}\oplus\Delta_{17}\oplus\Delta_{15}\oplus[14]$$
$$\Delta_{23}^2\oplus\Delta_{21,5}\oplus\Delta_{19}\oplus\Delta_{17}\oplus\Delta_{11}[5]\oplus[4]$$
$$\Delta_{23}^2\oplus\Delta_{21,9}\oplus\Delta_{19}\oplus\Delta_{15}[3]\oplus\Delta_{11}\oplus[8]$$
$$\Delta_{23}^2\oplus\Delta_{21,9}\oplus\Delta_{19,7}\oplus\Delta_{15}[3]\oplus\Delta_{11}\oplus[6]$$
$$\Delta_{23}^2\oplus\Delta_{21,13}\oplus\Delta_{19}\oplus\Delta_{17}\oplus\Delta_{15}\oplus[12]$$
$$\Delta_{23,7}\oplus\Delta_{21,9}\oplus\Delta_{19}\oplus\Delta_{15}[3]\oplus\Delta_{11}\oplus[6]$$
$$\Delta_{23,9}\oplus\Delta_{21}\oplus\Delta_{19}\oplus\Delta_{15}[3]\oplus\Delta_{11}\oplus[8]$$
$$\Delta_{23,9}\oplus\Delta_{21}\oplus\Delta_{19,7}\oplus\Delta_{15}[3]\oplus\Delta_{11}\oplus[6]$$
$$\Delta_{23,9}\oplus\Delta_{21,5}\oplus\Delta_{19,7}\oplus\Delta_{15}[3]\oplus\Delta_{11}\oplus[4]$$
$$\Delta_{23,13}\oplus\Delta_{21}\oplus\Delta_{19}\oplus\Delta_{17}\oplus\Delta_{15}\oplus[12]$$
$$\Delta_{23, 17, 5} \oplus\Delta_{21}\oplus\Delta_{19,7}\oplus\Delta_{15}\oplus\Delta_{11}[3]\oplus[4]$$
$$\Delta_{23, 17, 9} \oplus\Delta_{21,13}\oplus\Delta_{19}\oplus\Delta_{15}\oplus\Delta_{11}\oplus[8]$$
$$\Delta_{23, 17, 9} \oplus\Delta_{21,13}\oplus\Delta_{19,7}\oplus\Delta_{15}\oplus\Delta_{11}\oplus[6]$$
$$\Delta_{23, 15, 3}\oplus\Delta_{21,5}\oplus\Delta_{19,7}\oplus\Delta_{17}\oplus\Delta_{11}[3]\oplus[2]$$
$$\Delta_{23, 15, 7}\oplus\Delta_{21}\oplus\Delta_{19}\oplus\Delta_{17}\oplus\Delta_{11}[3]\oplus[6]$$
$$\Delta_{23, 15, 7}\oplus\Delta_{21,5}\oplus\Delta_{19}\oplus\Delta_{17}\oplus\Delta_{11}[3]\oplus[4]$$
$$\Delta_{23}^2\oplus\Delta_{21}\oplus\Delta_{19}\oplus\Delta_{17}\oplus\Delta_{15}\oplus\Delta_{11}[3]\oplus[8]$$
$$\Delta_{23}^2\oplus\Delta_{21}\oplus\Delta_{19,7}\oplus\Delta_{17}\oplus\Delta_{15}\oplus\Delta_{11}[3]\oplus[6]$$
$$\Delta_{23}^2\oplus\Delta_{21,5}\oplus\Delta_{19,7}\oplus\Delta_{17}\oplus\Delta_{15}\oplus\Delta_{11}[3]\oplus[4]$$
$$\Delta_{23}^2\oplus\Delta_{21,13}\oplus\Delta_{19}\oplus\Delta_{17}\oplus\Delta_{15}\oplus\Delta_{11}\oplus[10]$$
$$\Delta_{23,7}\oplus\Delta_{21}\oplus\Delta_{19}\oplus\Delta_{17}\oplus\Delta_{15}\oplus\Delta_{11}[3]\oplus[6]$$
$$\Delta_{23,7}\oplus\Delta_{21,5}\oplus\Delta_{19}\oplus\Delta_{17}\oplus\Delta_{15}\oplus\Delta_{11}[3]\oplus[4]$$
$$\Delta_{23,9}\oplus\Delta_{21,13}\oplus\Delta_{19}\oplus\Delta_{17}\oplus\Delta_{15}\oplus\Delta_{11}\oplus[8]$$
$$\Delta_{23,9}\oplus\Delta_{21,13}\oplus\Delta_{19,7}\oplus\Delta_{17}\oplus\Delta_{15}\oplus\Delta_{11}\oplus[6]$$
$$\Delta_{23,13}\oplus\Delta_{21}\oplus\Delta_{19}\oplus\Delta_{17}\oplus\Delta_{15}\oplus\Delta_{11}\oplus[10]$$
$$\Delta_{23,13}\oplus\Delta_{21,9}\oplus\Delta_{19}\oplus\Delta_{17}\oplus\Delta_{15}\oplus\Delta_{11}\oplus[8]$$
$$\Delta_{23,13}\oplus\Delta_{21,9}\oplus\Delta_{19,7}\oplus\Delta_{17}\oplus\Delta_{15}\oplus\Delta_{11}\oplus[6]$$
$$\Delta_{23, 13, 5}\oplus\Delta_{21,9}\oplus\Delta_{19,7}\oplus\Delta_{17}\oplus\Delta_{15}\oplus\Delta_{11}\oplus[4]$$

%% file: dim_form_fin.bbl
\newpage

\begin{thebibliography}{99}

{\footnotesize

\bibitem[\textsc{Ada}11]{adams} J. Adams, \textit{Discrete series and characters of the component group}, in
\textit{Stabilization of the Trace Formula, Shimura Varieties, and Arithmetic Applications}, Ed. L. Clozel, M. Harris, J.-P. Labesse \& B.-C. Ng\^o, International Press (2011).\ps

\bibitem[\textsc{ABV}92]{abv} J. Adams, D. Barbasch and  D. Vogan, \textit{The Langlands classification and irreducible characters for real reductive groups}, Progress in Mathematics 104, Birkha\"user, Boston-Basel-Berlin (1992).\ps

\bibitem[\textsc{AJ}87]{AJ} J. Adams \& J. Johnson, \textit{Endoscopic groups and packets of non-tempered representations}, Compositio Mathematica, 64 no. 3, 271--309 (1987).\ps  

\bibitem[\textsc{Art}89]{arthurunipotent} J. Arthur, \textit{Unipotent automorphic representations : Conjectures}, in {\it Orbites unipotentes et repr\'esentations II : Groupes $p$-adiques et r\'eels}, Ast\'erisque 171--172 (1989).\ps 

\bibitem[\textsc{Art}02]{arthurconjlan} J. Arthur, \textit{A note on the automorphic Langlands group}, Canad. Math. Bull. Vol. 45 (4), 466--482 (2002). \ps

\bibitem[\textsc{Art}04]{arthursp4} J. Arthur,\textit{Automorphic representations of ${\rm GSp}(4)$}, in {\it Contributions  to automorphic forms, geometry, and number theory}, 65/81, Johns Hopkins Univ. Press, Baltimore (2004). \ps

\bibitem[\textsc{Art}05]{arthurlivre} J. Arthur, {\it An introduction to the trace formula}, dans {\it Harmonic analysis, the trace formula and Shimura varieties}, A.M.S.\,, Clay Math. Institute (2005).\ps

\bibitem[\textsc{Art}11]{arthur} J. Arthur, \textit{The endoscopic classification of representations: orthogonal and symplectic groups}, available at the address  {\url{http://www.claymath.org/cw/arthur/}}.\ps

\bibitem[\textsc{AS}01]{asgarischmidt} M. Asgari \& R. Schmidt, \textit{Siegel modular forms and representations}, Manuscripta math. 104, 173--200 (2001). \ps

\bibitem[\textsc{AP}08]{ashpollack} A. Ash \& D. Pollack, \textit{Everywhere unramified automorphic cohomology for ${\rm GL}(3,\Z)$}, Int. Journal of Number Theory 4, 663--675 (2008).\ps

\bibitem[\textsc{BGGT}]{bggt} T. Barnet-Lamb, T. Gee, D. Geraghty \& R.
Taylor, \textit{Potential automorphy and change of weight}, to appear in
Annals of Math.\ps

\bibitem[\textsc{BC}11]{bchsigne} J. Bella\"iche \& G. Chenevier, \textit{The sign of Galois representations attached to automorphic forms for unitary groups}, Compositio Math. 147, 1337--1352 (2011).\ps

\bibitem[\textsc{BFG}11]{BFVdG} J. Bergstr\"om, C. Faber \& G. van der Geer, \textit{Siegel modular forms of degree three and the cohomology of local systems},  arXiv preprint \url{http://lanl.arxiv.org/abs/1108.3731v2} (2011). \ps

\bibitem[\textsc{BS}59]{blijspringer1} F. van der Blij \& T. A. Springer, \textit{The arithmetic of octaves and the groups ${\rm G}_2$}, Proc. Kon. Ak. Amsterdam 62 (= Ind. Math. 21), 406--418 (1959). \ps

\bibitem[\textsc{Bor}84]{borcherdsthese} R. Borcherds, \textit{The Leech lattice and other lattices}, Ph. D. Dissertation, University of Cambridge (1984).\ps

\bibitem[\textsc{Bor}63]{borelfini} A. Borel, \textit{Some finiteness theorems for adeles groups over number fields}, Publ. math. IHES 16, 101-126 (1963). \ps

\bibitem[\textsc{Bor}77]{borelcorvallis} A. Borel, \textit{Automorphic L-functions}, dans {\it Automorphic forms, representations, and L-functions}, P.S.P.M. 33 Part 2, conf\'erence de Corvallis (1977).\ps

\bibitem[\textsc{Bor}91]{borellivre} A. Borel, \textit{Linear algebraic groups}, Springer Verlag G. T. M. 126 (1991), 1st ed. 1969.\ps

\bibitem[\textsc{BJ}79]{boreljacquet} A. Borel \& H. Jacquet, \textit{Automorphic forms
and automorphic representation}, Corvallis, P.S.P.M. {\bf 33} vol. I (1979).\ps

\bibitem[\textsc{Bou}81]{bourbaki} N. Bourbaki, \textit{\'El\'ements de math\'ematique. Groupes et alg\`ebres de Lie. Chapitres IV, V, VI},  Masson, Paris  (1981). \ps

\bibitem[\textsc{Bum}96]{bump} D. Bump, \textit{Automorphic forms and representations}, Cambridge studies in adv. Math. 45 (1996).\ps

\bibitem[\textsc{BG}11]{buzzardgee} K. Buzzard \& T. Gee, \textit{The conjectural
connections between automorphic representations and Galois representations},
Proceedings of the LMS Durham Symposium 2011.\ps

\bibitem[\textsc{CG}11]{calegarigee} F. Calegari \& T. Gee, \textit{Irreducibility of   
automorphic Galois representations of $\GL(n)$, $n$ at most 5}, to appear in Annales de l'Institut Fourier.  \ps

\bibitem[\textsc{Car}12]{caraiani} A. Caraiani, \textit{Local-global compatibility and the action of monodromy on nearby cycles}, \`a para\^itre \`a Duke Math. J.\ps

\bibitem[\textsc{Car}72]{carter} R. W. Carter, \textit{Conjugacy classes in the Weyl group}, Compositio Math. 25, 1--59 (1972)  .\ps

\bibitem[\textsc{CL}12]{chaulau} P.-H. Chaudouard \& G. Laumon, 
\textit{Le lemme fondamental pond\'er\'e. II. \'Enonc\'es
cohomologiques}, Annals of Math.  176, 1647--1781 (2012). \ps

\bibitem[\textsc{Ch}13]{chhdr} G. Chenevier, {\it Repr\'esentations galoisiennes automorphes et cons\'equances arithm\'etiques des conjectures d'Arthur}, Habilitation thesis, available at the url~\url{http://www.math.polytechnique.fr/~chenevier/hdr/HDR.pdf} (2013).\ps

\bibitem[\textsc{CC}09]{chcl} G. Chenevier \& L. Clozel, {\it Corps de nombres peu ramifi\'es et formes automorphes autoduales}, Journal of the A.M.S. 22 Vol 2, 467-519 (2009). \ps

\bibitem[CH13]{chharris} G. Chenevier \& M. Harris, {\it Construction of
automorphic Galois representations II}, Cambridge Math. Journal 1 (2013).\ps

\bibitem[\textsc{CL}11]{cl} G. Chenevier \& J. Lannes, \textit{Kneser neighbours and orthogonal Galois representations in dimensions 16 and 24}, 
Algebraic Number Theory, Oberwolfach Report 31/2011. This paper announces some results of the forthcoming paper~\cite{cl2}. See the url~\url{http://www.math.polytechnique.fr/~chenevier/niemeier/niemeier.html} for a collection of tables. \ps

\bibitem[\textsc{CL}]{cl2} G. Chenevier \& J. Lannes, \textit{Formes automorphes et voisins de Kneser des r\'eseaux de Niemeier}, 
In preparation.  \ps

\bibitem[\textsc{CR}]{homepage} G. Chenevier \& D. Renard, {\it  The home page of level one algebraic cusp forms of classical groups}, \url{http://www.math.polytechnique.fr/~chenevier/levelone.html}. \ps

\bibitem[\textsc{Clo}90]{clozel} L. Clozel, {\it Motifs et formes automorphes : applications du principe de fonctorialit\'e}, Automorphic forms, Shimura varieties, and L-functions, Vol. I (Ann Arbor, MI, 1988), Perspect. Math., vol. 10, Academic Press, Boston, MA (1990).\ps

\bibitem[\textsc{Cog}04]{cogdell} J. Cogdell, {\it Lectures on ${\rm L}$-functions, converse theorems, and functoriality for ${\rm GL}_n$}, in {\it Lectures on automorphic ${\rm L}$-functions}, A.M.S. , Fields Institute Monographs (2004). \ps

\bibitem[\textsc{CNP}96]{cn} A. Cohen, G. Nebe \& W. Plesken, \textit{Cayley orders}, Compositio Math. 103, 63--74 (1996).\ps

\bibitem[\textsc{Con}11]{conradgp} B. Conrad, \textit{Reductive group schemes},  notes de l'\'ecole d'\'et\'e SGA3 (Luminy, 2011), notes disponibles \`a l'adresse~\url{http://math.stanford.edu/~conrad/papers/}. \ps

\bibitem[\textsc{CS}99]{conwaysloane} J. Conway \& N. Sloane,  {\it Sphere packings, lattices, and groups}, 3 ed, Springer  Verlag,  Grundlehren der math. wiss. 290 (1999).  \ps

\bibitem[\textsc{Cor}77]{corvallis}  {\it Automorphic forms, representations
and ${\rm L}$-functions, Part I \& II}, Proc.  Sympos.  Pure Math.  XXXIII, Oregon
State Univ., Corvallis, Ore.  (1977), Providence, R.I.: American
Mathematical Society.  \ps

\bibitem[\textsc{Cox}46]{coxeter} H. M. S Coxeter, \textit{Integral Cayley Numbers},
Duke Math. J. 13, 567--578 (1946).\ps

\bibitem[\textsc{Del}97]{delorme} P. Delorme, \textit{Infinitesimal
character and distribution character of representations of reductive lie
groups}, A.M.S., P.S.P.M. 61, Edinburgh conference, 73--81 (1997).

\bibitem[\textsc{Dix}74]{dixmier} J. Dixmier, \textit{Alg\`ebres enveloppantes}, Gauthier-Villars \'Ed. (1974).\ps

\bibitem[\textsc{FC}90]{faltings} G. Faltings \& C.-L. Chai, {\it Degeneration of abelian varieties}, 
Springer Verlag, Ergebnisse der Math. und ihrer Grenzgebiete (1990)
 
\bibitem[\textsc{Fer}96]{fermigier} S. Fermigier, {\it Annulation de la cohomologie cuspidale de sous-groupes de congruence de $\GL_n(\Z)$}, Math. Annalen 306, 247--256 (1996). \ps


\bibitem[\textsc{Fon}85]{fontaineabelienne} J.-M. Fontaine, {\it 
Il n'y a pas de vari\'et\'e ab\'elienne sur $\Z$}, 
Invent. Math. 81, 515--538 (1985). \ps

\bibitem[\textsc{Fon}93]{fontaineprojlisse} J.-M. Fontaine, {\it 
Sch\'emas propres et lisses sur $\Z$}, 
Proceedings of the Indo-French Conference on Geometry (Bombay, 1989), 43--56,
Hindustan Book Agency, Delhi (1993). \ps

\bibitem[\textsc{FM}95]{fontainemazur} J.-M. Fontaine, B. Mazur, {\it Geometric
Galois representations}, in Coates, John; Yau., S.-T., {\it Elliptic curves,
modular forms, \& Fermat's last theorem}, Series in Number Theory 1, Int.
Press, Cambridge, MA, 41-78 (1995).  \ps


\bibitem[\textsc{GGS}02]{gangrosssavin}  W. T. Gan, B. Gross \& G. Savin, {\it Fourier coefficients of modular forms on ${\rm G}_2$}, Duke Math. 115, 105-169 (2002). \ps

\bibitem[\textsc{GG}05]{gangurevich1} W. T. Gan \& N. Gurevich, {\it Non-tempered Arthur packets of ${\rm G}_2$}, proceedings of Rallis' 60th birthday conference: {\it Automorphic Representations, L-functions and Applications : Progress and Prospects}, 129-155 (2005). \ps

\bibitem[\textsc{GG}06]{gangurevich2} W. T. Gan \& N. Gurevich, {\it Non-tempered Arthur packets of ${\rm G}_2$: liftings from $\tilde{SL}_2$}, American Journal of Math 128, 1105-1185 (2006). \ps

\bibitem[\textsc{Gee}08]{vdg} G. van der Geer, \textit{Siegel modular forms and their applications}, The 1-2-3 of modular forms, 181--245, Universitext, Springer, Berlin (2008). \ps

\bibitem[\textsc{GJ}78]{GJ} S. Gelbart \& H. Jacquet, {\it A relation between automorphic representations of $\GL(2)$ and $\GL(3)$}, Annales Sci. \'E.N.S. 11, 1--73 (1978).\ps


\bibitem[\textsc{GGPS}66]{GGPS} I. M. Gel'fand, M. I. Graev \& I. I. Pyatetskii-Shapiro, {\it Representation theory and automorphic functions}, Academic Press (1990), 1st ed. 1966. \ps

\bibitem[GP]{GP} PARI/GP, version {\tt 2.5.0}, Bordeaux, 2011, \url{http://pari.math.u-bordeaux.fr/}.\ps

\bibitem[\textsc{Gol}12]{goldring} W. Goldring, \textit{Galois representations associated to holomorphic limits of discrete series I: Unitary Groups}, to appear in Compositio Math. \ps

\bibitem[\textsc{GW}98]{GW}R. Goodman \& N. Wallach  {\it Representations and Invariants of the Classical Groups},  Cambridge U. Press, (1998).\ps

\bibitem[\textsc{Gro}96]{grossinv} B. Gross, \textit{Reductive groups over $\Z$}, Invent. math.124, 263--279 (1996). \ps

\bibitem[\textsc{Gro}98]{grossatake} B. Gross, \textit{On the Satake isomorphism},
dans {\it Galois representations in arithmetic algebraic geometry}, A.
Scholl \& R.  Taylor Ed., Cambridge university press (1998).\ps

\bibitem[\textsc{Gro}99]{grossalgform} B. Gross, \textit{Algebraic modular forms}, Israel J. Math. 113, 61--93 (1999).\ps

\bibitem[\textsc{GP}05]{grosspollack} B. Gross \& D. Pollack, {\it On the Euler characteristic of the discrete spectrum}, Journal of Number Theory, 110 (2005) no. 1, 136-163. \ps

\bibitem[\textsc{GS}98]{grosssavin} B. Gross \& G. Savin, \textit{Motives with Galois
group of type ${\rm G}_2$ : an exceptional theta correspondence},  Compositio Math. 114, 153--217 (1998).\ps

\bibitem[\textsc{GRFA}11]{grfa} \textit{Stabilization of the Trace Formula, Shimura Varieties, and Arithmetic Applications}, \'edit\'e par L. Clozel, M. Harris, J.-P. Labesse et B.-C. Ng\^o, International Press (2011).\ps

\bibitem[\textsc{HC}68]{harishchandra} Harish-Chandra, \textit{Automorphic forms on semisimple Lie groups}, Springer Verlag, Lecture notes in Math. (1968). \ps

\bibitem[\textsc{Iwa}49]{iwasawa} K. Iwasawa, \textit{On some types of topological groups}, Annals of Math.  50, 507--558 (1949).\ps

\bibitem[\textsc{JS}81]{jasha} H. Jacquet \& J. Shalika, \textit{On Euler Products and the Classification of automorphic representations I,II}, American Journal of Mathematics 103, I. 499--558, II. 777--815 (1981).\ps

\bibitem[\textsc{Kha}06]{khare1} C. Khare, {\it 
Serre's modularity conjecture: the level one case}, 
Duke Math. J. 134 , 557--589 (2006).\ps

\bibitem[\textsc{Kha}07]{kh} C. Khare, {\it Modularity of Galois representations and motives with good reduction property}, J. Ramanujan math. soc. 22, 1--26 (2007). \ps

\bibitem[\textsc{KS}02]{kimsha} H. Kim \& F. Shahidi, {\it Cuspidality of symmetric powers with applications}, Duke math. journal 112, 177--197 (2002). \ps

\bibitem[\textsc{Kna}86]{knapplivre} A. Knapp, {\it Representation theory of semisimple Lie groups}, Princeton univ. Press (1986).\ps

\bibitem[\textsc{Kna}94]{knapp} A. Knapp, {\it Local Langlands correspondence: the archimedean case}, p. 393--410 in {\it Motives, Part 2} P.S.P.M. 55, AMS Providence, RI (1994).\ps

\bibitem[\textsc{KnV}95]{knappvogan} A. Knapp \& D. Vogan, {\it Cohomological induction and unitary representations}, Princeton Univ. Press (1995).\ps

\bibitem[\textsc{Kot}84]{kottwitzctt} R. Kottwitz, {\it Stable trace formula: Cuspidal tempered terms}, Duke Math. J. Volume 51, Number 3\,\, 611--650 (1984). \ps

\bibitem[\textsc{Kot}88]{kottwitz} R. Kottwitz, {\it Shimura varieties and $\lambda$-adic representations}, in {\it Automorphic forms, Shimura varieties, and L-functions}, Ann Arbor conference, Vol. I (1988).\ps

\bibitem[\textsc{LL}79]{labesselanglands} J.-P. Labesse \& R. Langlands, {\it ${\rm L}$-indistinguishability for {\rm SL}(2)}, Canadian Journal of Mathematics 31, 726Ð785 (1979).\ps

\bibitem[\textsc{LW}13]{labwals} J.-P. Labesse \& J.-L. Waldspurger, {\it 
La formule des traces tordue d'apr\`es le Friday
Morning Seminar}, with a foreword by Robert Langlands, CRM
Monograph Series, 31. American Mathematical Society, Providence, RI
(2013).\ps

\bibitem[\textsc{Lan}70]{langlandsyale} R. Langlands, {\it Problems in the theory of
automorphic forms}, Springer Lecture Notes 170 (1970).\ps

\bibitem[\textsc{Lan}73]{langlandsreel} R. Langlands, {\it The classification of representations of real reductive groups}, A.M.S. Math. Surveys and Monographs 31 (1973).\ps

\bibitem[\textsc{Lan}79]{Lgl} R.  Langlands, {\it Ein M\"archen}, in {\it Automorphic representations, Shimura varieties, and motives}, dans P.S.P.M. 33, AMS Providence, 205--246 (1979).\ps

\bibitem[\textsc{Lan}96]{langlandsedin} R, Langlands, {\it Where stands
functoriality today ?}, A.M.S.  P.S.P.M.  61, Edinburgh (1996).  \ps

\bibitem[\textsc{LP}02]{lanskypollack} J. Lansky \& D. Pollack, {\it Hecke algebras and automorphic forms}, Composition math. 130, 21-48 (2002). \ps

\bibitem[\textsc{Mes}86]{mestre} J.-F. Mestre, {\it Formules explicites et minorations de conducteurs des vari\'et\'es alg\'ebriques}, Compositio Math. 58, 209--232 (1986). \ps

\bibitem[\textsc{Mez}a]{mezo} P. Mezo, {\it Character identities in the twisted endoscopy of real reductive groups}, Memoirs of the A.M.S. (to appear).\ps

\bibitem[\textsc{Mez}b]{mezopreprint} P. Mezo, {\it Spectral transfer in the theory of twisted endoscopy 
of real groups}, available at the url~{\url{http://people.math.carleton.ca/~mezo/research.html}}.
\ps

\bibitem[\textsc{Mil}02]{miller} S. Miller, {\it The highest-lowest zero and other applications of positivity}, Duke math. J.  112,  83-116 (2002).\ps

\bibitem[\textsc{Mot}94]{motives} {\it Motives}, P.S.P.M.  55, AMS Providence, RI (1994).\ps

\bibitem[\textsc{Ng\^o}10]{ngo} B. C. Ng\^o, {\it 
Le lemme fondamental pour les alg\`ebres de Lie}, Publ. 
Math. Inst. Hautes \'Etudes Sci. No. 111, 1--169 (2010). \ps 

\bibitem[\textsc{PR}94]{pr} V. Platonov \& A. Rapinchuk, {\it Algebraic groups and number theory}, Pure and Applied Mathematics, 139. Academic Press, Inc., Boston (1994), 1st ed. 1991.\ps

\bibitem[\textsc{Pol}98]{pollack} D. Pollack, {\it Explicit Hecke actions on modular
forms}, PhD. Thesis Harvard univ. (1998).\ps

\bibitem[\textsc{RRST}]{rrst} M. Raum, N. C. Ryan, N.-P. Skoruppa, G. Tornar\`ia, {\it Explicit computations of Siegel modular forms of degree two}, ArXiV preprint \url{http://arxiv.org/pdf/1205.6255.pdf}. \ps

\bibitem[\textsc{Sat}63]{satake} I. Satake, {\it Theory of spherical functions on reductive algebraic groups over $p$-adic fields}, Publ. math. I.H.\'E.S. {\bf 18} (1963), 5--69.\ps

\bibitem[\textsc{Ser}68]{serreabelian} J.-P. Serre, {\it Abelian $\ell$-adic
representations and elliptic curves}, W. A. Benjamin, New York (1968).\ps

\bibitem[\textsc{Ser}70]{serre} J.-P. Serre, {\it Cours d'arithm\'etique}, P.U. F.\, , Paris (1970). \ps

\bibitem[\textsc{Ser}94]{serremotives} J.-P. Serre, {\it Propri\'et\'es conjecturales   
des groupes de Galois motiviques et des repr\'esentations $\ell$-adiques},
in~\cite{motives}. \ps


\bibitem[\textsc{Ser}97]{serregc} J.-P.Serre, {\it Cohomologie galoisienne}, Springer Lecture notes in Math.  5, edition 5 (1997). \ps



\bibitem[\textsc{SGA}3]{sga} M. Demazure \& A. Grothendieck, {\it  S\'eminaire de G\'eom\'etrie Alg\'ebrique du Bois Marie - 1962-64 - Sch\'emas en groupes - (SGA 3)} - tome 3, Soc. Math. France, 
Documents math\'ematiques 8, edited by P. Gille and P. Polo (2011). \ps


\bibitem[\textsc{She}82]{shelun}  D. Shelstad, {\it L-Indistinguishablility for Real Groups}, Math. Ann. 259, 385--430 (1982). \ps

\bibitem[\textsc{She}08]{shelstad} D. Shelstad, {\it Tempered endoscopy for real groups I: geometric transfer with canonical factors}, 
Contemporary Math,  Vol. 472 , 215--246 (2008).\ps

\bibitem[\textsc{Shi}11]{shin} S.-W. Shin, {\it Galois representations arising from some compact Shimura varieties}, Annals of Math.173, 1645--1741 (2011).\ps

\bibitem[\textsc{Sko}92]{sko} N. Skoruppa, {\it Computations of Siegel modular forms of genus 2},  Mathematics of Computation 58, (1992), 381--398. \ps

\bibitem[\textsc{Tai}12]{taibi} O. Ta\"ibi, {\it Eigenvarieties for classical groups and complex conjugations in Galois representations}, {\url{http://www.math.ens.fr/~taibi/}}. \ps

\bibitem[\textsc{Tat}79]{tate} J. Tate, {\it Number theoretic background}, in {\it Automorphic forms, representations and L-functions}, Proc. Sympos. Pure Math. XXXIII Part 2,  3--26 (1979).\ps

\bibitem[\textsc{Tit}79]{titscorvallis} J. Tits, {\it  Reductive groups over local fields}, dans~\cite{corvallis} Part I, 29-69 (1979). \ps

\bibitem[\textsc{Tsu}83]{tsushima} T. Tsushima, {\it An Explicit Dimension Formula for the Spaces of Generalized Automorphic Forms with Respect to ${\rm Sp}_2(\Z)$}, Proc. Japan Acad. 59 (1983).\ps

\bibitem[\textsc{Tsu}86]{tsuyumine} S. Tsuyumine, {\it Siegel modular forms of degree $3$}, American Journal of Mathematics 108 (1986), 755--862.\ps

\bibitem[\textsc{Wal}84]{wallach} N. Wallach, {\it On the constant term of a square integrable automorphic form}, Operator          
algebras and group representations, Vol. II (Neptun, 1980), 227-237, Monogr.    
Stud. Math., 18, Pitman, Boston, MA (1984). \ps

\bibitem[\textsc{Wal}09]{walspond} J.-L. Waldspurger, {\it 
Endoscopie et changement de caract\'eristique: int\'egrales orbitales
pond\'er\'ee}, Ann. Inst. Fourier (Grenoble) 59, 1753--1818
(2009).\ps 



}

\end{thebibliography}
